\documentclass[a4paper,reqno]{amsart}

\usepackage[utf8]{inputenc}
\usepackage{amssymb}
\usepackage{enumitem}
\usepackage{mathrsfs}
\usepackage{mathtools}
\setlist[enumerate]{label=\textnormal{(\roman*)}}

\usepackage{hyperref}
\usepackage{amsthm}
\usepackage{cleveref}

\usepackage{xcolor}

\newtheorem{theorem}{Theorem}[section]
\newtheorem{corollary}[theorem]{Corollary}
\newtheorem{lemma}[theorem]{Lemma}
\newtheorem{proposition}[theorem]{Proposition}

\theoremstyle{definition}
\newtheorem{definition}[theorem]{Definition}
\newtheorem{remark}[theorem]{Remark}

\numberwithin{equation}{section}
\newcommand\RR{\mathbb{R}}
\newcommand\ZZ{\mathbb{Z}}

\title[Decay of Schr\"odinger maps]
{Optimal decay and pointwise asymptotics for two dimensional Schr\"odinger maps}

\author[Z.~Li]{Ze Li}
\address{School of Mathematics and Statistics,
Ningbo University, Ningbo 315211, Zhejiang, China}
\email{rikudosennin@163.com}

\author[X.~Sun]{Xiaotong Sun}
\address{School of Mathematical Sciences,
University of Science and Technology of China, Hefei 230026, Anhui, China}
\email{sunxt@mail.ustc.edu.cn}

\author[L.~Zhao]{Lifeng Zhao}
\address{School of Mathematical Sciences,
University of Science and Technology of China, Hefei 230026, Anhui, China}
\email{zhaolf@ustc.edu.cn}

\thanks{}



\begin{document}

\begin{abstract}
The Schr\"odinger map equation describes the evolution of the unit magnetization in an idealized isotropic ferromagnet under the exchange interaction. In two space dimensions, the equation is energy critical. Beyond global well-posedness and scattering in the energy topology, a basic physical question is how rapidly the amplitude of a localized spin excitation decays and what waveform describes it as it spreads. 
For small, spatially localized initial perturbations of a constant state $Q$, we prove that both $\phi-Q$ and $\nabla\phi$ decay at the sharp rate $t^{-1}$ in $L_x^\infty$, assuming only $H^{0^+}$ regularity for $\nabla\phi_0$. We further show that a free Schr\"odinger evolution in the tangent plane at $Q$ gives the leading magnetization waveform with an error that is $o(t^{-1})$ uniformly in space. The proof combines the caloric gauge with the method of spacetime resonance. On the one hand, we exploit the negative derivative  structure of the caloric connection and the Fourier bound of the profile, to propagate frequency envelopes associated with the $\dot{H}^{-\frac{1}{2}}$ norm and the Galilean type field along the heat flow and obtain decay in physical time. On the other hand, in the spacetime resonance analysis, we use the structure of the connection to track  interactions between symbols and derivatives of the phase, where the hidden null structures are discovered. 

\end{abstract}

\maketitle

\tableofcontents

\section{Introduction}\label{sec:1}

We consider the Schr\"odinger map equation
\begin{equation}\label{eschrodinger-map}
\left\{
\begin{aligned}
&\partial_t\phi=\phi\times\Delta\phi, &&(t,x)\in\RR\times\RR^2,\\
&\phi(0,x)=\phi_0(x), &&x\in\RR^2,
\end{aligned}
\right.
\end{equation}
where $\phi:\RR\times\RR^2\to\mathbb S^2\subset\RR^3$ is the unit magnetization field. 
It describes the conservative evolution of the unit magnetization in an idealized isotropic ferromagnet when only the exchange interaction is considered, see \cite{LL,LP,SSB}.
The exchange interaction produces an effective field proportional to $\Delta\phi$, and the cross product describes the precession of the magnetization about this field while preserving the constraint $|\phi|=1$.

The Schr\"odinger map equation is also the Hamiltonian flow associated with the Dirichlet energy. For a K\"ahler target with complex structure $J$, its geometric form is $\partial_t\phi=J(\phi)\tau(\phi)$, where $\tau(\phi)$ is the tension field \cite{D2002,T2012,TU}.  For the sphere, $J(\phi)X=\phi\times X$ for $X\in T_\phi\mathbb S^2$, which gives \eqref{eschrodinger-map}.

For a fixed $Q\in\mathbb S^2$, sufficiently regular solutions of finite energy with $\phi_0-Q\in L_x^2$ satisfy the conservation laws for mass and energy:
\begin{align}
\|\phi(t)-Q\|_{L_x^2}^2&=\|\phi_0-Q\|_{L_x^2}^2,\label{emass-conservation}\\
\|\nabla\phi(t)\|_{L_x^2}^2&=\|\nabla\phi_0\|_{L_x^2}^2.\label{enerconser}
\end{align}
On $\RR^d$, the equation is invariant under the scaling $\phi(t,x)\mapsto\phi(\lambda^2t,\lambda x)$, which preserves the $\dot H^{\frac{d}{2}}$ seminorm of the initial map. In dimension two, this scaling preserves the exchange energy $\frac{1}{2}\|\nabla\phi\|_{L_x^2}^2$, so the equation is energy critical.

In the classical description of Lifshitz and Pitaevskii \cite[Section~69]{LP}, spin waves are ``oscillations of the magnetic moments, propagated through the lattice.'' The excitation, propagation, and attenuation of spin waves are fundamental topics in magnetization dynamics, as discussed in  Chapter~6 of Stancil and Prabhakar \cite{SP}.  Understanding how the amplitude and waveform of these oscillations evolve locally during propagation is therefore a basic physical question, directly addressed by time-resolved measurements of coherent spin waves \cite{BWSXC}. However, local magnetization dynamics can vary substantially across spatial regions and differ from the spatially averaged observations \cite{ZKHLSAR}, making a spatially resolved description essential to understanding the evolving magnetization field. 

A central mathematical problem is therefore to determine how rapidly the amplitude of the perturbation $\phi(t,x)-Q$ decays and to identify the leading waveform governing its long time evolution, uniformly in space. The natural amplitude scale is $t^{-1}$: linearization about the constant state $Q$ yields the free Schr\"odinger equation on $T_Q\mathbb S^2\simeq\mathbb C$, whose solutions with suitably localized initial data exhibit this dispersive decay in two space dimensions. For the nonlinear flow in two dimensions, Bejenaru, Ionescu, Kenig, and Tataru \cite{BIKT2011} established global well-posedness for sphere-valued maps with small energy, and Li \cite{L2023} proved both scattering in the energy topology and uniform convergence to $Q$ for smooth data of small energy. Scattering in $\dot H^1$ provides a fundamental description of the long time dynamics in the exchange energy topology, but does not by itself resolve this pointwise asymptotic problem. Uniform convergence to $Q$ shows that $\|\phi(t)-Q\|_{L_x^\infty}$ tends to zero, but does not determine how rapidly the perturbation decays. 

In this paper, we aim to establish the refined pointwise asymptotics of \eqref{eschrodinger-map} including the  
optimal decay and identifying the leading waveform with an error that is $o(t^{-1})$, uniformly in space.

\subsection{Historical background} 

The local well-posedness for data of arbitrary size was proved in \cite{DW,KK,KLPST,KN,M,SSB}. For maps into $\mathbb S^2$ near a constant, Ionescu and Kenig \cite{IK2006} used stereographic coordinates to prove local well-posedness in $H^s$ for $s>\frac{d+1}{2}$ in dimensions $d\geq2$, and Bejenaru \cite{B2008b} reached the full subcritical range $s>\frac{d}{2}$.  

Early global results include the one dimensional theory of Zhou, Guo, and Tan \cite{ZGT}. Chang, Shatah, and Uhlenbeck \cite{CSU} used a moving frame to obtain global results at finite energy in one dimension and at small energy under symmetry assumptions in two dimensions. At critical regularity, Ionescu and Kenig \cite{IK2007} and, independently, Bejenaru \cite{B2008a} proved global well-posedness for small data in $\dot B^{\frac{d}{2}}_{2,1}$ in dimensions $d\geq3$. In dimensions $d\geq4$, Bejenaru, Ionescu, and Kenig \cite{BIK} applied the Coulomb gauge to prove global regularity for smooth data small in $\dot H^{\frac{d}{2}}$. Using the caloric gauge, Bejenaru, Ionescu, Kenig, and Tataru \cite{BIKT2011} established global well-posedness for all $d\geq2$ under smallness in $\dot H^{\frac{d}{2}}$ and the low frequency condition $\phi_0-Q\in\dot H^{\frac{d}{2}-1}$, and showed that higher Sobolev norms remain uniformly bounded in time. Bejenaru \cite{B2026} developed a bilinear approach to global results for derivative nonlinearities with null structure, including Schr\"odinger maps. The small data theory for  compact K\"ahler targets was completed by Li \cite{L2021,L2023} in critical Sobolev spaces in $d\ge 2$  and by Dodson and Marzuola \cite{DM} in $\dot{B}^{1}_{2,\infty}$ in $d=2$.

Beyond the small data theory, progress has also been made in the study of large data dynamics under symmetry assumptions. For $1$-equivariant maps, Bejenaru, Ionescu,  Kenig, and Tataru  \cite{BIKT2013,BIKT2016} proved that global well-posedness and scattering of the differential fields hold below energy $4\pi$ for $\mathbb S^2$ and at arbitrary finite energy for $\mathbb H^2$. The ground state harmonic maps are proved to be orbitally stable by Gustafson, Kang, and Tsai \cite{GKT2007}, and asymptotically stable for $|m|\geq4$ by \cite{GKT2008} and $|m|=3$ by Gustafson, Nakanishi, and  Tsai \cite{GNT2010}. For $|m|=1$, the ground state is unstable in the energy topology but stable in a stronger topology, see Bejenaru and Tataru \cite{BT}. For $|m|=2$, Bejenaru, Pillai, and Tataru \cite{BPT} proved global existence near the ground state in the energy topology, dispersive and modulation bounds, and stability in a stronger Besov topology. Stability near harmonic maps on the hyperbolic plane was studied in Lawrie, L\"uhrmann,  Oh, and  Shahshahani \cite{LLOS} and Li \cite{L2022}. For 1-equivariant maps from $\RR^2$ into $\mathbb S^2$, finite time blowup solutions with various rates were constructed by \cite{JXYZ,MRR,P}, see also \cite{Y2026} for multi-bubbles in  $m\geq3$. Without the equivariant assumptions, 
Smith studied global regularity for dispersed energy data \cite{S2011,S2013}, while Dodson and Smith obtained a continuation criterion \cite{DS}.

These developments have revealed a rich range of global and singular dynamics for Schr\"odinger maps. However, determining the local amplitude and leading waveform of a localized spin excitation as it propagates in the plane remains an open problem. Even for small perturbations of a uniformly magnetized state $Q$ (constant map), neither the precise decay rate of the local oscillations nor their leading spatial profile has been established.  Resolving the leading spatial and temporal oscillations therefore requires an approximation to $\phi(t,x)-Q$ whose error is $o(t^{-1})$, uniformly for $x\in\RR^2$.

\subsection{Main results}

In this paper, we study two dimensional Schr\"odinger maps with small, localized initial data near a constant equilibrium $Q$, assuming only $H^\epsilon$ regularity for $\nabla\phi_0$. Our main result has two parts: optimal decay and spatially pointwise asymptotics. First, we prove that both $\|\phi(t)-Q\|_{L_x^\infty}$ and $\|\nabla\phi(t)\|_{L_x^\infty}$ decay at the sharp rate $t^{-1}$. Thus the deviation of the magnetization from equilibrium and its  spatial gradient both decay at the natural rate of the two dimensional free Schr\"odinger evolution. Second, we show that a free Schr\"odinger evolution in the tangent plane at $Q$ describes the leading magnetization waveform with an $O(t^{-2+\frac{\beta}{2}})$ error in $L_x^\infty$, where $0<\beta<\epsilon$. This determines the waveform uniformly in space with an error nearly one order smaller than its natural $t^{-1}$ amplitude scale, and we further identify the quadratic correction in the normal direction and obtain explicit stationary phase asymptotics.

The main result is as follows. 
\begin{theorem}\label{thm:main}
Fix $0<\beta<\epsilon\ll1$. Let $Q\in\mathbb S^2$ and let $(v_\infty,w_\infty)$ be a positively oriented orthonormal basis of $T_Q\mathbb S^2$, with $w_\infty=Q\times v_\infty$. Suppose that $\phi_0:\RR^2\to\mathbb S^2$ satisfies
\begin{equation}\label{ini}
\|\nabla\phi_0\|_{H^\epsilon}+\|\langle x\rangle^2\langle\nabla\rangle(\phi_0-Q)\|_{L_x^2}\leq\epsilon_0,
\end{equation}
where $\epsilon_0>0$ is sufficiently small. There exists a complex-valued scattering state $g_+\in H^{1+\epsilon}(\RR^2)$ satisfying
\begin{equation}\label{ethm1.1-scat}
\|g_+\|_{H^{1+\epsilon}}+\|\langle x\rangle^2g_+\|_{L_x^2}\lesssim\epsilon_0
\end{equation}
such that the following conclusions hold.

\begin{enumerate}
\item[(i)] \textit{Global existence, uniform bounds, and scattering.} Equation \eqref{eschrodinger-map} admits a unique global solution in the small data solution class of \cite{BIKT2011}, with
\begin{equation}\label{ethm1.1-glo}
\phi-Q\in C\bigl(\RR,H^{1+\epsilon}(\RR^2,\RR^3)\bigr).
\end{equation}
Moreover,
\begin{equation}\label{ethm1.1-sobolev}
\sup_{t\in\RR}\|\phi(t)-Q\|_{H^{1+\epsilon}}\lesssim\|\phi_0-Q\|_{H^{1+\epsilon}}.
\end{equation}
The solution scatters in $H^\epsilon$, with
\begin{equation}\label{ethm1.1-scattering}
\bigl\|\phi(t)-Q-v_\infty\operatorname{Re}(e^{it\Delta}g_+)-w_\infty\operatorname{Im}(e^{it\Delta}g_+)\bigr\|_{H^\epsilon}\lesssim\epsilon_0^2t^{-1},\qquad t\geq1.
\end{equation}

\item[(ii)] \textit{Optimal decay.} For every $t\geq1$,
\begin{align}
\|\nabla\phi(t)\|_{L_x^\infty}&\lesssim\epsilon_0t^{-1},\label{ethm1.1-decay1}\\
\|\phi(t)-Q\|_{L_x^\infty}&\lesssim\epsilon_0t^{-1}.\label{th1.1decay2}
\end{align}

\item[(iii)] \textit{Pointwise asymptotics.} For every $t\geq1$,
\begin{equation}\label{ethm1.1-uniform1}
\phi(t)=Q+v_\infty\operatorname{Re}(e^{it\Delta}g_+)+w_\infty\operatorname{Im}(e^{it\Delta}g_+)+O_{L_x^\infty}\bigl(\epsilon_0^2t^{-2+\frac{\beta}{2}}\bigr).
\end{equation}
\end{enumerate}
\end{theorem}

The leading radiation in \eqref{ethm1.1-uniform1} takes values in $T_Q\mathbb S^2$. The following corollary identifies the quadratic correction in the normal direction and gives explicit stationary phase formulas for both the map and its normal component.

\begin{corollary}\label{cmap asymptotics}
Under the assumptions of Theorem~\ref{thm:main}, let $g_+$ be its scattering state. For every $t\geq1$,
\begin{equation}\label{ethm1.1-uniform2}
\phi(t)=Q+v_\infty\operatorname{Re}(e^{it\Delta}g_+)+w_\infty\operatorname{Im}(e^{it\Delta}g_+)-\frac{1}{2}|e^{it\Delta}g_+|^2Q+O_{L_x^\infty}\bigl(\epsilon_0^3t^{-2+\frac{\beta}{2}}\bigr),
\end{equation}
and the normal component satisfies
\begin{equation}\label{ethm1.1-uniform3}
\left\|Q\cdot\phi(t)-1+\frac{1}{2}|e^{it\Delta}g_+|^2\right\|_{L_x^\infty}\lesssim\epsilon_0^4t^{-3+\frac{\beta}{2}}.
\end{equation}
The explicit stationary phase expansion is
\begin{equation}\label{ethm1.1-uniform4}
\begin{aligned}
\phi(t,x)=&Q+v_\infty\operatorname{Re}\left[\frac{e^{\frac{i|x|^2}{4t}}}{2it}\widehat g_+\left(\frac{x}{2t}\right)\right]\\
&+w_\infty\operatorname{Im}\left[\frac{e^{\frac{i|x|^2}{4t}}}{2it}\widehat g_+\left(\frac{x}{2t}\right)\right]+O_{L_x^\infty}\bigl(\epsilon_0t^{-\frac{3}{2}}\bigr),
\end{aligned}
\end{equation}
and the corresponding expansion for the normal component is
\begin{equation}\label{ethm1.1-uniform5}
\left\|Q\cdot\phi(t)-1+\frac{1}{8t^2}\left|\widehat g_+\left(\frac{x}{2t}\right)\right|^2\right\|_{L_x^\infty}\lesssim\epsilon_0^2t^{-\frac{5}{2}}.
\end{equation}
\end{corollary}

\begin{remark}[Localization of the initial map]
The smallness condition on the undifferentiated second moment in \eqref{ini} can be relaxed. Suppose that
\begin{equation}\label{thm1.1-2}
\|\nabla\phi_0\|_{H^\epsilon}+\|\langle x\rangle^2\nabla\phi_0\|_{L_x^2}\leq\epsilon_0
\end{equation}
and that
\begin{equation}\label{ethm1.1-init3}
\|\langle x\rangle^2(\phi_0-Q)\|_{L_x^2}<\infty.
\end{equation}
Integration by parts in $\RR^2$ gives
\begin{equation*}
\|\phi_0-Q\|_{L_x^2}\leq\||x|\nabla\phi_0\|_{L_x^2},\qquad \||x|(\phi_0-Q)\|_{L_x^2}\leq\frac{1}{2}\|\langle x\rangle^2\nabla\phi_0\|_{L_x^2}.
\end{equation*}
Hence,
\begin{equation*}
\|\phi_0-Q\|_{H^{1+\epsilon}}+\|\langle x\rangle(\phi_0-Q)\|_{L_x^2}\lesssim\epsilon_0.
\end{equation*}
The second commuted energy estimate for the stereographic coordinate requires only finiteness of the initial second moment to keep the second moment of its profile uniformly bounded in time. Thus the regularity and scattering conclusions, as well as the time exponents in the decay and asymptotic estimates, remain valid under \eqref{thm1.1-2} and \eqref{ethm1.1-init3}. The constants in estimates involving the second moment may now depend on its initial size, and in particular $\|\langle x\rangle^2g_+\|_{L_x^2}$ need no longer be $O(\epsilon_0)$.
\end{remark}

\begin{remark}[Regularity and moment assumptions]
The $\|\nabla\phi_0\|_{H^\epsilon}$ condition in \eqref{ini} only requires regularity slightly higher than the energy. 
And the second spatial moments assumption is natural for studying dispersive decay in two dimensional Schr\"odinger equations, see e.g. Germain, Masmoudi, and Shatah \cite{GMS2012a}. 
\end{remark}

\begin{remark}[Higher regularity]
Suppose that $\|\nabla\phi_0\|_{H^{1+\epsilon}}\leq\epsilon_0$ with the same localization as in \eqref{ini}. Since $H^{1+\epsilon}(\RR^2)\hookrightarrow L^\infty(\RR^2)$, the uniform Sobolev bounds also control $\nabla\phi$ in $L_x^\infty$ on $[0,1]$. Thus the decay estimates \eqref{ethm1.1-decay1} and \eqref{th1.1decay2} hold with $\langle t\rangle^{-1}$ for all $t\geq0$. Using this paper's arguments, one can show that the additional derivative also removes the loss $t^{\frac{\beta}{2}}$ in the comparison with the exact free evolution. For all $t\geq0$, the $L_x^\infty$ remainder is then $O(\epsilon_0^2\langle t\rangle^{-2})$, or $O(\epsilon_0^3\langle t\rangle^{-2})$ after the quadratic normal correction is included, while the remainder in the normal component is $O(\epsilon_0^4\langle t\rangle^{-3})$.
\end{remark}

\begin{remark}[Sharpness of the decay rates]
Estimate \eqref{th1.1decay2} measures the largest displacement of the magnetization from equilibrium, while \eqref{ethm1.1-decay1} controls its spatial variation. Both rates are attained for every nonconstant initial map. More precisely,
\begin{equation*}
\lim_{t\to\infty}t\|\phi(t)-Q\|_{L_x^\infty}=\frac{1}{2}\|\widehat g_+\|_{L_\xi^\infty},\qquad \lim_{t\to\infty}t\|\nabla\phi(t)\|_{L_x^\infty}=\frac{1}{2}\|\widehat{\nabla g_+}\|_{L_\xi^\infty}.
\end{equation*}
The first identity follows from \eqref{ethm1.1-uniform4} and the orthonormality of $(v_\infty,w_\infty)$. The second follows from $|\nabla\phi|=|\psi|$ and the stationary phase formula for the differential fields in Corollary~\ref{cpsi scattering}. By \eqref{emass-conservation} and \eqref{ethm1.1-scattering}, $\|g_+\|_{L_x^2}=\|\phi_0-Q\|_{L_x^2}$. Both limits are positive unless $\phi_0=Q$. Thus, as a localized spin excitation spreads, its peak amplitude decays like $t^{-1}$ and its peak exchange energy density $\frac{1}{2}\|\nabla\phi(t)\|_{L_x^\infty}^2$ decays like $t^{-2}$, while the total exchange energy remains conserved.
\end{remark}

\begin{remark}[Uniformity of the asymptotics]
The remainder in \eqref{ethm1.1-uniform1} is $O_{L_x^\infty}(\epsilon_0^2t^{-2+\frac{\beta}{2}})$ and hence decays strictly faster than the natural $t^{-1}$ amplitude scale of the radiation. The free evolution therefore identifies the leading magnetization waveform on this scale, with the same absolute error bound at a fixed observation point and along any path $x=x(t)$.
\end{remark}

\begin{remark}[Longitudinal magnetization]
The component $Q\cdot\phi$ is the magnetization measured along the equilibrium direction. Since $|\phi|=|Q|=1$, the exact identity $1-Q\cdot\phi=\frac{1}{2}|\phi-Q|^2$ shows why its first deviation from equilibrium is quadratic, even though the transverse oscillation is of order $t^{-1}$. In particular, \eqref{ethm1.1-uniform5} gives
\begin{equation*}
\sup_{\xi\in\RR^2}\left|t^2\bigl(1-Q\cdot\phi(t,2t\xi)\bigr)-\frac{1}{8}|\widehat g_+(\xi)|^2\right|\lesssim\epsilon_0^2t^{-\frac{1}{2}}.
\end{equation*}
Thus, along each ray, the rescaled longitudinal magnetization deficit recovers the squared amplitude of the leading radiation, up to the factor $\frac{1}{8}$, but carries no information about its phase. The sphere constraint makes this $t^{-2}$ contribution visible, whereas projecting \eqref{ethm1.1-uniform1} directly onto $Q$ gives only an $O_{L_x^\infty}(\epsilon_0^2t^{-2+\frac{\beta}{2}})$ bound.
\end{remark}

\begin{remark}[Choice of tangent frame]
Once $(v_\infty,w_\infty)$ is fixed, \eqref{ethm1.1-scattering} determines $g_+$ uniquely by unitarity in $L^2$. Rotating this basis through an angle $\theta$ replaces $g_+$ by $e^{-i\theta}g_+$. Thus the leading tangent wave in \eqref{ethm1.1-uniform1} is unchanged at every $(t,x)$ under a rotation of the basis, and the longitudinal profile $|\widehat g_+|^2$ in \eqref{ethm1.1-uniform5} is likewise independent of this choice.
\end{remark}

\begin{remark}[Decay and asymptotics for related models]
Global well-posedness alone does not describe the pointwise behavior of a solution at large times. For the massless Maxwell-Klein-Gordon system on $\RR^{1+3}$, Yang and Yu \cite{YY} proved pointwise peeling estimates for large data satisfying weighted Sobolev bounds, allowing nonzero charge. After the charge contribution is separated, these estimates distinguish the decay rates of the null components of the electromagnetic field. Wei, Yang, and Yu \cite{WYY} obtained pointwise decay estimates for Yang-Mills-Higgs fields in the future of a hyperboloid under weighted energy bounds on the hyperboloid. For a scalar quasilinear wave equation satisfying the weak null condition in three space dimensions, D.~Yu used the scattering data constructed in \cite{Y2024} to derive an asymptotic formula inside the light cone for small, smooth, compactly supported initial data \cite{Y2025}. A closer example is provided by one dimensional Schr\"odinger maps into $\mathbb S^2$, for which Li \cite{L2025} proved the sharp decay rate and modified scattering under a suitable choice of local coordinates. For the two dimensional maps considered here, \eqref{ethm1.1-uniform1} identifies the transverse waveform with an error of order $o(t^{-1})$ uniformly in space, and Corollary~\ref{cmap asymptotics} determines the leading term of the longitudinal magnetization deficit $1-Q\cdot\phi$ at order $t^{-2}$.
\end{remark}

\begin{remark}[Caloric gauge]
The proof of Theorem~\ref{thm:main} relies on the caloric gauge and the method of spacetime resonance.   The caloric gauge was introduced by Tao \cite{T2004} in wave maps and later extended to Schr\"odinger maps in \cite{T2006}. For a comparison with other frames, see Tataru \cite{T2012}. The construction extends the map by the harmonic map heat flow and, when this extension converges to $Q$, transports a limiting frame back parallel to the heat flow. The spatial connection is then recovered by integrating curvature in the heat parameter. This representation is essential in  the critical small data theory of \cite{BIKT2011}.   Analogous gauges based on the Yang-Mills heat flow were developed by Oh \cite{O} and Oh and Tataru  \cite{OT2019,OT2021} for hyperbolic Yang-Mills equations.
\end{remark}

 \begin{remark}[Spacetime resonance]
 The method of spacetime resonance was introduced by Germain, Masmoudi, and Shatah in  \cite{GMS2009,GMS2012a,GMS2012b}.  The method examines the oscillatory phase in the Fourier representation of the Duhamel formula. Away from time resonance, integration by parts in time is available, as in the normal form method of Shatah \cite{S1985}. Away from space resonance, integration by parts in frequency yields time decay at the cost of spatial weights on the profiles. Where both resonances occur, the geometry of the resonant set and the behavior of the nonlinear multiplier become important. For expositions, see Germain \cite{G2010}. The method has been widely applied to various models, see e.g.  \cite{DIP,GMS2012b,GMS2015,GIP,HPS,IP,PS}. 
\end{remark}

\subsection{Outline of the proof}

The main idea of the proof is to combine the method of spacetime resonance with the caloric gauge.  The key is to use Fourier bounds and spatial localization of the Schr\"odinger profile to improve the estimates for the differential fields and the caloric connection, and in turn, exploit the structure of the connection in the spacetime resonance analysis to control the profile and its spatial moments at low regularity.

We now describe the main ideas and steps of the proof in more detail.

\noindent\underline{\textit{1. The choice and structure of the caloric gauge.}} 

We formulate the equation for the differential fields in the caloric gauge. Let $\{A_x,A_t\}$ and $\{\psi_x,\psi_t\}$  be the associated induced connection and differential fields under the caloric gauge.  Its advantage over the Coulomb gauge can be seen from the following schematic quadratic expressions, without writing the component indices or conjugations explicitly, as in \cite{BIKT2011}:
\begin{align*}
A&\sim\sum_{j<k}2^{-k}P_j\psi P_k\psi+\sum_{j\leq k}2^{-j}P_j(P_k\psi P_k\psi),\qquad\text{Coulomb gauge},\\
A&\sim\sum_{j<k}2^{-k}P_j\psi P_k\psi+\sum_{j\leq k}2^{-k}P_j(P_k\psi P_k\psi),\qquad\text{caloric gauge}.
\end{align*}
In fact, for high-high interactions producing low frequencies, the negative derivative  in the caloric gauge is attached to the input frequency rather than the output frequency. This gives the additional factor $2^{j-k}$ when the output frequency $2^j$ is much smaller than the input frequency $2^k$.

In this work, we further exploit a balanced distribution of this negative derivative. Let $A^{(2)}$ denote the quadratic contribution obtained by replacing the differential fields in the caloric integral by their free heat evolutions, see \eqref{eAcom2}. Its bilinear symbol satisfies
\begin{equation*}
\frac{|\xi-\eta|+|\eta|}{|\xi-\eta|^2+|\eta|^2}\lesssim\min\left\{|\xi|^{-1},|\xi-\eta|^{-1},|\eta|^{-1}\right\},
\end{equation*}
and also
\begin{equation*}
\frac{|\xi-\eta|+|\eta|}{|\xi-\eta|^2+|\eta|^2}\lesssim\frac{1}{|\xi-\eta|^{\frac{1}{2}}|\eta|^{\frac{1}{2}}}.
\end{equation*}
The first estimate makes the improved high-high to low frequency factor explicit, since the multiplier is bounded by $2^{-k}$ when both inputs $\xi-\eta,\eta$  have frequency $2^k$, even if the output frequency $\xi$ is much smaller. The second estimate distributes half a negative derivative  to each input and motivates the propagation of a negative half order norm for $\psi$. The Fourier bound in the spacetime resonance framework, together with the energy bound, controls $\psi$ in $\dot H^{-\sigma}$ for every $0<\sigma<1$ in two dimensions, see Lemma~\ref{lqk}. The negative half order estimate is therefore compatible with our bootstrap and provides useful control of several norms of the connection. 

Let's compare the caloric gauge and the frame  $\{\frac{\partial}{\partial z},\frac{\partial}{\partial \bar{z}}\}$ associated with the stereographic coordinate. The map derivative $\partial_m\phi$ is represented by $\partial_m z$ in a fixed stereographic chart and by $\psi_m$ in the caloric frame. The first order terms in their evolution equations are, respectively,
$$
  \sum_{\ell=1}^2\frac{4\bar z\,\partial_\ell z}{1+|z|^2}\,\partial_\ell\partial_m z, \qquad\text{and} \qquad -2i\sum_{\ell=1}^2 A_\ell\,\partial_\ell\psi_m .
$$
To estimate the first and second spatial moments of the profiles of $\partial_m z$ and $\psi_m$, we commute their equations with $J_j=x_j+2it\partial_j$ and $J_iJ_j$. In stereographic coordinates, after taking $J_iJ_j$, the inhomogeneous term contains
$$
(2it)^2\sum_{\ell=1}^2\partial_i\partial_j\left(\frac{4\bar z\,\partial_\ell z}{1+|z|^2}\right)\partial_\ell\partial_m z.
$$
 A direct $L_x^2$ estimate of these terms by H\"older would use $\|z\|_{L_x^\infty}\|\nabla^2z\|_{L_x^2}$ and $\|z\|_{L_x^\infty}\|\nabla^3z\|_{L_x^2}$, neither of which is controlled by the assumed $H^{1+\epsilon}$ regularity.
However, 
 the caloric gauge suffices to close the bootstrap of spatial moments of profiles to $\psi$ by exploiting the structure of caloric connections and null structures, see Step 7 for details.   

These properties motivate our choice of the caloric gauge to reformulate the equation.

\noindent\underline{\textit{2. The dispersive bootstrap at low regularity.}} 

For the profile $f=e^{-it\Delta}\psi$, the linear estimate
\begin{equation*}
\|\psi(t)\|_{L_x^\infty}\lesssim t^{-1}\|\widehat f(t)\|_{L_\xi^\infty}+t^{-\frac{3}{2}+}\|\langle x\rangle^2f(t)\|_{L_x^2},\qquad t\geq1,
\end{equation*}
motivates the simultaneous control of the Fourier profile and its spatial moments. With $\epsilon_1=\epsilon_0^{\frac{1}{2}}$ and $0<2\gamma<\beta<\epsilon\ll1$, we bootstrap the $H^\epsilon$ bound, the $t^{-1}$ dispersive decay, the uniform Fourier bound, and the first and second moments with permitted growth $t^\gamma$ and $t^\beta$, respectively. See Proposition~\ref{pmain bootstrap} for the bootstrap setting. 

We address the roles of the  Fourier bounds and spatial moment bounds in the bootstrap setup. 
The Fourier bound gives an extra factor of $2^k$ at low frequencies compared with energy control. This gain compensates for the negative derivative  in the connection, see Remark~\ref{rmk:Alinfty} and Remark~\ref{rmk:Al2}, and gives the initial bound needed to propagate the negative half order envelope. As a bound on the initial data for the heat flow, it also yields an additional factor of $s^{-\frac{1}{2}}$ in the large $s$ decay of the free heat, see Lemma~\ref{lpsiAdecay}. Second spatial moment bound interpolated with the $H^\epsilon$ bound, by Proposition~\ref{pmain bootstrap}, gives  at high frequencies that 
\begin{equation*}
\|P_kf(t)\|_{L_x^1}\lesssim\epsilon_1t^{\frac{\beta}{2}}2^{-\frac{\epsilon}{2}k},\qquad k>0.
\end{equation*}
Furthermore, we have
\begin{equation}\label{op0}
\|P_k\psi(t)\|_{L_x^\infty}\lesssim\epsilon_1\min\left\{2^{2k},t^{-1},t^{-1+\frac{\beta}{2}}2^{-\frac{\epsilon}{2}k}\right\},\qquad k\in\mathbb Z.
\end{equation}
Note that the $2^{-\frac{\epsilon}{2}k} $ term yields 
\begin{equation*}
\sum_k\|P_k\psi(t)\|_{L_x^\infty}\lesssim\epsilon_1t^{-1}(1+\log t).
\end{equation*}
See Corollary~\ref{cbootstrap}.
Thus spatial localization provides the high frequency summability  with assuming only additional small $\epsilon$ derivative of $\nabla \phi_0$. Improving the Fourier and weighted bounds themselves requires the more precise resonance analysis described in the subsequent steps. 

Since $H^\epsilon$ does not control the initial $L_x^\infty$ norm of $\psi$, we first propagate the second spatial moment on $[0,1]$, see Proposition~\ref{locest}, and initialize the large time bootstrap at $t=1$. 

\noindent\underline{\textit{3. Propagation of the energy and negative half order envelopes along the heat direction.}} 

The differential fields $\psi$ satisfy the nonlinear heat equation \eqref{ecaloric4}, where the connection $A$  is given by \eqref{ecaloric6}. 
In Section~3, we establish in Proposition~\ref{propenv} that 
\begin{equation*}
\|P_k\psi(s,t)\|_{L_x^2}\lesssim\epsilon_1 2^{-\sigma k}b_k(\sigma)(1+s2^{2k})^{-20},\qquad 0\leq\sigma\leq\epsilon,
\end{equation*}
and in Proposition~\ref{prohaboo} that 
\begin{equation*}
\|P_k\psi(s,t)\|_{L_x^2}\lesssim\epsilon_1 2^{\frac{k}{2}}q_k(1+s2^{2k})^{-12}.
\end{equation*}
The envelopes $b_k(\sigma)$ and $q_k$ measure the energy and negative half order norms of the initial heat field.

The proof uses the heat Duhamel formula and paraproduct decomposition. Using the bootstrap assumptions of $\psi$, the caloric representation \eqref{ecaloric6} yields corresponding frequency envelope estimates for $A$.  One of the key observations in this step is  the following inequality.  For fixed $M>1$,
\begin{equation*}
\int_0^s\bigl(1+(s-r)2^{2k}\bigr)^{-M}\bigl(1+r2^{2j}\bigr)^{-M}\,\mathrm{d}r\lesssim2^{-2j}\bigl(1+s2^{2k}\bigr)^{-M},\qquad j\geq k-4.
\end{equation*}
This keeps the gain from the higher input frequency (i.e. $2^{-2j}$) while preserving the output heat weight $-M$.

The envelopes and their propagation estimates are formulated at each fixed physical time, preserving the time dependence needed in the subsequent dispersive argument. This heat flow analysis requires only the propagation of basic $L_x^2$ based energy and negative half order norms, without propagating the Strichartz, local smoothing, or maximal function norms used in \cite{BIKT2011}.

\noindent\underline{\textit{4. Estimates for the connection.}} 

Section~4 converts the heat flow bounds into Fourier, decay, and weighted estimates for the connection. 

The Fourier bound of $A$ is proved in Proposition~\ref{pAhat}: 
$$\sup_{s\geq0}\|\widehat A(s,t)\|_{L_\xi^\infty}\lesssim\epsilon_1^2.$$ To explain the main idea, we take $A^{(2)}$ as the candidate model. Indeed, the negative half order norm gives 
\begin{equation*}
\|\widehat{A^{(2)}}(0,t)\|_{L_\xi^\infty}\lesssim\bigl\||\nabla|^{-\frac{1}{2}}\psi(0,t)\bigr\|_{L_x^2}^2.
\end{equation*}
This follows from the balanced symbol bound and Cauchy-Schwarz in the convolution variable.

The principal Sobolev and $L^{\infty}$ bounds for the connection are given in Lemma~\ref{Ainfty}, Lemma~\ref{lAL2}, Proposition~\ref{Aderi2}, Proposition~\ref{proinfty}:
\begin{align*}
\|A(s,t)\|_{L_x^\infty}&\lesssim\epsilon_1^2(t+s)^{-\frac{3}{2}},\\
\|A(s,t)\|_{L_x^2}&\lesssim\epsilon_1^2(t+s)^{-1}(1+\log t),
\end{align*}
and 
\begin{align}\label{71}
\|A(s,t)\|_{\dot H^{\sigma+1}}&\lesssim\epsilon_1(t+s)^{-1}\|\psi(0,t)\|_{\dot H^\sigma}\lesssim\epsilon_1^2(t+s)^{-1},\qquad 0\leq\sigma\leq\epsilon\nonumber\\
\sup_{s\geq0}\|\nabla A(s,t)\|_{L_x^\infty}&\lesssim\epsilon_1^2t^{-2}(1+\log t).
\end{align}
In order to prove the above estimates, we need the $\|P_k \psi\|_{L^{\infty}_x}$ bounds. 
In fact, the paraproduct decomposition along the heat flow,  together with the high frequency estimate \eqref{op0} from spatial localization, gives
\begin{equation}\label{h00}
\|P_k\psi(s,t)\|_{L_x^\infty}\lesssim\epsilon_1\min\left\{2^{2k},(t+s)^{-1},t^{-1+\frac{\beta}{2}}2^{-\frac{\epsilon}{2}k}\right\}(1+s2^{2k})^{-10},
\end{equation}
see Proposition~\ref{psiinfty} and Proposition~\ref{highpsi}. The same argument establishes $$\|\psi(s,t)\|_{L_x^\infty}\lesssim\epsilon_1(t+s)^{-1}. $$
Note that \eqref{h00} captures both dispersive decay in $t$ and parabolic decay in $s$, where the $s^{-1}$ decay in \eqref{h00}  essentially relies on the Fourier bound of $\psi$ at $s=0$. 

To understand the estimates for $A$ listed above, take $A^{(2)}$ as the candidate, for which there are two main ingredients. First, the Fourier bound on the heat flow initial data gives an extra factor of $2^k$ at low frequencies compared with energy control, see Remark~\ref{rmk:Alinfty} and Remark~\ref{rmk:Al2}. Second, differentiating $A^{(2)}$ produces a zero order symbol, see Remark~\ref{conde}.

For the first spatial weight, we aim to prove
$$
\|xA(0,t)\|_{L_x^2}\lesssim\epsilon_1^2.
$$
See Proposition~\ref{pxA estimate}.
To prove this estimate, inspired by the Galilean vector field, we introduce a $(s,t)$ dependent vector field commuting with the linear heat operator, namely $J_s=x+2(s+it)\nabla$. A weighted frequency envelope $E_k$ defined in \eqref{eEk}, incorporating the bound for $J_0\psi(0,t)=e^{it\Delta}(xf(t))$ and the unweighted energy envelope, allows us to estimate the commuted heat equation. Using a bootstrap argument and the negative half order estimates, we obtain in Proposition~\ref{pJbootstrap} that 
\begin{equation*}
\|P_kJ_s\psi(s,t)\|_{L_x^2}\lesssim\epsilon_1t^\gamma E_k(1+s2^{2k})^{-10}.
\end{equation*}
The identity used to recover the spatial weight has the schematic form
\begin{equation}\label{xpsi}
x\nabla\psi\sim\nabla J_s\psi-2(s+it)\nabla^2\psi-\psi.
\end{equation}
Inside the caloric integral of $A$ in \eqref{ecaloric6}, the factor $s+it$ is compensated by $\|\psi(s,t)\|_{L_x^\infty}\lesssim\epsilon_1(t+s)^{-1}$, which  removes the apparent time growth factor in \eqref{xpsi}.

Recall the Schr\"odinger profile $f(t)=e^{-it\Delta}\psi$. 
The bootstrap of $\|\widehat{f}\|_{L^{\infty}_{\xi}}$ and $\|\langle x\rangle^2f\|_{L^2_{x}}$  requires two integrations by parts in frequency, which further reduces to dominating $\|x^2 A\|_{L^2_x}$. The $J_s$ commutator argument above   would involve a twice-commuted heat equation and lead to factors of $t^2$, which is beyond  control. We instead split the connection into quadratic part and the high order parts, i.e., $A=A^{(2)}+A^{(\geq4)}$. 
The quadratic part is estimated by the  refined spacetime resonance analysis, while the higher order part has stronger decay and does not require a second frequency integration by parts. Section 4 concludes by collecting the bounds for $A^{(\geq4)}$ and for $A_t$, whose quadratic part is a zero order bilinear expression in $\psi$ and its conjugate.

\noindent\underline{\textit{5. The energy estimate.}}

Uniform Sobolev bounds are already available from the small data theory of \cite{BIKT2011}. To keep the proof within the fixed time framework used here, Section~5 gives a direct energy argument, see Proposition~\ref{p5.1}. A paraproduct decomposition, integration by parts, and commutator estimates transfer the derivative in the magnetic term to the connection, while its one derivative gain \eqref{71} controls the remaining frequency interactions. The connection estimates in Step 4  give a time integrable coefficient in the energy inequality. Gronwall inequality, together with conservation of $\|\psi(t)\|_{L_x^2}$, then yields the uniform $H^\epsilon$ bound.

\noindent\underline{\textit{6. spacetime resonance: Fourier estimates for the profile.}}

In section 6 we prove Proposition~\ref{phatf improve}. Duhamel formula gives 
\begin{equation*}
\widehat f(t,\xi)=\widehat f(1,\xi)-i\int_1^t e^{i\tau|\xi|^2}\mathcal{F}\bigl[-2iA\cdot\nabla\psi+A_3\psi+F\psi\bigr](\tau,\xi)\,\mathrm{d}\tau.
\end{equation*}
The main difficulty of estimating $\|\widehat{f}\|_{L^{\infty}_{\xi}}$ is to dominate the magnetic term $iA\cdot\nabla\psi$, i.e.,
\begin{equation}\label{ddf}
\int_1^t\int_{\RR^2}e^{i\tau(|\xi|^2-|\xi-\eta|^2)}(\xi-\eta)\cdot\widehat A(\tau,\eta)\widehat f(\tau,\xi-\eta)\,\mathrm{d}\eta\,\mathrm{d}\tau.
\end{equation}
To estimate \eqref{ddf}, we split the outer frequency $|\xi-\eta|$ according to time size. For $A^{(\geq4)}$, we split at $\tau^{\frac{1}{2}}$ and integrate by parts at most once in $\eta$. For $A^{(2)}$, we split at $\tau^{-\frac{1}{3}-\delta}$ and use the identity $\nabla_\eta(|\xi|^2-|\xi-\eta|^2)=2(\xi-\eta)$, which is a null structure. In the low frequency scale $|\xi-\eta|\lesssim \tau^{-\frac{1}{3}-\delta}$, the Fourier bounds for $A^{(2)}$ and $f$ give $\epsilon_1^3\tau^{-1-3\delta}$, which is admissible. In the high frequency regime $|\xi-\eta|\gtrsim\tau^{-\frac{1}{3}-\delta}$, we integrate the quadratic contribution by parts twice in $\eta$. We explain the idea of bounding the most difficult term involved with $A^{(2)}$, i.e., \eqref{ehatf7}, which is the case when two derivatives hit $A^{(2)}$ in the two times of integration by parts in $\eta$. Inserting \eqref{eAcom2} into \eqref{ehatf7} gives cubic space time integrals of $\psi\overline{\psi}\psi$ type, 
\begin{align}\label{Hust}
\int_1^t\iint e^{i\tau \Phi_3}m(\xi,\eta,\zeta) \overline{\widehat{x^{a_1}f}(\tau,\zeta-\eta)}\widehat{x^{a_2}f}(\tau,\zeta)\widehat{x^{a_3}f}(\tau,\xi-\eta)\,\mathrm{d}\zeta\,\mathrm{d}\eta\,\mathrm{d}\tau,
\end{align}
where integers $a_1,a_2,a_3\ge 0$ and $ a_1+a_2+a_3\leq 2$, and $m$ denotes a multiplier of degree zero or negative degree, see \eqref{twde}, and phase
\begin{equation*}
\Phi_3(\xi,\eta,\zeta)=|\xi|^2-|\xi-\eta|^2+|\eta-\zeta|^2-|\zeta|^2.
\end{equation*}
The most complicated term among \eqref{Hust} is \eqref{complexterm}, with a degree zero multiplier $m_0$ satisfying
\begin{equation*}
|m_0(\xi,\eta,\zeta)|\lesssim\frac{|\eta-\zeta|^2|\zeta|}{|\xi-\eta|(|\eta-\zeta|^2+|\zeta|^2)}.
\end{equation*}
Near the space resonant set of $\Phi_3$, the small frequency support and the local integrability of $|\xi-\eta|^{-1}$ give
\begin{align}
&\biggl|\iint e^{i\tau\Phi_3}\chi_h\bigl(\tau^{\frac{1}{3}+\delta}|\xi-\eta|\bigr)\chi_l\bigl(\tau^{\frac{1}{2}-\delta}\sqrt{|\eta|^2+|\xi-\zeta|^2}\bigr)m_0(\xi,\eta,\zeta)\nonumber\\
&\qquad\cdot\overline{\widehat f(\zeta-\eta)}\widehat f(\zeta)\widehat f(\xi-\eta)\,\mathrm{d}\zeta\,\mathrm{d}\eta\biggr|\lesssim\epsilon_1^3\tau^{-2+4\delta},\label{yy9}
\end{align}
where $\chi_{l}(r)$ is a smooth even function  such that
$\operatorname{supp}\chi_l\subset[-2,2]$ and $\chi_{l}(r)=1 $ for $|r|\leq1$,  and $\chi_h(r):=1-\chi_l(r)$.

On the complementary region, i.e., $\chi_l$ in \eqref{yy9} is replaced by $\chi_{h}$, we integrate by parts twice using $\nabla_\eta\Phi_3=2(\xi-\zeta)$ and $\nabla_\zeta\Phi_3=-2\eta$. We localize
\begin{equation*}
|\xi-\eta|\simeq2^k,\qquad \sqrt{|\eta-\zeta|^2+|\zeta|^2}\simeq2^l,\qquad \sqrt{|\eta|^2+|\xi-\zeta|^2}\simeq2^j.
\end{equation*}
If $k\leq l-10$, then $j=l+O(1)$. If $l\leq k-10$, then $j=k+O(1)$. If $|k-l|\leq10$, then $j\leq k+O(1)$. The multiplier has size at most $2^{l-k}$.
Summing over the three frequency cases gives, for $\delta>0$ sufficiently small,
\begin{equation*}
\begin{aligned}
&\biggl|\iint e^{i\tau\Phi_3}\chi_h\bigl(\tau^{\frac{1}{3}+\delta}|\xi-\eta|\bigr)\chi_h\bigl(\tau^{\frac{1}{2}-\delta}\sqrt{|\eta|^2+|\xi-\zeta|^2}\bigr)m_0(\xi,\eta,\zeta)\\
&\qquad\cdot\overline{\widehat f(\zeta-\eta)}\widehat f(\zeta)\widehat f(\xi-\eta)\,\mathrm{d}\zeta\,\mathrm{d}\eta\biggr|\lesssim\epsilon_1^3\tau^{-1-2\delta}.
\end{aligned}
\end{equation*}
Combining the estimates from all regions and integrating in time yields the uniform $L_\xi^\infty$ bound for $\widehat f$.

\noindent\underline{\textit{7. Spacetime resonance: first and second spatial moments of the profile.} }

In Sections 7 and 8, we estimate the spatial moment $\|\langle x\rangle^2 f\|_{L^2}$, see Proposition~\ref{xfimp} and~\ref{xx2}.

For the first moment, by commuting with $J=x+2it\nabla$, $J_i \psi_m$ satisfies the equation \eqref{exf3}, where  the most tough inhomogeneous term is the  magnetic commutator  $\tau(\nabla A^{(2)})\nabla\psi$. Using $\nabla_\zeta\Phi_3=-2\eta$, which is a null structure, to integrate by parts in $\zeta$, we obtain the key connection estimate
\begin{equation*}
\|x\nabla A^{(2)}(\tau)\|_{L_x^2}+\tau\|\nabla^2A^{(2)}(\tau)\|_{L_x^2}\lesssim\epsilon_1^2\tau^{-\frac{1}{2}}.
\end{equation*}
Together with the decay of $\psi$, this controls the  magnetic commutator. See Proposition~\ref{xfimp}. 

For the second moment, $J_iJ_j\psi_m$ fulfills the equation \eqref{ex2f2}, where the principal inhomogeneous terms are $\tau(\nabla A^{(2)})\nabla e^{i\tau\Delta}(xf)$ and $\tau^2(\nabla^2A^{(2)})\nabla\psi$. For the former term, arguments for the first moment suffice  to give a time integrable bound. For the second term, we insert formula~\eqref{eAcom2} for $A^{(2)}$ and split the  frequencies according to the size  $\tau \sqrt{|\eta-\zeta|^2+|\zeta|^2}$. Below this scale, one integration by parts in $\eta$ and the small frequency support give an $L_\xi^2$ bound of $\epsilon_1^3\tau^{-3+\gamma}$. 
On the complementary frequency  region,  to remove $\tau^2\eta_i\eta_j$, we  can integrate by parts twice, since $\nabla_\xi\Phi_3=-\nabla_\zeta\Phi_3=2\eta$. 

Now, we dominate the remainder terms in \eqref{ex2f8} after twice integration by parts. The symbol  in \eqref{ex2f8} reads as  $\{(\xi-\eta)\nabla^j_{\zeta}m_{-1}\}^2_{j=0}$, and the high order derivatives of $m_{-1}$ are denoted by $m_{-2},m_{-3},m_{-4}$, see  \eqref{99}. Heuristically speaking,  the symbol $m_{-\alpha}$
 acts as $|\nabla|^{-\alpha}$, thus the Sobolev embeddings  and Bernstein inequalities hold for the corresponding symbols $m_{-\alpha}$. See Lemma \ref{luu1} and Lemma \ref{luu2}.

The two most difficult cases in  \eqref{ex2f8} are those in which both frequency derivatives fall on the same $\widehat{f}$  or both fall on $m_{-1}$.

When both derivatives fall on the same $\widehat{f}$, the more difficult region is $|\xi-\eta|\gg|\eta-\zeta|+|\zeta|$. We integrate by parts along $\nabla_\eta+\nabla_\zeta$ if $\widehat{x^2f}$ is evaluated at $\zeta-\eta$, and along $\nabla_\eta$ if it is evaluated at $\zeta$. In this way, no $\widehat{x^kf}$ with $k\geq 3$ occurs. Note also that the corresponding phase gradient has size comparable to $|\xi-\eta|$, which is a null structure. Thus
\begin{equation*}
\left\|\iint_{|\xi-\eta|\gg|\eta-\zeta|+|\zeta|}e^{i\tau\Phi_3}(\xi-\eta)m_{-1}(\eta,\zeta)\overline{\widehat{x^2f}(\zeta-\eta)}\widehat f(\zeta)\widehat f(\xi-\eta)\,\mathrm{d}\zeta\,\mathrm{d}\eta\right\|_{L_\xi^2}\lesssim\epsilon_1^3\tau^{-1+\beta}.
\end{equation*}

When both derivatives fall on $m_{-1}$, we use the phase identity $(\xi-\eta)e^{i\tau\Phi_3}=\frac{1}{2i\tau}\nabla_\eta e^{i\tau\Phi_3}-(\eta-\zeta)e^{i\tau\Phi_3}$. The factor $\eta-\zeta$ turns $m_{-3}$ into a multiplier of order $-2$, and  the first term gains $\tau^{-1}$ after integration by parts. Therefore,
\begin{equation*}
\left\|\iint e^{i\tau\Phi_3}(\xi-\eta)m_{-3}(\eta,\zeta)\overline{\widehat f(\zeta-\eta)}\widehat f(\zeta)\widehat f(\xi-\eta)\,\mathrm{d}\zeta\,\mathrm{d}\eta\right\|_{L_\xi^2}\lesssim\epsilon_1^3\tau^{-1+\gamma}(1+\log\tau).
\end{equation*}
After integration in time, these contributions lie within the allowed $t^\beta$ growth of the second moment and close the weighted estimates, see  Proposition \ref{xx2}.

\noindent\underline{\textit{8. Returning to the original map: stereographic projection in a fixed frame.}}

At the regularity $\nabla\phi\in H^\epsilon$, direct weighted estimates for the stereographic equation encounter a loss of derivatives. In Section~9, we pass to the stereographic coordinate $z$ in the fixed tangent frame $(v_\infty,w_\infty)$ only after establishing decay and weighted estimates for the caloric fields. The link between the two descriptions is the identity
\begin{equation*}
\partial_l z=\Theta \psi_l, \qquad\Theta:=\frac{(v_\infty+iw_\infty-zQ)\cdot v}{1+Q\cdot\phi}
\end{equation*}
The caloric connection controls the spatial derivative of $\Theta$, so this identity allows us to estimate $\nabla Jz$ using $J\psi$ without requiring an additional derivative of $\psi$. In the twice commuted equation for $z$, the energy weight $(1+|z|^2)^{-2}$ cancels the term that would otherwise require control of $\nabla J^2z$.

These commuted energy estimates give a uniform second moment for the profile of $z$. Interpolating the first and second moments of $f$ and using the same coordinate identity gives, for some $\sigma>1$ sufficiently close to one,
\begin{equation*}
\|\langle x\rangle^2e^{-it\Delta}z(t)\|_{L_x^2}+t^{-\frac{\beta}{2}}\|\langle x\rangle^\sigma e^{-it\Delta}\nabla z(t)\|_{L_x^2}\lesssim\epsilon_0,\qquad t\geq1.
\end{equation*}
The identity $|\nabla\phi|=|\psi|$ gives the optimal gradient decay in \eqref{ethm1.1-decay1}, while the second moment of the coordinate profile gives the optimal decay of $\phi-Q$ in \eqref{th1.1decay2}. Finally, the inverse stereographic projection yields the pointwise asymptotics in \eqref{ethm1.1-uniform1}.

\subsection{Organization of the paper}

Section~\ref{sec:2} introduces the caloric gauge formulation and states the main bootstrap proposition. In Sections~\ref{sec:heat envelope} and~\ref{sec:connection estimates}, we propagate bounds for the differential fields along the heat flow and use them to estimate the connection coefficients. In Sections~\ref{sec:energy estimate}-\ref{sec:decayestimate3}, we close the bootstrap by estimating the energy norm of the differential fields, the Fourier transform of the profile, and its first two spatial moments. In Section~\ref{sec:9}, we return to the original map and prove the main results on optimal decay and pointwise asymptotics. Appendix~\ref{app:A} treats the higher order terms in the connection, and Appendix~\ref{app:B} gives the local estimates needed to start the bootstrap.

\subsection{Notation}

We write $X\lesssim Y$ if $X\leq CY$ for some universal constant $C>0$, and $X\simeq Y$ if both $X\lesssim Y$ and $Y\lesssim X$. We write $a+$ for $a$ plus an arbitrarily small positive number. Throughout the paper, we fix a universal constant $0<\delta\ll1$, for example $\delta=10^{-10}$.

We sometimes omit the summation sign for repeated spatial indices, in which case summation over $\{1,2\}$ is understood. We denote the Galilean vector field by $J=x+2it\nabla$ and its extension along the heat flow by $J_s=x+2(s+it)\nabla$. At heat time $s=0$, we sometimes write $J_0$ for $J$.

\subsection{Acknowledgments} Z. Li was supported by the National Natural Science Foundation of China under grants (No.12422110 and No.12371244). X. Sun was supported by National Natural Science Foundation of China (No. 12271497).  L. Zhao was supported by National Natural Science Foundation of China (No. 12341102 and No. 12271497).

\section{Gauge formulation and bootstrap estimates}\label{sec:2}

This section is divided into three parts. In the first part, we introduce the differential fields and connection coefficients, rewrite the Schr\"odinger map equation in the caloric gauge, and present the local estimates needed to start the bootstrap at $t=1$. In the second part, we collect several analytical tools used throughout the proof. In the third part, we state the main bootstrap proposition, derive estimates of the bootstrap bounds, and state the scattering and uniform asymptotic formulas for the differential fields.

\subsection{differential fields and gauge formulations}

We set
$$
  \partial_3=\partial_t.
$$
Let $(v,w)$ be a positively oriented orthonormal frame of $T_\phi\mathbb{S}^2$, so that
\[
  \phi\times v=w, \qquad\phi\times w=-v.
\]
For $m=1,2,3$, define the differential fields and the connection coefficients by
\begin{equation}\label{efield1}
\psi_m=v\cdot\partial_m\phi+i\,w\cdot\partial_m\phi, \qquad A_m=w\cdot\partial_m v, \qquad D_m=\partial_m+iA_m.
\end{equation}
The derivatives of $\phi$ and of the frame can then be recovered from
\begin{equation}\label{efield2}
\left\{
  \begin{aligned}
  \partial_m\phi&=v\operatorname{Re}\psi_m+w\operatorname{Im}\psi_m,\\
  \partial_m v&=-\phi\operatorname{Re}\psi_m+wA_m,\\
  \partial_m w&=-\phi\operatorname{Im}\psi_m-vA_m.
\end{aligned}
\right.
\end{equation}
Commuting the derivatives in \eqref{efield2} gives the compatibility and curvature identities
\begin{equation}\label{efield3}
  D_l\psi_m=D_m\psi_l, \qquad\partial_lA_m-\partial_mA_l=\operatorname{Im}\bigl(\psi_l\overline{\psi_m}\bigr), \qquad l,m\in\{1,2,3\}.
\end{equation}
In particular,
\begin{equation*}
  [D_l,D_m]=i\operatorname{Im}\bigl(\psi_l\overline{\psi_m}\bigr).
\end{equation*}

Taking the components of \eqref{eschrodinger-map} in the frame $(v,w)$ and using $\phi\times v=w$ and $\phi\times w=-v$, we obtain
\begin{equation}\label{efield4}
  \psi_3=i\sum_{l=1}^2D_l\psi_l.
\end{equation}
Hence, for $m=1,2$,
\begin{equation*}
  D_3\psi_m =D_m\psi_3 =i\sum_{l=1}^2D_mD_l\psi_l =i\sum_{l=1}^2D_lD_l\psi_m+\sum_{l=1}^2\operatorname{Im}\bigl(\psi_l\overline{\psi_m}\bigr)\psi_l.
\end{equation*}

\subsubsection{The Coulomb gauge}

The Coulomb gauge is determined by the condition
\begin{equation*}
  \sum_{l=1}^2\partial_lA_l=0.
\end{equation*}
After fixing the residual gauge at spatial infinity, taking a spatial divergence in \eqref{efield3} gives
\begin{equation}\label{coumb}
  A_m=\Delta^{-1}\sum_{l=1}^2\partial_l\operatorname{Im}\bigl(\psi_l\overline{\psi_m}\bigr), \qquad m=1,2,3.
\end{equation}
Thus, in the Coulomb gauge, the differential fields satisfy
\begin{equation*}
  (i\partial_t+\Delta)\psi_m=-2i\sum_{l=1}^2A_l\partial_l\psi_m+\left(A_3+\sum_{l=1}^2A_l^2\right)\psi_m-i\sum_{l=1}^2\psi_l\operatorname{Im}\bigl(\overline{\psi_l}\psi_m\bigr), \qquad m=1,2,
\end{equation*}
where $\psi_3$ and $A_m$ are determined by \eqref{efield4} and \eqref{coumb}, respectively.

As discussed in the Introduction, when two high frequencies interact to produce a low frequency, the negative derivative  in the Coulomb connection acts on the low output frequency. In two dimensions, this does not provide the low frequency control needed for our decay estimates. We therefore work in the caloric gauge and study the Schr\"odinger equation \eqref{ecaloric7} for the differential fields, with the connection coefficients given by \eqref{ecaloric6}.

\subsubsection{The caloric gauge}

For every fixed physical time $t$, we extend $\phi(t,x)$ in the parabolic variable $s\geq0$ by solving the harmonic map heat flow
\begin{equation}\label{ecaloric1}
\left\{
\begin{aligned}
  &\partial_s\tilde\phi=\Delta\tilde\phi+\tilde\phi\sum_{l=1}^2|\partial_l\tilde\phi|^2, &&s>0,\\
  &\tilde\phi(0,t,x)=\phi(t,x).
\end{aligned}
\right.
\end{equation}
The frame $(v,w)$ is extended along the heat flow by parallel transport from $s=\infty$. Setting $\partial_0=\partial_s$, this is equivalent to the caloric gauge condition
\begin{equation}\label{ecaloric2}
  A_0=w\cdot\partial_s v=0.
\end{equation}
For $m=0,1,2,3$, the fields $\psi_m$, the coefficients $A_m$, and the covariant derivatives $D_m$ are defined as in \eqref{efield1}. The identities \eqref{efield3} remain valid for $l,m\in\{0,1,2,3\}$.

Taking the frame components of \eqref{ecaloric1} gives
\begin{equation}\label{ecaloric3}
  \psi_0=\sum_{l=1}^2D_l\psi_l.
\end{equation}
Using \eqref{efield3} and \eqref{ecaloric2}, for $m=1,2,3$ we obtain
\begin{equation*}
  \partial_s\psi_m=D_m\psi_0 =\sum_{l=1}^2D_mD_l\psi_l=\sum_{l=1}^2D_lD_l\psi_m+i\sum_{l=1}^2\operatorname{Im}\bigl(\psi_m\overline{\psi_l}\bigr)\psi_l.
\end{equation*}
Equivalently,
\begin{equation}
  (\partial_s-\Delta)\psi_m=2i\sum_{l=1}^2\partial_l(A_l\psi_m)-\sum_{l=1}^2\bigl(A_l^2+i\partial_lA_l\bigr)\psi_m+i\sum_{l=1}^2\operatorname{Im}\bigl(\psi_m\overline{\psi_l}\bigr)\psi_l, \qquad m=1,2,3. \label{ecaloric4}
\end{equation}

The curvature identity with $l=0$, together with $A_0=0$, implies
\begin{equation}\label{ecaloric5}
  \partial_sA_m=\operatorname{Im}\bigl(\psi_0\overline{\psi_m}\bigr), \qquad m=1,2,3.
\end{equation}
Since the connection coefficients vanish as $s\to\infty$, integration of \eqref{ecaloric5} gives
\begin{equation}
  A_m(s)=-\int_s^\infty\operatorname{Im}\bigl(\psi_0(r)\overline{\psi_m(r)}\bigr)\,\mathrm{d}r=-\sum_{l=1}^2\int_s^\infty\operatorname{Im}\bigl(D_l\psi_l(r)\overline{\psi_m(r)} \bigr)\,\mathrm{d}r, \qquad m=1,2,3. \label{ecaloric6}
\end{equation}
At $s=0$, the physical differential fields therefore satisfy
\begin{equation}(i\partial_t+\Delta)\psi_m=-2i\sum_{l=1}^2A_l\partial_l\psi_m+\left(A_3+\sum_{l=1}^2\bigl(A_l^2-i\partial_lA_l\bigr)\right)\psi_m-i\sum_{l=1}^2\psi_l\operatorname{Im}
  \bigl(\overline{\psi_l}\psi_m\bigr), \quad m=1,2, \label{ecaloric7}
\end{equation}
where the connection coefficients are given by \eqref{ecaloric6}.

\begin{proposition}[Construction of the caloric gauge]
\label{pcaloric gauge}
Let $\phi$ solve \eqref{eschrodinger-map} on a time interval $I$ containing $0$, with $\phi-Q\in C(I;H^{1+\epsilon})$ and initial data satisfying \eqref{ini}. For every fixed $t\in I$, the heat flow \eqref{ecaloric1} exists uniquely for all $s\geq0$. It admits a unique caloric frame $(v,w)$ satisfying \eqref{ecaloric2} and
\[
  \lim_{s\to\infty}\bigl(\|\tilde\phi(s,t)-Q\|_{L_x^\infty}+\|v(s,t)-v_\infty\|_{L_x^\infty}+\|w(s,t)-w_\infty\|_{L_x^\infty}\bigr)=0.
\]
Fix $0<\delta\ll1$. At each fixed physical time $t$, define
\begin{equation*}
  \gamma_k(\sigma):=\sup_{j\in\ZZ}2^{-\delta|k-j|}2^{\sigma j}\|P_j(\phi(t)-Q)\|_{L_x^2},\qquad 0\leq\sigma\leq1+\epsilon.
\end{equation*}
Then, uniformly in $t\in I$, $s\geq0$, and $k\in\ZZ$,
\begin{equation}\label{ecaloric8}
\begin{aligned}
&\|P_k(\tilde\phi(s,t)-Q)\|_{L_x^2}+\|P_k(v(s,t)-v_\infty)\|_{L_x^2}+\|P_k(w(s,t)-w_\infty)\|_{L_x^2}\\
&\qquad\lesssim2^{-\sigma k}\gamma_k(\sigma)(1+s2^{2k})^{-20},\qquad 0\leq\sigma\leq1+\epsilon.
\end{aligned}
\end{equation}
\end{proposition}

\begin{proof}
By \eqref{enerconser} and \eqref{ini},
$$
\|\nabla\phi(t)\|_{L_x^2}=\|\nabla\phi_0\|_{L_x^2}\leq\epsilon_0, \qquad t\in I.
$$
We apply the argument of \cite[Proposition~4.2]{BIKT2011} to $\phi(t)$ at each fixed physical time, using $L_x^2$ in place of $L_t^\infty L_x^2$. This gives global existence of the heat flow, the normalized caloric frame, and \eqref{ecaloric8} for the map when $0\leq\sigma\leq 1+\epsilon$ and for the frame when $1\leq\sigma\leq 1+\epsilon$. For the frame, a more precise frequency decomposition, keeping the high-high interactions inside the heat time integral, extends the estimate to $0\leq\sigma<1$. We omit the details.
\end{proof}

Now we turn to initial estimates for the differential fields. Write $\psi=(\psi_1,\psi_2)$ and, at $s=0$, abbreviate $\psi(t)=\psi(0,t)$.
\begin{proposition}\label{pinitial differentiate fields}
Under \eqref{ini}, the differential fields in the caloric gauge satisfy
\begin{equation}\label{efield-init}
\|\psi(0)\|_{H^\epsilon}+\|\langle x\rangle^2\psi(0)\|_{L_x^2}+\|\widehat\psi(0)\|_{L_\xi^\infty}\lesssim\epsilon_0.
\end{equation}
\end{proposition}
\begin{proof}
For every physical time $t$, the vectors $\partial_m\phi(t)$ are tangent to $\mathbb S^2$, and $(v,w)$ is an orthonormal basis of the tangent plane. Thus \eqref{efield1} gives
\begin{equation}\label{efield-init2}
|\psi(t)|^2=\sum_{m=1}^2\bigl(|v\cdot\partial_m\phi(t)|^2+|w\cdot\partial_m\phi(t)|^2\bigr)=|\nabla\phi(t)|^2.
\end{equation}
In particular,
\begin{equation*}
\|\widehat\psi(0)\|_{L_\xi^\infty}\lesssim\|\psi(0)\|_{L_x^1}\lesssim\|\langle x\rangle^2\psi(0)\|_{L_x^2}=\|\langle x\rangle^2\nabla\phi_0\|_{L_x^2}\lesssim\epsilon_0.
\end{equation*}
By \eqref{ecaloric8} with $s=t=0$ and $\sigma=1$,
\begin{equation*}
\|\nabla v(0,0)\|_{L_x^2}+\|\nabla w(0,0)\|_{L_x^2}\lesssim\|\nabla\phi_0\|_{L_x^2}\lesssim\epsilon_0.
\end{equation*}
Since $0<\epsilon<1$, the fractional Sobolev inequalities give
\begin{align*}
\bigl\||\nabla|^\epsilon(v(0,0)-v_\infty)\bigr\|_{L_x^{\frac{2}{\epsilon}}}+\bigl\||\nabla|^\epsilon(w(0,0)-w_\infty)\bigr\|_{L_x^{\frac{2}{\epsilon}}}&\lesssim\|\nabla v(0,0)\|_{L_x^2}+\|\nabla w(0,0)\|_{L_x^2},\\
\|\nabla\phi_0\|_{L_x^{\frac{2}{1-\epsilon}}}&\lesssim\|\nabla\phi_0\|_{H^\epsilon}.
\end{align*}
Using $|v|=|w|=1$ and the fractional Leibniz inequality, we obtain
\begin{align*}
\|\psi(0)\|_{H^\epsilon}&\lesssim\|\nabla\phi_0\|_{H^\epsilon}+\|\nabla\phi_0\|_{L_x^{\frac{2}{1-\epsilon}}}\left(\bigl\||\nabla|^\epsilon(v(0,0)-v_\infty)\bigr\|_{L_x^{\frac{2}{\epsilon}}}+\bigl\||\nabla|^\epsilon(w(0,0)-w_\infty)\bigr\|_{L_x^{\frac{2}{\epsilon}}}\right)\\
&\lesssim(1+\epsilon_0)\|\nabla\phi_0\|_{H^\epsilon}\lesssim\epsilon_0.
\end{align*}
This proves \eqref{efield-init}.
\end{proof}

Define the Schr\"odinger profile by
\begin{equation}\label{edef-profile}
f(t)=e^{-it\Delta}\psi(t).
\end{equation}
The following local estimate initializes the large time bootstrap.
\begin{proposition}[Local estimates]\label{locest}
Under \eqref{ini}, the differential fields satisfy
\begin{equation}\label{elocalestimate}
\begin{aligned}
&\sup_{0\leq t\leq1}\left(\|\psi(t)\|_{H^\epsilon}+\|xf(t)\|_{L_x^2}+\|x^2f(t)\|_{L_x^2}+\|\widehat f(t)\|_{L_\xi^\infty}\right)\\
&\qquad+\left(\int_0^1\left(\sum_k\|P_k\psi(t)\|_{L_x^\infty}\right)^2\,\mathrm dt\right)^{\frac{1}{2}}\lesssim\epsilon_0.
\end{aligned}
\end{equation}
In particular,
\begin{equation}\label{eboot-init}
\|\psi(1)\|_{H^\epsilon}+\|\psi(1)\|_{L_x^\infty}+\|\widehat f(1)\|_{L_\xi^\infty}+\|xf(1)\|_{L_x^2}+\|x^2f(1)\|_{L_x^2}\lesssim\epsilon_0.
\end{equation}
\end{proposition}
\begin{proof}
The proof is given in Appendix~\ref{app:B}.
\end{proof}

\subsection{Some analysis tools}

For $f\in\mathcal S(\RR^2)$, we use the Fourier transform convention
\[
  \widehat f(\xi)=\mathcal Ff(\xi):=\frac{1}{2\pi}\int_{\RR^2}e^{-ix\cdot\xi}f(x)\,\mathrm{d}x, \qquad\mathcal F^{-1}f(x):=\frac{1}{2\pi}\int_{\RR^2}e^{ix\cdot\xi}f(\xi)\,\mathrm{d} \xi.
\]
Fix an even function $\chi\in C_c^\infty(\RR;[0,1])$ such that
$$
  \operatorname{supp}\chi\subset[-2,2], \qquad\chi(r)=1\quad\text{for }|r|\leq1,
$$
and set
$$
  \chi_l(r):=\chi(r), \qquad\chi_h(r):=1-\chi(r).
$$
For each $k\in\mathbb Z$, define
$$
  \varphi_k(\xi):=\chi(2^{-k}|\xi|)-\chi(2^{-k+1}|\xi|).
$$
Then $\varphi_k$ is supported in $\{2^{k-1}\leq|\xi|\leq2^{k+1}\}$, and
$$
  \sum_{k\in\mathbb Z}\varphi_k(\xi)=1, \qquad\sum_{j\leq k}\varphi_j(\xi)=\chi(2^{-k}|\xi|), \qquad \xi\neq0.
$$
We define the Littlewood-Paley projections by
$$
  \widehat{P_kf}(\xi):=\varphi_k(\xi)\widehat f(\xi), \qquad\widehat{P_{\leq k}f}(\xi):=\chi(2^{-k}|\xi|)\widehat f(\xi),
$$
and write
$$
  P_{<k}:=P_{\leq k-1}, \qquad P_{>k}:=I-P_{\leq k}, \qquad P_{\geq k}:=I-P_{<k}.
$$

\begin{lemma}[Linear dispersive estimate, \cite{HN}]
\label{llinear dispersive estimate}
Let $0\leq\alpha_1\leq1$ and $\alpha_2>1+2\alpha_1$. For every $f:\mathbb R^2\to\mathbb C$ with $\langle x\rangle^{\alpha_2}f\in L^2(\mathbb R^2)$ and every $t\neq0$, one has
\[
  \|e^{it\Delta}f\|_{L_x^\infty(\mathbb R^2)}\lesssim |t|^{-1}\|\widehat f\|_{L_\xi^\infty(\mathbb R^2)}+|t|^{-1-\alpha_1}\|\langle x\rangle^{\alpha_2}f\|_{L_x^2(\mathbb R^2)},
\]
where the implicit constant depends only on $\alpha_1,\alpha_2$.
\end{lemma}

\begin{proof}
The Schr\"odinger kernel formula and $|e^{\frac{i|x|^2}{4t}}-1|\lesssim |t|^{-\alpha_1}|x|^{2\alpha_1}$ yield
\begin{align*}
  \|e^{it\Delta}f\|_{L_x^\infty}&\lesssim |t|^{-1}\|\widehat f\|_{L_\xi^\infty}+|t|^{-1-\alpha_1}\|\langle x\rangle^{2\alpha_1}f\|_{L_x^1}\\
  &\lesssim|t|^{-1}\|\widehat f\|_{L_\xi^\infty}+|t|^{-1-\alpha_1}\|\langle x\rangle^{2\alpha_1-\alpha_2}\|_{L_x^2}\|\langle x\rangle^{\alpha_2}f\|_{L_x^2}.
\end{align*}
Since $\alpha_2-2\alpha_1>1$, the estimate admits.
\end{proof}

\begin{lemma}[Decay estimate for the heat operator, \cite{MYZ}]\label{ldecay of heat operator}
  For any Schwartz function $f\in\mathcal{S}(\mathbb{R}^2)$, $s>0$, $l\geq 0$, there holds
  \begin{equation}\label{eheatoperator}
    \|e^{s\Delta}|\nabla|^lf\|_{L^q}\lesssim s^{-(\frac{1}{p}-\frac{1}{q})-\frac{l}{2}}\|f\|_{L^p},\qquad 1\leq p\leq q\leq\infty.
  \end{equation}
\end{lemma}

\begin{proof}
By Young inequality and scaling, we have
\begin{equation*}
  \|e^{s\Delta}|\nabla|^l f\|_{L^q}\lesssim\|\mathcal F^{-1}(e^{-s|\xi|^2}|\xi|^l)\|_{L^r}\|f\|_{L^p}\lesssim s^{-(\frac1p-\frac1q)-\frac l2}\|f\|_{L^p},
\end{equation*}
where $1+\frac{1}{q}=\frac{1}{r}+\frac{1}{p}$.
\end{proof}

\noindent\textit{A polynomial heat convolution lemma}
The following integral is frequently employed in heat estimates. In particular, it reveals that sufficient frequency control factors persist at high-high interaction. Furthermore, this estimate involves no loss of polynomial weight, which motivates our use of polynomial weights instead of exponential weights.
\begin{lemma}\label{lpolynomialconvolution}
Let $\bar k\geq k-c$, where $c$ is a universal constant, and let $M>1$. Then
\begin{equation}\label{epolynomialconvolution}
  \int_0^s \bigl(1+(s-r)2^{2k}\bigr)^{-M}\bigl(1+r2^{2\bar k}\bigr)^{-M}\,\mathrm{d}r\lesssim 2^{-2\bar k}\bigl(1+s2^{2k}\bigr)^{-M}.
\end{equation}
\end{lemma}

\begin{proof}
The harmless difference between $\bar k\geq k-c$ and $\bar k\geq k$ can be absorbed into the implicit constant, so it suffices to consider $\bar k\geq k$.  

If $s2^{2\bar k}\leq1$, then $s\lesssim2^{-2\bar k}$ and $s2^{2k}\lesssim 1$. Hence
\[
  \int_0^s\bigl(1+(s-r)2^{2k}\bigr)^{-M}\bigl(1+r2^{2\bar k}\bigr)^{-M}\,\mathrm{d}r\leq s\lesssim2^{-2\bar k},
\]
which proves \eqref{epolynomialconvolution} in this case.

Suppose now that $s2^{2\bar k}>1$. On $0\leq r\leq \frac{s}{2}$, one has $s-r\geq \frac{s}{2}$, and therefore
\begin{align}\label{epolynomialconvolution1}
  \int_0^{\frac{s}{2}}\bigl(1+(s-r)2^{2k}\bigr)^{-M}\bigl(1+r2^{2\bar k}\bigr)^{-M}\,\mathrm{d}r&\lesssim (1+s2^{2k})^{-M}\int_0^\infty(1+r2^{2\bar k})^{-M}\,\mathrm{d}r\notag\\
  &\lesssim 2^{-2\bar k}(1+s2^{2k})^{-M}.
\end{align}
On $\frac{s}{2}\leq r\leq s$, we have $1+r2^{2\bar k}\gtrsim1+s2^{2\bar k}$. Hence,
\begin{equation*}
  \int_{\frac{s}{2}}^s\bigl(1+(s-r)2^{2k}\bigr)^{-M}\bigl(1+r2^{2\bar k}\bigr)^{-M}\mathrm{d}r\lesssim (1+s2^{2\bar k})^{-M}\min\{s,2^{-2k}\}.
\end{equation*}
Since $\bar k\geq k$, a direct distinction between $s2^{2k}\leq1$ and $s2^{2k}>1$ gives
\[
  (1+s2^{2\bar k})^{-M}\min\{s,2^{-2k}\}\lesssim 2^{-2\bar k}(1+s2^{2k})^{-M}.
\]
Together with \eqref{epolynomialconvolution1}, this proves the lemma.
\end{proof}

We will also repeatedly use the elementary tail estimate
\begin{equation}\label{epolynomialconvolution2}
  \int_s^\infty (1+r2^{2k})^{-M}\,\mathrm{d}r\lesssim 2^{-2k}(1+s2^{2k})^{-M+1},\qquad M>1.
\end{equation}

Next we recall the classical Coifman-Meyer estimates.
\begin{lemma}[Coifman-Meyer estimates, \cite{CM,NT}]\label{lcoifman-meyer estimates}
Let $m$ be smooth away from $(\xi,\eta)=(0,0)$ and satisfy
\[
  |\nabla_\xi^\mu\nabla_\eta^\nu m(\xi,\eta)|\lesssim_{\mu,\nu}(|\xi-\eta|+|\eta|)^{-|\mu|-|\nu|}
\]
for all multi-indices $\mu,\nu$. If $1<p,q\leq\infty$, $1<r<\infty$, and $\frac1r=\frac1p+\frac1q$, then
\[
  \left\|\mathcal F_\xi^{-1}\left[\int_{\RR^2}m(\xi,\eta)\widehat f(\xi-\eta)\widehat g(\eta)\,\mathrm{d}\eta\right]\right\|_{L_x^r}\lesssim\|f\|_{L_x^p}\|g\|_{L_x^q}.
\]
Moreover, for every $\sigma\geq0$,
\[
  \left\|\mathcal F_\xi^{-1}\left[\int_{\RR^2}m(\xi,\eta)\widehat f(\xi-\eta)\widehat g(\eta)\,\mathrm{d}\eta\right]\right\|_{H^\sigma}\lesssim\|f\|_{H^\sigma}\|g\|_{L_x^\infty} +\|f\|_{L_x^\infty}\|g\|_{H^\sigma}.
\]
The implicit constants depend on finitely many symbol seminorms.
\end{lemma}

\begin{proof}
The first estimate is the bilinear Coifman-Meyer theorem. The second is its fractional Leibniz estimate, obtained by paraproduct decomposition and Littlewood-Paley square function estimates.
\end{proof}

We also state multilinear estimates for multipliers with integrable kernels. For a bilinear multiplier on \(\RR^2\times\RR^2\) and a trilinear multiplier on \(\RR^2\times\RR^2\times\RR^2\), respectively, we set
\[
  \|m(\xi,\eta)\|_{\mathcal S^\infty}:=\|\mathcal F^{-1}_{\xi,\eta}m\|_{L^1(\RR^4)}, \qquad\|m(\xi,\eta,\zeta)\|_{\mathcal S^\infty}:=\|\mathcal F^{-1}_{\xi,\eta,\zeta}m \|_{L^1(\RR^6)}.
\]

\begin{lemma}[Multilinear estimates, \cite{IP}]\label{l multilinear estimates}
Let \(1\leq p,q,r\leq\infty\) and \(\frac1r=\frac1p+\frac1q\), then
\[
  \left\|\mathcal F^{-1}_\xi\left[\int_{\RR^2}m(\xi,\eta)\widehat f(\xi-\eta)\widehat g(\eta)\,\mathrm{d}\eta\right]\right\|_{L_x^r}\lesssim\|m(\xi,\eta)\|_{\mathcal S^\infty}
  \|f\|_{L_x^p}\|g\|_{L_x^q}.
\]
If \(1\leq p_1,p_2,p_3,r\leq\infty\) and \(\frac{1}{r}=\frac{1}{p_1}+\frac{1}{p_2}+\frac{1}{p_3}\), then
\[
\begin{aligned}
  &\left\|\mathcal F^{-1}_\xi\left[\iint m(\xi,\eta,\zeta)\widehat f(\xi-\eta)\widehat g(\eta-\zeta)\widehat h(\zeta)\,\mathrm{d}\zeta\,\mathrm{d}\eta\right]\right\|_{L_x^r}\\
  &\lesssim\|m(\xi,\eta,\zeta)\|_{\mathcal S^\infty}\|f\|_{L_x^{p_1}}\|g\|_{L_x^{p_2}}\|h\|_{L_x^{p_3}}.
\end{aligned}
\]
The same estimates hold with different functions in the two or three input positions.
\end{lemma}

\begin{proof}
    The proof is standard, which follows by writing the multiplier as its inverse Fourier transform and then applying translation invariance, Minkowski, and H\"older.
\end{proof}

We shall also use the following scale invariant criterion.
\begin{lemma}[Localized multiplier estimate, \cite{DIP}]\label{lmet}
For a smooth bilinear symbol \(m(\xi,\eta)\), define
\[
    m_{k_1,k_2}(\xi,\eta):=m(\xi,\eta)\varphi_{k_1}(\xi)\varphi_{k_2}(\eta).
\]
Then
\[
\begin{aligned}
    \|m_{k_1,k_2}\|_{\mathcal S^\infty}\lesssim\sum_{\substack{|\mu|\leq2\\|\nu|\leq2}}2^{k_1|\mu|}2^{k_2|\nu|}\left\|\nabla_\xi^\mu\nabla_\eta^\nu m_{k_1,k_2}(\xi,\eta)\right\|_{L^\infty_{\xi,\eta}}.
\end{aligned}
\]
Similarly, for a smooth trilinear symbol \(m(\xi,\eta,\zeta)\), define
\[
    m_{k_1,k_2,k_3}(\xi,\eta,\zeta):=m(\xi,\eta,\zeta)\varphi_{k_1}(\xi)\varphi_{k_2}(\eta)\varphi_{k_3}(\zeta).
\]
Then
\[
\begin{aligned}
    \|m_{k_1,k_2,k_3}\|_{\mathcal S^\infty}\lesssim\sum_{\substack{|\mu|\leq2,\ |\nu|\leq2\\|\kappa|\leq2}}2^{k_1|\mu|}2^{k_2|\nu|}2^{k_3|\kappa|}\left\|\nabla_\xi^\mu\nabla_\eta^\nu\nabla_\zeta^\kappa m_{k_1,k_2,k_3}(\xi,\eta,\zeta)\right\|_{L^\infty_{\xi,\eta,\zeta}}.
\end{aligned}
\]
\end{lemma}

\begin{remark}
The same estimates hold after any fixed nonsingular linear change of the frequency variables. In particular, when the multiplier is written in output input variables, all the relevant Littlewood-Paley factors, such as
\[
    \varphi_k(\xi)\varphi_{k_1}(\xi-\eta)\varphi_{k_2}(\eta)
\]
in the bilinear case and
\[
    \varphi_k(\xi)\varphi_{k_1}(\xi-\eta)\varphi_{k_2}(\eta-\zeta)\varphi_{k_3}(\zeta)
\]
in the trilinear case, are included in the localized symbol before the derivative estimates above are applied.
\end{remark}

\begin{proof}
We only prove the bilinear estimate, since the trilinear estimate is identical. Let
\[
  K(x,y):=\mathcal F_{\xi,\eta}^{-1}[m_{k_1,k_2}](x,y).
\]
By Cauchy-Schwarz and Plancherel equality,
\[
\begin{aligned}
  \|K\|_{L_{x,y}^1}&\lesssim 2^{-k_1-k_2}\left\|\langle2^{k_1}x\rangle^2\langle2^{k_2}y\rangle^2K\right\|_{L_{x,y}^2}\\
  &\lesssim 2^{-k_1-k_2}\sum_{\substack{|\mu|\leq2\\|\nu|\leq2}}2^{k_1|\mu|}2^{k_2|\nu|}\left\|\nabla_\xi^\mu\nabla_\eta^\nu m_{k_1,k_2}\right\|_{L_{\xi,\eta}^2}\\
  &\lesssim\sum_{\substack{|\mu|\leq2\\|\nu|\leq2}}2^{k_1|\mu|}2^{k_2|\nu|}\left\|\nabla_\xi^\mu\nabla_\eta^\nu m_{k_1,k_2}\right\|_{L_{\xi,\eta}^\infty},
\end{aligned}
\]
this proves the desired estimate.
\end{proof}

\subsection{The main bootstrap proposition}
We use the profile defined in \eqref{edef-profile} and fix
\begin{equation}\label{bop}
0<2\gamma<\beta<\epsilon\ll\frac{1}{2}.
\end{equation}
By \eqref{eboot-init}, the following bootstrap assumption is initialized at $t=1$ when $\epsilon_0$ is sufficiently small.
\begin{proposition}[Bootstrap improvement]\label{pmain bootstrap}
Set $\epsilon_1=\epsilon_0^{\frac{1}{2}}$. Suppose that, for every $t\in[1,T]$, the solution satisfies
\begin{equation}\label{booassum}
\|\psi(t)\|_{H^\epsilon}+t\|\psi(t)\|_{L_x^\infty}+\|\widehat f(t)\|_{L_\xi^\infty}+t^{-\gamma}\|xf(t)\|_{L_x^2}+t^{-\beta}\|x^2f(t)\|_{L_x^2}\leq2\epsilon_1.
\end{equation}
Then, for every $t\in[1,T]$,
\begin{equation}\label{etimebootstrapimprovement}
\|\psi(t)\|_{H^\epsilon}+t\|\psi(t)\|_{L_x^\infty}+\|\widehat f(t)\|_{L_\xi^\infty}+t^{-\gamma}\|xf(t)\|_{L_x^2}+t^{-\beta}\|x^2f(t)\|_{L_x^2}\leq\epsilon_1.
\end{equation}
\end{proposition}
\begin{proof}
The proof is divided into Sections~\ref{sec:energy estimate}, \ref{sec:decayestimate1}, \ref{sec:decayestimate2}, and~\ref{sec:decayestimate3}, where we improve the bounds for $\|\psi(t)\|_{H^\epsilon}$, $\|\widehat f(t)\|_{L_\xi^\infty}$, $\|xf(t)\|_{L_x^2}$, and $\|x^2f(t)\|_{L_x^2}$, respectively. The $L_x^\infty$ decay estimate follows from Lemma~\ref{llinear dispersive estimate}, with $\alpha_2=2$ and $\beta<\alpha_1<\frac12$. Combining these estimates and taking $\epsilon_0$ sufficiently small yields \eqref{etimebootstrapimprovement}.
\end{proof}

\begin{corollary}\label{cbootstrap}
Under the bootstrap assumption \eqref{booassum}, for every $t\geq1$ and $k\in\mathbb Z$, one has
\begin{equation}\label{bspco1}
\|\widehat{P_kf}(t)\|_{L_\xi^1}+\|\widehat{P_{\leq k}f}(t)\|_{L_\xi^1}\lesssim\epsilon_1
\begin{cases}
2^{2k},&k\leq0,\\
2^{(1-\epsilon)k},&k>0,
\end{cases}
\end{equation}
hence,
\begin{equation}\label{bspco2}
\sum_k2^{-k}\|\widehat{P_kf}(t)\|_{L_\xi^1}\lesssim\epsilon_1.
\end{equation}

For $k>0$, one also has
\begin{equation}\label{ebootstrapcor3}
\|P_kf(t)\|_{L_x^1}\lesssim\epsilon_1t^{\frac\beta2}2^{-\frac\epsilon2k}.
\end{equation}
In particular, for every $k\in\mathbb Z$,
\begin{equation}\label{ebootstrapcor4}
\|P_k\psi(t)\|_{L_x^\infty}\lesssim\epsilon_1\min\left\{2^{2k},t^{-1},t^{-1+\frac\beta2}2^{-\frac\epsilon2k}\right\},
\end{equation}
and
\begin{align}
\sum_k\|P_k\psi(t)\|_{L_x^\infty}&\lesssim\epsilon_1t^{-1}(1+\log t),\label{bs5}\\
\sum_k\|P_k\psi(t)\|_{L_x^\infty}^2&\lesssim\epsilon_1^2t^{-2}(1+\log t),\label{ebootstrapcor6}\\
\sum_k2^{-k}\|P_k\psi(t)\|_{L_x^\infty}&\lesssim\epsilon_1t^{-\frac12}.\label{bso7}
\end{align}
For every $2\leq p<\infty$,
\begin{equation}\label{bsc8}
\sum_k\|P_k\psi(t)\|_{L_x^p}\lesssim\epsilon_1t^{-1+\frac2p}.
\end{equation}
Moreover, for every $2\leq p<\infty$,
\begin{equation}\label{bsc9}
\|e^{it\Delta}(xf(t))\|_{L_x^p}\lesssim\epsilon_1t^{-1+\beta+\frac2p(1+\gamma-\beta)},
\end{equation}
and, in the particular case $p=4$,
\begin{equation}\label{bsc10}
\|e^{it\Delta}(xf(t))\|_{L_x^4}\lesssim\epsilon_1t^{-\frac12+\frac\beta2}.
\end{equation}
For every $k\in\mathbb Z$,
\begin{equation}\label{bsc11}
\|P_k(xf(t))\|_{L_x^2}\lesssim\epsilon_1 2^{\frac k2}t^{\frac\beta2},
\end{equation}
whereas, for $k>0$,
\begin{align}
\|P_k(xf(t))\|_{L_x^2}&\lesssim\epsilon_1t^{\frac\beta2}2^{-\frac\epsilon2k},\label{bsc12}\\
\|P_ke^{it\Delta}(xf(t))\|_{L_x^4}&\lesssim\epsilon_1t^{-\frac12+\frac{3\beta}{4}}2^{-\frac\epsilon4k}.\label{bsc13}
\end{align}
Hence,
\begin{equation}\label{bsc14}
\sum_k\|P_ke^{it\Delta}(xf(t))\|_{L_x^4}\lesssim\epsilon_1t^{-\frac12+\frac\beta2}(1+\log t).
\end{equation}
The implicit constants may depend on $\epsilon$ and, where applicable, on $p$.
\end{corollary}

\begin{proof}
For $k\leq0$, Plancherel and \eqref{booassum} give $\|P_kf(t)\|_{L_x^2}\lesssim2^k\|\widehat f(t)\|_{L_\xi^\infty}\lesssim\epsilon_12^k$, for $k>0$, $\|P_kf(t)\|_{L_x^2}\lesssim2^{-\epsilon k}\|f(t)\|_{H^\epsilon}\lesssim\epsilon_12^{-\epsilon k}$. Using
\begin{equation*}
\|\widehat{P_kf}(t)\|_{L_\xi^1}\lesssim2^k\|P_kf(t)\|_{L_x^2},
\end{equation*}
and summing over frequencies below $2^k$, we obtain \eqref{bspco1}, and
\begin{equation*}
\sum_k2^{-k}\|\widehat{P_kf}(t)\|_{L_\xi^1}\lesssim\epsilon_1\left(\sum_{k\leq0}2^k+\sum_{k>0}2^{-\epsilon k}\right)\lesssim\epsilon_1,
\end{equation*}
which proves \eqref{bspco2}.

For $k>0$, commuting $x^2$ with $P_k$ gives
\begin{equation*}
\|x^2P_kf(t)\|_{L_x^2}\lesssim\|x^2f(t)\|_{L_x^2}+2^{-k}\|xf(t)\|_{L_x^2}+2^{-2k}\|f(t)\|_{L_x^2}\lesssim\epsilon_1t^\beta.
\end{equation*}
Then by \eqref{booassum}, we obtain
\begin{equation*}
\|P_kf(t)\|_{L_x^1}\lesssim\|P_kf(t)\|_{L_x^2}^{\frac12}\|x^2P_kf(t)\|_{L_x^2}^{\frac12}\lesssim\epsilon_1t^{\frac\beta2}2^{-\frac\epsilon2k}.
\end{equation*}
This proves \eqref{ebootstrapcor3}. The dispersive estimate gives the third bound in \eqref{ebootstrapcor4} for $k>0$. The first two bounds follow from $\|\widehat f(t)\|_{L_\xi^\infty}\lesssim\epsilon_1$ and $\|\psi(t)\|_{L_x^\infty}\lesssim\epsilon_1t^{-1}$. For $k\leq0$, the third bound is weaker than the second, so \eqref{ebootstrapcor4} holds for every $k$.

For \eqref{bs5}, splitting the sum at $2^k=t^{-\frac12}$ and $2^k=t^{\frac\beta\epsilon}$, we obtain
\begin{align*}
\sum_k\|P_k\psi(t)\|_{L_x^\infty}&\lesssim\epsilon_1\left(\sum_{2^k\leq t^{-\frac12}}2^{2k}+t^{-1}\sum_{t^{-\frac12}<2^k\leq t^{\frac\beta\epsilon}}1+t^{-1+\frac\beta2}\sum_{2^k>t^{\frac\beta\epsilon}}2^{-\frac\epsilon2k}\right)\\
&\lesssim\epsilon_1t^{-1}(1+\log t).
\end{align*}
Since $\sup_k\|P_k\psi(t)\|_{L_x^\infty}\lesssim\epsilon_1t^{-1}$, multiplying this estimate by the sum above proves \eqref{ebootstrapcor6}.
Moreover,
\begin{equation*}
\sum_k2^{-k}\|P_k\psi(t)\|_{L_x^\infty}\lesssim\epsilon_1\left(\sum_{2^k\leq t^{-\frac12}}2^k+t^{-1}\sum_{2^k>t^{-\frac12}}2^{-k}\right)\lesssim\epsilon_1t^{-\frac12}.
\end{equation*}
Thus \eqref{bso7} follows. 

For $2\leq p<\infty$, interpolation and \eqref{booassum} give
\begin{equation*}
\sum_k\|P_k\psi(t)\|_{L_x^p}\lesssim\|\psi(t)\|_{L_x^\infty}^{1-\frac2p}\sum_k\|P_kf(t)\|_{L_x^2}^{\frac2p}\lesssim\epsilon_1t^{-1+\frac2p}\left(\sum_{k\leq0}2^{\frac{2k}{p}}+\sum_{k>0}2^{-\frac{2\epsilon k}{p}}\right).
\end{equation*}
This proves \eqref{bsc8}.

For $2<p<\infty$, H\"older inequality gives, for every $M>0$,
\begin{align*}
  \|xf(t)\|_{L_x^{\frac{p}{p-1}}}&\lesssim\|xf(t)\|_{L^{\frac{p}{p-1}}(|x|\leq M)}+\|xf(t)\|_{L^{\frac{p}{p-1}}(|x|>M)}\\
  &\lesssim M^{1-\frac2p}\|xf(t)\|_{L_x^2}+M^{-\frac2p}\|x^2f(t)\|_{L_x^2},
\end{align*}
optimizing in $M$ and using the dispersive estimate, we obtain
\begin{equation*}
\|e^{it\Delta}(xf(t))\|_{L_x^p}\lesssim t^{-1+\frac2p}\|xf(t)\|_{L_x^2}^{\frac2p}\|x^2f(t)\|_{L_x^2}^{1-\frac2p}\lesssim\epsilon_1t^{-1+\beta+\frac2p(1+\gamma-\beta)},
\end{equation*}
which proves \eqref{bsc9}. For $p=4$, recalling that $J=x+2it\nabla$ and \(J\psi(t)=e^{it\Delta}(xf(t))\),
\[
  \bigl\|e^{it\Delta}(xf(t))\bigr\|_{L_x^4}^2=\|J\psi(t)\|_{L_x^4}^2\lesssim\|\psi(t)\|_{L_x^\infty}\|J^2\psi(t)\|_{L_x^2}\lesssim\epsilon_1^2t^{-1+\beta},
\]
which proves \eqref{bsc10}.

By H\"older and Gagliardo-Nirenberg inequality in Fourier space,
\begin{align*}
\|P_k(xf(t))\|_{L_x^2}\lesssim2^{\frac k2}\|\nabla_\xi\widehat f(t)\|_{L_\xi^4}\lesssim2^{\frac k2}\|\widehat f(t)\|_{L_\xi^\infty}^{\frac12}\|\nabla_\xi^2\widehat f(t)\|_{L_\xi^2}^{\frac12}\lesssim\epsilon_12^{\frac k2}t^{\frac\beta2}.
\end{align*}
This proves \eqref{bsc11}. For $k>0$, Fourier integration by parts and Cauchy-Schwarz give
\begin{align*}
\|P_k(xf)\|_{L_x^2}^2&=-\operatorname{Re}\int\varphi_k^2\overline{\widehat f}\,\Delta_\xi\widehat f\,\mathrm d\xi-2\operatorname{Re}\int\varphi_k\nabla_\xi\varphi_k\cdot\nabla_\xi\widehat f\,\overline{\widehat f}\,\mathrm d\xi\\
&\leq\frac12\|P_k(xf)\|_{L_x^2}^2+C\|P_{[k-2,k+2]}f\|_{L_x^2}\|x^2f\|_{L_x^2}+C2^{-2k}\|P_{[k-2,k+2]}f\|_{L_x^2}^2.
\end{align*}
Hence,
\begin{equation*}
\|P_k(xf(t))\|_{L_x^2}\lesssim\|P_{[k-2,k+2]}f(t)\|_{L_x^2}^{\frac12}\|x^2f(t)\|_{L_x^2}^{\frac12}+2^{-k}\|P_{[k-2,k+2]}f(t)\|_{L_x^2}\lesssim\epsilon_1t^{\frac\beta2}2^{-\frac\epsilon2k},
\end{equation*}
which proves \eqref{bsc12}. For $k>0$, commuting $x$ with $P_k$ also gives
\begin{equation*}
\|xP_k(xf(t))\|_{L_x^2}\lesssim\|x^2f(t)\|_{L_x^2}+2^{-k}\|xf(t)\|_{L_x^2}\lesssim\epsilon_1t^\beta.
\end{equation*}
Using the weighted interpolation established above with $p=4$, followed by the dispersive estimate, we obtain
\begin{align*}
\|P_ke^{it\Delta}(xf(t))\|_{L_x^4}&\lesssim t^{-\frac12}\|P_k(xf(t))\|_{L_x^{\frac43}}\lesssim t^{-\frac12}\|P_k(xf(t))\|_{L_x^2}^{\frac12}\|xP_k(xf(t))\|_{L_x^2}^{\frac12}\\
&\lesssim\epsilon_1t^{-\frac12+\frac{3\beta}{4}}2^{-\frac\epsilon4k}.
\end{align*}
This proves \eqref{bsc13}.

Finally, Bernstein and \eqref{booassum} imply $\|P_ke^{it\Delta}(xf(t))\|_{L_x^4}\lesssim\epsilon_12^{\frac k2}t^\gamma$. Splitting at $2^k=t^{-1}$ and $2^k=t^{\frac\beta\epsilon}$, and using \eqref{bsc10} on the middle frequencies and \eqref{bsc13} on the high frequencies, we obtain
\begin{align*}
\sum_k\|P_ke^{it\Delta}(xf(t))\|_{L_x^4}&\lesssim\epsilon_1t^\gamma\sum_{2^k\leq t^{-1}}2^{\frac k2}+(1+\log t)\|e^{it\Delta}(xf(t))\|_{L_x^4}+\epsilon_1t^{-\frac12+\frac{3\beta}{4}}\sum_{2^k>t^{\frac\beta\epsilon}}2^{-\frac\epsilon4k}\\
&\lesssim\epsilon_1t^{-\frac12+\frac\beta2}(1+\log t),
\end{align*}
which proves \eqref{bsc14}.
\end{proof}

\begin{corollary}[Scattering and pointwise asymptotics of the differential fields]\label{cpsi scattering}
Under \eqref{ini}, let $\psi$ be the differential fields in the caloric gauge and let $g_+$ be the scattering state in Theorem~\ref{thm:main}. For every $t\geq1$,
\begin{equation}\label{epsiscat1}
\|\psi(t)-e^{it\Delta}\nabla g_+\|_{H^{\epsilon-1}}\lesssim\epsilon_0^3t^{-1}.
\end{equation}
Moreover,
\begin{equation}\label{epsiscat2}
\psi(t)=e^{it\Delta}\nabla g_++O_{L_x^\infty}\left(\epsilon_0^{1+\epsilon}t^{-1-\frac{\epsilon-\beta}{2}}\right),
\end{equation}
and
\begin{equation}\label{epsiscat3}
\psi(t,x)=\frac{e^{\frac{i|x|^2}{4t}}}{2it}\widehat{\nabla g_+}\left(\frac{x}{2t}\right)+O_{L_x^\infty}\left(\epsilon_0t^{-1-\frac{\epsilon-\beta}{2}}\right).
\end{equation}
\end{corollary}
\begin{proof}
See Subsection~\ref{subsec:psi-scatter}.
\end{proof}

\section{Estimates of differential fields along the heat flow: energy norms and negative half order norms}\label{sec:heat envelope}

In this section we establish decay estimates along the heat flow for the differential fields $\psi$ and the connection coefficients $A$. These estimates serve two purposes. On the one hand, they describe quantitatively how the differential fields propagate and decay along the heat parameter $s$. On the other hand, they will be used later to derive a collection of estimates for the connection coefficients $A$ at the initial heat time $s=0$. We first control the energy norms of the differential fields and then turn to the corresponding negative half order norms.

Throughout this section we fix the physical time $t$ and omit it from the notation, all estimates are uniform in $t$. The frequency localized heat Duhamel formula is
\begin{equation}
\label{eheatduhamelpk}
\begin{aligned}
P_k\psi_m(s)=&e^{s\Delta}P_k\psi_m(0)+2i\int_0^s e^{(s-r)\Delta}P_k\nabla\cdot(A\psi_m)(r)\,\mathrm{d}r\\
&-\int_0^s e^{(s-r)\Delta}P_k\bigl((A_lA_l+i\partial_lA_l)\psi_m\bigr)(r)\,\mathrm{d}r+i\int_0^s e^{(s-r)\Delta}P_k\bigl(\operatorname{Im}(\psi_m\overline{\psi_l})\psi_l\bigr)(r)\,\mathrm{d}r.
\end{aligned}
\end{equation}
The main tool used to estimate the Duhamel term is the Littlewood-Paley paraproduct decomposition. More precisely, 
\begin{equation}\label{eparaproductdecomposition}
\begin{aligned}
P_k(fg)=&\sum_{\substack{k_1\leq k-5\\ |k_2-k|\leq 4}}P_k\bigl(P_{k_1}f\,P_{k_2}g\bigr)+\sum_{\substack{|k_1-k|\leq 4\\ k_2\leq k-5}}P_k\bigl(P_{k_1}f\,P_{k_2}g\bigr)\\
&\quad+\sum_{\substack{k_1,k_2\geq k-4\\ |k_1-k_2|\leq 8}}P_k\bigl(P_{k_1}f\,P_{k_2}g\bigr).
\end{aligned}
\end{equation}
The three sums in \eqref{eparaproductdecomposition} will be referred to, respectively, as the Low-High, High-Low, and High-High interactions. We estimate these three interactions separately. We emphasize two points in the argument. First, for the low-high interaction in which \(A\) carries the high frequency, the polynomial decay weights must be treated more carefully. Indeed, owing to the integral representation of \(A\), the corresponding tail integral loses one power of polynomial decay, and the refined frequency dependent estimate for \(A\) is needed to compensate for this loss. Second, for the high-high interaction, we rely on a simple polynomial heat kernel convolution lemma, which allows us to sum efficiently over high frequencies.

We formulate all the estimates below in terms of frequency envelopes. We recall the definition that will be used throughout this section.

\begin{definition}[Frequency envelope]
\label{def:frequency envelope}
Fix a sufficiently small constant $\delta>0$. A positive sequence $\{b_k\}_{k\in\ZZ}$ is called a frequency envelope if
\begin{equation*}
\sum_{k\in\ZZ}b_k^2<\infty
\end{equation*}
and if it is slowly varying in the sense that
\begin{equation*}
b_k\leq 2^{\delta|k-j|}b_j, \qquad k,j\in\ZZ.
\end{equation*}
\end{definition}

\begin{remark}
\label{rem:fixed time envelope}
In contrast with the spacetime frequency envelopes used in \emph{Bejenaru--Ionescu--Kenig--Tataru}, every frequency envelope in this section is defined after fixing the physical time
$t$. In particular, no supremum over $t$ is included in its definition. All estimates below are therefore fixed $t$ estimates along the heat flow, with constants uniform in $t$. This formulation preserves the dependence on the physical time and is adapted to the decay estimates established later in the paper.
\end{remark}

\begin{remark}
In the course of the frequency localized estimates, a finite enlargement of a frequency envelope frequently occurs:
\begin{equation*}
\widetilde b_k:=\sup_{|k'-k|\leq 4}b_{k'}.
\end{equation*}
The slow variation property immediately gives $b_k\leq\widetilde b_k\leq2^{4\delta}b_k.$ Hence, $\widetilde b_k\sim b_k$, with constants depending only on $\delta$. The same observation applies to $b_k(\sigma)$ and to all the other frequency envelopes used below. To avoid unnecessary notation, all such finite enlargements will henceforth be denoted simply by $b_k$, or by $b_k(\sigma)$ when the regularity index must be displayed.
\end{remark}

\subsection{The energy frequency envelope and relative heat estimates}

In this subsection the physical time $t$ is fixed and is therefore omitted from the notation. In particular, no spacetime norm in the physical time variable is involved. After expanding the imaginary part in \eqref{ecaloric4}, all terms other than $2i\nabla\cdot(A\psi_m)$ are finite linear combinations of $F\psi$, where
\begin{equation}\label{eFdefinition}
F\in\bigl\{\nabla\cdot A,\ A_lA_m,\ \psi_l\overline{\psi_m}: l,m\in\{1,2\}\bigr\}.
\end{equation}
We omit component indices in the estimates below. The caloric representation of the connection coefficients takes the form
\begin{equation}\label{ecaloricconnection}
  A(s)=-\int_s^\infty \operatorname{Im}(\overline\psi\nabla\psi)(r)\,\mathrm{d}r-\int_s^\infty \operatorname{Im}(iA\overline\psi\psi)(r)\,\mathrm{d}r.
\end{equation}
The former is the quadratic term (to be distinguished from the definition $A^{(2)}$ introduced later), and the latter is the remainder term.

For $0\leq\sigma\leq\epsilon$, define the frequency envelope at the fixed physical time $t$ by
\begin{equation}\label{eenergyenvelope}
b_k(\sigma):=\epsilon_1^{-1}\sup_{j\in\ZZ}2^{-\delta|k-j|}2^{\sigma j}\|P_j\psi(0)\|_{L_x^2},\qquad b_k:=b_k(0).
\end{equation}
It follows immediately that
\begin{equation}\label{eenergyenvelope2}
\|P_k\psi(0)\|_{L_x^2}\leq\epsilon_1 2^{-\sigma k}b_k(\sigma),\qquad b_j(\sigma)\leq2^{\delta|j-k|}b_k(\sigma).
\end{equation}
Moreover, the bootstrap assumption \eqref{booassum} gives, for every $0\leq\sigma\leq\epsilon$,
\begin{equation}\label{eenergyenvelope3}
\sum_{k\in\ZZ}b_k^2(\sigma)\lesssim\epsilon_1^{-2}\sum_{k\in\ZZ}2^{2\sigma k}\|P_k\psi(0)\|_{L_x^2}^2\lesssim1.
\end{equation}

We shall close the following bootstrap estimate.

\begin{proposition}[Propagation of the energy envelope]\label{propenv}
Under the bootstrap assumption \eqref{booassum}, assume that, for every $0\leq\sigma\leq\epsilon$,
\begin{equation}\label{eheatenergybootstrap}
\|P_k\psi(s)\|_{L_x^2}\leq2C\epsilon_1 2^{-\sigma k}b_k(\sigma)(1+s2^{2k})^{-20},
\end{equation}
then, for every $0\leq\sigma\leq\epsilon$,
\begin{equation}\label{eheatenergyimproved2}
\|P_k\psi(s)\|_{L_x^2}\leq\frac{3}{2}C\epsilon_1 2^{-\sigma k}b_k(\sigma)(1+s2^{2k})^{-20}.
\end{equation}
\end{proposition}

\begin{remark}
At each fixed physical time, the argument in the proof of \cite[Corollary~4.4]{BIKT2011} and (4.11) therein give \eqref{eheatenergybootstrap} for all $s\geq0$ once $C$ is sufficiently large. For smooth approximations, (4.12) therein and slow variation ensure that the best constants in the energy envelope estimates and in the negative half order estimates \eqref{enegativeenve1} and \eqref{erk3} below are finite. For the weighted envelope in \eqref{eJbootstrap} below, we iterate the linear commuted heat equation \eqref{eJs6} on finite heat intervals, with bounds independent of the endpoint. Thus, in the proofs of Propositions~\ref{propenv}, \ref{prohaboo}, and~\ref{pJbootstrap}, it suffices to establish the stated bootstrap improvements.
\end{remark}

The main ingredient in the proof of Proposition~\ref{propenv} is the following pair of estimates for the connection coefficients.

\begin{proposition}[Connection estimates in the energy norm]\label{pconnection energy estimates}
Under \eqref{eheatenergybootstrap}, for $0\leq\sigma\leq\epsilon$, define
\begin{equation}\label{edkdef}
d_k(\sigma):=b_k(\sigma)\sum_{j\leq k-5}2^{j-k}b_j+\sum_{j\geq k-4}2^{k-j}b_jb_j(\sigma).
\end{equation}
Then
\begin{equation}\label{edkdef2}
d_k(\sigma)\approx b_kb_k(\sigma),\qquad d_k:=d_k(0)\approx b_k^2,
\end{equation}
and the following estimates hold:
\begin{equation}\label{econnectionenve1}
\|P_kA(s)\|_{L_x^2}\lesssim\epsilon_1^2 2^{-\sigma k}d_k(\sigma)(1+s2^{2k})^{-19},
\end{equation}
\begin{equation}\label{econnectionenve2}
\|P_kA(s)\|_{L_x^2}\lesssim\epsilon_1^2 2^{-\sigma k}b_k(\sigma)\left[\sum_{j\leq k-5}2^{j-k}b_j(1+s2^{2j})^{-20}(1+s2^{2k})^{-19}+b_k(1+s2^{2k})^{-20}\right].
\end{equation}
\end{proposition}

\begin{remark}
Estimate~\eqref{econnectionenve1} is a convenient closed expression for subsequent bounds on the connection along the heat flow, and it is also the form used later in the propagation of the negative half order norm.  The refined estimate~\eqref{econnectionenve2} has a different purpose.  In the heat Duhamel formula, when $A$ has the output high frequency and $\psi$ has a much lower frequency, the factor $(1+s2^{2k})^{-19}$ in \eqref{econnectionenve1} is not sufficient by itself to recover the
desired power $(1+s2^{2k})^{-20}$. The first term in \eqref{econnectionenve2} preserves both the low frequency gain $2^{j-k}$ and the low frequency heat weight, while its second term reflects the fact that the High-High contribution to $A_k$ already has the full power $(1+s2^{2k})^{-20}$. More precisely, for the first term in \eqref{econnectionenve2}, we rely on inequality
\[
\big[\sum_{j\leq k-5}2^{j-k}b_j(1+s2^{2j})^{-20}\big]^2\lesssim(1+s2^{2k})^{-1}
\]
to gain an extra polynomial heat decay weight.
\end{remark}

\begin{remark}
The estimates above are simpler than the corresponding statements in Lemmas 5.2 and 5.3 of \cite{BIKT2011}. There the connection is estimated in spaces containing local smoothing and maximal function components, and the heat variable is divided into dyadic intervals to compare $2^k$ with the thermal frequency $s^{-\frac{1}{2}}$. This leads to a heat time dependent frequency envelope. In the present fixed time $L_x^2$ argument, Bernstein inequalities give the explicit factors, while \eqref{epolynomialconvolution} and \eqref{epolynomialconvolution2} perform the entire heat time summation. Thus neither a dyadic decomposition of the heat variable nor a more complicated
heat time dependent envelope is required.
\end{remark}

\begin{proof}[Proof of Proposition~\ref{pconnection energy estimates}]
By the definition \eqref{edkdef}, we note from the slow variation of $b_k$ and $b_k(\sigma)$ that
\[
\begin{aligned}
  \sum_{j\leq k-5}2^{j-k}b_j&\lesssim b_k\sum_{j\leq k-5}2^{-(1-\delta)(k-j)}\lesssim b_k,\\
  \sum_{j\geq k-4}2^{k-j}b_jb_j(\sigma)&\lesssim b_kb_k(\sigma)\sum_{j\geq k-4}2^{-(1-2\delta)(j-k)}\lesssim b_kb_k(\sigma).
\end{aligned}
\]
The reverse bound in~\eqref{edkdef2} is obvious.

Since \eqref{econnectionenve1} is a direct result of \eqref{econnectionenve2}, it suffices to prove the refined estimate \eqref{econnectionenve2}. We first consider the quadratic term in
\eqref{ecaloricconnection}.

\smallskip
\noindent\textit{Quadratic Low-High interactions: $j\leq k-5$ and the differentiated factor has frequency $k+O(1)$.} By Bernstein inequality and \eqref{eheatenergybootstrap},
\begin{align}
  &\left\|P_k\int_s^\infty P_j\psi(r)\nabla P_{k+O(1)}\psi(r)\,\mathrm{d}r\right\|_{L_x^2}\notag\\
  &\quad\lesssim 2^{j+k}\int_s^\infty\|P_j\psi(r)\|_{L_x^2}\|P_{k+O(1)}\psi(r)\|_{L_x^2}\,\mathrm{d}r\notag\\
  &\quad\lesssim\epsilon_1^2 2^{j+k}2^{-\sigma k}b_jb_k(\sigma)\int_s^\infty(1+r2^{2j})^{-20}(1+r2^{2k})^{-20}\,\mathrm{d}r\notag\\
  &\quad\lesssim\epsilon_1^2 2^{j-k}2^{-\sigma k}b_jb_k(\sigma)(1+s2^{2j})^{-20}(1+s2^{2k})^{-19}. \label{eAquad1}
\end{align}
In the last line we used $r\geq s$ to take $(1+s2^{2j})^{-20}$ outside the integral and then applied \eqref{epolynomialconvolution2} at frequency $2^k$.

\smallskip
\noindent\textit{Quadratic High-Low interactions: $j\leq k-5$ and the differentiated factor has frequency $2^j$.}
In this case,
\begin{align}
  \left\|P_k\int_s^\infty P_{k+O(1)}\psi(r)\nabla P_j\psi(r)\,\mathrm{d}r\right\|_{L_x^2}&\lesssim 2^{2j}\int_s^\infty\|P_{k+O(1)}\psi(r)\|_{L_x^2}\|P_j\psi(r)\|_{L_x^2}\,\mathrm{d}r\notag\\
  &\lesssim\epsilon_1^2 2^{2(j-k)}2^{-\sigma k}b_jb_k(\sigma)(1+s2^{2j})^{-20}(1+s2^{2k})^{-19}. \label{eAquad2}
\end{align}
Since $2^{2(j-k)}\leq 2^{j-k}$ for $j\leq k-5$, this contribution is dominated by the Low-High term in~\eqref{econnectionenve2}.

\smallskip
\noindent\textit{Quadratic High-High interactions: the two input frequencies are $j+O(1)$ and $j\geq k-4$.}
Using Bernstein inequality again, we have
\[
\begin{aligned}
  \|P_k(P_j\psi\nabla P_{j+O(1)}\psi)(r)\|_{L_x^2}&\lesssim 2^k\|P_j\psi(r)\nabla P_{j+O(1)}\psi(r)\|_{L_x^1}\notag\\&\lesssim 2^{k+j}\|P_j\psi(r)\|_{L_x^2}\|P_{j+O(1)} \psi(r)\|_{L_x^2}.
\end{aligned}
\]
Hence,
\[
\begin{aligned}
  \left\|P_k\int_s^\infty P_j\psi(r)\nabla P_{j+O(1)}\psi(r)\,\mathrm{d}r\right\|_{L_x^2}&\lesssim\epsilon_1^22^{k+j}2^{-\sigma j}b_jb_j(\sigma)\int_s^\infty(1+r2^{2j})^{-40}\,\mathrm{d}r\notag\\
  &\lesssim\epsilon_1^2 2^{k-j}2^{-\sigma j}b_jb_j(\sigma)(1+s2^{2j})^{-39}.
\end{aligned}
\]
Since $j\geq k-4$,$(1+s2^{2j})^{-39}\lesssim (1+s2^{2k})^{-20}$. Moreover,
\begin{equation*}
  \sum_{j\geq k-4}2^{k-j}2^{-\sigma j}b_jb_j(\sigma)\lesssim 2^{-\sigma k}b_kb_k(\sigma)\sum_{j\geq k-4}2^{-(1+\sigma-2\delta)(j-k)}\lesssim 2^{-\sigma k}b_kb_k(\sigma).
\end{equation*}
Thus the High-High part is bounded directly by
\begin{equation}
  \epsilon_1^2 2^{-\sigma k}b_kb_k(\sigma)(1+s2^{2k})^{-20}. \label{eAquad3}
\end{equation}
Combining~\eqref{eAquad1}, \eqref{eAquad2}, and \eqref{eAquad3} proves both desired estimates for the quadratic term in~\eqref{ecaloricconnection}.

We turn to the remainder term.  We use the following two elementary trilinear inequalities.  After arranging three frequencies as $k_1\leq k_2\leq k_3$, if $k=k_3+O(1)$, then
\begin{equation}
  \|P_k(P_{k_1}fP_{k_2}gP_{k_3}h)\|_{L_x^2}\lesssim 2^{k_1+k_2}\|P_{k_1}f\|_{L_x^2}\|P_{k_2}g\|_{L_x^2}\|P_{k_3}h\|_{L_x^2}. \label{etrilinear1}
\end{equation}
If $k\leq k_3-5$ and $k_2=k_3+O(1)$, then
\begin{equation}
  \|P_k(P_{k_1}fP_{k_2}gP_{k_3}h)\|_{L_x^2}\lesssim 2^{k+k_1}\|P_{k_1}f\|_{L_x^2}\|P_{k_2}g\|_{L_x^2}\|P_{k_3}h\|_{L_x^2}. \label{etrilinear2}
\end{equation}
The first estimate follows by placing the two lower frequencies in $L_x^\infty$. For the second one, we first use $L_x^1\to L_x^2$ Bernstein at the output and then place the lowest frequency in $L_x^\infty$.

We now close~\eqref{econnectionenve2}. Let $C$ be the best constant in that estimate, then also \eqref{econnectionenve1} holds under the constant $2C$. The constant in \eqref{eheatenergybootstrap} is fixed, and all implicit constants below are independent of the connection constant $C$. In the remainder term we write the two $\psi$ frequencies as $k_1\leq k_2$, the frequency of $A$ as $k_3$, and $\bar k=\max\{k_1,k_2,k_3\}$.

\smallskip
\noindent\textit{Refined remainder case 1: $k_2=\bar k=k+O(1)$.}
We place the $\sigma$ derivatives on the high frequency $\psi$ factor. First assume $k_1\leq k-5$. We use the rough estimate \eqref{econnectionenve1} for $P_{k_3}A$ and keep the heat weight of the low frequency $P_{k_1}\psi$ factor. Since $r\geq s$,
\begin{equation*}
  \int_s^\infty (1+r2^{2k_1})^{-20}(1+r2^{2k})^{-20}(1+r2^{2k_3})^{-19}\,\mathrm{d}r\lesssim 2^{-2k}(1+s2^{2k_1})^{-20}(1+s2^{2k})^{-19}.
\end{equation*}
Hence, \eqref{etrilinear1}, \eqref{eheatenergybootstrap}, the rough bootstrap bound \eqref{econnectionenve1}, \eqref{edkdef2}, and \eqref{eenergyenvelope3} give
\[
\begin{aligned}
  &\sum_{\substack{k_1\leq k-5\\k_3\leq k+4}}\left\|P_k\int_s^\infty P_{k_1}\psi P_{k+O(1)}\psi P_{k_3}A(r)\,\mathrm{d}r\right\|_{L_x^2}\notag\\
  &\quad\lesssim C\epsilon_1^4 2^{-\sigma k}b_k(\sigma)(1+s2^{2k})^{-19}\notag\\
  &\hspace{15mm}\times\left[\sum_{k_1\leq k-5}2^{k_1-k}b_{k_1}(1+s2^{2k_1})^{-20}\right]\left[\sum_{k_3\leq k+4}2^{k_3-k}d_{k_3}\right]\notag\\
  &\quad\lesssim C\epsilon_1^4 2^{-\sigma k}b_k(\sigma)\sum_{k_1\leq k-5}2^{k_1-k}b_{k_1}(1+s2^{2k_1})^{-20}(1+s2^{2k})^{-19}.
\end{aligned}
\]
If $k_1=k+O(1)$, then the two $\psi$ factors have frequency $k+O(1)$, and
\[
\begin{aligned}
  &\sum_{k_3\leq k+4}\left\|P_k\int_s^\infty P_{k+O(1)}\psi P_{k+O(1)}\psi P_{k_3}A(r)\,\mathrm{d}r\right\|_{L_x^2}\notag\\
  &\quad\lesssim C\epsilon_1^4 2^{-\sigma k}b_kb_k(\sigma)\sum_{k_3\leq k+4}2^{k+k_3}d_{k_3}\int_s^\infty(1+r2^{2k})^{-40}(1+r2^{2k_3})^{-19}\,\mathrm{d}r\notag\\
  &\quad\lesssim C\epsilon_1^4 2^{-\sigma k}b_kb_k(\sigma)(1+s2^{2k})^{-39}\left[\sum_{k_3\leq k+4}2^{k_3-k}d_{k_3}\right]\notag\\
  &\quad\lesssim C\epsilon_1^4 2^{-\sigma k}b_k(\sigma)\left[b_k(1+s2^{2k})^{-20}\right].
\end{aligned}
\]

\smallskip
\noindent\textit{Refined remainder case 2: $k_3=\bar k=k+O(1)$ and $k_2\leq k-5$.}
This is the most complicated case, since here $A$ is the only output high frequency, we must counteract the loss of one order polynomial weight. Using \eqref{etrilinear1} and we insert the full refined bootstrap estimate for the high frequency $A$ factor, and we place the $\sigma$ derivatives on $A$:
\begin{align}
  &\sum_{k_1\leq k_2\leq k-5}\left\|P_k\int_s^\infty P_{k_1}\psi P_{k_2}\psi P_{k+O(1)}A(r)\,\mathrm{d}r\right\|_{L_x^2}\notag\\
  &\quad\lesssim\epsilon_1^2\sum_{k_1\leq k_2\leq k-5}\int_s^\infty 2^{k_1+k_2}b_{k_1}(1+r2^{2k_1})^{-20}b_{k_2}(1+r2^{2k_2})^{-20}\|P_{k+O(1)}A(r)\|_{L_x^2}\,\mathrm{d}r\notag\\
  &\quad\lesssim C\epsilon_1^4 2^{-\sigma k}b_k(\sigma)\int_s^\infty 2^{2k}\left[\sum_{j\leq k-5}2^{j-k}b_j(1+r2^{2j})^{-20}\right]^2\notag\\
  &\hspace{15mm}\times\Bigg[\sum_{j\leq k-5}2^{j-k}b_j(1+r2^{2j})^{-20}(1+r2^{2k})^{-19}+b_k(1+r2^{2k})^{-20}\Bigg]\,\mathrm{d}r. \label{eAremainder1}
\end{align}
For the first term in the last bracket, we use the following relation: by Cauchy-Schwarz and~\eqref{eenergyenvelope3},
\begin{align}
  \left[\sum_{j\leq k-5}2^{j-k}b_j(1+r2^{2j})^{-20}\right]^2&\leq\left(\sum_{j\leq k-5}b_j^2\right)\left[\sum_{j\leq k-5}2^{2(j-k)}(1+r2^{2j})^{-40}\right]\notag\\
  &\lesssim (1+r2^{2k})^{-1}.\label{eAremainder2}
\end{align}
Indeed, if $r2^{2k}\leq 1$, the last frequency sum is bounded by a geometric series. If $r2^{2k}\geq 1$, then on the part $r2^{2j}\leq 1$ one has
\begin{equation*}
  \sum_{\substack{j\leq k-5\\r2^{2j}\leq 1}}2^{2(j-k)}\lesssim (r2^{2k})^{-1},
\end{equation*}
whereas on the part $r2^{2j}\geq 1$,
\begin{equation*}
  \sum_{\substack{j\leq k-5\\r2^{2j}\geq 1}}2^{2(j-k)}(1+r2^{2j})^{-40}\lesssim(r2^{2k})^{-1}\sum_{\substack{j\leq k-5\\r2^{2j}\geq 1}}(r2^{2j})^{-39}\lesssim (r2^{2k})^{-1}.
\end{equation*}

Hence, keeping one low frequency heat weight and applying \eqref{epolynomialconvolution2}, we obtain
\begin{align}
  &\int_s^\infty 2^{2k}(1+r2^{2k})^{-1}\sum_{j\leq k-5}2^{j-k}b_j(1+r2^{2j})^{-20}(1+r2^{2k})^{-19}\,\mathrm{d}r\notag\\
  &\quad\leq\sum_{j\leq k-5}2^{j-k}b_j(1+s2^{2j})^{-20}\int_s^\infty 2^{2k}(1+r2^{2k})^{-20}\,\mathrm{d}r\notag\\
  &\quad\lesssim\sum_{j\leq k-5}2^{j-k}b_j(1+s2^{2j})^{-20}(1+s2^{2k})^{-19}. \label{eAremainder3}
\end{align}
For the second term in the bracket of \eqref{eAremainder1},
\begin{equation}
  b_k\int_s^\infty 2^{2k}(1+r2^{2k})^{-1}(1+r2^{2k})^{-20}\,\mathrm{d}r\lesssim b_k(1+s2^{2k})^{-20}. \label{eAremainder4}
\end{equation}
Combining~\eqref{eAremainder3} and \eqref{eAremainder4} proves the desired refined estimate in case 2.

\smallskip
\noindent\textit{Refined remainder case 3: $k_2=k_3+O(1)=\bar k$ and $\bar k\geq k+5$.}
Here one $\psi$ factor and $A$ form the high frequency pair. We place the $\sigma$ derivatives on $A$. By \eqref{etrilinear2}, \eqref{eheatenergybootstrap}, \eqref{econnectionenve1}, \eqref{epolynomialconvolution2}, and \eqref{eenergyenvelope2}, we obtain
\[
\begin{aligned}
  &\sum_{\bar k\geq k+5}\sum_{k_1\leq \bar k+4}\left\|P_k\int_s^\infty P_{k_1}\psi P_{\bar k+O(1)}\psi P_{\bar k+O(1)}A(r)\,\mathrm{d}r\right\|_{L_x^2}\\
  &\quad\lesssim C\epsilon_1^4\sum_{\bar k\geq k+5}\sum_{k_1\leq \bar k+4}2^{k+k_1}2^{-\sigma\bar k}b_{k_1}b_{\bar k}d_{\bar k}(\sigma)\int_s^\infty(1+r2^{2k_1})^{-20}(1+r2^{2\bar k})^{-39}\,\mathrm{d}r\\
  &\quad\lesssim C\epsilon_1^4\sum_{\bar k\geq k+5}2^{k-\bar k}2^{-\sigma\bar k}b_{\bar k}^2d_{\bar k}(\sigma)(1+s2^{2\bar k})^{-38}\\
  &\quad\lesssim C\epsilon_1^4 2^{-\sigma k}b_kb_k(\sigma)(1+s2^{2k})^{-20}.
\end{aligned}
\]

\smallskip
\noindent\textit{Refined remainder case 4: $k_1=k_2+O(1)=\bar k$ and $\bar k\geq k+5$.}
The two $\psi$ factors form the high-frequency pair. We put the $\sigma$ derivatives on one of them and  again use the rough estimate \eqref{econnectionenve1} for $A$. By \eqref{etrilinear2},
\[
\begin{aligned}
  &\sum_{\bar k\geq k+5}\sum_{k_3\leq \bar k+4}\left\|P_k\int_s^\infty P_{\bar k+O(1)}\psi P_{\bar k+O(1)}\psi P_{k_3}A(r)\,\mathrm{d}r\right\|_{L_x^2}\\
  &\quad\lesssim C\epsilon_1^4\sum_{\bar k\geq k+5}\sum_{k_3\leq \bar k+4}2^{k+k_3}2^{-\sigma\bar k}b_{\bar k}b_{\bar k}(\sigma)d_{k_3}\int_s^\infty(1+r2^{2\bar k})^{-40}(1+r2^{2k_3})^{-19}\,\mathrm{d}r\\
  &\quad\lesssim C\epsilon_1^4\sum_{\bar k\geq k+5}2^{k-\bar k}2^{-\sigma\bar k}b_{\bar k}b_{\bar k}(\sigma)d_{\bar k}(1+s2^{2\bar k})^{-39}\\
  &\quad\lesssim C\epsilon_1^4 2^{-\sigma k}b_kb_k(\sigma)(1+s2^{2k})^{-20}.
\end{aligned}
\]

The four cases show that the remainder term is bounded by $C\epsilon_1^2$ times the right hand side of \eqref{econnectionenve2}.  The quadratic term is bounded by an absolute constant times the same right hand side.  Hence the best constant in both estimates satisfies
\begin{equation*}
  C\lesssim 1+C\epsilon_1^2.
\end{equation*}
Choosing $\epsilon_1$ sufficiently small completes the proof of Proposition~\ref{pconnection energy estimates}.
\end{proof}

\begin{lemma}
Under \eqref{eheatenergybootstrap}, for $F$ as in \eqref{eFdefinition}, $s>0$, and $0\leq\sigma\leq\epsilon$, one has
\begin{equation}\label{eFenve1}
\|P_kF(s)\|_{L_x^2}\lesssim\epsilon_1^2s^{-\frac12}2^{-\sigma k}b_k(\sigma)(1+s2^{2k})^{-19}.
\end{equation}
\end{lemma}

\begin{proof}
For $F=\nabla\cdot A$, the estimate follows from Bernstein inequality, \eqref{econnectionenve2}, and
\[
\sum_{j\leq k-5}2^jb_j(1+s2^{2j})^{-20}+2^kb_k(1+s2^{2k})^{-1}\lesssim s^{-\frac{1}{2}}.
\]
For the two products, we use the paraproduct decomposition, \eqref{eheatenergybootstrap}, \eqref{econnectionenve1}, and \eqref{edkdef2}. The Low-High and High-Low interactions follow from
\[
\|\psi(s)\|_{L_x^\infty}+\epsilon_1^{-1}\|A(s)\|_{L_x^\infty}\lesssim\epsilon_1s^{-\frac{1}{2}}.
\]
The High-High contributions are bounded by
\[
\epsilon_1^2 2^k\sum_{j\geq k-4}2^{-\sigma j}b_jb_j(\sigma)(1+s2^{2j})^{-38}\lesssim\epsilon_1^2s^{-\frac{1}{2}}2^{-\sigma k}b_k(\sigma)(1+s2^{2k})^{-19},
\]
the additional factor $\epsilon_1^2$ for $AA$ is harmless.
\end{proof}

We now prove Proposition~\ref{propenv}.
\begin{proof}[Proof of Proposition~\ref{propenv}]
By \eqref{eheatduhamelpk}, the linear term satisfies
\begin{equation}\label{eenergyenveproof1}
  \|e^{s\Delta}P_k\psi(0)\|_{L_x^2}\lesssim e^{-cs2^{2k}}\epsilon_1 2^{-\sigma k}b_k(\sigma)\lesssim \epsilon_1 2^{-\sigma k}b_k(\sigma)(1+s2^{2k})^{-20}.
\end{equation}
We first estimate the term containing $\nabla\cdot(A\psi)$ by the paraproduct decomposition and Lemma~\ref{lpolynomialconvolution}.

\smallskip
\noindent\textit{Duhamel Low-High interactions: $A$ has frequency $j\leq k-5$ and $\psi$ has frequency $k+O(1)$.}
Using Bernstein, \eqref{epolynomialconvolution} and the rough form \eqref{econnectionenve1} for $A$, we obtain
\begin{align}
  &\sum_{j\leq k-5}\left\|\int_0^s e^{(s-r)\Delta}P_k\nabla\cdot(P_jA P_{k+O(1)}\psi)(r)\,\mathrm{d}r\right\|_{L_x^2}\notag\\
  &\quad\lesssim\sum_{j\leq k-5}2^{k+j}\int_0^s e^{-c(s-r)2^{2k}}\|P_jA(r)\|_{L_x^2}\|P_{k+O(1)}\psi(r)\|_{L_x^2}\,\mathrm{d}r\notag\\
  &\quad\lesssim\epsilon_1^3 2^{-\sigma k}b_k(\sigma)\sum_{j\leq k-5}2^{k+j}d_j\int_0^s e^{-c(s-r)2^{2k}}(1+r2^{2j})^{-19}(1+r2^{2k})^{-20}\,\mathrm{d}r\notag\\
  &\quad\lesssim\epsilon_1^3 2^{-\sigma k}b_k(\sigma)(1+s2^{2k})^{-20}\sum_{j\leq k-5}2^{j-k}d_j\notag\\
  &\quad\lesssim\epsilon_1^3 2^{-\sigma k}b_k(\sigma)(1+s2^{2k})^{-20}. \label{eheatLH}
\end{align}

\smallskip
\noindent\textit{Duhamel High-Low interactions: $A$ has frequency $k+O(1)$ and $\psi$ has frequency $j\leq k-5$.}
This is the interaction for which the refined connection estimate \eqref{econnectionenve2} is essential. By Bernstein,
\[
\begin{aligned}
  &\sum_{j\leq k-5}\left\|\int_0^s e^{(s-r)\Delta}P_k\nabla\cdot(P_{k+O(1)}A P_j\psi)(r)\,\mathrm{d}r\right\|_{L_x^2}\notag\\
  &\quad\lesssim\sum_{j\leq k-5}2^{k+j}\int_0^s e^{-c(s-r)2^{2k}}\|P_{k+O(1)}A(r)\|_{L_x^2}\|P_j\psi(r)\|_{L_x^2}\,\mathrm{d}r.
\end{aligned}
\]
Inserting~\eqref{econnectionenve2} gives
\begin{align}
  &\sum_{j\leq k-5}\left\|\int_0^s e^{(s-r)\Delta}P_k\nabla\cdot(P_{k+O(1)}A P_j\psi)(r)\,\mathrm{d}r\right\|_{L_x^2}\notag\\
  &\quad\lesssim\epsilon_1^32^{-\sigma k}b_k(\sigma)\int_0^s e^{-c(s-r)2^{2k}}2^{2k}\left[\sum_{j\leq k-5}2^{j-k}b_j(1+r2^{2j})^{-20}\right]\notag\\
  &\hspace{15mm}\times\Bigg[\sum_{j\leq k-5}2^{j-k}b_j(1+r2^{2j})^{-20}(1+r2^{2k})^{-19}+b_k(1+r2^{2k})^{-20}\Bigg]\,\mathrm{d}r. \label{eenergyenveproof2}
\end{align}
For the product of the two low frequency sums, using \eqref{eAremainder2} to recover one polynomial heat decay power:
\begin{align}
  &\left[\sum_{j\leq k-5}2^{j-k}b_j(1+r2^{2j})^{-20}\right]\left[\sum_{j\leq k-5}2^{j-k}b_j(1+r2^{2j})^{-20}(1+r2^{2k})^{-19}\right]\notag\\
  &\quad=(1+r2^{2k})^{-19}\left[\sum_{j\leq k-5}2^{j-k}b_j(1+r2^{2j})^{-20}\right]^2\lesssim (1+r2^{2k})^{-20}. \label{eenergyenveproof3}
\end{align}
For the remaining term,
\begin{equation}
  b_k(1+r2^{2k})^{-20}\sum_{j\leq k-5}2^{j-k}b_j(1+r2^{2j})^{-20}\lesssim b_k^2(1+r2^{2k})^{-20}\lesssim (1+r2^{2k})^{-20}. \label{eenergyenveproof4}
\end{equation}
Substituting~\eqref{eenergyenveproof3} and \eqref{eenergyenveproof4} into \eqref{eenergyenveproof2}, and by \eqref{epolynomialconvolution}, we conclude that
\begin{align}
  &\sum_{j\leq k-5}\left\|\int_0^s e^{(s-r)\Delta}P_k\nabla\cdot(P_{k+O(1)}A P_j\psi)(r)\,\mathrm{d}r\right\|_{L_x^2}\notag\\
  &\quad\lesssim\epsilon_1^3 2^{-\sigma k}b_k(\sigma)\int_0^s e^{-c(s-r)2^{2k}}2^{2k}(1+r2^{2k})^{-20}\,\mathrm{d}r\notag\\
  &\quad\lesssim\epsilon_1^3 2^{-\sigma k}b_k(\sigma)(1+s2^{2k})^{-20}.\label{eenergyenveproof5}
\end{align}

\smallskip
\noindent\textit{Duhamel High-High interactions: the two input frequencies are $j+O(1)$ and $j\geq k-4$.}
At this stage the polynomial heat convolution estimate\eqref{epolynomialconvolution} yields the crucial factor $2^{-2j}$. For the output frequency, Bernstein gives
\begin{equation*}
  \|P_k\nabla\cdot(P_jA P_{j+O(1)}\psi)(r)\|_{L_x^2}\lesssim 2^{2k}\|P_jA(r)\|_{L_x^2}\|P_{j+O(1)}\psi(r)\|_{L_x^2}.
\end{equation*}
Using~\eqref{econnectionenve1} with the $\sigma$ derivatives on $A$, \eqref{epolynomialconvolution}, \eqref{edkdef2}, we obtain
\begin{align}
  &\sum_{j\geq k-4}\left\|\int_0^s e^{(s-r)\Delta}P_k\nabla\cdot(P_jA P_{j+O(1)}\psi)(r)\,\mathrm{d}r\right\|_{L_x^2}\notag\\
  &\quad\lesssim\epsilon_1^3\sum_{j\geq k-4}2^{2k}2^{-\sigma j}d_j(\sigma)b_j\int_0^s e^{-c(s-r)2^{2k}}(1+r2^{2j})^{-39}\,\mathrm{d}r\notag\\
  &\quad\lesssim\epsilon_1^3\sum_{j\geq k-4}2^{2(k-j)}2^{-\sigma j}d_j(\sigma)b_j(1+s2^{2k})^{-20}\notag\\
  &\quad\lesssim\epsilon_1^3 2^{-\sigma k}b_k(\sigma)(1+s2^{2k})^{-20}. \label{eenergyenveproof6}
\end{align}

For the remaining nonlinearities, the paraproduct decomposition, \eqref{eFenve1}, and \eqref{eAremainder2} give
\[
\|P_k(F\psi)(r)\|_{L_x^2}\lesssim\epsilon_1^3r^{-\frac{1}{2}}2^{(1-\sigma)k}b_k(\sigma)(1+r2^{2k})^{-\frac{39}{2}}+\epsilon_1^3r^{-\frac{1}{2}}2^k\sum_{j\geq k-4}2^{-\sigma j}b_jb_j(\sigma)(1+r2^{2j})^{-39}.
\]
Using
\begin{equation}\label{esingularheatconvolution}
\int_0^s e^{-c(s-r)2^{2k}}r^{-\frac{1}{2}}(1+r2^{2j})^{-M}\,\mathrm{d}r\lesssim 2^{-j}(1+s2^{2k})^{-M-\frac{1}{2}},
\end{equation}
valid for $s>0$, $j\geq k-4$, and $M>\frac{1}{2}$, we obtain
\begin{equation}\label{eenergyenveproof7}
\left\|\int_0^s e^{(s-r)\Delta}P_k(F\psi)(r)\,\mathrm{d}r\right\|_{L_x^2}\lesssim\epsilon_1^3 2^{-\sigma k}b_k(\sigma)(1+s2^{2k})^{-20}.
\end{equation}

Combining \eqref{eenergyenveproof1}, \eqref{eheatLH}, \eqref{eenergyenveproof5}, \eqref{eenergyenveproof6}, and \eqref{eenergyenveproof7}, we obtain
\begin{equation*}
  \|P_k\psi(s)\|_{L_x^2}\leq C\epsilon_1 2^{-\sigma k}b_k(\sigma)(1+s2^{2k})^{-20}+\tilde{C}\epsilon_1^3 2^{-\sigma k}b_k(\sigma)(1+s2^{2k})^{-20}.
\end{equation*}
After increasing $C$ once to dominate the linear estimate and then choosing $\epsilon_1$ sufficiently small, the second term on the right-hand side is at most $\frac12 C\epsilon_1 2^{-\sigma k}b_k(\sigma)(1+s2^{2k})^{-20}$.  This proves \eqref{eheatenergyimproved2} and closes the bootstrap argument.
\end{proof}

\subsection{The negative half order frequency envelope and relative heat estimates}
In this subsection, we define the negative half order frequency envelope $q_k$ at fixed time $t$ and derive its propagation along the heat flow, which is crucial for the subsequent estimates of the connection coefficients $A$.

\begin{definition}[Negative half order envelope]\label{def:qk}
  For every fixed physical time $t$, define the fixed time negative half order frequency envelope
\begin{equation}\label{eqk1}
  q_k:=\epsilon_1^{-1}\sup_{j\in\ZZ}2^{-\delta|k-j|}2^{-\frac{j}{2}}\|P_j\psi(0,t)\|_{L_x^2}.
\end{equation}
\end{definition}
\begin{lemma}\label{lqk}
For every fixed physical time $t$, one has
\begin{equation}\label{eqk2}
  \|\psi(0,t)\|_{\dot H_x^{-\frac{1}{2}}}\lesssim\epsilon_1.
\end{equation}
Hence, the envelope $q_k$ satisfies
\begin{equation}\label{eqk3}
  \|P_k\psi(0,t)\|_{L_x^2}\leq\epsilon_1 2^{\frac{k}{2}}q_k, \qquad\sum_{k\in\ZZ}q_k^2\lesssim 1, \qquad q_k\lesssim2^{\delta|k-j|}q_j.
\end{equation}
\end{lemma}

\begin{proof}
For $k\leq0$, Plancherel theorem and the size of the annulus $\{|\xi|\sim2^k\}$ imply
\begin{equation*}
  \|P_k\psi(0,t)\|_{L_x^2}\lesssim 2^k\|\widehat{\psi}(0,t)\|_{L_\xi^\infty}=2^k\|\hat{f}(t)\|_{L_\xi^\infty}\lesssim \epsilon_1 2^k.
\end{equation*}
Therefore,
\begin{equation}\label{eqk4}
  \sum_{k\leq0}2^{-k}\|P_k\psi(0,t)\|_{L_x^2}^2\lesssim\epsilon_1^2\sum_{k\leq0}2^k\lesssim\epsilon_1^2.
\end{equation}
For $k>0$, the $L_x^2$ bound in \eqref{booassum} gives
\begin{equation}\label{eqk5}
\sum_{k>0}2^{-k}\|P_k\psi(0,t)\|_{L_x^2}^2\leq\|\psi(0,t)\|_{L_x^2}^2\lesssim\epsilon_1^2.
\end{equation}
Combining \eqref{eqk4} and \eqref{eqk5} proves \eqref{eqk2}.

The estimate \eqref{eqk3} follows directly from \eqref{eqk2} and the definition \eqref{eqk1}.
\end{proof}

\begin{remark}[Motivation for the negative half order estimate]
The choice of the negative half order norm is suggested directly by the quadratic part of the connection $A^{(2)}$. Indeed, after replacing both heat fields by their free evolutions and integrating in the heat variable, then
\begin{equation*}
  A^{(2)}(0,t)\sim\int_0^\infty \overline{e^{s\Delta}\psi(0,t)}\nabla e^{s\Delta}\psi(0,t)\,\mathrm{d}s,
\end{equation*}
the relevant bilinear symbol is bounded by
\begin{equation}\label{eqk6}
 \frac{|\xi-\eta|}{|\xi-\eta|^2+|\eta|^2}+\frac{|\eta|}{|\xi-\eta|^2+|\eta|^2}\lesssim\frac{1}{|\xi-\eta|^{\frac12}|\eta|^{\frac12}}.
\end{equation}
Thus, instead of using only the order $-1$ estimate
\[
 \frac{|\xi-\eta|+|\eta|}{|\xi-\eta|^2+|\eta|^2}\lesssim\min\left\{|\xi|^{-1},|\xi-\eta|^{-1},|\eta|^{-1}\right\},
\]
one may distribute the full smoothing effect evenly between the two input frequencies. By Cauchy-Schwarz in the convolution variable, \eqref{eqk6} gives
\begin{equation*}
 \|\widehat{A^{(2)}}(0,t)\|_{L_\xi^\infty}\lesssim\bigl\||\nabla|^{-\frac12}\psi(0,t)\bigr\|_{L_x^2}^2.
\end{equation*}
This balanced estimate is particularly useful in high-low interactions: it avoids assigning the entire negative derivative  to one input and treats the two heat fields symmetrically.

It is therefore natural to propagate the negative half order norm of $\psi$ along the heat flow. The free heat evolution does not increase this norm, while the heat Duhamel term can be controlled perturbatively. At the dyadic level, the resulting bound contains the factor $2^{\frac{k}{2}}$ at low frequencies, which provides the summability needed in the estimates for the connection. At high frequencies, the polynomial heat kernel factor converts positive powers of $2^k$ into decay in the heat variable. Hence the negative half order estimate simultaneously captures the balanced order $-1$ structure of the quadratic connection, resolves the low frequency summation, and yields the heat time decay required for differentiated estimates of both $\psi$ and $A$.
\end{remark}

The main focus of this subsection is to prove the following proposition, which provides negative half order heat estimates for the differential field $\psi$.
\begin{proposition}[Propagation of the negative half order envelope]\label{prohaboo}
Under the bootstrap assumption \eqref{booassum}, assume that
\begin{equation}\label{enegativeenve1}
  \|P_k\psi(s,t)\|_{L_x^2}\leq 2C\epsilon_12^{\frac{k}{2}}q_k(1+s2^{2k})^{-12},
\end{equation}
then 
\begin{equation}\label{enegativeenve2}
  \|P_k\psi(s,t)\|_{L_x^2}\leq\frac{3}{2}C\epsilon_12^{\frac{k}{2}}q_k(1+s2^{2k})^{-12}.
\end{equation}
\end{proposition}

Analogously to Proposition\ref{pconnection energy estimates}, we also establish the negative half order estimate for the connection coefficients.

\begin{proposition}[Connection estimates in the negative half order norm]\label{pconnection negative estimates}
Under \eqref{enegativeenve1}, define
\begin{equation}\label{erk1}
  r_k:=b_k\sum_{j\leq k-5}2^{\frac{3}{2}(j-k)}q_j+\sum_{j\geq k-4}2^{\frac{1}{2}(k-j)}b_jq_j.
\end{equation}
Then 
\begin{equation}\label{erk2}
  r_k\sim b_kq_k,
\end{equation}
the following estimate holds:
\begin{equation}\label{erk3}
  \|P_kA(s)\|_{L_x^2}\lesssim\epsilon_1^2 2^{\frac{k}{2}}r_k(1+s2^{2k})^{-18}.
\end{equation}
\end{proposition}

\begin{remark}
When the connection coefficient \(A\) in the higher order term \(A\psi\psi\) carries a frequency much higher than those of the two differential fields, the tail integration in the heat variable loses one power of polynomial decay. As in Proposition~\ref{pconnection energy estimates}, this loss could be avoided by establishing a more refined frequency-dependent estimate for \(A\). We shall not pursue such a refinement here, since the stronger decay rate \((1+s2^{2k})^{-19}\) is not needed for our argument. Indeed, the bound \eqref{erk3} already provides more than sufficient decay to close the propagation estimate for the negative half order norm of \(\psi\), for which only the weight \((1+s2^{2k})^{-12}\) is required. 
\end{remark}
\begin{proof}
Estimate \eqref{erk2} follows directly from \eqref{erk1}. The proof of \eqref{erk3} is analogous to \eqref{econnectionenve1}. We prove the quadratic part and clarify the correction to the remainder term estimates. We emphasize here that since the negative half order framework mainly handles the low frequency regime, we employ energy estimates \eqref{eheatenergyimproved2} or \eqref{econnectionenve1} for high frequencies. Therefore, the polynomial decay of $A$ is better than that of $\psi$ in this setting.

Consider the quadratic term in \eqref{ecaloricconnection}, for the low-high interaction, $j\leq k-5$, 
\begin{equation*}
\|P_k(P_{j}\psi(r)\nabla P_{k}\psi(r))\|_{L_x^2}\lesssim\|P_{j}\psi(r)\|_{L_x^\infty}2^{k}\|P_{k}\psi(r)\|_{L_x^2}\lesssim\epsilon_1^22^{\frac{3}{2}j+k}q_{j}b_k(1+r2^{2j})^{-12}(1+r2^{2k})^{-20}.
\end{equation*}
After integration in $r$, we have
\begin{equation*}
  \int_s^\infty\|P_k(P_j\psi(r)\nabla P_k\psi(r))\|_{L_x^2}\,\mathrm{d}r\lesssim\epsilon_1^2 2^{\frac{k}{2}}2^{\frac{3}{2}(j-k)}b_kq_{j}(1+s2^{2k})^{-18}.
\end{equation*}
When the derivative falls on the low frequency factor, the resulting factor is $2^{\frac{5}{2}(j-k)}$, which is smaller than the low-high case. 

For the high-high interaction, we use one  energy estimate and one negative half order estimate. For $j\geq k-4$, 
\begin{equation*}
  \|P_k(P_j\psi(r)\nabla P_{\sim j}\psi(r))\|_{L_x^2}\lesssim 2^k\|P_j\psi(r)\|_{L_x^2}\|\nabla P_{\sim j}\psi(r)\|_{L_x^2}\lesssim\epsilon_1^22^{k+\frac{3}{2}j}b_jq_j(1+r2^{2j})^{-32}.
\end{equation*}
Thus
\begin{equation*}
  \int_s^\infty\|P_k(P_j\psi(r)\nabla P_{\sim j}\psi(r))\|_{L_x^2}\,\mathrm{d}r\lesssim\epsilon_1^2 2^{\frac{k}{2}}2^{\frac{1}{2}(k-j)}b_jq_j(1+s2^{2k})^{-18}.
\end{equation*}
Summing the low-high and high-high contributions gives
\begin{align*}
\left\|P_k\int_s^\infty\psi(r)\nabla\psi(r)\,\mathrm{d}r\right\|_{L_x^2}&\lesssim\epsilon_1^2 2^{\frac{k}{2}}\left(b_k\sum_{j\leq k}2^{\frac{3}{2}(j-k)}q_j+\sum_{j\geq k-4}2^{\frac{1}{2}(k-j)}b_jq_j\right)(1+s2^{2k})^{-18}\notag\\
&\lesssim\epsilon_1^2 2^{\frac{k}{2}}r_k(1+s2^{2k})^{-18}.
\end{align*}

Now in order to close \eqref{erk3}, let $C$ be the best constant in that estimate, denote $\bar{k}=\max\{k_1, k_2, k_3\}$, assume $k_1\leq k_2$, and estimate
\begin{equation*}
\sum_{k_1,k_2,k_3}\int_s^\infty\left\|P_k\bigl(P_{k_1}\psi(r)P_{k_2}\psi(r)P_{k_3}A(r)\bigr)\right\|_{L_x^2}\,\mathrm{d}r.
\end{equation*}

In each frequency interaction, the lowest frequency factor is placed in the corresponding negative half order norm, while the other two factors are placed in their energy norms. 
\smallskip
\noindent\textit{Negative half order remainder case 1: \(k_2=\bar k=k+O(1)\).}

If \(k_1\leq k_3\), we place \(P_{k_1}\psi\) in the negative half order norm and use the energy norms for \(P_{k+O(1)}\psi\) and \(P_{k_3}A\). If \(k_3<k_1\), we instead place \(P_{k_3}A\) in the negative half order norm \eqref{erk3} and use the energy norms for both differential fields. By \eqref{etrilinear1} and \eqref{erk2},
\[
\begin{aligned}
&\sum_{k_1,k_3\leq k+4}\left\|P_k\int_s^\infty P_{k_1}\psi(r)P_{k+O(1)}\psi(r)P_{k_3}A(r)\,\mathrm{d}r\right\|_{L_x^2}\\
&\quad\lesssim\epsilon_1^4(1+s2^{2k})^{-19}\sum_{\substack{k_1,k_3\leq k+4\\ k_1\leq k_3}}2^{\frac{3k_1}{2}+k_3-2k}q_{k_1}b_kd_{k_3}\\
&\qquad\quad+C\epsilon_1^4(1+s2^{2k})^{-19}\sum_{\substack{k_1,k_3\leq k+4\\ k_3<k_1}}2^{k_1+\frac{3k_3}{2}-2k}b_{k_1}b_kr_{k_3}\\
&\quad\lesssim C\epsilon_1^4 2^{\frac{k}{2}}r_k(1+s2^{2k})^{-18}.
\end{aligned}
\]

\smallskip

\noindent\textit{Negative half order remainder case 2: \(k_3=\bar k=k+O(1)\) and \(k_2\leq k-5\).}

Here \(A\) is the unique highest frequency factor. This is the only case in which we do not insert the negative half order bootstrap bound for \(A\). Instead, we use the energy estimate \eqref{econnectionenve1} for \(A\). Thus by \eqref{etrilinear1} and \eqref{erk2},
\[
\begin{aligned}
&\sum_{k_1\leq k_2\leq k-5}\left\|P_k\int_s^\infty P_{k_1}\psi(r)P_{k_2}\psi(r)P_{k+O(1)}A(r)\,\mathrm{d}r\right\|_{L_x^2}\\
&\quad\lesssim\epsilon_1^4\sum_{k_1\leq k_2\leq k-5}2^{k_1+k_2}2^{\frac{k_1}{2}}q_{k_1}b_{k_2}d_k\int_s^\infty(1+r2^{2k_1})^{-12}(1+r2^{2k_2})^{-20}(1+r2^{2k})^{-19}\,\mathrm{d}r\\
&\quad\lesssim\epsilon_1^4\sum_{k_1\leq k_2\leq k-5}2^{\frac{3k_1}{2}+k_2-2k}q_{k_1}b_{k_2}d_k(1+s2^{2k})^{-18}\\
&\quad\lesssim\epsilon_1^42^{\frac{k}{2}}r_k(1+s2^{2k})^{-18}.
\end{aligned}
\]

\smallskip

\noindent\textit{Negative half order remainder case 3: \(k_2=k_3+O(1)=\bar k\) and \(\bar k\geq k+5\).}

We place \(P_{k_1}\psi\) in the negative half order norm and use the energy estimates for the two high frequency factors. By \eqref{etrilinear2} and \eqref{erk2}, 
\[
\begin{aligned}
&\sum_{\bar k\geq k+5}\sum_{k_1\leq\bar k+4}\left\|P_k\int_s^\infty P_{k_1}\psi(r)P_{\bar k+O(1)}\psi(r)P_{\bar k+O(1)}A(r)\,\mathrm{d}r\right\|_{L_x^2}\\
&\quad\lesssim\epsilon_1^4\sum_{\bar k\geq k+5}\sum_{k_1\leq\bar k+4}2^{k+k_1}2^{\frac{k_1}{2}}q_{k_1}b_{\bar k}d_{\bar k}\int_s^\infty(1+r2^{2k_1})^{-12}(1+r2^{2\bar k})^{-39}\,\mathrm{d}r\\
&\quad\lesssim\epsilon_1^4 2^{\frac{k}{2}}\sum_{\bar k\geq k+5}\sum_{k_1\leq\bar k+4}2^{\frac12(k-\bar k)}2^{\frac32(k_1-\bar k)}q_{k_1}b_{\bar k}d_{\bar k}(1+s2^{2\bar k})^{-38}\\
&\quad\lesssim\epsilon_1^4 2^{\frac{k}{2}}r_k(1+s2^{2k})^{-18}.
\end{aligned}
\]

\smallskip

\noindent\textit{Negative half order remainder case 4: \(k_1=k_2+O(1)=\bar k\) and \(\bar k\geq k+5\).}

We simply place $A$ in the negative half order norm \eqref{erk3} and place both differential fields in their energy norms. By \eqref{etrilinear2} and \eqref{erk2}, 
\[
\begin{aligned}
&\sum_{\bar k\geq k+5}\sum_{k_3\leq\bar k+4}\left\|P_k\int_s^\infty P_{\bar k+O(1)}\psi(r)P_{\bar k+O(1)}\psi(r)P_{k_3}A(r)\,\mathrm{d}r\right\|_{L_x^2}\\
&\quad\lesssim C\epsilon_1^4\sum_{\bar k\geq k+5}\sum_{k_3\leq\bar k+4}2^{k+k_3}2^{\frac{k_3}{2}}b_{\bar k}^2r_{k_3}\int_s^\infty(1+r2^{2\bar k})^{-40}(1+r2^{2k_3})^{-18}\,\mathrm{d}r\\
&\quad\lesssim C\epsilon_1^4 2^{\frac{k}{2}}\sum_{\bar k\geq k+5}\sum_{k_3\leq\bar k+4}2^{\frac12(k-\bar k)}2^{\frac32(k_3-\bar k)}b_{\bar k}^2r_{k_3}(1+s2^{2\bar k})^{-39}\\
&\quad\lesssim C\epsilon_1^4 2^{\frac{k}{2}}r_k(1+s2^{2k})^{-18}.
\end{aligned}
\]

Finally, applying the same argument like Proposition~\ref{pconnection energy estimates}, the best constant satisfies $C\lesssim 1+C\epsilon_1^2$, choosing $\epsilon_1$ sufficiently small yields the proof.
\end{proof}

Finally we prove Proposition~\ref{prohaboo}.

\begin{proof}[Proof of Proposition~\ref{prohaboo}]
By \eqref{eheatduhamelpk}, the linear term is bounded by
\begin{equation}\label{erk4}
  \|e^{s\Delta}P_k\psi(0)\|_{L_x^2}\lesssim e^{-cs2^{2k}}\|P_k\psi(0)\|_{L_x^2}\lesssim\epsilon_1 2^{\frac{k}{2}}q_k(1+s2^{2k})^{-12}.
\end{equation}

For $j\leq k-5$, we combine the negative half order estimate for the low frequency factor with the energy estimate for the high frequency factor. This gives
\[
\begin{aligned}
  \|P_k(P_jA\,P_{k+O(1)}\psi)(r)\|_{L_x^2}&\lesssim\epsilon_1^32^{\frac{3j}{2}}r_jb_k(1+r2^{2j})^{-18}(1+r2^{2k})^{-20},\\
  \|P_k(P_{k+O(1)}A\,P_j\psi)(r)\|_{L_x^2}&\lesssim\epsilon_1^32^{\frac{3j}{2}}d_kq_j(1+r2^{2k})^{-19}(1+r2^{2j})^{-12}.
\end{aligned}
\]
Applying Lemma \ref{lpolynomialconvolution}, there holds 
\begin{align}
  &\sum_{j\leq k-5}\bigg(\left\|\int_0^s e^{(s-r)\Delta}P_k\nabla\cdot(P_jA P_{k+O(1)}\psi)(r)\,\mathrm{d}r\right\|_{L_x^2}+\left\|\int_0^s e^{(s-r)\Delta}P_k\nabla\cdot(P_{k+O(1)}A P_j\psi)(r)\,\mathrm{d}r\right\|_{L_x^2}\bigg)\notag\\
  &\quad\lesssim\epsilon_1^3 2^{\frac{k}{2}}q_k(1+s2^{2k})^{-12}. \label{erk5}
\end{align}

For high-high interactions, $j\geq k-4$, we place $A$ in its energy norm and $\psi$ in its negative half order norm. Thus
\begin{equation*}
  \|P_k(P_jA\,P_{j+O(1)}\psi)(r)\|_{L_x^2}\lesssim\epsilon_1^32^{k+\frac{j}{2}}d_jq_j(1+r2^{2j})^{-31}.
\end{equation*}
By Lemma~\ref{lpolynomialconvolution},
\begin{align}
  &\sum_{j\geq k-4}\left\|\int_0^s e^{(s-r)\Delta}P_k\nabla\cdot(P_jA P_{j+O(1)}\psi)(r)\,\mathrm{d}r\right\|_{L_x^2}\notag\\
  &\quad\lesssim\epsilon_1^3 2^{\frac{k}{2}}(1+s2^{2k})^{-12}\sum_{j\geq k-4}2^{\frac32(k-j)}d_jq_j\lesssim \epsilon_1^32^{\frac{k}{2}}q_k(1+s2^{2k})^{-12}. \label{erk6}
\end{align}

For the remaining nonlinearities, we combine \eqref{eFenve1} with $\sigma=0$ and \eqref{enegativeenve1}. The paraproduct decomposition gives
\[
\|P_k(F\psi)(r)\|_{L_x^2}\lesssim\epsilon_1^3r^{-\frac{1}{2}}2^{\frac{3k}{2}}q_k(1+r2^{2k})^{-12}+\epsilon_1^3r^{-\frac{1}{2}}2^k\sum_{j\geq k-4}2^{\frac{j}{2}}b_jq_j (1+r2^{2j})^{-31}.
\]
By \eqref{esingularheatconvolution},
\begin{equation}\label{erk7}
\left\|\int_0^s e^{(s-r)\Delta}P_k(F\psi)(r)\,\mathrm{d}r\right\|_{L_x^2}\lesssim\epsilon_1^32^{\frac{k}{2}}q_k(1+s2^{2k})^{-12}.
\end{equation}

Combining \eqref{erk4}, \eqref{erk5}, \eqref{erk6}, and \eqref{erk7}, we obtain
\begin{equation*}
  \|P_k\psi(s)\|_{L_x^2}\leq C\epsilon_12^{\frac{k}{2}}q_k(1+s2^{2k})^{-12}+\tilde{C}\epsilon_1^3 2^{\frac{k}{2}}q_k(1+s2^{2k})^{-12}.
\end{equation*}
For $\epsilon_1$ sufficiently small, this improves the constant $2C$ in \eqref{enegativeenve1} to $\frac{3}{2}C$ and closes the bootstrap.
\end{proof}

\section{Estimates for the connection coefficients}
\label{sec:connection estimates}

This section is devoted to establishing a collection of estimates for the connection coefficients $A$, which will be used in the energy and dispersive estimates developed in the subsequent sections. The discussion is divided into six parts. First, we establish the Fourier estimate $\bigl\|\widehat{A}\bigr\|_{L_\xi^\infty}.$ Second, we derive preliminary decay estimates for $\psi$ and $A$ with respect to both the heat parameter $s$ and the physical time $t$. Combined with the bootstrap argument on $L^\infty$ norms, we present some more refined estimates. In the third part, we prove decay estimates for $A$ and its derivatives in $L_x^2$ and $L_x^\infty$. The fourth part is concerned with the first order weighted estimate $\|xA\|_{L_x^2}$. In the fifth part, we establish refined estimates for the higher order part $A^{(\geq4)}$ of the spatial connection coefficients. Finally, we estimate the temporal connection coefficient $A_3$ and its higher order part. The proofs rely primarily on the propagation estimates along the heat flow for the energy norm and the negative half order norm obtained in the previous section, together with the Galilean type heat modified vector field $J_s=x+2(s+it)\nabla$. In this section, all estimates are uniform with respect to the physical time $t$.

\subsection{Uniform Fourier control for $\|\hat{A}\|_{L_\xi^\infty}$ along the heat time}

We first estimate $\|\hat{A}\|_{L_\xi^\infty}$, which relies heavily on the negative half order estimates established in the previous section to control the low frequency. 
\begin{proposition}\label{ppsihat} 
Under the bootstrap assumption \eqref{booassum}, for every fixed physical time $t$,
\begin{equation}\label{epsihat1}
  \sup_{s\geq0}\|(1+s|\xi|^2)^{12}\widehat{\psi}(s,t)\|_{L_\xi^\infty}\lesssim\epsilon_1.
\end{equation}
In particular,
\begin{equation}\label{epsihat2}
  \|\widehat{P_k\psi}(s,t)\|_{L_\xi^\infty}\lesssim\epsilon_1(1+s2^{2k})^{-12}.
\end{equation}
\end{proposition}

\begin{proof}
We first estimate the contribution of $\nabla\cdot(A\psi)$ in \eqref{eheatduhamelpk}. By Young inequality, \eqref{enegativeenve2}, and \eqref{erk3}, for the low-high interaction,
\begin{equation*}
  \left\|\widehat{P_k(P_jA(r)P_k\psi(r))}\right\|_{L_\xi^\infty}\lesssim\epsilon_1^3 2^{\frac{j+k}{2}}r_jq_k(1+r2^{2j})^{-18}(1+r2^{2k})^{-12}, \qquad j\leq k-5.
\end{equation*}
For the high-low interaction,
\begin{equation*}
  \left\|\widehat{P_k(P_kA(r)P_j\psi(r))}\right\|_{L_\xi^\infty}\lesssim\epsilon_1^3 2^{\frac{j+k}{2}}r_kq_j(1+r2^{2k})^{-18}(1+r2^{2j})^{-12}, \qquad j\leq k-5.
\end{equation*}
For the high-high interaction,
\begin{equation*}
  \left\|\widehat{P_k(P_jA(r)P_j\psi(r))}\right\|_{L_\xi^\infty}\lesssim\epsilon_1^3 2^jr_jq_j(1+r2^{2j})^{-30}, \qquad j\geq k-4.
\end{equation*}
By \eqref{epolynomialconvolution}, for every $s\geq0$, we obtain
\begin{equation}\label{epsihat3}
\begin{aligned}
  \left\|\widehat{\int_0^s e^{(s-r)\Delta}P_k\nabla\cdot(A\psi)(r)\,\mathrm{d}r}\right\|_{L_\xi^\infty}&\lesssim 2^k\int_0^s e^{-c(s-r)2^{2k}}\|\widehat{P_k(A\psi)(r)}\|_{L_\xi^\infty}\,\mathrm{d}r\\
  &\lesssim\epsilon_1^3(1+s2^{2k})^{-12}\Bigg[\sum_{j\leq k-5}2^{\frac12(j-k)}(r_jq_k+q_jr_k)+\sum_{j\geq k-4}2^{k-j}r_jq_j\Bigg]\\
  &\lesssim\epsilon_1^3r_kq_k(1+s2^{2k})^{-12}.
\end{aligned}
\end{equation}
Since $\sum_k r_kq_k\lesssim1$, \eqref{epsihat3} gives
\[
  \left\|(1+s|\xi|^2)^{12}\widehat{\int_0^s e^{(s-r)\Delta}\nabla\cdot(A\psi)(r)\,\mathrm{d}r}\right\|_{L_\xi^\infty}\lesssim\epsilon_1^3.
\]
For $F$ as in \eqref{eFdefinition}, Bernstein inequality, Young inequality, \eqref{eheatenergyimproved2}, \eqref{econnectionenve1}, and \eqref{edkdef2} give
\[
  \int_0^\infty\|(1+r|\xi|^2)^{12}\widehat F(r)\|_{L_\xi^1}\,\mathrm{d}r\lesssim\epsilon_1^2\sum_kd_k+\epsilon_1^2\sum_{j,k}2^{-|j-k|}b_jb_k\lesssim\epsilon_1^2.
\]
For $0\leq r\leq s$, the polynomial weight satisfies
\[
  e^{-(s-r)|\xi|^2}(1+s|\xi|^2)^{12}\lesssim e^{-(s-r)|\xi|^2}\bigl(1+(s-r)|\xi|^2\bigr)^{12}(1+r|\xi|^2)^{12}\lesssim(1+r|\xi|^2)^{12}.
\]
Together with
\[
  (1+r|\xi|^2)^{12}\lesssim(1+r|\eta|^2)^{12}(1+r|\xi-\eta|^2)^{12},
\]
the weighted Young inequality gives
\[
  \left\|(1+s|\xi|^2)^{12}\widehat{\int_0^s e^{(s-r)\Delta}(F\psi)(r)\,\mathrm{d}r}\right\|_{L_\xi^\infty}\lesssim\int_0^s\|(1+r|\xi|^2)^{12}\widehat F(r)\|_{L_\xi^1}\|(1+r|\xi|^2)^{12}\widehat\psi(r)\|_{L_\xi^\infty}\,\mathrm{d}r.
\]
The linear term satisfies
\[
  \|(1+s|\xi|^2)^{12}\widehat{e^{s\Delta}\psi(0)}\|_{L_\xi^\infty}\lesssim\|\widehat\psi(0)\|_{L_\xi^\infty}\lesssim\epsilon_1.
\]
Hence Duhamel formula implies
\[
  \|(1+s|\xi|^2)^{12}\widehat\psi(s)\|_{L_\xi^\infty}\lesssim\epsilon_1+\int_0^s\max_F\|(1+r|\xi|^2)^{12}\widehat F(r)\|_{L_\xi^1}\|(1+r|\xi|^2)^{12}\widehat\psi(r)\|_{L_\xi^\infty} \,\mathrm{d}r,
\]
where the maximum is taken over the finite collection in \eqref{eFdefinition}. By the integral bound above and Gronwall inequality,
\[
  \|(1+s|\xi|^2)^{12}\widehat\psi(s)\|_{L_\xi^\infty}\lesssim\epsilon_1e^{C\epsilon_1^2}\lesssim\epsilon_1,
\]
uniformly for $s\geq0$, proving \eqref{epsihat1}. Substituting this bound into the weighted convolution estimate and localizing to frequency $2^k$, we obtain
\begin{equation}\label{epsihat4}
  \left\|\widehat{\int_0^s e^{(s-r)\Delta}P_k(F\psi)(r)\,\mathrm{d}r}\right\|_{L_\xi^\infty}\lesssim\epsilon_1^3(1+s2^{2k})^{-12}.
\end{equation}
\end{proof}

\begin{lemma}\label{lApsiL1}
Under the bootstrap assumption \eqref{booassum}, one has
\begin{equation*}
  \int_0^\infty\|\widehat{A\psi}(s)\|_{L_\xi^1}\,\mathrm{d}s\lesssim\epsilon_1^3.
\end{equation*}
\end{lemma}

\begin{proof}
A direct paraproduct estimate gives
\begin{equation*}
  \|P_k(A\psi)\|_{L_x^2}\lesssim \epsilon_1^3 2^k\big[(b_k\sum_{j\leq k-5}2^{j-k}d_j+d_k\sum_{j\leq k-5}2^{j-k}b_j)(1+s2^{2k})^{-19}+\sum_{j\geq k-4}d_jb_j(1+s2^{2j})^{-39}\big],
\end{equation*}
hence
\[
\begin{aligned}
  &\int_0^\infty\|\widehat{A\psi}(s)\|_{L_\xi^1}\,\mathrm{d}s\\
  &\qquad\lesssim\sum_k2^k\int_0^\infty\|P_k(A\psi)(s)\|_{L_x^2}\,\mathrm{d}s\\
  &\qquad\lesssim\sum_k 2^k\int_0^\infty\epsilon_1^3 2^kd_kb_k(1+s2^{2k})^{-19}\,\mathrm{d}s+\sum_k2^k\int_0^\infty\epsilon_1^3 2^k\sum_{j\geq k-4}d_jb_j(1+s2^{2j})^{-39} \,\mathrm{d}s\\
  &\qquad\lesssim\epsilon_1^3\sum_kd_kb_k\lesssim\epsilon_1^3.
\end{aligned}
\]
\end{proof}

We now prove the central result of this section. The key is to employ the negative half order estimates to control the low frequency.
\begin{proposition}\label{pAhat}
Under the bootstrap assumption \eqref{booassum}, uniformly in the physical time $t$, we have
\begin{equation*}
  \sup_{s\geq0}\|\widehat A(s,t)\|_{L_\xi^\infty}\lesssim\epsilon_1^2.
\end{equation*}
In particular,
\begin{equation}\label{eAhat1}
  \|\widehat{A}(0,t)\|_{L_\xi^\infty}\lesssim\epsilon_1^2.
\end{equation}
\end{proposition}

\begin{proof}
We again use the formula \eqref{ecaloricconnection}. For the quadratic contribution, the negative half order estimates imply
\[
\begin{aligned}
  &\int_s^\infty\left\|\widehat{P_{k_1}\psi(r)\nabla P_{k_2}\psi(r)}\right\|_{L_\xi^\infty}\,\mathrm{d}r\\
  &\lesssim\epsilon_1^22^{\frac{k_1}{2}}2^{\frac{3k_2}{2}}q_{k_1}q_{k_2}\int_s^\infty(1+r2^{2k_1})^{-12}(1+r2^{2k_2})^{-12}\,\mathrm{d}r.
\end{aligned}
\]
The time integral is bounded by $2^{-2\max{\{k_1,k_2\}}}$, hence
\[
\begin{aligned}
  \sup_{s\geq0}\left\|\int_s^\infty\widehat{\psi(r)\nabla\psi(r)}\,\mathrm{d}r\right\|_{L_\xi^\infty}&\lesssim\epsilon_1^2\sum_{k_1\leq k_2}2^{\frac12(k_1-k_2)}q_{k_1}q_{k_2}+\epsilon_1^2\sum_{k_2<k_1}2^{-\frac32(k_1-k_2)}q_{k_1}q_{k_2}\\
  &\lesssim\epsilon_1^2\sum_kq_k^2\lesssim\epsilon_1^2.
\end{aligned}
\]

For the term containing $A\psi\psi$, Young inequality, Proposition~\ref{ppsihat}, and Lemma~\ref{lApsiL1} give
\begin{equation*}
  \sup_{s\geq0}\left\|\int_s^\infty\widehat{A\psi\psi}(r)\,\mathrm{d}r\right\|_{L_\xi^\infty}\lesssim\sup_{r\geq0}\|\widehat{\psi}(r)\|_{L_\xi^\infty}\int_0^\infty\|\widehat{A\psi}(r)\|_{L_\xi^1}\,\mathrm{d}r\lesssim\epsilon_1^4.
\end{equation*}
We finish the proof.
\end{proof}

\subsection{Decay estimates for $\psi$ and $A$}

\subsubsection{The \(L^\infty\) bootstrap along the heat flow}
To obtain the decay of the connection coefficients $A$ in physical time, we require some control over the $L^\infty$ norm of the differential field $\psi$.

We first derive preliminary heat decay estimates using Fourier control.
\begin{lemma}\label{lpsiAdecay} 
For every $0\leq a\leq\epsilon$, one has
\begin{align}
  \|\nabla^a\psi(s,t)\|_{L_x^2}&\lesssim\epsilon_1\langle s\rangle^{-\frac a2-\frac12},\label{epsiAdecay1}\\
  \|\nabla^aA(s,t)\|_{L_x^2}&\lesssim\epsilon_1^2\langle s\rangle^{-\frac a2-1}.\label{epsiAdecay2}
\end{align}
Moreover, \eqref{epsiAdecay1} and \eqref{epsiAdecay2} hold for $s\geq1$ and $0\leq a\leq3$.
\end{lemma}

\begin{proof} 
The case $0\leq s\leq1$ follows from \eqref{eheatenergyimproved2} and \eqref{econnectionenve1}. By \eqref{epsihat2} and Plancherel theorem,
\[
\|P_k\psi(s,t)\|_{L_x^2}\lesssim2^k\|\widehat{P_k\psi}(s,t)\|_{L_\xi^\infty}\lesssim\epsilon_1 2^k(1+s2^{2k})^{-12}.
\]
Hence, for $s\geq1$ and $0\leq a\leq4$,
\[
\|\nabla^a\psi(s,t)\|_{L_x^2}^2\lesssim\epsilon_1^2\sum_k2^{2(a+1)k}(1+s2^{2k})^{-24}\lesssim\epsilon_1^2s^{-a-1}.
\]
The same Fourier bound gives
\[
\|\psi(r,t)\|_{L_x^\infty}+r^{\frac12}\|\nabla\psi(r,t)\|_{L_x^\infty}\lesssim\epsilon_1\sum_k2^{2k}(1+r^{\frac12}2^k)(1+r2^{2k})^{-12}\lesssim\epsilon_1r^{-1}.
\]
By \eqref{econnectionenve1}, \eqref{edkdef2}, and \eqref{eenergyenvelope3},
\begin{equation}\label{epsiAdecay3}
\|A(r,t)\|_{L_x^\infty}\lesssim\epsilon_1^2\sum_k2^kd_k(1+r2^{2k})^{-19}\lesssim\epsilon_1^2r^{-\frac12}.
\end{equation}
Together with $\|\nabla^aA(r,t)\|_{L_x^2}\lesssim\epsilon_1^2r^{-\frac a2}$, the fractional Leibniz rule gives, for $r\geq1$ and $0\leq a\leq3$,
\begin{align*}
  \|\nabla^a(\overline\psi\nabla\psi)(r,t)\|_{L_x^2}&\lesssim\|\psi(r,t)\|_{L_x^\infty}\|\nabla^{a+1}\psi(r,t)\|_{L_x^2}+\|\nabla^a\psi(r,t)\|_{L_x^2}\|\nabla\psi(r,t)\|_{L_x^\infty}\lesssim\epsilon_1^2r^{-\frac a2-2},\\
  \|\nabla^a(\overline\psi A\psi)(r,t)\|_{L_x^2}&\lesssim\|\nabla^aA(r,t)\|_{L_x^2}\|\psi(r,t)\|_{L_x^\infty}^2+\|A(r,t)\|_{L_x^\infty}\|\psi(r,t)\|_{L_x^\infty}\|\nabla^a\psi(r,t)\|_{L_x^2}\lesssim\epsilon_1^4r^{-\frac a2-2}.
\end{align*}
Thus \eqref{ecaloricconnection} implies
\[
\|\nabla^aA(s,t)\|_{L_x^2}\lesssim\epsilon_1^2\int_s^\infty r^{-\frac a2-2}\,\mathrm{d}r\lesssim\epsilon_1^2s^{-\frac a2-1},\qquad s\geq1.
\]
This proves both estimates.
\end{proof}

The following conclusion is the main $L^\infty$ estimates for the differential field $\psi$. 
\begin{proposition}\label{psiinfty}
Under the bootstrap assumption \eqref{booassum}, uniformly for $s\geq0$ and $k\in\ZZ$, one has
\begin{align}
  \|\psi(s,t)\|_{L_x^\infty}&\lesssim\epsilon_1(t+s)^{-1}, \label{epsidecay1}\\
  \|P_{\leq k}\psi(s,t)\|_{L_x^\infty}&\lesssim\epsilon_1\min\left\{2^{2k},\frac{1}{t+s}\right\}, \label{epsidecay2}\\
  \|P_k\psi(s,t)\|_{L_x^\infty}&\lesssim\epsilon_1\min\left\{2^{2k},\frac{1}{t+s}\right\}(1+s2^{2k})^{-10}. \label{epsidecay3}
\end{align}
\end{proposition}

We again proceed with a bootstrap argument. For every fixed physical time $t$, define
\begin{align*}
  U_1&:=\sup_{s\geq0}(t+s)\|\psi(s,t)\|_{L_x^\infty},\\
  U_2&:=\sup_{s\geq0}\sup_{k\in\ZZ}2^{-2k}\|P_{\leq k}\psi(s,t)\|_{L_x^\infty},\\
  U_3&:=\sup_{s\geq0}\sup_{k\in\ZZ}\big[\min\left\{2^{2k},\frac{1}{t+s}\right\}\big]^{-1}(1+s2^{2k})^{10}\|P_k\psi(s,t)\|_{L_x^\infty}.
\end{align*}
It suffices to show that the bootstrap assumption
\begin{equation}\label{eU-bootstrap}
  U_1+U_2+U_3\leq C\epsilon_1
\end{equation}
improves to
\begin{equation}\label{eU-improvement}
  U_1+U_2+U_3\leq \frac{3}{4}C\epsilon_1,
\end{equation}
where $C$ is first chosen sufficiently large and $\epsilon_1$ is then chosen sufficiently small.

\begin{remark}
In the subsequent estimates, we omit the dependence on the fixed bootstrap constant $C$. Before closing \eqref{eU-bootstrap}, references to \eqref{epsidecay1}, \eqref{epsidecay2}, and \eqref{epsidecay3} always mean their bootstrap versions, which follow directly from \eqref{eU-bootstrap}. The dependence on $C$ is kept when closing the bootstrap, while the Duhamel terms contain at least an additional factor $\epsilon_1^2$.
\end{remark}

\begin{corollary}\label{cpsiintegrability}
Under \eqref{eU-bootstrap}, uniformly for $s\geq0$, one has
\begin{equation}\label{epsiinte1}
\int_s^\infty\|\nabla\psi(r,t)\|_{L_x^\infty}\,\mathrm{d}r\lesssim\epsilon_1(t+s)^{-\frac12}.
\end{equation}
In particular,
\begin{equation}\label{epsiinte2}
\int_0^\infty\|\nabla\psi(s,t)\|_{L_x^\infty}\,\mathrm{d}s\lesssim\epsilon_1t^{-\frac12}.
\end{equation}
Moreover,
\begin{equation}\label{epsiinte3}
\int_s^\infty\|\psi(r,t)\|_{L_x^\infty}^2\,\mathrm{d}r\lesssim\epsilon_1^2(t+s)^{-1}.
\end{equation}
In particular,
\begin{equation}\label{epsiinte4}
\int_0^\infty\|\psi(r,t)\|_{L_x^\infty}^2\,\mathrm{d}r\lesssim\epsilon_1^2t^{-1}.
\end{equation}
Finally,
\begin{equation}\label{epsiinte5}
\int_0^\infty\|\nabla^2\psi(s,t)\|_{L_x^2}\,\mathrm{d}s\lesssim\epsilon_1.
\end{equation}
\end{corollary}
\begin{proof}
By Bernstein inequality and \eqref{epsidecay3}, for $r>0$,
\begin{equation}\label{epsiinte6}
\|\nabla\psi(r,t)\|_{L_x^\infty}\lesssim\epsilon_1(t+r)^{-1}\sum_k2^k(1+r2^{2k})^{-10}\lesssim\epsilon_1(t+r)^{-1}r^{-\frac12}.
\end{equation}
If $0\leq s\leq t$, then
\begin{equation*}
\int_s^\infty r^{-\frac12}(t+r)^{-1}\,\mathrm{d}r\leq t^{-1}\int_0^t r^{-\frac12}\,\mathrm{d}r+\int_t^\infty r^{-\frac32}\,\mathrm{d}r\lesssim t^{-\frac12}.
\end{equation*}
If $s>t$, then
\begin{equation*}
\int_s^\infty r^{-\frac12}(t+r)^{-1}\,\mathrm{d}r\leq\int_s^\infty r^{-\frac32}\,\mathrm{d}r\lesssim s^{-\frac12}.
\end{equation*}
These estimates prove \eqref{epsiinte1}. Taking $s=0$ gives \eqref{epsiinte2}. By \eqref{epsidecay1},
\begin{equation*}
\int_s^\infty\|\psi(r,t)\|_{L_x^\infty}^2\,\mathrm{d}r\lesssim\epsilon_1^2\int_s^\infty(t+r)^{-2}\,\mathrm{d}r\lesssim\epsilon_1^2(t+s)^{-1},
\end{equation*}
which proves \eqref{epsiinte3}. Taking $s=0$ gives \eqref{epsiinte4}. For $0<s\leq1$, \eqref{eheatenergyimproved2} with $\sigma=\epsilon$ and \eqref{eenergyenvelope3} give
\begin{equation*}
\|\nabla^2\psi(s,t)\|_{L_x^2}^2\lesssim\epsilon_1^2\sum_k2^{2(2-\epsilon)k}b_k(\epsilon)^2(1+s2^{2k})^{-40}\lesssim\epsilon_1^2s^{-2+\epsilon}.
\end{equation*}
For $s\geq1$, \eqref{epsiAdecay1} with $a=2$ gives $\|\nabla^2\psi(s,t)\|_{L_x^2}\lesssim\epsilon_1s^{-\frac32}$. Therefore,
\begin{equation*}
\int_0^\infty\|\nabla^2\psi(s,t)\|_{L_x^2}\,\mathrm{d}s\lesssim\epsilon_1\int_0^1s^{-1+\frac\epsilon2}\,\mathrm{d}s+\epsilon_1\int_1^\infty s^{-\frac32}\,\mathrm{d}s\lesssim\epsilon_1.
\end{equation*}
This proves \eqref{epsiinte5}.
\end{proof}

\subsubsection{Estimates for the connection coefficients}

\begin{lemma}\label{Ainfty}
Under \eqref{eU-bootstrap}, one has
\begin{align}
  \|A(s,t)\|_{L_x^\infty}&\lesssim\epsilon_1^2(t+s)^{-\frac32}, \label{eAinfty1}\\
  \|P_{\leq k}A(s,t)\|_{L_x^\infty}&\lesssim\epsilon_1^2\min\left\{2^{2k},(t+s)^{-\frac32}\right\}. \label{eAinfty2}
\end{align}
\end{lemma}

\begin{proof}
To prove \eqref{eAinfty1}, we first consider the quadratic part. From \eqref{epsiinte1} and \eqref{epsidecay1},
\begin{equation*}
\left\|\int_s^\infty\psi(r,t)\nabla\psi(r,t)\,\mathrm{d}r\right\|_{L_x^\infty}\lesssim\sup_{r\geq s}\|\psi(r,t)\|_{L_x^\infty}\int_s^\infty\|\nabla\psi(r,t)\|_{L_x^\infty} \,\mathrm{d}r\lesssim\epsilon_1^2(t+s)^{-\frac32}.
\end{equation*}
By \eqref{epsiAdecay3}, the remainder term satisfies
\begin{align*}
  \left\|\int_s^\infty A(r,t)\psi(r,t)\psi(r,t)\,\mathrm{d}r\right\|_{L_x^\infty}&\lesssim\int_s^\infty\|A(r,t)\|_{L_x^\infty}\|\psi(r,t)\|_{L_x^\infty}^2\,\mathrm{d}r\\
  &\lesssim\epsilon_1^4\int_s^\infty r^{-\frac12}(t+r)^{-2}\,\mathrm{d}r.
\end{align*}
For $0\leq s\leq t$,
\[
  \int_s^\infty r^{-\frac12}(t+r)^{-2}\,\mathrm{d}r\leq t^{-2}\int_0^tr^{-\frac12}\,\mathrm{d}r+\int_t^\infty r^{-\frac52}\,\mathrm{d}r\lesssim t^{-\frac32},
\]
whereas, for $s>t$,
\[
  \int_s^\infty r^{-\frac12}(t+r)^{-2}\,\mathrm{d}r\leq\int_s^\infty r^{-\frac52}\,\mathrm{d}r\lesssim s^{-\frac32}.
\]
Combining these estimates proves \eqref{eAinfty1}. On the other hand, Proposition \ref{pAhat} gives
\begin{equation*}
  \|P_{\leq k}A(s,t)\|_{L_x^\infty}\lesssim 2^{2k}\|\hat{A}(s,t)\|_{L_\xi^\infty}\lesssim\epsilon_1^2 2^{2k}.
\end{equation*}
Taking the minimum of it and \eqref{eAinfty1} proves \eqref{eAinfty2}.
\end{proof}

\begin{remark}\label{rmk:Alinfty}
The decay in \eqref{eAinfty1} can be understood from the quadratic caloric structure, which has the schematic dyadic size
\begin{align*}
  \sum_k2^{-k}\|P_k\psi(s,t)\|_{L_x^\infty}\|P_{\leq k+4}\psi(s,t)\|_{L_x^\infty}&\lesssim\epsilon_1^2\sum_k2^{-k}\min\left\{2^{2k},(t+s)^{-1}\right\}^2\\
  &\lesssim\epsilon_1^2(t+s)^{-\frac32}.
\end{align*}
Here we used \eqref{epsidecay2} and \eqref{epsidecay3}, discarding the additional heat weights. The improved low frequency factor originates from the Fourier bound at the initial heat time:
\[
  \|P_{\leq k}\psi(0,t)\|_{L_x^\infty}\lesssim 2^{2k}\|\widehat{\psi}(0,t)\|_{L_\xi^\infty}\lesssim\epsilon_1 2^{2k},
\]
whereas the energy bound alone gives
\[
  \|P_{\leq k}\psi(0,t)\|_{L_x^\infty}\lesssim2^k\|\psi(0,t)\|_{L_x^2}\lesssim\epsilon_1 2^k.
\]
Using only the latter frequency factor in the same calculation yields
\[
  \epsilon_1^2\sum_k2^{-k}\min\left\{2^k,(t+s)^{-1}\right\}^2\lesssim\epsilon_1^2(t+s)^{-1}.
\]
Thus the Fourier control gives an additional factor $2^k$ at low frequencies. Its propagation along the heat flow improves the summation against the negative derivative  $2^{-k}$, yielding the better decay $(t+s)^{-\frac32}$ in place of $(t+s)^{-1}$.
\end{remark}

\begin{proposition}\label{pALinftypk}
Under \eqref{eU-bootstrap}, one has
\begin{equation}
  \|P_kA(s,t)\|_{L_x^\infty}\lesssim\epsilon_1^2\min\left\{2^{2k},(t+s)^{-\frac32}\right\}(1+s2^{2k})^{-9}. \label{eAinfty3}
\end{equation}
\end{proposition}

\begin{proof}
We first establish the bound with decay $(t+s)^{-\frac{3}{2}}$ and heat weight $(1+s2^{2k})^{-9}$. For the quadratic term, using \eqref{epsidecay2} and \eqref{epsidecay3} in the low-high, high-low, and high-high interactions gives
\begin{align*}
  &\left\|P_k\int_s^\infty\psi(r,t)\nabla\psi(r,t)\,\mathrm{d}r\right\|_{L_x^\infty}\\
  &\quad\lesssim\epsilon_1^2\int_s^\infty(t+r)^{-2}\left[2^k(1+r2^{2k})^{-10}+\sum_{j\geq k-4}2^j(1+r2^{2j})^{-20}\right]\,\mathrm{d}r\\
  &\quad\lesssim\epsilon_1^2\int_s^\infty r^{-\frac12}(t+r)^{-2}(1+r2^{2k})^{-9}\,\mathrm{d}r\\
  &\quad\lesssim\epsilon_1^2(t+s)^{-\frac32}(1+s2^{2k})^{-9},
\end{align*}
the heat integral was estimated in the proof of \eqref{eAinfty1}.

For the remainder term, first assume $2^{2k}\leq(t+s)^{-1}$. Since $s2^{2k}\leq1$, by \eqref{eAinfty1} and \eqref{epsidecay1},
\[
  \left\|P_k\int_s^\infty A\psi\psi(r,t)\,\mathrm{d}r\right\|_{L_x^\infty}\lesssim\epsilon_1^4\int_s^\infty(t+r)^{-\frac72}\,\mathrm{d}r\lesssim\epsilon_1^4(t+s)^{-\frac32}(1+s2^{2k})^{-9}.
\]
Now assume $2^{2k}>(t+s)^{-1}$, we distinguish the following four cases.

\smallskip
\noindent\textit{One highest frequency $\psi$ has frequency $2^{k+O(1)}$.}
By \eqref{eAinfty1}, \eqref{epsidecay2}, and \eqref{epsidecay3},
\begin{align*}
  &\left\|P_k\int_s^\infty P_{k+O(1)}\psi\,P_{\leq k+4}\psi\,P_{\leq k+4}A(r,t)\,\mathrm{d}r\right\|_{L_x^\infty}\\
  &\quad\lesssim\epsilon_1^4\int_s^\infty(t+r)^{-\frac72}(1+r2^{2k})^{-10}\,\mathrm{d}r\\
  &\quad\lesssim\epsilon_1^4(t+s)^{-\frac52}(1+s2^{2k})^{-10}\\
  &\quad\lesssim\epsilon_1^4(t+s)^{-\frac32}(1+s2^{2k})^{-9}.
\end{align*}

\smallskip
\noindent\textit{$A$ has frequency $2^{k+O(1)}$ and both $\psi$ factors have lower frequencies.}
By Bernstein inequality, \eqref{econnectionenve1}, \eqref{epsidecay2}, and $2^{-k}<(t+s)^{\frac{1}{2}}$,
\begin{align*}
  &\left\|P_k\int_s^\infty P_{k+O(1)}A\,P_{\leq k-5}\psi\,P_{\leq k-5}\psi(r,t)\,\mathrm{d}r\right\|_{L_x^\infty}\\
  &\quad\lesssim\epsilon_1^4 2^kd_k\int_s^\infty(t+r)^{-2}(1+r2^{2k})^{-19}\,\mathrm{d}r\\
  &\quad\lesssim\epsilon_1^4(t+s)^{-2}2^{-k}d_k(1+s2^{2k})^{-18}\\
  &\quad\lesssim\epsilon_1^4(t+s)^{-\frac32}(1+s2^{2k})^{-9}.
\end{align*}

\smallskip
\noindent\textit{$A$ and one $\psi$ factor form the high frequency pair.}
By Bernstein inequality, \eqref{econnectionenve1}, \eqref{epsidecay2}, and \eqref{epsidecay3},
\begin{align*}
  &\left\|P_k\int_s^\infty\sum_{j\geq k-4}P_jA\,P_{j+O(1)}\psi\,P_{\leq j+4}\psi(r,t)\,\mathrm{d}r\right\|_{L_x^\infty}\\
  &\quad\lesssim\epsilon_1^4 2^k\sum_{j\geq k-4}d_j\int_s^\infty(t+r)^{-2}(1+r2^{2j})^{-29}\,\mathrm{d}r\\
  &\quad\lesssim\epsilon_1^4(t+s)^{-2}\sum_{j\geq k-4}2^{k-2j}d_j(1+s2^{2j})^{-28}\\
  &\quad\lesssim\epsilon_1^4(t+s)^{-2}2^{-k}(1+s2^{2k})^{-28}\\
  &\quad\lesssim\epsilon_1^4(t+s)^{-\frac32}(1+s2^{2k})^{-9}.
\end{align*}

\smallskip
\noindent\textit{The two $\psi$ factors form the high frequency pair.}
By Bernstein inequality, \eqref{eheatenergyimproved2}, \eqref{epsidecay3}, and \eqref{eAinfty1},
\begin{align*}
  &\left\|P_k\int_s^\infty\sum_{j\geq k-4}P_j\psi\,P_{j+O(1)}\psi\,P_{\leq j+4}A(r,t)\,\mathrm{d}r\right\|_{L_x^\infty}\\
  &\quad\lesssim\epsilon_1^42^k\sum_{j\geq k-4}b_j\int_s^\infty(t+r)^{-\frac52}(1+r2^{2j})^{-30}\,\mathrm{d}r\\
  &\quad\lesssim\epsilon_1^4(t+s)^{-\frac52}\sum_{j\geq k-4}2^{k-2j}b_j(1+s2^{2j})^{-29}\\
  &\quad\lesssim\epsilon_1^4(t+s)^{-\frac52}2^{-k}(1+s2^{2k})^{-29}\\
  &\quad\lesssim\epsilon_1^4(t+s)^{-\frac32}(1+s2^{2k})^{-9}.
\end{align*}

Combining the quadratic and remainder estimates with \eqref{ecaloricconnection}, we obtain
\[
  \|P_kA(s,t)\|_{L_x^\infty}\lesssim\epsilon_1^2(t+s)^{-\frac32}(1+s2^{2k})^{-9}.
\]
Finally, if $2^{2k}\leq(t+s)^{-\frac32}$, then \(s2^{2k}\leq 1\), \eqref{eAinfty2} gives
\[
  \|P_kA(s,t)\|_{L_x^\infty}\lesssim\epsilon_1^22^{2k}\lesssim\epsilon_1^2 2^{2k}(1+s2^{2k})^{-9}.
\]
If $2^{2k}>(t+s)^{-\frac32}$, the time decay bound above applies directly. These two cases prove \eqref{eAinfty3}.
\end{proof}

We also need the following \(L^2\) estimate.

\begin{lemma}\label{lAL2}
Under \eqref{eU-bootstrap},
\begin{equation}\label{eAL2}
  \|A(s,t)\|_{L_x^2}\lesssim\epsilon_1^2(t+s)^{-1}(1+\log t).
\end{equation}
\end{lemma}

\begin{proof}
By \eqref{epsiAdecay1} and \eqref{epsiinte6}:
\[
  \left\|\int_s^\infty\psi(r,t)\nabla\psi(r,t)\,\mathrm{d}r\right\|_{L_x^2}\lesssim\epsilon_1^2\int_s^\infty\langle r\rangle^{-\frac12}r^{-\frac12}(t+r)^{-1} \,\mathrm{d}r.
\]
For $0\leq s\leq t$,
\[
  \int_s^\infty \langle r\rangle^{-\frac12}r^{-\frac12}(t+r)^{-1}\,\mathrm{d}r\lesssim t^{-1}\int_0^1r^{-\frac12}\,\mathrm{d}r+t^{-1} \int_1^{t}r^{-1}\,\mathrm{d}r+\int_{t}^\infty r^{-2}\,\mathrm{d}r\lesssim t^{-1}(1+\log t).
\]
For $s>t$,
\[
  \int_s^\infty \langle r\rangle^{-\frac12}r^{-\frac12}(t+r)^{-1}\,\mathrm{d}r\lesssim\int_s^\infty r^{-2}\,\mathrm{d}r\lesssim s^{-1}.
\]
Hence,
\begin{equation}\label{eAL21}
\left\|\int_s^\infty\psi(r,t)\nabla\psi(r,t)\,\mathrm{d}r\right\|_{L_x^2}\lesssim\epsilon_1^2(t+s)^{-1}(1+\log t).
\end{equation}
For the remainder term, by \eqref{econnectionenve1} and \eqref{epsiinte3},
\begin{align}
  &\left\|\int_s^\infty A(r,t)\psi(r,t)\psi(r,t)\,\mathrm{d}r\right\|_{L_x^2}\lesssim\int_s^\infty\|A(r,t)\|_{L_x^2}\|\psi(r,t)\|_{L_x^\infty}^2\,\mathrm{d}r\notag\\
  &\quad\lesssim\sup_{r\geq 0}\|A(r,t)\|_{L_x^2}\int_s^\infty\|\psi(r,t)\|_{L_x^\infty}^2\,\mathrm{d}r\lesssim\epsilon_1^4(t+s)^{-1}. \label{eAL22}
\end{align}
Combining \eqref{eAL21} and \eqref{eAL22} proves \eqref{eAL2}.
\end{proof}

\begin{remark}\label{rmk:Al2}
The quadratic caloric structure gives the bound
\begin{align*}
  \sum_k2^{-k}\|P_k\psi(s,t)\|_{L_x^2}\|P_{\leq k+4}\psi(s,t)\|_{L_x^\infty}&\lesssim\epsilon_1^2\sum_k2^{-k}\min\left\{2^k,\langle s\rangle^{-\frac12}\right\}\min\left\{2^{2k},(t+s)^{-1}\right\}\\
  &\lesssim\epsilon_1^2\sum_k
  \begin{cases}
    2^{2k}, & 2^k\leq(t+s)^{-\frac12},\\
    (t+s)^{-1}, &(t+s)^{-\frac12}<2^k\leq\langle s\rangle^{-\frac12},\\
    2^{-k}\langle s\rangle^{-\frac12}(t+s)^{-1}, &2^k>\langle s\rangle^{-\frac12}
  \end{cases}\\
  &\lesssim\epsilon_1^2(t+s)^{-1}(1+\log t),
\end{align*}
consistently with \eqref{eAL2}.
\end{remark}

\subsubsection{The free heat evolution}
Recall the bootstrap assumption \eqref{booassum}
\[
  \|\psi(0,t)\|_{L_x^\infty}\lesssim\epsilon_1t^{-1}, \qquad\|\widehat{\psi}(0,t)\|_{L_\xi^\infty}\lesssim\epsilon_1.
\]
The \(L^\infty\)-contractivity of the heat semigroup implies
\[
  \|e^{s\Delta}\psi(0,t)\|_{L_x^\infty}\lesssim\epsilon_1t^{-1}.
\]
For $s>0$, the Fourier estimate gives
\[
  \|e^{s\Delta}\psi(0,t)\|_{L_x^\infty}\lesssim\epsilon_1\int_{\RR^2}e^{-s|\xi|^2}\,\mathrm{d}\xi\lesssim\epsilon_1s^{-1}.
\]
Taking the minimum,
\begin{equation}\label{eUfree1}
  \|e^{s\Delta}\psi(0,t)\|_{L_x^\infty}\lesssim\frac{\epsilon_1}{t+s}.
\end{equation}
Similarly,
\begin{equation*}
  \|P_{\leq k}e^{s\Delta}\psi(0,t)\|_{L_x^\infty}\lesssim 2^{2k}\|\hat\psi(0,t)\|_{L_\xi^\infty}\lesssim\epsilon_1 2^{2k}.
\end{equation*}
At a fixed frequency,
\begin{align}
  \|P_ke^{s\Delta}\psi(0,t)\|_{L_x^\infty}&\lesssim e^{-cs2^{2k}}\min\left\{\epsilon_1 2^{2k},\epsilon_1t^{-1}\right\}\notag\\
  &\lesssim\epsilon_1\min\left\{2^{2k},\frac{1}{t+s}\right\}(1+s2^{2k})^{-10}, \label{eUfree2}
\end{align}
where the last inequality follows by separating the cases
\[
  (t+s)2^{2k}\leq1 \qquad\text{and}\qquad (t+s)2^{2k}>1,
\]
and using the exponential decay to dominate any fixed polynomial when \(s2^{2k}\gtrsim 1\).

\subsubsection{Improvement of \(U_1\)}

From \eqref{epsidecay1} and \eqref{eAinfty1},
\begin{equation}\label{eUimpro1}
  \|A(r,t)\psi(r,t)\|_{L_x^\infty}\lesssim\epsilon_1^3(t+r)^{-\frac52}.
\end{equation}
Moreover, by \eqref{eAL2} and \eqref{epsiAdecay1},
\begin{equation}\label{eUimpro2}
  \|A(r,t)\psi(r,t)\|_{L_x^1}\leq\|A(r,t)\|_{L_x^2}\|\psi(r,t)\|_{L_x^2}\lesssim\epsilon_1^3(t+r)^{-1}\langle r\rangle^{-\frac12}(1+\log t).
\end{equation}
If $s\leq t$, by \eqref{eheatoperator} $L_x^\infty\to L_x^\infty$ and \eqref{eUimpro1},
\[
  \left\|\int_0^s e^{(s-r)\Delta}\nabla\cdot(A\psi)(r,t)\,\mathrm{d}r\right\|_{L_x^\infty}\lesssim\epsilon_1^3\int_0^s(s-r)^{-\frac12}(t+r)^{-\frac52}\,\mathrm{d}r\lesssim\frac{\epsilon_1^3}{t+s}.
\]
If $s\geq t$, on $\frac{s}{2}\leq r\leq s$, we again use \eqref{eUimpro1}:
\begin{equation*}
  \int_{\frac{s}{2}}^s(s-r)^{-\frac12}\|A(r,t)\psi(r,t)\|_{L_x^\infty}\,\mathrm{d}r\lesssim\epsilon_1^3s^{\frac12}(t+s)^{-\frac52}\lesssim\frac{\epsilon_1^3}{t+s},
\end{equation*}
on $0\leq r\leq \frac{s}{2}$, we use \eqref{eheatoperator} $L_x^1\to L_x^\infty$ and \eqref{eUimpro2}:
\[
  \int_0^{\frac{s}{2}}(s-r)^{-\frac32}\|A(r,t)\psi(r,t)\|_{L_x^1}\,\mathrm{d}r\lesssim\epsilon_1^3s^{-\frac32}(1+\log t)\int_0^{\frac{s}{2}}\frac{\langle r\rangle^{-\frac12}}{t+r}\,\mathrm{d}r\lesssim\frac{\epsilon_1^3}{t+s}.
\]
Thus
\begin{equation*}
  \left\|\int_0^s e^{(s-r)\Delta}\nabla\cdot(A\psi)(r,t)\,\mathrm{d}r\right\|_{L_x^\infty}\lesssim\frac{\epsilon_1^3}{t+s}.
\end{equation*}
For the remaining nonlinearities, \eqref{eFenve1} with $\sigma=0$ and \eqref{eenergyenvelope3} give $\max_F\|F(r,t)\|_{L_x^2}\lesssim\epsilon_1^2r^{-\frac12}$ for $0<r\leq1$. For $r\geq1$, Gagliardo-Nirenberg inequality, \eqref{epsiAdecay1}, and \eqref{epsiAdecay2} give
\[
  \max_F\|F(r,t)\|_{L_x^2}\lesssim\|\nabla A(r,t)\|_{L_x^2}+\|\psi(r,t)\|_{L_x^4}^2+\|A(r,t)\|_{L_x^4}^2\lesssim\epsilon_1^2r^{-\frac32}.
\]
Hence, by \eqref{epsiAdecay1} and \eqref{epsidecay1},
\begin{align}
  \|(F\psi)(r,t)\|_{L_x^1}&\lesssim\epsilon_1^3r^{-\frac12}\langle r\rangle^{-\frac32},\notag\\
  \|(F\psi)(r,t)\|_{L_x^2}&\lesssim\epsilon_1^3r^{-\frac12}\langle r\rangle^{-1}(t+r)^{-1}.\label{eUimpro3}
\end{align}
For $s\leq t$, \eqref{eheatoperator} gives
\[
  \left\|\int_0^s e^{(s-r)\Delta}(F\psi)(r,t)\,\mathrm{d}r\right\|_{L_x^\infty}\lesssim\epsilon_1^3t^{-1}\int_0^s(s-r)^{-\frac12}r^{-\frac12}\langle r\rangle^{-1}\,\mathrm{d}r\lesssim\frac{\epsilon_1^3}{t+s}.
\]
For $s\geq t$, applying \eqref{eheatoperator} $L_x^1\to L_x^\infty$ on $[0,\frac{s}{2}]$ and $L_x^2\to L_x^\infty$ on $[\frac{s}{2},s]$, we obtain
\[
  \left\|\int_0^s e^{(s-r)\Delta}(F\psi)(r,t)\,\mathrm{d}r\right\|_{L_x^\infty}\lesssim\epsilon_1^3s^{-1}\int_0^\infty r^{-\frac12}\langle r\rangle^{-\frac32}\,\mathrm{d}r+\epsilon_1^3s^{-1}(t+s)^{-1}\lesssim \frac{\epsilon_1^3}{t+s}.
\]
Together with \eqref{eUfree1}, and restoring the dependence on $C$, this gives
\[
  U_1\lesssim\epsilon_1+C^3\epsilon_1^3.
\]

\subsubsection{Improvement of \(U_2\)}

By \eqref{epsihat1}, uniformly for $s\geq0$ and $k\in\ZZ$,
\[
  \|P_{\leq k}\psi(s,t)\|_{L_x^\infty}\lesssim 2^{2k}\|\widehat\psi(s,t)\|_{L_\xi^\infty}\lesssim\epsilon_12^{2k},
\]
then we have
\[
  U_2\lesssim\epsilon_1.
\]

\subsubsection{Improvement of \(U_3\)}
We first estimate the term containing $\nabla\cdot(A\psi)$ by the paraproduct decomposition, as in the proof of Proposition~\ref{propenv}.

\smallskip
\noindent\textit{Low-High interaction:}
By \eqref{eAinfty2} and \eqref{epsidecay3},
\begin{align}
    &\left\|\int_0^sP_ke^{(s-r)\Delta}\nabla\cdot\bigl(P_{\leq k-5}A(r,t)P_{k+O(1)}\psi(r,t)\bigr)\,\mathrm{d}r\right\|_{L_x^\infty}\notag\\
    &\quad\lesssim 2^k\int_0^se^{-c(s-r)2^{2k}}\|P_{\leq k-5}A(r,t)\|_{L_x^\infty}\|P_{k+O(1)}\psi(r,t)\|_{L_x^\infty}\,\mathrm{d}r\notag\\
    &\quad\lesssim\epsilon_1^3 2^k\int_0^se^{-c(s-r)2^{2k}}\min\left\{2^{2k},\frac{1}{t+r}\right\}^{2}(1+r2^{2k})^{-10}\,\mathrm{d}r\notag\\
    &\quad\lesssim\epsilon_1^3 2^{2k}\int_0^se^{-c(s-r)2^{2k}}\min\left\{2^{2k},\frac{1}{t+r}\right\}(1+r2^{2k})^{-10}\,\mathrm{d}r\notag\\
    &\quad\lesssim\epsilon_1^3\min\left\{2^{2k},\frac{1}{t+s}\right\}(1+s2^{2k})^{-10}, \label{eUimpro4}
\end{align}
Here we used a refined heat convolution estimate
\begin{equation}\label{eUimpro5}
    2^{2k}\int_0^se^{-c(s-r)2^{2k}}\min\left\{2^{2k},\frac{1}{t+r}\right\}(1+r2^{2k})^{-10}\,\mathrm{d}r\lesssim\min\left\{2^{2k},\frac{1}{t+s}\right\}(1+s2^{2k})^{-10}. 
\end{equation}
To verify this estimate, observe that, for \(0\leq r\leq s\),
\[
    \min\left\{2^{2k},\frac{1}{t+r}\right\}\lesssim\bigl(1+(s-r)2^{2k}\bigr)\min\left\{2^{2k},\frac{1}{t+s}\right\},
\]
and absorb the additional factor into the exponential kernel. Then \eqref{eUimpro5} follows from Lemma~\ref{lpolynomialconvolution}.

\smallskip
\noindent\textit{High-Low interaction:}
By Bernstein inequality, \eqref{econnectionenve1}, and \eqref{epsidecay2},
\begin{align}
    &\left\|\int_0^sP_ke^{(s-r)\Delta}\nabla\cdot\bigl(P_{k+O(1)}A(r,t)P_{\leq k-5}\psi(r,t)\bigr)\,\mathrm{d}r\right\|_{L_x^\infty}\notag\\
    &\quad\lesssim 2^k\int_0^s e^{-c(s-r)2^{2k}}\|P_{k+O(1)}A(r,t)\|_{L_x^\infty}\|P_{\leq k-5}\psi(r,t)\|_{L_x^\infty}\,\mathrm{d}r\notag\\
    &\quad\lesssim\epsilon_1^32^{2k}d_k\int_0^s e^{-c(s-r)2^{2k}}\min\left\{2^{2k},\frac{1}{t+r}\right\}(1+r2^{2k})^{-19}\,\mathrm{d}r\notag\\
    &\quad\lesssim\epsilon_1^3\min\left\{2^{2k},\frac{1}{t+s}\right\}(1+s2^{2k})^{-10}, \label{eUimpro6}
\end{align}
here $d_k\lesssim 1$ and the last line follows from \eqref{eUimpro5}.

\smallskip
\noindent\textit{High-High interaction:}
By Bernstein inequality, \eqref{econnectionenve1}, and \eqref{epsidecay3}, it follows that
\[
\begin{aligned}
    &\sum_{j\geq k-4}\left\|\int_0^sP_ke^{(s-r)\Delta}\nabla\cdot\bigl(P_jA(r,t)P_{j+O(1)}\psi(r,t)\bigr)\,\mathrm{d}r\right\|_{L_x^\infty}\\
    &\quad\lesssim\epsilon_1^3 2^{2k}\sum_{j\geq k-4}d_j\int_0^se^{-c(s-r)2^{2k}}\min\left\{2^{2j},\frac{1}{t+r}\right\}(1+r2^{2j})^{-29}\,\mathrm{d}r.
\end{aligned}
\]
For \(j\geq k-4\),
\[
  2^{2(k-j)}\min\left\{2^{2j},\frac{1}{t+r}\right\}\lesssim\min\left\{2^{2k},\frac{1}{t+r}\right\}.
\]
As in the proof of \eqref{eUimpro5}, for $j\geq k-4$, Lemma~\ref{lpolynomialconvolution} gives
\[
    2^{2j}\int_0^s e^{-c(s-r)2^{2k}}\min\left\{2^{2k},\frac{1}{t+r}\right\}(1+r2^{2j})^{-10}\,\mathrm{d}r\lesssim\min\left\{2^{2k},\frac{1}{t+s}\right\}(1+s2^{2k})^{-10}.
\]
Hence, by \eqref{edkdef2} and \eqref{eenergyenvelope3},
\begin{align}
    &\sum_{j\geq k-4}\left\|\int_0^sP_ke^{(s-r)\Delta}\nabla\cdot\bigl(P_jA(r,t)P_{j+O(1)}\psi(r,t)\bigr)\,\mathrm{d}r\right\|_{L_x^\infty}\notag\\
    &\quad\lesssim\epsilon_1^3\min\left\{2^{2k},\frac{1}{t+s}\right\}(1+s2^{2k})^{-10}\sum_jd_j\notag\\
    &\quad\lesssim\epsilon_1^3\min\left\{2^{2k},\frac{1}{t+s}\right\}(1+s2^{2k})^{-10}. \label{eUimpro7}
\end{align}

Combining \eqref{eUimpro4}, \eqref{eUimpro6}, and \eqref{eUimpro7}, we conclude that
\begin{equation*}
    \left\|\int_0^sP_ke^{(s-r)\Delta}\nabla\cdot\bigl(A(r,t)\psi(r,t)\bigr)\,\mathrm{d}r\right\|_{L_x^\infty}\lesssim\epsilon_1^3\min\left\{2^{2k},\frac{1}{t+s}\right\}(1+s2^{2k})^{-10}.
\end{equation*}
For the remaining nonlinearities, the paraproduct decomposition, \eqref{eFenve1} with $\sigma=0$, \eqref{epsidecay2} and \eqref{epsidecay3} give
\[
  \|P_k(F\psi)(r,t)\|_{L_x^\infty}\lesssim\epsilon_1^3r^{-\frac{1}{2}}(t+r)^{-1}2^k(1+r2^{2k})^{-10}+\epsilon_1^3r^{-\frac{1}{2}}(t+r)^{-1}2^k\sum_{j\geq k-4} b_j(1+r2^{2j})^{-29}.
\]
For $s\leq t$, \eqref{esingularheatconvolution} gives the required bound with $(t+s)^{-1}$. For $s\geq t$, the same argument applies on $\frac{s}{2}\leq r\leq s$, while \eqref{eUimpro3} gives
\[
  \left\|\int_0^{\frac{s}{2}}e^{(s-r)\Delta}P_k(F\psi)(r,t)\,\mathrm{d}r\right\|_{L_x^\infty}\lesssim\epsilon_1^3 2^{2k}e^{-\frac{c}{2}s2^{2k}}\lesssim\frac{\epsilon_1^3}{t+s}(1+s2^{2k})^{-10}.
\]
Combining this with \eqref{epsihat4}, we obtain
\[
  \left\|\int_0^s e^{(s-r)\Delta}P_k(F\psi)(r,t)\,\mathrm{d}r\right\|_{L_x^\infty}\lesssim\epsilon_1^3\min\left\{2^{2k},\frac{1}{t+s}\right\}(1+s2^{2k})^{-10}.
\]
Together with \eqref{eUfree2}, these estimates give
\[
    U_3\lesssim\epsilon_1+C^3\epsilon_1^3.
\]

\begin{proof}[Proof of Proposition~\ref{psiinfty}]
Combining the bounds for $U_1$, $U_2$, and $U_3$, we obtain
\[
    U_1+U_2+U_3\lesssim\epsilon_1+C^3\epsilon_1^3.
\]
We first choose $C$ sufficiently large to dominate the implicit absolute constant and then choose $\epsilon_1$ sufficiently small depending on this fixed $C$. It follows that
\[
    U_1+U_2+U_3\leq\frac{3}{4}C\epsilon_1,
\]
which proves \eqref{eU-improvement}.
\end{proof}

\subsubsection{High frequency heat estimates for $\psi$}
We now combine the heat flow bounds above with the $H^\epsilon$ control and the second weighted profile estimate to obtain additional high frequency decay and the frequency summation bounds used below.
\begin{proposition}\label{highpsi}
Under \eqref{booassum}, for every $s\geq0$ and $k\in\ZZ$, one has
\begin{equation}\label{epsihighfre1}
  \|P_k\psi(s,t)\|_{L_x^\infty}\lesssim\epsilon_1t^{-1+\frac\beta2}2^{-\frac\epsilon2k}(1+s2^{2k})^{-10}.
\end{equation}
Together with \eqref{epsidecay3}, this gives
\begin{align}
  \sum_k\sup_{s\geq0}(1+s2^{2k})^{10}\|P_k\psi(s,t)\|_{L_x^\infty}&\lesssim\epsilon_1t^{-1}(1+\log t),\label{epsihighfre2}\\
  \sum_k\sup_{s\geq0}(1+s2^{2k})^{20}\|P_k\psi(s,t)\|_{L_x^\infty}^2&\lesssim\epsilon_1^2t^{-2}(1+\log t).\label{epsihighfre3}
\end{align}
\end{proposition}
\begin{proof}
For $k\leq0$, \eqref{epsidecay3} implies \eqref{epsihighfre1}. We then assume $k>0$. By \eqref{ebootstrapcor3} and the dispersive estimate,
\begin{equation*}
  \|P_ke^{s\Delta}\psi(0,t)\|_{L_x^\infty}\lesssim e^{-cs2^{2k}}t^{-1}\|P_kf(t)\|_{L_x^1}\lesssim\epsilon_1t^{-1+\frac\beta2}2^{-\frac\epsilon2k}(1+s2^{2k})^{-10}.
\end{equation*}
For the nonlinear terms, \eqref{eAinfty1}, \eqref{epsidecay1}, and \eqref{eFdefinition} imply
\begin{equation*}
  \max_F\|P_{\leq k-5}F(r,t)\|_{L_x^\infty}\lesssim2^k\|A(r,t)\|_{L_x^\infty}+\|A(r,t)\|_{L_x^\infty}^2+\|\psi(r,t)\|_{L_x^\infty}^2\lesssim\epsilon_1^22^k(t+r)^{-\frac32}.
\end{equation*}
Using the same paraproduct decomposition as in \eqref{eUimpro4}, \eqref{eUimpro6}, and \eqref{eUimpro7}, we apply \eqref{eheatenergyimproved2}, \eqref{econnectionenve1}, and \eqref{eFenve1} with $\sigma=\epsilon$ to the high frequency $L_x^2$ factors. Together with \eqref{epsidecay2}, \eqref{epsidecay3}, \eqref{edkdef2}, and \eqref{eenergyenvelope3}, this gives
\begin{align*}
  &\|P_k\nabla\cdot(A\psi)(r,t)\|_{L_x^2}+\max_F\|P_k(F\psi)(r,t)\|_{L_x^2}\\
  &\quad\lesssim\epsilon_1^3t^{-1}(2^k+r^{-\frac12})\left[2^{-\epsilon k}(1+r2^{2k})^{-19}+\sum_{j\geq k-4}2^{-\epsilon j}(1+r2^{2j})^{-29}\right].
\end{align*}
By Bernstein inequality, \eqref{epolynomialconvolution}, and \eqref{esingularheatconvolution},
\begin{align*}
  &(1+s2^{2k})^{10}\|P_k(\psi(s,t)-e^{s\Delta}\psi(0,t))\|_{L_x^\infty}\\
  &\quad\lesssim\epsilon_1^3t^{-1}\left[2^{-\epsilon k}+\sum_{j\geq k-4}\left(2^{2(k-j)}+2^{k-j}\right)2^{-\epsilon j}\right]\lesssim\epsilon_1^3t^{-1}2^{-\epsilon k}.
\end{align*}
Combining this with the free heat estimate proves \eqref{epsihighfre1}. By \eqref{epsidecay3} and \eqref{epsihighfre1},
\begin{align*}
  &\sum_k\sup_{s\geq0}(1+s2^{2k})^{10}\|P_k\psi(s,t)\|_{L_x^\infty}\\
  &\quad\lesssim\epsilon_1\sum_{2^k\leq t^{-\frac12}}2^{2k}+\epsilon_1t^{-1}\sum_{t^{-\frac12}<2^k<t^{\frac\beta\epsilon}}1+\epsilon_1t^{-1+\frac\beta2}\sum_{2^k\geq t^{\frac\beta\epsilon}}2^{-\frac\epsilon2k}\lesssim\epsilon_1t^{-1}(1+\log t).
\end{align*}
Similarly,
\begin{align*}
  &\sum_k\sup_{s\geq0}(1+s2^{2k})^{20}\|P_k\psi(s,t)\|_{L_x^\infty}^2\\
  &\quad\lesssim\epsilon_1^2\left[\sum_{2^k\leq t^{-\frac12}}2^{4k}+t^{-2}\sum_{t^{-\frac12}<2^k<t^{\frac\beta\epsilon}}1+t^{-2+\beta}\sum_{2^k\geq t^{\frac\beta\epsilon}}2^{-\epsilon k}\right]\lesssim\epsilon_1^2t^{-2}(1+\log t).
\end{align*}
This proves both frequency summation estimates.
\end{proof}

\begin{corollary}\label{cpsiheatintegrability}
Under \eqref{booassum}, one has
\begin{align}
  \int_0^\infty\|\nabla\psi(s,t)\|_{L_x^\infty}^2\,\mathrm{d}s&\lesssim\epsilon_1^2t^{-2}(1+\log t),\label{epsihighfre4}\\
  \int_0^\infty\|\nabla^2\psi(s,t)\|_{L_x^\infty}\,\mathrm{d}s&\lesssim\epsilon_1t^{-1}(1+\log t).\label{epsihighfre5}
\end{align}
\end{corollary}
\begin{proof}
Since
\begin{equation*}
  2^{j+k}\int_0^\infty(1+s2^{2j})^{-10}(1+s2^{2k})^{-10}\,\mathrm{d}s\lesssim2^{-|j-k|},
\end{equation*}
Bernstein inequality, Young inequality for sequences, and \eqref{epsihighfre3} give
\begin{align*}
  \int_0^\infty\|\nabla\psi(s,t)\|_{L_x^\infty}^2\,\mathrm{d}s&\lesssim\sum_{j,k}2^{-|j-k|}\left[\sup_{s\geq0}(1+s2^{2j})^{10}\|P_j\psi(s,t)\|_{L_x^\infty}\right]\left[\sup_{s\geq0}(1+s2^{2k})^{10}\|P_k\psi(s,t)\|_{L_x^\infty}\right]\\
  &\lesssim\sum_k\sup_{s\geq0}(1+s2^{2k})^{20}\|P_k\psi(s,t)\|_{L_x^\infty}^2\lesssim\epsilon_1^2t^{-2}(1+\log t).
\end{align*}
Similarly, by \eqref{epsihighfre2},
\begin{align*}
  \int_0^\infty\|\nabla^2\psi(s,t)\|_{L_x^\infty}\,\mathrm{d}s&\lesssim\sum_k\left[\sup_{s\geq0}(1+s2^{2k})^{10}\|P_k\psi(s,t)\|_{L_x^\infty}\right]\int_0^\infty2^{2k}(1+s2^{2k})^{-10}\,\mathrm{d}s\\
  &\lesssim\sum_k\sup_{s\geq0}(1+s2^{2k})^{10}\|P_k\psi(s,t)\|_{L_x^\infty}\lesssim\epsilon_1t^{-1}(1+\log t).
\end{align*}
\end{proof}

\subsection{Decay estimates for $A$ at the initial heat time}
\begin{corollary}\label{cA time decay}
Under the bootstrap assumption \eqref{booassum}, uniformly in the physical time $t$,
\begin{align}
  \|A(0,t)\|_{L_x^\infty}&\lesssim\epsilon_1^2t^{-\frac32}, \label{eAlinftydecay}\\
  \|A(0,t)\|_{L_x^2}&\lesssim\epsilon_1^2t^{-1}(1+\log t).
\label{eAl2decay}
\end{align}
\end{corollary}

\begin{proof}
The bound \eqref{eAlinftydecay} follows immediately from \eqref{eAinfty1} with $s=0$, while \eqref{eAl2decay} follows from \eqref{eAL2} with $s=0$.
\end{proof}

\begin{proposition}\label{Aderi2}
Under the bootstrap assumption \eqref{booassum}, for $0\leq\sigma\leq\epsilon$ and every $s\geq0$, one has
\begin{equation}\label{eAonederivativel2}
  2^{(\sigma+1)k}\|P_kA(s,t)\|_{L_x^2}\lesssim\epsilon_1^2(t+s)^{-1}b_k(\sigma)(1+s2^{2k})^{-19}.
\end{equation}
Hence,
\begin{equation}\label{eAonederivativel21}
  \|A(s,t)\|_{\dot H^{\sigma+1}}\lesssim\epsilon_1(t+s)^{-1}\|\psi(0,t)\|_{\dot H^\sigma}\lesssim\epsilon_1^2(t+s)^{-1}.
\end{equation}
In particular,
\[
  \|A(0,t)\|_{\dot H^{\sigma+1}}\lesssim\epsilon_1t^{-1}\|\psi(0,t)\|_{\dot H^\sigma}\lesssim\epsilon_1^2t^{-1}.
\]
\end{proposition}

\begin{proof}
Recall from \eqref{epsidecay1} that, for every \(r\geq s\),
\[
  \|\psi(r,t)\|_{L_x^\infty}\lesssim\epsilon_1(t+r)^{-1}\lesssim\epsilon_1(t+s)^{-1}.
\]
For the quadratic Low-High interaction, using \eqref{eheatenergyimproved2}, we obtain
\begin{align*}
&2^{(\sigma+1)k}\left\|P_k\int_s^\infty P_{\leq k-5}\psi(r,t)\nabla P_{k+O(1)}\psi(r,t)\,\mathrm{d}r\right\|_{L_x^2}\\
&\quad\lesssim 2^{(\sigma+1)k}\int_s^\infty\|P_{\leq k-5}\psi(r,t)\|_{L_x^\infty}2^k\|P_{k+O(1)}\psi(r,t)\|_{L_x^2}\,\mathrm{d}r\\
&\quad\lesssim\epsilon_1^2(t+s)^{-1}2^{(\sigma+2)k}2^{-\sigma k}b_k(\sigma)\int_s^\infty(1+r2^{2k})^{-20}\,\mathrm{d}r\\
&\quad\lesssim\epsilon_1^2(t+s)^{-1}b_k(\sigma)(1+s2^{2k})^{-19}.
\end{align*}
The quadratic High-Low interaction is estimated in the same way and admits the same bound. For the quadratic High-High interaction,
\begin{align*}
&2^{(\sigma+1)k}\sum_{j\geq k-4}\left\|P_k\int_s^\infty P_{j+O(1)}\psi(r,t)\nabla P_j\psi(r,t)\,\mathrm{d}r\right\|_{L_x^2}\\
&\quad\lesssim 2^{(\sigma+1)k}\sum_{j\geq k-4}2^j\int_s^\infty\|P_{j+O(1)}\psi(r,t)\|_{L_x^\infty}\|P_j\psi(r,t)\|_{L_x^2}\,\mathrm{d}r\\
&\quad\lesssim\epsilon_1^2(t+s)^{-1}\sum_{j\geq k-4}2^{(\sigma+1)(k-j)}b_j(\sigma)(1+s2^{2j})^{-19}\\
&\quad\lesssim\epsilon_1^2(t+s)^{-1}b_k(\sigma)(1+s2^{2k})^{-19}.
\end{align*}
It remains to estimate the remainder in \eqref{ecaloricconnection}. Let $C$ be the best constant in
\[
  2^{(\sigma+1)k}\|P_kA(s,t)\|_{L_x^2}\leq C\epsilon_1^2(t+s)^{-1}b_k(\sigma)(1+s2^{2k})^{-19}, \qquad s\geq0.
\]
We divide the remainder interactions into two cases according to the factor carrying the highest frequency.

\smallskip
\noindent\textit{Remainder case 1: $A$ is the unique highest frequency factor.}

In this case, the frequency of $A$ is $k+O(1)$, while both $\psi$ factors have frequencies at most $k-5$. By the bootstrap estimate above and \eqref{epsiinte3},
\[
\begin{aligned}
&2^{(\sigma+1)k}\left\|P_k\int_s^\infty P_{k+O(1)}A(r,t)P_{\leq k-5}\psi(r,t)P_{\leq k-5}\psi(r,t)\,\mathrm{d}r\right\|_{L_x^2}\\
&\quad\lesssim C\epsilon_1^2b_k(\sigma)(1+s2^{2k})^{-19}\int_s^\infty(t+r)^{-1}\|\psi(r,t)\|_{L_x^\infty}^2\,\mathrm{d}r\\
&\quad\lesssim C\epsilon_1^4(t+s)^{-2}b_k(\sigma)(1+s2^{2k})^{-19}.
\end{aligned}
\]

\smallskip
\noindent\textit{Remainder case 2: at least one $\psi$ factor carries the highest frequency.}

By \eqref{eheatenergyimproved2}, \eqref{eAinfty1}, and \eqref{epsidecay1},
\[
\begin{aligned}
&2^{(\sigma+1)k}\sum_{j\geq k-8}\left\|P_k\int_s^\infty P_{j+O(1)}\psi(r,t)P_{\leq j+4}A(r,t)P_{\leq j+4}\psi(r,t)\,\mathrm{d}r\right\|_{L_x^2}\\
&\quad\lesssim\epsilon_1^4\sum_{j\geq k-8}2^{(\sigma+1)(k-j)}b_j(\sigma)2^j\int_s^\infty(t+r)^{-\frac{5}{2}}(1+r2^{2j})^{-20}\,\mathrm{d}r\\
&\quad\lesssim\epsilon_1^4(t+s)^{-2}\sum_{j\geq k-8}2^{(\sigma+1)(k-j)}b_j(\sigma)2^j\int_s^\infty r^{-\frac{1}{2}}(1+r2^{2j})^{-20}\,\mathrm{d}r.
\end{aligned}
\]
Using
\[
  2^j\int_s^\infty r^{-\frac{1}{2}}(1+r2^{2j})^{-20}\,\mathrm{d}r\lesssim (1+s2^{2j})^{-\frac{39}{2}},
\]
the last expression is bounded by
\[
  \epsilon_1^4(t+s)^{-2}\sum_{j\geq k-8}2^{(\sigma+1)(k-j)}b_j(\sigma)(1+s2^{2j})^{-\frac{39}{2}}\lesssim\epsilon_1^4(t+s)^{-2}b_k(\sigma)(1+s2^{2k})^{-19}.
\]
Combining the quadratic and remainder estimates, the best constant satisfies
\[
  C\lesssim1+C\epsilon_1^2.
\]
Choosing $\epsilon_1$ sufficiently small proves \eqref{eAonederivativel2}. Finally, by \eqref{eenergyenvelope3},
\[
  \|A(s,t)\|_{\dot H^{\sigma+1}}\lesssim\epsilon_1(t+s)^{-1}\|\psi(0,t)\|_{\dot H^\sigma}\lesssim\epsilon_1^2(t+s)^{-1}.
\]
This proves \eqref{eAonederivativel21}, taking \(s=0\) gives the corresponding estimate at the initial heat time.
\end{proof}

\begin{proposition}\label{proinfty}
Under the bootstrap assumption \eqref{booassum}, one has
\begin{equation}\label{eAonederivativelinfty}
  \sup_{s\geq0}\|\nabla A(s,t)\|_{L_x^\infty}\lesssim\epsilon_1^2t^{-2}(1+\log t).
\end{equation}
In particular,
\[
  \|\nabla A(0,t)\|_{L_x^\infty}\lesssim\epsilon_1^2t^{-2}(1+\log t).
\]
\end{proposition}

\begin{proof} 
It suffices to consider $s=0$, since \eqref{ecaloricconnection} gives the same integral bounds for $s>0$ after restricting the heat integral to $[s,\infty)$. By Bernstein inequality and \eqref{eAinfty3}, for $s>0$,
\[
\|\nabla A(s,t)\|_{L_x^\infty}\lesssim\epsilon_1^2(t+s)^{-\frac32}\sum_k2^k(1+s2^{2k})^{-9}\lesssim\epsilon_1^2s^{-\frac12}(t+s)^{-\frac32}.
\]
Differentiating \eqref{ecaloricconnection}, and using \eqref{epsihighfre4}, \eqref{epsihighfre5}, \eqref{epsidecay1}, \eqref{eAinfty1}, and \eqref{epsiinte2}, we obtain
\begin{align*}
  \|\nabla A(0,t)\|_{L_x^\infty}&\lesssim\int_0^\infty\left(\|\nabla\psi(s,t)\|_{L_x^\infty}^2+\|\psi(s,t)\|_{L_x^\infty}\|\nabla^2\psi(s,t)\|_{L_x^\infty}\right)\,\mathrm{d}s\\
  &\quad+\int_0^\infty\left(\|\nabla A(s,t)\|_{L_x^\infty}\|\psi(s,t)\|_{L_x^\infty}^2+\|A(s,t)\|_{L_x^\infty}\|\psi(s,t)\|_{L_x^\infty}\|\nabla\psi(s,t)\|_{L_x^\infty}\right)\,\mathrm{d}s\\
  &\lesssim\epsilon_1^2t^{-2}(1+\log t)+\epsilon_1^4\int_0^\infty s^{-\frac12}(t+s)^{-\frac72}\,\mathrm{d}s+\epsilon_1^3t^{-\frac52}\int_0^\infty\|\nabla\psi(s,t)\|_{L_x^\infty}\,\mathrm{d}s\\
  &\lesssim\epsilon_1^2t^{-2}(1+\log t)+\epsilon_1^4t^{-3}\lesssim\epsilon_1^2t^{-2}(1+\log t).
\end{align*}
This proves \eqref{eAonederivativelinfty}.
\end{proof}

\begin{remark}\label{conde} 
It is useful to compare this heat time argument with the quadratic contribution $A^{(2)}$ defined in Definition~\ref{def:connectiondecom} below. The frequency estimate \eqref{epsihighfre1} and its summation \eqref{epsihighfre2} account for the logarithmic loss.

After integrating the two free heat evolutions in \eqref{eAcom2}, the multiplier for $\partial_iA_m^{(2)}$, with $i,m\in\{1,2\}$, has the symbol
\begin{equation}
  \frac{\xi_i(\xi-\eta)_\ell}{|\eta|^2+|\xi-\eta|^2}. \label{eAonederivativelinfty1}
\end{equation}
The numerator and denominator in \eqref{eAonederivativelinfty1} have the same degree. The resulting zero order symbol can therefore be estimated on each localized frequency interaction without losing derivatives. From \eqref{eAcom2}, we have
\[
  \partial_iA_m^{(2)}(0,t)=\frac{1}{2\pi}\operatorname{Im}\mathcal{F}^{-1}\left[\int_{\RR^2}\frac{\xi_i(\xi-\eta)_\ell}{|\eta|^2+|\xi-\eta|^2}\widehat{\overline{\psi_m(0,t)}}(\eta)\widehat{\psi_\ell(0,t)}(\xi-\eta)
  \,\mathrm{d}\eta\right].
\]
We estimate the associated bilinear multiplier after localizing the output and the two input frequencies. On every fixed frequency interaction, all the rescaled derivatives of the localized rational symbol satisfy the same bounds as the symbol itself. Since we do not exploit the oscillatory structure here, we omit the complex conjugate notation in what follows. For the Low-High interaction, Lemma~\ref{lmet} gives
\[
\left\|\varphi_k(\xi)\varphi_{\leq k-5}(\eta)\varphi_{k+O(1)}(\xi-\eta)\frac{\xi_i(\xi-\eta)_\ell}{|\eta|^2+|\xi-\eta|^2}\right\|_{\mathcal S^\infty}\lesssim1.
\]
Hence, by Lemma~\ref{l multilinear estimates},
\[
\begin{aligned}
    \sum_k\left\|P_k\nabla A^{(2)}\bigl(P_{\leq k-5}\psi,P_{k+O(1)}\psi\bigr)\right\|_{L_x^\infty}&\lesssim\sum_k\|P_{\leq k-5}\psi(0,t)\|_{L_x^\infty}\|P_{k+O(1)}\psi(0,t)\|_{L_x^\infty}\\
    &\lesssim\|\psi(0,t)\|_{L_x^\infty}\sum_k\|P_k\psi(0,t)\|_{L_x^\infty}.
\end{aligned}
\]
The High-Low interaction admits the same bound. For the High-High interactions, Lemma~\ref{lmet} gives
\[
\left\|\varphi_k(\xi)\varphi_{j+O(1)}(\eta)\varphi_j(\xi-\eta)\frac{\xi_i(\xi-\eta)_\ell}{|\eta|^2+|\xi-\eta|^2}\right\|_{\mathcal S^\infty}\lesssim 2^{k-j}.
\]
It follows from Lemma~\ref{l multilinear estimates} that
\[
\begin{aligned}
    \sum_j\sum_{k\leq j+4}\left\|P_k\nabla A^{(2)}\bigl(P_j\psi,P_{j+O(1)}\psi\bigr)\right\|_{L_x^\infty}&\lesssim\sum_j\sum_{k\leq j+4}2^{k-j}\|P_j\psi(0,t)\|_{L_x^\infty}\|P_{j+O(1)}\psi(0,t)\|_{L_x^\infty}\\
    &\lesssim\|\psi(0,t)\|_{L_x^\infty}\sum_j\|P_j\psi(0,t)\|_{L_x^\infty}.
\end{aligned}
\]
Combining the above estimates with \(\sup_{s\geq0}\|\psi(s,t)\|_{L_x^\infty}\lesssim\epsilon_1t^{-1}\) and \eqref{epsihighfre2}, we conclude that
\begin{equation}
    \|\nabla A^{(2)}(0,t)\|_{L_x^\infty}\lesssim\epsilon_1^2t^{-2}(1+\log t). \label{eAonederivativelinfty2}
\end{equation}

The localized $\mathcal S^\infty$ bounds above give control by
\[
  \|\psi(0,t)\|_{L_x^\infty}\sum_k\|P_k\psi(0,t)\|_{L_x^\infty}.
\]
Thus the endpoint estimate uses a Besov type sum over one input. The zero order symbol of \eqref{eAonederivativelinfty1} alone does not justify an $L_x^\infty\times L_x^\infty\to L_x^\infty$ bound for the full multiplier. The logarithmic factor in \eqref{eAonederivativelinfty2} comes from summing the intermediate frequency scales in \eqref{epsihighfre2}.
\end{remark}

\subsection{One order spatial moment of the connection}
This subsection is devoted to proving the one order spatial moment estimate of the connection coefficients. We have the following proposition:
\begin{proposition}\label{pxA estimate}
Under the bootstrap assumption \eqref{booassum}, there holds
\begin{equation}
  \|xA(0,t)\|_{L_x^2}\lesssim \epsilon_1^2. \label{exA1}
\end{equation}
\end{proposition}
The main ingredient is a weighted estimate for the differential field along the heat flow. Since multiplication by \(x\) does not commute with the heat operator, inspired by the Galilean vector field associated to the Schr\"odinger equation, we introduce a heat modified vector field which also incorporates the physical Schr\"odinger evolution.

\subsubsection{The heat modified vector field}
For \(m=1,2\), define
\begin{equation}
  J_{s,m}:=x_m+2(s+it)\partial_m, \qquad J_s:=x+2(s+it)\nabla. \label{eJs1}
\end{equation}
A direct computation gives
\[
  [\partial_s-\Delta,x_m]=-2\partial_m, \qquad [\partial_s-\Delta,2(s+it)\partial_m]=2\partial_m.
\]
Hence,
\begin{equation}
  [\partial_s-\Delta,J_{s,m}]=0. \label{eJs2}
\end{equation}
The commutator with a spatial derivative is
\begin{equation}
  [J_{s,m},\partial_l]=-\delta_{ml}. \label{eJs3}
\end{equation}
It also follows from \eqref{eJs2} that, for \(0\leq r\leq s\),
\[
  J_s e^{(s-r)\Delta}=e^{(s-r)\Delta}J_r.
\]
Recall from \eqref{edef-profile} that \(\psi(0,t)=e^{it\Delta}f(t)\). The standard Schr\"odinger commutation identity
\[
  (x+2it\nabla)e^{it\Delta}=e^{it\Delta}x
\]
therefore gives
\begin{equation}
  J_0\psi(0,t)=(x+2it\nabla)e^{it\Delta}f(t)=e^{it\Delta}(xf(t)). \label{eJs4}
\end{equation}
In particular, by \eqref{eJs4},
\begin{equation}
  \|P_kJ_0\psi(0,t)\|_{L_x^2}=\|P_k(xf(t))\|_{L_x^2}. \label{eJs5}
\end{equation}
Using \eqref{eJs2}, \eqref{eJs3}, and the product rule, the full heat equation \eqref{ecaloric4} gives
\begin{equation}
\begin{aligned}
  (\partial_s-\Delta)J_{s,n}\psi_m&=2i\partial_l(A_lJ_{s,n}\psi_m)+4i(s+it)\partial_l((\partial_nA_l)\psi_m)-2iA_n\psi_m\\
  &\quad-(A_lA_l+i\partial_lA_l)J_{s,n}\psi_m-2(s+it)\partial_n(A_lA_l+i\partial_lA_l)\psi_m\\
  &\quad+i\operatorname{Im}(\psi_m\overline{\psi_l})J_{s,n}\psi_l+2i(s+it)\partial_n\operatorname{Im}(\psi_m\overline{\psi_l})\psi_l.
\end{aligned} \label{eJs6}
\end{equation}
The last four terms are finite linear combinations of $FJ_s\psi$ and $(s+it)(\nabla F)\psi$, with $F$ as in \eqref{eFdefinition}. We omit component indices in the estimates below.

\subsubsection{The weighted frequency envelope}
For every fixed physical time \(t\), define
\begin{equation}
  e_k(t):=\epsilon_1^{-1}t^{-\gamma}\sup_{j\in\mathbb Z}2^{-\delta|j-k|}\|P_j(xf(t))\|_{L_x^2}, \label{eek1}
\end{equation}
and
\begin{equation}\label{eEk}
  E_k(t):=e_k(t)+t^{-\gamma}b_k(t).
\end{equation}
The bootstrap assumption \eqref{booassum} implies
\begin{equation}
  \sum_{k\in\mathbb Z}e_k(t)^2+\sum_{k\in\mathbb Z}E_k(t)^2\lesssim 1. \label{eek2}
\end{equation}
Moreover, both sequences are slowly varying:
\begin{equation}
  e_k(t)\lesssim 2^{\delta|k-j|}e_j(t), \qquad E_k(t)\lesssim 2^{\delta|k-j|}E_j(t). \label{eek3}
\end{equation}
By \eqref{eJs5} and \eqref{eek1},
\begin{equation}
  \|P_kJ_0\psi(0,t)\|_{L_x^2}\lesssim\epsilon_1t^\gamma e_k(t). \label{eek4}
\end{equation}
\begin{proposition}[Propagation of the weighted envelope]
\label{pJbootstrap}
Under the bootstrap assumption \eqref{booassum}, assume that
\begin{equation}
  \|P_kJ_s\psi(s,t)\|_{L_x^2}\leq 2C\epsilon_1t^\gamma E_k(1+s2^{2k})^{-10}, \label{eJbootstrap}
\end{equation}
then
\begin{equation}
  \|P_kJ_s\psi(s,t)\|_{L_x^2}\leq\frac{3}{2}C\epsilon_1t^\gamma E_k(1+s2^{2k})^{-10}. \label{eJimprove}
\end{equation}
\end{proposition}

\begin{proof}
By the Duhamel formula associated with \eqref{eJs6} and \eqref{eek4}, the free heat evolution satisfies
\[
  \|P_ke^{s\Delta_x}J_0\psi(0,t)\|_{L_x^2}\lesssim e^{-cs2^{2k}}\|P_kJ_0\psi(0,t)\|_{L_x^2}\lesssim\epsilon_1t^\gamma e_k(1+s2^{2k})^{-10}.
\]

\smallskip
\noindent\textit{Estimate of the term containing \(\nabla\cdot(AJ_s\psi)\).}

For the Low-High interaction, by \eqref{eAinfty2} and \eqref{eJbootstrap},
\[
  \|P_k\nabla\cdot(P_{\leq k-5}A\,P_{k+O(1)}J_r\psi)(r,t)\|_{L_x^2}\lesssim C\epsilon_1^3t^\gamma 2^k\min\left\{2^{2k},(t+r)^{-\frac32}\right\}E_k(1+r2^{2k})^{-10}.
\]
Since
\[
  2^{-k}\min\left\{2^{2k}, (t+r)^{-\frac32}\right\}\lesssim (t+r)^{-\frac34}\lesssim 1,
\]
Lemma~\ref{lpolynomialconvolution} gives
\[
  \left\|P_k\int_0^se^{(s-r)\Delta}\nabla\cdot(P_{\leq k-5}A\,P_{k+O(1)}J_r\psi)(r,t)\,\mathrm{d}r\right\|_{L_x^2}\lesssim C\epsilon_1^3t^\gamma E_k(1+s2^{2k})^{-10}.
\]
For the High-Low interaction, Bernstein inequality, \eqref{eek3} and \eqref{eJbootstrap} imply
\[
  \|P_{\leq k-5}J_r\psi(r,t)\|_{L_x^\infty}\lesssim\sum_{j\leq k-5}2^j\|P_jJ_r\psi(r,t)\|_{L_x^2}\lesssim C\epsilon_1t^\gamma 2^kE_k.
\]
Using \eqref{econnectionenve1}, we obtain
\[
  \|P_k\nabla\cdot(P_{k+O(1)}A\,P_{\leq k-5}J_r\psi)(r,t)\|_{L_x^2}\lesssim C\epsilon_1^3t^\gamma 2^{2k}d_kE_k(1+r2^{2k})^{-19}.
\]
Hence Lemma~\ref{lpolynomialconvolution} gives
\[
  \left\|P_k\int_0^se^{(s-r)\Delta}\nabla\cdot(P_{k+O(1)}A\,P_{\leq k-5}J_r\psi)(r,t)\,\mathrm{d}r\right\|_{L_x^2}\lesssim C\epsilon_1^3t^\gamma E_k(1+s2^{2k})^{-10}.
\]
For the High-High interaction, by Bernstein inequality, \eqref{econnectionenve1}, and \eqref{eJbootstrap}, we have
\[
\begin{aligned}
  \|P_k\nabla\cdot(P_jA\,P_{j+O(1)}J_r\psi)(r,t)\|_{L_x^2}&\lesssim 2^{2k}\|P_jA(r,t)\|_{L_x^2}\|P_{j+O(1)}J_r\psi(r,t)\|_{L_x^2}\\
  &\lesssim C\epsilon_1^3t^\gamma 2^{2k}d_jE_j(1+r2^{2j})^{-29}.
\end{aligned}
\]
Again by Lemma~\ref{lpolynomialconvolution},
\[
\begin{aligned}
  &\sum_{j\geq k-4}\left\|P_k\int_0^se^{(s-r)\Delta}\nabla\cdot(P_jA\,P_{j+O(1)}J_r\psi)(r,t)\,\mathrm{d}r\right\|_{L_x^2}\\
  &\quad\lesssim C\epsilon_1^3t^\gamma(1+s2^{2k})^{-10}\sum_{j\geq k-4}2^{2(k-j)}d_jE_j\\
  &\quad\lesssim C\epsilon_1^3t^\gamma E_k(1+s2^{2k})^{-10}.
\end{aligned}
\]
Combining the Low-High, High-Low, and High-High estimates,
\[
  \left\|P_k\int_0^se^{(s-r)\Delta}\nabla\cdot(AJ_r\psi)(r,t)\,\mathrm{d}r\right\|_{L_x^2}\lesssim C\epsilon_1^3t^\gamma E_k(1+s2^{2k})^{-10}.
\]

\smallskip
\noindent\textit{Estimate of the term containing \(\nabla\cdot((r+it)(\nabla A)\psi)\).}

It follows from \eqref{eAinfty2} and Bernstein inequality that
\[
  \|\nabla P_{\leq k}A(r,t)\|_{L_x^\infty}\lesssim\epsilon_1^2 2^k\min\left\{2^{2k}, (t+r)^{-\frac32}\right\}.
\]
Since \(|r+it|\lesssim t+r\), we have
\[
  |r+it|\|\nabla P_{\leq k}A(r,t)\|_{L_x^\infty}\lesssim\epsilon_1^22^k(t+r)^{-\frac12}\lesssim\epsilon_1^22^k.
\]
For the Low-High interaction, by \eqref{eheatenergyimproved2},
\[
  \left\|P_k\int_0^se^{(s-r)\Delta}\nabla\cdot\bigl((r+it)(\nabla P_{\leq k-5}A)P_{k+O(1)}\psi\bigr)\,\mathrm{d}r\right\|_{L_x^2}\lesssim\epsilon_1^3b_k(1+s2^{2k})^{-10}.
\]
For the High-Low interaction, by \eqref{econnectionenve1} and \eqref{epsidecay2},
\[
  \|P_k\nabla\cdot\bigl((r+it)(\nabla P_{k+O(1)}A)P_{\leq k-5}\psi\bigr)(r,t)\|_{L_x^2}\lesssim\epsilon_1^32^{2k}d_k(1+r2^{2k})^{-19},
\]
this yields the same estimate as the low-high case. For the High-High interaction, using \eqref{econnectionenve1} and \eqref{epsidecay3}, we obtain
\[
  \|P_k\nabla\cdot\bigl((r+it)(\nabla P_jA)P_{j+O(1)}\psi\bigr)(r,t)\|_{L_x^2}\lesssim\epsilon_1^32^{k+j}d_j(1+r2^{2j})^{-19}.
\]
Applying Lemma~\ref{lpolynomialconvolution},
\[
\begin{aligned}
  &\sum_{j\geq k-4}\left\|P_k\int_0^se^{(s-r)\Delta}\nabla\cdot\bigl((r+it)(\nabla P_jA)P_{j+O(1)}\psi\bigr)\,\mathrm{d}r\right\|_{L_x^2}\\
  &\quad\lesssim\epsilon_1^3(1+s2^{2k})^{-10}\sum_{j\geq k-4}2^{k-j}d_j\\
  &\quad\lesssim\epsilon_1^3b_k(1+s2^{2k})^{-10}.
\end{aligned}
\]
Combining the three frequency interactions and using \(t^{-\gamma}b_k\leq E_k\), we conclude that
\[
  \left\|P_k\int_0^se^{(s-r)\Delta}\nabla\cdot\bigl((r+it)(\nabla A)\psi\bigr)(r,t)\,\mathrm{d}r\right\|_{L_x^2}\lesssim\epsilon_1^3t^\gamma E_k(1+s2^{2k})^{-10}.
\]

\smallskip
\noindent\textit{Estimate of the term containing \(A\psi\).}

For the Low-High interaction, by \eqref{eAinfty2} and \eqref{eheatenergyimproved2},
\[
  \|P_k(P_{\leq k-5}A\,P_{k+O(1)}\psi)(r,t)\|_{L_x^2}\lesssim\epsilon_1^3\min\left\{2^{2k}, (t+r)^{-\frac32}\right\}b_k(1+r2^{2k})^{-20}.
\]
Thus
\[
  \left\|P_k\int_0^se^{(s-r)\Delta}(P_{\leq k-5}A\,P_{k+O(1)}\psi)(r,t)\,\mathrm{d}r\right\|_{L_x^2}\lesssim\epsilon_1^3b_k(1+s2^{2k})^{-10}.
\]
For the High-Low interaction, using \eqref{econnectionenve1} and \eqref{epsidecay2}, we have
\[
  \|P_k(P_{k+O(1)}A\,P_{\leq k-5}\psi)(r,t)\|_{L_x^2}\lesssim\epsilon_1^3d_k\min\left\{2^{2k}, \frac{1}{t+r}\right\}(1+r2^{2k})^{-19}.
\]
The same convolution estimate and \(d_k\lesssim b_k\) give
\[
  \left\|P_k\int_0^se^{(s-r)\Delta}(P_{k+O(1)}A\,P_{\leq k-5}\psi)(r,t)\,\mathrm{d}r\right\|_{L_x^2}\lesssim\epsilon_1^3b_k(1+s2^{2k})^{-10}.
\]
For the High-High interaction, we use the negative half order bounds \eqref{erk3} and \eqref{enegativeenve2}. Combining with Bernstein inequality we have
\[
\begin{aligned}
  \|P_k(P_jA\,P_{j+O(1)}\psi)(r,t)\|_{L_x^2}&\lesssim2^k\|P_jA(r,t)\|_{L_x^2}\|P_{j+O(1)}\psi(r,t)\|_{L_x^2}\\
  &\lesssim\epsilon_1^32^{k+j}r_jq_j(1+r2^{2j})^{-30}.
\end{aligned}
\]
Lemma~\ref{lpolynomialconvolution} therefore gives
\[
\begin{aligned}
  &\sum_{j\geq k-4}\left\|P_k\int_0^se^{(s-r)\Delta}(P_jA\,P_{j+O(1)}\psi)(r,t)\,\mathrm{d}r\right\|_{L_x^2}\\
  &\quad\lesssim\epsilon_1^3(1+s2^{2k})^{-10}\sum_{j\geq k-4}2^{k-j}r_jq_j\\
  &\quad\lesssim\epsilon_1^3b_k(1+s2^{2k})^{-10}.
\end{aligned}
\]
Combining the Low-High, High-Low, and High-High interactions,
\[
  \left\|P_k\int_0^se^{(s-r)\Delta}(A\psi)(r,t)\,\mathrm{d}r\right\|_{L_x^2}\lesssim\epsilon_1^3t^\gamma E_k(1+s2^{2k})^{-10}.
\]

\smallskip
\noindent\textit{Estimate of the remaining nonlinearities.}

By \eqref{eFenve1} with $\sigma=0$, \eqref{eJbootstrap}, and \eqref{esingularheatconvolution}, the same paraproduct estimates give
\[
  \left\|\int_0^s e^{(s-r)\Delta}P_k(FJ_r\psi)(r,t)\,\mathrm{d}r\right\|_{L_x^2}\lesssim C\epsilon_1^3t^\gamma E_k(1+s2^{2k})^{-10}.
\]
For the differentiated terms, the paraproduct decomposition, \eqref{eAonederivativel2} with $\sigma=0$, \eqref{eheatenergyimproved2}, \eqref{econnectionenve1}, \eqref{epsidecay3}, and \eqref{eAinfty3} imply, for $r>0$,
\[
\begin{aligned}
  \|P_kF(r,t)\|_{L_x^2}&\lesssim\epsilon_1^2(t+r)^{-1}\sum_{j\geq k-4}b_j(1+r2^{2j})^{-19},\\
  \|F(r,t)\|_{L_x^\infty}&\lesssim\epsilon_1^2r^{-\frac{1}{2}}(t+r)^{-\frac{3}{2}}.
\end{aligned}
\]
Using these bounds, $|r+it|\lesssim t+r$, and $\|P_{\leq k}\psi(r,t)\|_{L_x^\infty}\lesssim\epsilon_1 2^k$, we apply \eqref{esingularheatconvolution} to the Low-High interaction and Lemma~\ref{lpolynomialconvolution} to the other interactions. Thus
\[
\begin{aligned}
&\left\|\int_0^s e^{(s-r)\Delta}P_k\bigl((r+it)(\nabla F)\psi\bigr)(r,t)\,\mathrm{d}r\right\|_{L_x^2}\\
&\quad\lesssim\epsilon_1^3(1+s2^{2k})^{-10}\left[b_k+\sum_{j\geq k-4}2^{k-j}b_j\sum_{l\geq j-4}2^{2(j-l)}b_l\right]\\
&\quad\lesssim\epsilon_1^3t^\gamma E_k(1+s2^{2k})^{-10}.
\end{aligned}
\]

Combining the free evolution and the five nonlinear estimates, we obtain
\[
  \|P_kJ_s\psi(s,t)\|_{L_x^2}\lesssim\epsilon_1t^\gamma e_k(1+s2^{2k})^{-10}+C\epsilon_1^3t^\gamma E_k(1+s2^{2k})^{-10}.
\]
After choosing \(C\) sufficiently large to dominate the linear term and then choosing \(\epsilon_1\) sufficiently small depending on this fixed \(C\), this improves \eqref{eJbootstrap} to \eqref{eJimprove}.
\end{proof}

As an immediate result of \eqref{eJimprove} and \eqref{eek2}, for every \(s>0\),
\[
  \|\nabla J_s\psi(s,t)\|_{L_x^2}^2\lesssim\epsilon_1^2t^{2\gamma}\sum_{k\in\mathbb{Z}}2^{2k}E_k^2(1+s2^{2k})^{-20}\lesssim\epsilon_1^2t^{2\gamma}s^{-1}\sum_{k\in\mathbb Z}E_k^2.
\]
Therefore,
\begin{equation}
  \|\nabla J_s\psi(s,t)\|_{L_x^2}\lesssim\epsilon_1t^\gamma s^{-\frac12}. \label{exA2}
\end{equation}

\subsubsection{Estimate of \(xA\)}
We first estimate the quadratic term and then absorb the remainder. For the quadratic contribution
\[
  -\int_0^\infty\operatorname{Im}\bigl(\partial_l\psi_l(s,t)\overline{\psi_m(s,t)}\bigr)\,\mathrm{d}s,
\]
in the connection formula \eqref{ecaloricconnection}, we emphasize that this is the quadratic term in the original caloric representation of \(A\), and is not the separately defined \(A^{(2)}\) in the sequel. Since   \(x_n\) is a multiplication operator, we may place the entire weight on the differentiated factor:
\[
  x_n\int_0^\infty\operatorname{Im}\bigl(\partial_l\psi_l(s,t)\overline{\psi_m(s,t)}\bigr)\,\mathrm{d}s=\int_0^\infty\operatorname{Im}\bigl(x_n\partial_l\psi_l(s,t)\overline{\psi_m(s,t)}\bigr)\,\mathrm{d}s.
\]
By \eqref{eJs1} and \eqref{eJs3},
\begin{equation}
  x_n\partial_l\psi_l=\partial_l(J_{s,n}\psi_l)-2(s+it)\partial_l\partial_n\psi_l-\delta_{nl}\psi_l. \label{exA3}
\end{equation}
We estimate the three terms in \eqref{exA3} separately, again omitting unnecessary subscripts.

\begin{proof}[Proof of Proposition~\ref{pxA estimate}] 
For the term containing $\nabla J_s\psi$, \eqref{exA2} and \eqref{epsidecay1} give
\[
\int_0^\infty\|\nabla J_s\psi(s,t)\|_{L_x^2}\|\psi(s,t)\|_{L_x^\infty}\,\mathrm{d}s\lesssim\epsilon_1^2t^\gamma\int_0^\infty\frac{s^{-\frac12}}{t+s}\,\mathrm{d}s\lesssim\epsilon_1^2t^{\gamma-\frac12}\lesssim\epsilon_1^2.
\]
For the term containing $(s+it)\nabla^2\psi$, by \eqref{epsiinte5}, \eqref{epsidecay1}, and $|s+it|\lesssim t+s$,
\[
\int_0^\infty|s+it|\|\nabla^2\psi(s,t)\|_{L_x^2}\|\psi(s,t)\|_{L_x^\infty}\,\mathrm{d}s\lesssim\epsilon_1\int_0^\infty\|\nabla^2\psi(s,t)\|_{L_x^2}\,\mathrm{d}s\lesssim\epsilon_1^2.
\]
For the last term in \eqref{exA3}, \eqref{epsiAdecay1} and \eqref{epsidecay1} imply
\[
\int_0^\infty\|\psi(s,t)\|_{L_x^2}\|\psi(s,t)\|_{L_x^\infty}\,\mathrm{d}s\lesssim\epsilon_1^2\int_0^\infty\frac{\langle s\rangle^{-\frac12}}{t+s}\,\mathrm{d}s\lesssim\epsilon_1^2t^{-\frac12}.
\]
Combining these estimates with \eqref{exA3}, we obtain
\begin{equation}\label{exA4}
  \sup_{s_0\geq0}\left\|x\int_{s_0}^\infty\operatorname{Im}(\overline\psi\nabla\psi)(s,t)\,\mathrm{d}s\right\|_{L_x^2}\lesssim\epsilon_1^2.
\end{equation}
For the remainder, we place the spatial weight on $A$. By \eqref{ecaloricconnection}, \eqref{epsiinte4}, and \eqref{exA4},
\begin{align*}
\sup_{s_0\geq0}\|xA(s_0,t)\|_{L_x^2}&\lesssim\epsilon_1^2+\sup_{s_0\geq0}\|xA(s_0,t)\|_{L_x^2}\int_0^\infty\|\psi(s,t)\|_{L_x^\infty}^2\,\mathrm{d}s\\
&\lesssim\epsilon_1^2+\epsilon_1^2t^{-1}\sup_{s_0\geq0}\|xA(s_0,t)\|_{L_x^2}.
\end{align*}
Absorbing the last term gives
\[
  \sup_{s_0\geq0}\|xA(s_0,t)\|_{L_x^2}\lesssim\epsilon_1^2.
\]
Taking $s_0=0$ proves \eqref{exA1}.
\end{proof}

\subsection{Improved estimates for the higher order connection}
We now introduce the decomposition used in the subsequent refined estimates.

\begin{definition}\label{def:connectiondecom}
For $m\in\{1,2\}$, we decompose the connection coefficient at $s=0$ as
\begin{equation}\label{eAcom1}
A_m(0,t)=A_m^{(2)}(0,t)+A_m^{(\geq4)}(0,t),
\end{equation}
where
\begin{equation}\label{eAcom2}
A_m^{(2)}(0,t):=-\int_0^\infty\operatorname{Im}\left(\overline{e^{s\Delta}\psi_m(0,t)}\,\partial_\ell e^{s\Delta}\psi_\ell(0,t)\right)\,\mathrm{d}s,
\end{equation}
and
\begin{equation*}
A_m^{(\geq4)}(0,t):=A_m(0,t)-A_m^{(2)}(0,t).
\end{equation*}
\end{definition}
For $A_m^{(\geq4)}$, by \eqref{ecaloric4}, the heat correction $\psi(s,t)-e^{s\Delta}\psi(0,t)$ starts at cubic order. Thus replacing at least one factor in the quadratic heat integral by this correction, together with the remainder term in \eqref{ecaloricconnection}, produces terms of order at least four.

In this subsection, we state improved estimates for the higher order contribution $A_m^{(\geq4)}(0,t)$. These bounds will be used in the decay estimates developed in the subsequent sections. The quadratic contribution $A_m^{(2)}(0,t)$ will be treated differently. More precisely, in the later analysis we insert its explicit representation \eqref{eAcom2} into the relevant Schr\"odinger Duhamel terms. This produces cubic oscillatory integrals with phase
\begin{equation*}
\Phi_3(\xi,\eta,\zeta):=|\xi|^2-|\xi-\eta|^2+|\eta-\zeta|^2-|\zeta|^2.
\end{equation*}
We shall analyze the spacetime resonant structure of $\Phi_3$ and combine it with the null structures of the equation to obtain dispersive estimates. For the higher order contribution, the following bounds allow us to work directly with its physical space norms.

\begin{proposition}[Improved estimates for the higher order connection]\label{phigher order connection}
Under the bootstrap assumption \eqref{booassum}, uniformly for $m\in\{1,2\}$ and $t\geq1$, one has
\begin{align}
\|A_m^{(\geq4)}(0,t)\|_{L_x^2}&\lesssim\epsilon_1^4t^{-2},\label{eAhigh1}\\
\|\nabla^2A_m^{(\geq4)}(0,t)\|_{H^1}&\lesssim\epsilon_1^4t^{-\frac83+\frac{2\gamma}{3}}(1+\log t)^{\frac{2}{3}}\lesssim\epsilon_1^4t^{-\frac52},\label{eAhigh2}\\
\|\nabla^3A_m^{(\geq4)}(0,t)\|_{L_x^2}&\lesssim\epsilon_1^4t^{-3+\gamma}(1+\log t),\label{eAhigh3}\\
\|x\nabla A_m^{(\geq4)}(0,t)\|_{L_x^\infty}&\lesssim\epsilon_1^4t^{-\frac52+\gamma}(1+\log t),\label{eAhigh4}\\
\|x\nabla^2A_m^{(\geq4)}(0,t)\|_{L_x^2}&\lesssim\epsilon_1^4t^{-2+\gamma}(1+\log t).\label{eAhigh5}
\end{align}
\end{proposition}
\begin{proof}
To avoid interrupting the main line of the argument, we postpone the proof of Proposition~\ref{phigher order connection} to Appendix~\ref{app:A}.
\end{proof}

\subsection{Estimates for the temporal connection coefficient}

We now consider $A_3$. Recall that $\psi_3(0,t)=i\partial_l\psi_l(0,t)-A_l(0,t)\psi_l(0,t)$ and
\begin{equation}\label{epsi3equation}
  (\partial_s-\Delta)\psi_3=2i\partial_j(A_j\psi_3)-(A_jA_j+i\partial_jA_j)\psi_3+i\operatorname{Im}(\psi_3\overline{\psi_j})\psi_j.
\end{equation}
Since $A_3(\infty,t)=0$, we have
\begin{equation}\label{eA31}
  A_3(s,t)=-\int_s^\infty\operatorname{Im}\left[(\partial_l\psi_l+iA_l\psi_l)(r,t)\overline{\psi_3(r,t)}\right]\,\mathrm{d}r.
\end{equation}

We first establish a positive homogeneous heat estimate for the temporal field.

\begin{lemma}\label{lpsi3sobolev}
Under the bootstrap assumption \eqref{booassum}, for every fixed $0<\sigma\leq\epsilon$, one has
\begin{equation}\label{epsi3sobolev}
\begin{aligned}
  &\left(\sum_k2^{2(\sigma-1)k}\sup_{s\geq0}(1+s2^{2k})^{16}\left[\|P_k\psi_3(s,t)\|_{L_x^2}^2+\|P_k(\partial_l\psi_l+iA_l\psi_l)(s,t)\|_{L_x^2}^2\right]\right)^{\frac12}\\
  &\qquad\lesssim\|\psi(0,t)\|_{\dot H^\sigma}.
\end{aligned}
\end{equation}
\end{lemma}

\begin{proof}
We postpone the proof to Appendix~\ref{app:psi3sobolev}.
\end{proof}

We next compare $\psi_3$ with its leading free heat evolution.
\begin{lemma}\label{lpsi3estimate}
Under the bootstrap assumption \eqref{booassum}, one has
\begin{align}
\sup_{s\geq0}\|\psi_3(s,t)-i\partial_l e^{s\Delta}\psi_l(0,t)\|_{L_x^2}&\lesssim\epsilon_1^3t^{-\frac{3}{2}},\label{epsi31}\\
\sup_{s\geq0}\|\psi_3(s,t)-i\partial_l e^{s\Delta}\psi_l(0,t)\|_{L_x^\infty}&\lesssim\epsilon_1^3t^{-\frac{5}{2}}(1+\log t),\label{epsi32}\\
\left(\int_0^\infty\left\|\nabla\left(\psi_3(s,t)-i\partial_l e^{s\Delta}\psi_l(0,t)\right)\right\|_{L_x^2}^2\,\mathrm{d}s\right)^{\frac{1}{2}}&\lesssim\epsilon_1^3t^{-\frac{3}{2}}.\label{epsi33}
\end{align}
\end{lemma}

\begin{proof}
We fix $t$ and omit it in the heat flow expressions. The initial correction is $-A_l(0)\psi_l(0)$. By \eqref{eAL2}, \eqref{eAinfty1}, and \eqref{booassum},
\[
  \|A(0)\psi(0)\|_{L_x^2}\lesssim\epsilon_1^3t^{-2}(1+\log t), \qquad \|A(0)\psi(0)\|_{L_x^\infty}\lesssim\epsilon_1^3t^{-\frac52}.
\]
By Bernstein inequality and \eqref{eAinfty3},
\[
  \|\nabla A(r)\|_{L_x^\infty}\lesssim\epsilon_1^2r^{-\frac12}(t+r)^{-\frac32}.
\]
Together with \eqref{eAinfty1} and \eqref{epsidecay1}, this gives
\[
  \sup_{s>0}\int_0^s(s-r)^{-\frac12}\|A(r)\|_{L_x^\infty}\,\mathrm{d}r+\int_0^\infty\left(\|A(r)\|_{L_x^\infty}^2+\|\nabla A(r)\|_{L_x^\infty}+\|\psi(r)\|_{L_x^\infty}^2\right)\,\mathrm{d}r\lesssim\epsilon_1^2.
\]

Subtracting the free heat equation from \eqref{epsi3equation}, and integrating the divergence term by parts, the $L_x^2$ energy estimate gives
\[
\begin{aligned}
  &\frac12\frac{\mathrm{d}}{\mathrm{d}s}\|\psi_3(s)-i\partial_l e^{s\Delta}\psi_l(0)\|_{L_x^2}^2+\frac12\left\|\nabla\bigl(\psi_3(s)-i\partial_l e^{s\Delta}\psi_l(0)\bigr)\right\|_{L_x^2}^2\\
  &\lesssim\bigl(\|A(s)\|_{L_x^\infty}^2+\|\psi(s)\|_{L_x^\infty}^2\bigr)\|\psi_3(s)-i\partial_l e^{s\Delta}\psi_l(0)\|_{L_x^2}^2\\
  &\quad+\|\psi_3(s)-i\partial_l e^{s\Delta}\psi_l(0)\|_{L_x^2}\Bigl[\|A(s)\|_{L_x^\infty}\|\nabla^2e^{s\Delta}\psi(0)\|_{L_x^2}+\|\nabla A(s)\,\nabla e^{s\Delta}\psi(0)\|_{L_x^2}\\
  &\hspace{40mm}+\bigl(\|A(s)\|_{L_x^\infty}^2+\|\psi(s)\|_{L_x^\infty}^2\bigr)\|\nabla e^{s\Delta}\psi(0)\|_{L_x^2}\Bigr].
\end{aligned}
\]

For the $L_x^\infty$ estimate, we use the Duhamel formula and \eqref{eheatoperator}. The terms containing the correction are bounded by
\[
\begin{aligned}
  &\int_0^s\Bigl[(s-r)^{-\frac12}\|A(r)\|_{L_x^\infty}+\|A(r)\|_{L_x^\infty}^2+\|\nabla A(r)\|_{L_x^\infty}+\|\psi(r)\|_{L_x^\infty}^2\Bigr]\|\psi_3(r)-i\partial_l e^{r\Delta}\psi_l(0)\|_{L_x^\infty}\,\mathrm{d}r\\
  &\lesssim\epsilon_1^2\sup_{0\leq r\leq s}\|\psi_3(r)-i\partial_l e^{r\Delta}\psi_l(0)\|_{L_x^\infty}.
\end{aligned}
\]
Integrating the energy inequality on $[0,s]$ and absorbing the correction terms in both estimates, we obtain, for $p=2,\infty$,
\[
\begin{aligned}
  &\|\psi_3(s)-i\partial_l e^{s\Delta}\psi_l(0)\|_{L_x^p}\\
  &\lesssim\|A(0)\psi(0)\|_{L_x^p}+\int_0^s\Bigl[\|A(r)\|_{L_x^\infty}\|\nabla^2e^{r\Delta}\psi(0)\|_{L_x^p}+\|\nabla A(r)\,\nabla e^{r\Delta}\psi(0)\|_{L_x^p}\\
  &\hspace{32mm}+\bigl(\|A(r)\|_{L_x^\infty}^2+\|\psi(r)\|_{L_x^\infty}^2\bigr)\|\nabla e^{r\Delta}\psi(0)\|_{L_x^p}\Bigr]\,\mathrm{d}r.
\end{aligned}
\]

By \eqref{eheatoperator}, \eqref{booassum}, \eqref{epsihat2}, and \eqref{bs5},
\begin{equation}\label{epsi34}
\begin{gathered}
\|\nabla e^{r\Delta}\psi(0)\|_{L_x^2}\lesssim\epsilon_1r^{-\frac{1}{2}},\qquad \|\nabla e^{r\Delta}\psi(0)\|_{L_x^\infty}\lesssim\epsilon_1r^{-\frac{1}{2}}(t+r)^{-1},\\
\int_0^\infty\|\nabla^2e^{r\Delta}\psi(0)\|_{L_x^2}\,\mathrm{d}r\lesssim\epsilon_1,\qquad \int_0^\infty\|\nabla^2e^{r\Delta}\psi(0)\|_{L_x^\infty}\,\mathrm{d}r\lesssim\epsilon_1t^{-1}(1+\log t).
\end{gathered}
\end{equation}
For $p=2$, \eqref{eAonederivativel21}, \eqref{eAinfty1}, \eqref{epsidecay1}, and \eqref{epsi34} bound this expression by
\begin{equation*}
\epsilon_1^3\left[t^{-\frac{3}{2}}+\int_0^\infty r^{-\frac{1}{2}}(t+r)^{-2}\,\mathrm{d}r\right]+\epsilon_1^5\int_0^\infty r^{-\frac{1}{2}}(t+r)^{-3}\,\mathrm{d}r\lesssim\epsilon_1^3t^{-\frac{3}{2}}.
\end{equation*}
For $p=\infty$, use \eqref{eAonederivativelinfty} and the free heat bounds above to obtain
\[
  \epsilon_1^3t^{-\frac52}(1+\log t)+\epsilon_1^3\int_0^\infty\frac{r^{-\frac12}}{(t+r)^3}\,\mathrm{d}r+\epsilon_1^5t^{-\frac72}\lesssim\epsilon_1^3t^{-\frac52}(1+\log t).
\]
Finally, integrating the energy inequality and using \eqref{epsi31}, \eqref{epsi34}, and the estimate above for $p=2$, we obtain
\begin{equation*}
\int_0^\infty\left\|\nabla\left(\psi_3(s)-i\partial_l e^{s\Delta}\psi_l(0)\right)\right\|_{L_x^2}^2\,\mathrm{d}s\lesssim\|A(0)\psi(0)\|_{L_x^2}^2+\epsilon_1^6t^{-3}\lesssim\epsilon_1^6t^{-3}.
\end{equation*}
This proves \eqref{epsi33}.
\end{proof}

We now use these estimates to control the full temporal connection in homogeneous Sobolev spaces.
\begin{proposition}\label{pA3sobolev}
Under the bootstrap assumption \eqref{booassum}, for $0<\sigma\leq\epsilon$, one has
\begin{equation}\label{eA32}
  \|A_3(0,t)\|_{\dot H^\sigma}\lesssim\epsilon_1t^{-1}\|\psi(0,t)\|_{\dot H^\sigma}\lesssim\epsilon_1^2t^{-1}.
\end{equation}
\end{proposition}

\begin{proof}
We fix $t$ below. By Bernstein inequality, \eqref{epsi31}, and \eqref{booassum},
\[
  \|P_{\leq k}\psi_3(s)\|_{L_x^\infty}\lesssim2^k\|\psi(0)\|_{L_x^\infty}+2^k\|\psi_3(s)-i\partial_l e^{s\Delta}\psi_l(0)\|_{L_x^2}\lesssim\epsilon_1t^{-1}2^k.
\]
Likewise, \eqref{econnectionenve1} and \eqref{epsidecay1} give
\[
  \|P_{\leq k}(\partial_l\psi_l+iA_l\psi_l)(s)\|_{L_x^\infty}\lesssim2^k\bigl(1+\|A(s)\|_{L_x^2}\bigr)\|\psi(s)\|_{L_x^\infty}\lesssim\epsilon_1t^{-1}2^k.
\]
We apply the Low-High, High-Low, and High-High decomposition to \eqref{eA31}. Using the bounds above for the $L_x^\infty$ factor and integrating the heat weight of the $L_x^2$ factor gives
\[
\begin{aligned}
  &2^{\sigma k}\|P_kA_3(0,t)\|_{L_x^2}\\
  &\lesssim\epsilon_1t^{-1}\sum_{j\geq k-4}2^{\sigma(k-j)}2^{(\sigma-1)j}\sup_{s\geq0}\Bigl[(1+s2^{2j})^8\left(\|P_j\psi_3(s)\|_{L_x^2}+\|P_j(\partial_l\psi_l+iA_l\psi_l)(s)\|_{L_x^2}\right)\Bigr].
\end{aligned}
\]
Taking the $\ell_k^2$ norm, using $\sigma>0$ and \eqref{epsi3sobolev}, proves \eqref{eA32}.
\end{proof}

We next separate the quadratic contribution. Define
\begin{equation}\label{eA33}
  A_3(0,t)=A_3^{(2)}(0,t)+A_3^{(\geq4)}(0,t), \qquad A_3^{(2)}(0,t):=\int_0^\infty|\partial_l e^{s\Delta}\psi_l(0,t)|^2\,\mathrm{d}s.
\end{equation}
Then
\begin{equation}\label{eA34}
  \widehat{A_3^{(2)}(0,t)}(\eta)=\frac{1}{2\pi}\int_{\RR^2}\frac{\zeta_l(\zeta-\eta)_j}{|\zeta|^2+|\zeta-\eta|^2}\widehat{\psi_l(0,t)}(\zeta)\overline{\widehat{\psi_j(0,t)}(\zeta-\eta)}\,\mathrm{d}\zeta.
\end{equation}
This is a zero order bilinear symbol satisfying the hypotheses of Lemma~\ref{lcoifman-meyer estimates}.

Subtracting the quadratic contribution from \eqref{eA31} gives
\begin{equation}\label{eA3higherrep}
\begin{aligned}
  A_3^{(\geq4)}(0,t)=-\int_0^\infty\operatorname{Im}\Bigl[&\left\{\partial_l\bigl(\psi_l(s,t)-e^{s\Delta}\psi_l(0,t)\bigr)+iA_l(s,t)\psi_l(s,t)\right\}\overline{i\partial_j e^{s\Delta}\psi_j(0,t)}\\
  &+(\partial_l\psi_l+iA_l\psi_l)(s,t)\overline{\psi_3(s,t)-i\partial_j e^{s\Delta}\psi_j(0,t)}\Bigr]\,\mathrm{d}s.
\end{aligned}
\end{equation}
Each summand contains a cubic correction and is therefore of order at least four. We state the estimate needed below.

\begin{lemma}[higher order temporal connection]\label{lA3higherorder}
Under the bootstrap assumption \eqref{booassum}, one has
\begin{equation}\label{eA35}
  \|A_3^{(\geq4)}(0,t)\|_{H^2}\lesssim\epsilon_1^4t^{-2}.
\end{equation}
\end{lemma}

\begin{proof}
We postpone the proof to Appendix~\ref{app:A3remainder}.
\end{proof}

We now establish the bounds for the full temporal connection.

\begin{proposition}\label{pA3decay}
Under the bootstrap assumption \eqref{booassum}, one has
\begin{align}
\sup_{s\geq0}\|A_3(s,t)\|_{L_x^\infty}&\lesssim\epsilon_1^2t^{-2}(1+\log t),\label{eA36}\\
\|A_3(0,t)\|_{H^\epsilon}&\lesssim\epsilon_1^2t^{-1}.\label{eA37}
\end{align}
\end{proposition}

\begin{proof}
By \eqref{epsiinte2}, \eqref{epsihighfre4}, \eqref{eAinfty1}, and \eqref{epsidecay1},
\[
  \int_0^\infty\|(\partial_l\psi_l+iA_l\psi_l)(r,t)\|_{L_x^\infty}\,\mathrm{d}r\lesssim\epsilon_1t^{-\frac12},
\]
\[
  \int_0^\infty\|(\partial_l\psi_l+iA_l\psi_l)(r,t)\|_{L_x^\infty}^2\,\mathrm{d}r\lesssim\epsilon_1^2t^{-2}(1+\log t).
\]
The same square integral bound holds for $\partial_l e^{r\Delta}\psi_l(0,t)$, by the free heat bounds and the same splitting as \eqref{epsihighfre4}. Therefore, \eqref{eA31}, \eqref{epsi32}, and Cauchy-Schwarz inequality imply
\[
\begin{aligned}
  \sup_{s\geq0}\|A_3(s,t)\|_{L_x^\infty}&\lesssim\int_0^\infty\|(\partial_l\psi_l+iA_l\psi_l)(r,t)\|_{L_x^\infty}\|\partial_j e^{r\Delta}\psi_j(0,t)\|_{L_x^\infty}\,\mathrm{d}r\\
  &\quad+\epsilon_1^3t^{-\frac52}(1+\log t)\int_0^\infty\|(\partial_l\psi_l+iA_l\psi_l)(r,t)\|_{L_x^\infty}\,\mathrm{d}r\lesssim\epsilon_1^2t^{-2}(1+\log t).
\end{aligned}
\]
For \eqref{eA37}, Lemma~\ref{lcoifman-meyer estimates}, \eqref{eA34}, and \eqref{booassum} give
\begin{equation*}
\|A_3^{(2)}(0,t)\|_{H^\epsilon}\lesssim\|\psi(0,t)\|_{L_x^\infty}\|\psi(0,t)\|_{H^\epsilon}\lesssim\epsilon_1^2t^{-1}.
\end{equation*}
Combining this with \eqref{eA35} proves \eqref{eA37}.
\end{proof}

\section{Energy estimate}\label{sec:energy estimate}

\begin{proposition}\label{p5.1}
Under the bootstrap assumption \eqref{booassum}, there holds
\begin{equation*}
\frac{\mathrm{d}}{\mathrm{d}t}\|\psi(t)\|_{L_x^2}^2=0.
\end{equation*}
For every fixed $0<\sigma\leq\epsilon$ and every $t\in[1,T]$,
\begin{equation*}
\frac{\mathrm{d}}{\mathrm{d}t}\sum_k\|P_k|\nabla|^\sigma\psi(t)\|_{L_x^2}^2\lesssim\epsilon_1^2t^{-2}(1+\log t)\sum_k\|P_k|\nabla|^\sigma\psi(t)\|_{L_x^2}^2.
\end{equation*}
Hence, by \eqref{eboot-init},
\begin{equation*}
\sup_{1\leq t\leq T}\|\psi(t)\|_{H_x^\epsilon}\lesssim\|\psi(1)\|_{H_x^\epsilon}\lesssim\epsilon_0.
\end{equation*}
\end{proposition}

\begin{proof}
We first estimate the $L_x^2$ energy, integration by parts gives
\[
  \frac{1}{2}\frac{\mathrm{d}}{\mathrm{d}t}\|\psi(t)\|_{L_x^2}^2=-\operatorname{Re}\int_{\RR^2}\overline{\psi_m}\left(2A_l\partial_l\psi_m+(\partial_lA_l)\psi_m\right)\,\mathrm{d}x-\operatorname{Re}\int_{\RR^2}\overline{\psi_m}\psi_l\operatorname{Im}(\overline{\psi_l}\psi_m)\,\mathrm{d}x=0.
\]

We next estimate the homogeneous energy for $0<\sigma\leq\epsilon$. Applying $P_k|\nabla|^\sigma$ to \eqref{ecaloric7}, Cauchy-Schwarz gives
\[
\begin{aligned}
  \frac{1}{2}\frac{\mathrm{d}}{\mathrm{d}t}\sum_k\|P_k|\nabla|^\sigma\psi\|_{L_x^2}^2&\lesssim\left|\sum_k\operatorname{Re}\int_{\RR^2}\overline{P_k|\nabla|^\sigma\psi_m}\,P_k|\nabla|^\sigma(A\cdot\nabla\psi_m)\,\mathrm{d}x\right|\\
  &\quad+\|\psi\|_{\dot H_x^\sigma}\max_F\bigl\||\nabla|^\sigma(F\psi)\bigr\|_{L_x^2}+\|\psi\|_{\dot H_x^\sigma}\||\nabla|^\sigma(A_3\psi)\|_{L_x^2}.
\end{aligned}
\]
Here $F$ is as in \eqref{eFdefinition}. For the transport term, paraproduct decomposition reduces estimates to the following three interactions,
\begin{align*}
&\sum_{k\in\mathbb Z}\operatorname{Re}\int_{\mathbb R^2}\overline{P_k|\nabla|^\sigma\psi}\,P_k|\nabla|^\sigma\left(P_{\leq k-5}A\cdot\nabla P_{[k-4,k+4]}\psi\right)\,\mathrm{d}x\\
&\quad+\sum_{k\in\mathbb Z}\operatorname{Re}\int_{\mathbb R^2}\overline{P_k|\nabla|^\sigma\psi}\,P_k|\nabla|^\sigma\left(P_{[k-4,k+4]}A\cdot\nabla P_{\leq k-5}\psi\right)\,\mathrm{d}x\\
&\quad+\sum_{k\in\mathbb Z}\operatorname{Re}\int_{\mathbb R^2}\overline{P_k|\nabla|^\sigma\psi}\,\sum_{j\geq k-4}P_k|\nabla|^\sigma\left(P_jA\cdot\nabla P_{[j-8,j+8]}\psi\right)\,\mathrm{d}x.
\end{align*}

\noindent
\textit{Low-High interaction.}
For the first term, we use
\[
P_k|\nabla|^\sigma\left(P_{\leq k-5}A\cdot\nabla P_{[k-4,k+4]}\psi\right)=\left[P_k|\nabla|^\sigma, P_{\leq k-5}A\cdot\nabla\right]P_{[k-4,k+4]}\psi+P_{\leq k-5}A\cdot\nabla P_k|\nabla|^\sigma\psi.
\]
The contribution of the second term on the right satisfies
\begin{align*}
\left|\sum_{k\in\mathbb Z}\operatorname{Re}\int_{\mathbb R^2}\overline{P_k|\nabla|^\sigma\psi}\,P_{\leq k-5}A\cdot\nabla P_k|\nabla|^\sigma\psi\,\mathrm{d}x\right|&=\frac{1}{2}\left|\sum_{k\in\mathbb Z}\int_{\mathbb R^2}\nabla\cdot P_{\leq k-5}A\left|P_k|\nabla|^\sigma\psi\right|^2\,\mathrm{d}x\right|\\
&\lesssim\|\nabla A\|_{L_x^\infty}\sum_{k\in\mathbb Z}\|P_k|\nabla|^\sigma\psi\|_{L_x^2}^2\\
&\lesssim\|\nabla A\|_{L_x^\infty}\|\psi\|_{\dot H_x^\sigma}^2.
\end{align*}

We now estimate the commutator. Its physical space expression is
\begin{align*}
&\left[P_k|\nabla|^\sigma, P_{\leq k-5}A\cdot\nabla\right]P_{[k-4,k+4]}\psi(x)\\
&\quad=\frac{1}{2\pi}\int_{\mathbb R^2}\big(P_{\leq k-5}A(y)-P_{\leq k-5}A(x)\big)\cdot\mathcal F^{-1}\big(|\cdot|^\sigma\varphi_k\big)(x-y)\nabla P_{[k-4,k+4]}\psi(y)\,\mathrm{d}y\\
&\quad=\frac{1}{2\pi}\int_0^1\int_{\mathbb R^2}\nabla P_{\leq k-5}A\big(x+\theta(y-x)\big)(y-x)\cdot\mathcal F^{-1}\big(|\cdot|^\sigma\varphi_k\big)(x-y)\nabla P_{[k-4,k+4]}\psi(y)\,\mathrm{d}y\,\mathrm{d}\theta.
\end{align*}
Here the last equality follows from the fundamental theorem of calculus. By H\"older inequality, Young inequality, and Bernstein inequality,
\begin{align*}
&\left\|\left[P_k|\nabla|^\sigma, P_{\leq k-5}A\cdot\nabla\right]P_{[k-4,k+4]}\psi\right\|_{L_x^2}\\
&\quad\lesssim\|\nabla A\|_{L_x^\infty}\left\|2^kx\,\mathcal F^{-1}\big(|\cdot|^\sigma\varphi_k\big)(x)\right\|_{L_x^1}\left\|2^{-k}\nabla P_{[k-4,k+4]}\psi\right\|_{L_x^2}.
\end{align*}
The scaling of the Littlewood-Paley kernel gives \(\left\|2^kx\mathcal F^{-1}\big(|\cdot|^\sigma\varphi_k\big)(x)\right\|_{L_x^1}\lesssim2^{\sigma k}\), hence, by Bernstein inequality,
\begin{align*}
&\left|\sum_{k\in\mathbb Z}\operatorname{Re}\int_{\mathbb R^2}\overline{P_k|\nabla|^\sigma\psi}\,\left[P_k|\nabla|^\sigma, P_{\leq k-5}A\cdot\nabla\right]P_{[k-4,k+4]}\psi\,\mathrm{d}x\right|\\
&\quad\lesssim\|\nabla A\|_{L_x^\infty}\sum_{k\in\mathbb Z}\|P_k|\nabla|^\sigma\psi\|_{L_x^2}2^{\sigma k}\|P_{[k-4,k+4]}\psi\|_{L_x^2}\\
&\quad\lesssim\|\nabla A\|_{L_x^\infty}\|\psi\|_{\dot H_x^\sigma}^2.
\end{align*}
Combining the two parts, the entire Low-High interaction is bounded by \(\|\nabla A\|_{L_x^\infty}\|\psi\|_{\dot H_x^\sigma}^2\).

\medskip
\noindent\textit{High-Low interaction.}
By Cauchy-Schwarz inequality,
\begin{align*}
&\left|\sum_{k\in\mathbb Z}\operatorname{Re}\int_{\mathbb R^2}\overline{P_k|\nabla|^\sigma\psi}\,P_k|\nabla|^\sigma\left(P_{[k-4,k+4]}A\cdot\nabla P_{\leq k-5}\psi\right)\,\mathrm{d}x\right|\\
&\quad\lesssim\|\psi\|_{\dot H_x^\sigma}\left(\sum_{k\in\mathbb Z}2^{2\sigma k}\left\|P_{[k-4,k+4]}A\cdot\nabla P_{\leq k-5}\psi\right\|_{L_x^2}^2\right)^{\frac12}.
\end{align*}
For every \(k\in\mathbb Z\),
\begin{equation*}
2^{\sigma k}\left\|P_{[k-4,k+4]}A\cdot\nabla P_{\leq k-5}\psi\right\|_{L_x^2}\lesssim 2^{(\sigma+1)k}\|P_{[k-4,k+4]}A\|_{L_x^2}\|\psi\|_{L_x^\infty},
\end{equation*}
hence
\begin{align*}
&\left|\sum_{k\in\mathbb Z}\operatorname{Re}\int_{\mathbb R^2}\overline{P_k|\nabla|^\sigma\psi}\,P_k|\nabla|^\sigma\left(P_{[k-4,k+4]}A\cdot\nabla P_{\leq k-5}\psi\right)\,\mathrm{d}x\right|\\
&\quad\lesssim\|\psi\|_{\dot H_x^\sigma}\|\psi\|_{L_x^\infty}\left(\sum_{k\in\mathbb Z}2^{2(\sigma+1)k}\|P_{[k-4,k+4]}A\|_{L_x^2}^2\right)^{\frac12}\\
&\quad\lesssim\|\psi\|_{\dot H_x^\sigma}\|\psi\|_{L_x^\infty}\|A\|_{\dot H_x^{\sigma+1}}.
\end{align*}

\noindent\textit{High-High interaction.}
By Cauchy-Schwarz inequality in \(k\), this contribution is bounded by
\begin{align*}
\|\psi\|_{\dot H_x^\sigma}\left(\sum_{k\in\mathbb Z}\left[2^{\sigma k}\left\|\sum_{j\geq k-4}P_k\left(P_jA\cdot\nabla P_{[j-8,j+8]}\psi\right)\right\|_{L_x^2}\right]^2\right)^{\frac12}.
\end{align*}
For the second factor, H\"older inequality and Bernstein inequality give
\begin{align*}
&\left(\sum_{k\in\mathbb Z}\left[2^{\sigma k}\left\|\sum_{j\geq k-4}P_k\left(P_jA\cdot\nabla P_{[j-8,j+8]}\psi\right)\right\|_{L_x^2}\right]^2\right)^{\frac12}\\
&\quad\lesssim\left(\sum_{k\in\mathbb Z}\left[\sum_{j\geq k-4}2^{\sigma(k-j)}2^{(\sigma+1)j}\|P_jA\|_{L_x^2}2^{-j}\|\nabla P_{[j-8,j+8]}\psi\|_{L_x^\infty}\right]^2\right)^{\frac12}\\
&\quad\lesssim\|\psi\|_{L_x^\infty}\left(\sum_{k\in\mathbb Z}\left[\sum_{j\geq k-4}2^{\sigma(k-j)}2^{(\sigma+1)j}\|P_jA\|_{L_x^2}\right]^2\right)^{\frac12}\\
&\quad\lesssim \|\psi\|_{L_x^\infty}\|A\|_{\dot H_x^{\sigma+1}}.
\end{align*}
Hence the entire High-High interaction satisfies
\begin{equation*}
\left|\sum_{k\in\mathbb Z}\operatorname{Re}\int_{\mathbb R^2}\overline{P_k|\nabla|^\sigma\psi}\,\sum_{j\geq k-4}P_k|\nabla|^\sigma\left(P_jA\cdot\nabla P_{[j-8,j+8]}\psi\right)\,\mathrm{d}x\right| \lesssim\|\psi\|_{\dot H_x^\sigma}\|\psi\|_{L_x^\infty}\|A\|_{\dot H_x^{\sigma+1}}.
\end{equation*}
Combining the Low-High, High-Low, and High-High estimates, \eqref{booassum}, \eqref{eAonederivativel21}, and \eqref{eAonederivativelinfty}, we obtain
\[
\begin{aligned}
  \left|\sum_k\operatorname{Re}\int_{\RR^2}\overline{P_k|\nabla|^\sigma\psi_m}\,P_k|\nabla|^\sigma(A\cdot\nabla\psi_m)\,\mathrm{d}x\right|&\lesssim\|\nabla A\|_{L_x^\infty}\|\psi\|_{\dot H_x^\sigma}^2+\|A\|_{\dot H_x^{\sigma+1}}\|\psi\|_{L_x^\infty}\|\psi\|_{\dot H_x^\sigma}\\
  &\lesssim\epsilon_1^2t^{-2}(1+\log t)\|\psi(t)\|_{\dot H_x^\sigma}^2.
\end{aligned}
\]

For the $F\psi$ terms, the fractional Leibniz rule, \eqref{econnectionenve1}, \eqref{edkdef2}, \eqref{eenergyenvelope3}, \eqref{eAlinftydecay}, \eqref{eAonederivativel21}, \eqref{eAonederivativelinfty}, and \eqref{booassum} give
\[
\begin{aligned}
  \max_F\bigl\||\nabla|^\sigma(F\psi)\bigr\|_{L_x^2}&\lesssim\left(\|\nabla A\|_{L_x^\infty}+\|A\|_{L_x^\infty}^2+\|\psi\|_{L_x^\infty}^2\right)\|\psi\|_{\dot H_x^\sigma}\\&\quad+\|\psi\|_{L_x^\infty}\|A\|_{\dot H_x^{\sigma+1}}+\|A\|_{L_x^\infty}\|\psi\|_{L_x^\infty}\|A\|_{\dot H_x^\sigma}\\
  &\lesssim\epsilon_1^2t^{-2}(1+\log t)\|\psi(t)\|_{\dot H_x^\sigma}.
\end{aligned}
\]
For the $A_3\psi$ term, the fractional Leibniz rule, \eqref{eA32}, \eqref{eA36}, and \eqref{booassum} give
\begin{equation*}
\bigl\||\nabla|^\sigma(A_3\psi)\bigr\|_{L_x^2}\lesssim\|A_3\|_{L_x^\infty}\|\psi\|_{\dot H_x^\sigma}+\|A_3\|_{\dot H_x^\sigma}\|\psi\|_{L_x^\infty}\lesssim\epsilon_1^2t^{-2}(1+\log t)\|\psi(t)\|_{\dot H_x^\sigma}.
\end{equation*}

Combining these estimates and using Littlewood-Paley theory proves the differential estimate for the dyadic energy. Applying Gronwall inequality on $[1,t]$ with $\sigma=\epsilon$ and using conservation of the $L_x^2$ norm, we obtain
\begin{equation*}
\|\psi(t)\|_{H_x^\epsilon}\lesssim\exp\left(C\epsilon_1^2\int_1^t\tau^{-2}(1+\log\tau)\,\mathrm{d}\tau\right)\|\psi(1)\|_{H_x^\epsilon}\lesssim\|\psi(1)\|_{H_x^\epsilon},\qquad 1\leq t\leq T.
\end{equation*}
This completes the proof.
\end{proof}

\section{Decay estimate $\mathrm{I}$: improvement of $\|\hat{f}\|_{L_\xi^\infty}$}
\label{sec:decayestimate1}

In this section we improve the Fourier bootstrap bound \eqref{booassum} for $\|\widehat f\|_{L_\xi^\infty}$. We omit the initial heat parameter $s=0$ and indices.
\begin{proposition}\label{phatf improve}
Under the bootstrap assumption \eqref{booassum}, there holds
\begin{equation}\label{ehatfimprove}
\|\widehat f(t)\|_{L_\xi^\infty}\leq\|\widehat f(1)\|_{L_\xi^\infty}+C\epsilon_1^3.
\end{equation}
\end{proposition}

\begin{proof}
By \eqref{ecaloric7} and \eqref{eFdefinition}, the Duhamel formula gives
\begin{equation}\label{efhat}
\widehat f(t,\xi)=\widehat f(1,\xi)-i\int_1^t e^{i\tau|\xi|^2}\mathcal{F}\bigl[-2iA\cdot\nabla\psi+A_3\psi+F\psi\bigr](\tau,\xi)\,\mathrm{d}\tau.
\end{equation}

We first estimate the $F\psi$ terms. By \eqref{eAlinftydecay}, \eqref{eAl2decay}, \eqref{eAhigh1}, \eqref{eAhigh2}, and \eqref{booassum},
\begin{align*}
\|\widehat{AA\psi}(\tau)\|_{L_\xi^\infty}&\lesssim\|A(\tau)\|_{L_x^\infty}\|A(\tau)\|_{L_x^2}\|\psi(\tau)\|_{L_x^2}\lesssim\epsilon_1^5\tau^{-\frac52}(1+\log\tau),\\
\|\widehat{(\nabla\cdot A^{(\geq4)})\psi}(\tau)\|_{L_\xi^\infty}&\lesssim\|A^{(\geq4)}(\tau)\|_{L_x^2}^{\frac12}\|\nabla^2A^{(\geq4)}(\tau)\|_{L_x^2}^{\frac12}\|\psi(\tau)\|_{L_x^2}\lesssim\epsilon_1^5\tau^{-\frac94}.
\end{align*}

For the cubic terms, we use the identity
\begin{equation}\label{eJproduct}
  J_{0,i}(\psi\overline\psi\psi)=(J_{0,i}\psi)\overline\psi\psi-\psi\overline{J_{0,i}\psi}\psi+\psi\overline\psi J_{0,i}\psi.
\end{equation}
Applying this identity twice, \eqref{booassum} and \eqref{bsc10} give
\[
\begin{aligned}
  \||\psi|^2\psi\|_{L_x^2}&\lesssim\|\psi\|_{L_x^\infty}^2\|\psi\|_{L_x^2}\lesssim\epsilon_1^3\tau^{-2},\\
  \|J_0^2(|\psi|^2\psi)\|_{L_x^2}&\lesssim\|\psi\|_{L_x^\infty}^2\|J_0^2\psi\|_{L_x^2}+\|\psi\|_{L_x^\infty}\|J_0\psi\|_{L_x^4}^2\lesssim\epsilon_1^3\tau^{-2+\beta},
\end{aligned}
\]
weighted interpolation and unitarity yield
\begin{equation}
  \|\widehat{|\psi|^2\psi}(\tau)\|_{L_\xi^\infty}\lesssim\|e^{-i\tau\Delta}(|\psi|^2\psi)\|_{L_x^1}\lesssim\||\psi|^2\psi\|_{L_x^2}^{\frac{1}{2}}\|J_0^2(|\psi|^2\psi)\|_{L_x^2}^{\frac{1}{2}}
  \lesssim\epsilon_1^3\tau^{-2+\frac{\beta}{2}}. \label{epsipsipsiestimate}
\end{equation}

To estimate $(\nabla\cdot A^{(2)})\psi$, we use \eqref{eAcom2}. For $a=1,2$ and $0<s\leq1$, \eqref{booassum} gives
\begin{equation*}
\|\nabla^ae^{s\Delta}\psi(\tau)\|_{L_x^2}\lesssim s^{-\frac{a-\epsilon}{2}}\|\psi(\tau)\|_{H^\epsilon}\lesssim\epsilon_1s^{-\frac{a-\epsilon}{2}}.
\end{equation*}
For $a=0,1,2$ and $s>0$, \eqref{eheatoperator}, \eqref{eJs4}, and \eqref{booassum} also give
\begin{align*}
\|\nabla^ae^{s\Delta}\psi(\tau)\|_{L_x^\infty}&\lesssim\epsilon_1\tau^{-1}s^{-\frac a2},\\
\|J_0^2\nabla^ae^{s\Delta}\psi(\tau)\|_{L_x^2}&=\|x^2\nabla^ae^{s\Delta}f(\tau)\|_{L_x^2}\lesssim\epsilon_1(\tau^\beta+s)s^{-\frac a2},
\end{align*}
where the weighted estimate follows by distributing $x^2$ over the heat convolution. By \eqref{eJproduct}, H\"older inequality, and Gagliardo-Nirenberg inequality,
\begin{align*}
\left\|\overline{\nabla^ae^{s\Delta}\psi}\,\nabla^{2-a}e^{s\Delta}\psi\,\psi\right\|_{L_x^2}&\lesssim\epsilon_1^3\tau^{-2}\min\left\{s^{-1+\frac\epsilon2},s^{-1}\right\},\\
\left\|J_0^2\left(\overline{\nabla^ae^{s\Delta}\psi}\,\nabla^{2-a}e^{s\Delta}\psi\,\psi\right)\right\|_{L_x^2}&\lesssim\epsilon_1^3\tau^{-2}(\tau^\beta+s)s^{-1}.
\end{align*}
Applying the weighted interpolation in the proof of \eqref{epsipsipsiestimate} on $0<s\leq\tau$, and using $\|\nabla^ae^{s\Delta}\psi(\tau)\|_{L_x^2}\lesssim\epsilon_1s^{-\frac{a+1}{2}}$ on $s\geq\tau$, we obtain
\begin{align*}
\|\widehat{(\nabla\cdot A^{(2)})\psi}(\tau)\|_{L_\xi^\infty}&\lesssim\epsilon_1^3\tau^{-2}\int_0^\tau(\tau^\beta+s)^{\frac12}\min\left\{s^{-1+\frac\epsilon4},s^{-1}\right\}\,\mathrm{d}s+\epsilon_1^3\tau^{-1}\int_\tau^\infty s^{-2}\,\mathrm{d}s\\
&\lesssim\epsilon_1^3\tau^{-2}\left(\tau^{\frac\beta2}(1+\log\tau)+\tau^{\frac12}\right)\lesssim\epsilon_1^3\tau^{-\frac32}.
\end{align*}

For $A_3^{(2)}\psi$, the same heat integral argument applies with $a=1$, by \eqref{eA33}. Hence
\begin{equation*}
\|\widehat{A_3^{(2)}\psi}(\tau)\|_{L_\xi^\infty}\lesssim\epsilon_1^3\tau^{-2}\left(\tau^{\frac{\beta}{2}}(1+\log\tau)+\tau^{\frac{1}{2}}\right)\lesssim\epsilon_1^3\tau^{-\frac{3}{2}}.
\end{equation*}
By \eqref{eA35} and \eqref{booassum},
\begin{equation*}
\|\widehat{A_3^{(\geq4)}\psi}(\tau)\|_{L_\xi^\infty}\lesssim\|A_3^{(\geq4)}(\tau)\|_{L_x^2}\|\psi(\tau)\|_{L_x^2}\lesssim\epsilon_1^5\tau^{-2}.
\end{equation*}

Combining these estimates, with $F$ as in \eqref{eFdefinition},
\begin{equation}\label{ehatf1}
  \int_1^t\left(\max_F\|\widehat{F\psi}(\tau)\|_{L_\xi^\infty}+\|\widehat{A_3\psi}(\tau)\|_{L_\xi^\infty}\right)\,\mathrm{d}\tau\lesssim\epsilon_1^3\int_1^t\tau^{-\frac32}\,\mathrm{d}\tau\lesssim\epsilon_1^3.
\end{equation}

It remains to estimate the magnetic term $A\cdot\nabla\psi$. Up to a fixed factor, its contribution to the profile is
\begin{equation}\label{frofile}
\int_1^t\int_{\RR^2}e^{i\tau(|\xi|^2-|\xi-\eta|^2)}(\xi-\eta)\cdot\widehat A(\tau,\eta)\widehat f(\tau,\xi-\eta)\,\mathrm{d}\eta\,\mathrm{d}\tau.
\end{equation}
The quadratic phase is
\begin{equation*}
\Phi_2(\xi,\eta):=|\xi|^2-|\xi-\eta|^2=2\xi\cdot\eta-|\eta|^2.
\end{equation*}
After inserting \eqref{eAcom2} into the contribution of $A^{(2)}$, the resulting cubic phase is
\begin{equation*}
\Phi_3(\xi,\eta,\zeta):=|\xi|^2-|\xi-\eta|^2+|\eta-\zeta|^2-|\zeta|^2=2\eta\cdot(\xi-\zeta).
\end{equation*}
For the quadratic phase, the time resonance and space resonance sets are
\begin{align*}
\{(\xi,\eta):\Phi_2(\xi,\eta)=0\}&=\{(\xi,\eta):2\xi\cdot\eta=|\eta|^2\},\\
\{(\xi,\eta):\nabla_\eta\Phi_2(\xi,\eta)=0\}&=\{(\xi,\eta):\eta=\xi\},
\end{align*}
and their intersection is $\{(0,0)\}$. For the cubic phase,
\begin{equation*}
\{(\xi,\eta,\zeta):\Phi_3(\xi,\eta,\zeta)=0\}=\{(\xi,\eta,\zeta):\eta\cdot(\xi-\zeta)=0\},
\end{equation*}
whereas
\begin{equation*}
\nabla_\eta\Phi_3=2(\xi-\zeta),\qquad\nabla_\zeta\Phi_3=-2\eta.
\end{equation*}
Thus its space resonance and spacetime resonance sets are both $\{(\xi,\eta,\zeta):\eta=0,\ \zeta=\xi\}$. We do not integrate by parts in time. For the quadratic connection contribution, we estimate the region where $|\xi-\eta|$ is small by support volume and integrate twice in $\eta$ on its complement. When both resulting $\eta$ derivatives fall on the inner phase, we split into a neighborhood of the space resonance set of $\Phi_3$ and its complement, on the complement we integrate twice in the frequency direction generated by $(\nabla_\eta\Phi_3,\nabla_\zeta\Phi_3)$.

We first estimate the higher order magnetic contribution. To this end, we establish the spatial moment bound for $A^{(2)}$. Placing $x$ on the undifferentiated factor in \eqref{eAcom2}, we use
\begin{equation*}
xe^{s\Delta}\psi(\tau)=e^{s\Delta}J_0\psi(\tau)-2(s+i\tau)\nabla e^{s\Delta}\psi(\tau).
\end{equation*}
By \eqref{eUfree2}, \eqref{eJs4}, and the free heat estimates following from \eqref{booassum}, we obtain
\begin{equation}\label{ehatf2}
\|xA^{(2)}(\tau)\|_{L_x^2}\lesssim\epsilon_1^2\left[\tau^\gamma\int_0^\infty\frac{s^{-\frac12}}{\tau+s}\,\mathrm{d}s+\int_0^1s^{-1+\frac\epsilon2}\,\mathrm{d}s+\int_1^\infty s^{-\frac32}\,\mathrm{d}s\right]\lesssim\epsilon_1^2.
\end{equation}
Therefore, \eqref{exA1} and \eqref{eAcom1} imply
\begin{equation*}
\|xA^{(\geq4)}(\tau)\|_{L_x^2}\leq\|xA(\tau)\|_{L_x^2}+\|xA^{(2)}(\tau)\|_{L_x^2}\lesssim\epsilon_1^2.
\end{equation*}
We split the outer profile frequency at $|\xi-\eta|=\tau^{\frac12}$. By Cauchy-Schwarz inequality, \eqref{eAhigh1}, and \eqref{booassum}, the low frequency contribution satisfies
\begin{align*}
&\int_1^t\left\|\int_{\RR^2}e^{i\tau\Phi_2(\xi,\eta)}\chi_l(\tau^{-\frac12}|\xi-\eta|)(\xi-\eta)\cdot\widehat{A^{(\geq4)}}(\tau,\eta)\widehat f(\tau,\xi-\eta)\,\mathrm{d}\eta\right\|_{L_\xi^\infty}\,\mathrm{d}\tau\\
&\quad\lesssim\int_1^t\tau^{\frac{1-\epsilon}{2}}\|A^{(\geq4)}(\tau)\|_{L_x^2}\|f(\tau)\|_{H^\epsilon}\,\mathrm{d}\tau\lesssim\epsilon_1^5\int_1^t\tau^{-\frac32-\frac\epsilon2}\,\mathrm{d}\tau\lesssim\epsilon_1^5.
\end{align*}
For the high frequency contribution, using $\nabla_\eta e^{i\tau\Phi_2(\xi,\eta)}=2i\tau(\xi-\eta)e^{i\tau\Phi_2(\xi,\eta)}$, one integration by parts in $\eta$ gives
\begin{align*}
&\int_1^t\left\|\int_{\RR^2}e^{i\tau\Phi_2(\xi,\eta)}\chi_h(\tau^{-\frac12}|\xi-\eta|)(\xi-\eta)\cdot\widehat{A^{(\geq4)}}(\tau,\eta)\widehat f(\tau,\xi-\eta)\,\mathrm{d}\eta\right\|_{L_\xi^\infty}\,\mathrm{d}\tau\\
&\quad\lesssim\int_1^t\tau^{-1}\left(\tau^{-\frac\epsilon2}\|xA^{(\geq4)}\|_{L_x^2}\|f\|_{H^\epsilon}+\tau^{-\frac12}\|A^{(\geq4)}\|_{L_x^2}\|f\|_{L_x^2}+\|A^{(\geq4)}\|_{L_x^2}\|xf\|_{L_x^2}\right)\,\mathrm{d}\tau\lesssim\epsilon_1^3.
\end{align*}
Combining \eqref{efhat} and \eqref{ehatf1} with the above bound for the contribution of $A^{(\geq4)}$ in \eqref{frofile}, we obtain
\begin{equation}\label{ehatf3}
\begin{aligned}
\|\widehat f(t)-\widehat f(1)\|_{L_\xi^\infty}&\lesssim\epsilon_1^3+\left\|\int_1^t\int_{\RR^2}e^{i\tau(|\xi|^2-|\xi-\eta|^2)}(\xi-\eta)\cdot\widehat{A^{(2)}}(\tau,\eta)\widehat f(\tau,\xi-\eta)\,\mathrm{d}\eta\,\mathrm{d}\tau\right\|_{L_\xi^\infty}.
\end{aligned}
\end{equation}
It remains to estimate the quadratic magnetic contribution.

\subsection{The outer low frequency region}

The quadratic connection admits the same Fourier estimate as \eqref{eAhat1}. Indeed, by \eqref{eAcom2}, using the $H^\epsilon$ bound for $0<s\leq1$ and the Fourier $L^\infty$ bound for $s\geq1$ in \eqref{booassum}, we obtain
\begin{equation*}
\|\widehat{A^{(2)}}(\tau)\|_{L_\xi^\infty}\lesssim\int_0^\infty\|e^{s\Delta}\psi(\tau)\|_{L_x^2}\|\nabla e^{s\Delta}\psi(\tau)\|_{L_x^2}\,\mathrm{d}s\lesssim\epsilon_1^2\left(\int_0^1s^{-\frac{1-\epsilon}{2}}\,\mathrm{d}s+\int_1^\infty s^{-\frac32}\,\mathrm{d}s\right)\lesssim\epsilon_1^2.
\end{equation*}

\noindent\textit{The outer profile has frequency $|\xi-\eta|\lesssim\tau^{-\frac13-\delta}$.}
We insert $\chi_l(\tau^{\frac13+\delta}|\xi-\eta|)$ in \eqref{ehatf3}. By this bound and \eqref{booassum}, uniformly in $\xi$,
\begin{equation*}
\int_1^t\int_{\RR^2}\chi_l(\tau^{\frac13+\delta}|\xi-\eta|)|\xi-\eta|\,|\widehat{A^{(2)}}(\tau,\eta)|\,|\widehat f(\tau,\xi-\eta)|\,\mathrm{d}\eta\,\mathrm{d}\tau\lesssim\epsilon_1^3\int_1^t\tau^{-1-3\delta}\,\mathrm{d}\tau\lesssim\epsilon_1^3.
\end{equation*}

\subsection{The first integration by parts in the outer high frequency region}

We insert $\chi_h(\tau^{\frac13+\delta}|\xi-\eta|)$, which is supported where $|\xi-\eta|\gtrsim\tau^{-\frac13-\delta}$. Since $\nabla_\eta\Phi_2(\xi,\eta)=2(\xi-\eta)$, the derivative in the magnetic nonlinearity gives
\begin{align}
&\int_{\RR^2}e^{i\tau\Phi_2(\xi,\eta)}\chi_h(\tau^{\frac13+\delta}|\xi-\eta|)(\xi-\eta)\cdot\widehat{A^{(2)}}(\tau,\eta)\widehat f(\tau,\xi-\eta)\,\mathrm{d}\eta\notag\\
&\quad=-\frac{1}{2i\tau}\int_{\RR^2}e^{i\tau\Phi_2(\xi,\eta)}\nabla_\eta\cdot\Bigl[\chi_h(\tau^{\frac13+\delta}|\xi-\eta|)\widehat{A^{(2)}}(\tau,\eta)\widehat f(\tau,\xi-\eta)\Bigr]\,\mathrm{d}\eta.\label{ehatf4}
\end{align}
By \eqref{eAl2decay} and \eqref{eAhigh1},
\begin{equation}\label{ehatf5}
\|A^{(2)}(\tau)\|_{L_x^2}\leq\|A(\tau)\|_{L_x^2}+\|A^{(\geq4)}(\tau)\|_{L_x^2}\lesssim\epsilon_1^2\tau^{-1}(1+\log\tau).
\end{equation}

\smallskip
\noindent\textit{The frequency derivative falls on the cutoff.}
By \eqref{booassum} and \eqref{ehatf5},
\begin{align*}
&\int_1^t\tau^{-1}\int_{\RR^2}\left|\nabla_\eta\chi_h(\tau^{\frac13+\delta}|\xi-\eta|)\right|\,|\widehat{A^{(2)}}(\tau,\eta)|\,|\widehat f(\tau,\xi-\eta)|\,\mathrm{d}\eta\,\mathrm{d}\tau\\
&\quad\lesssim\int_1^t\tau^{-1+\frac13+\delta}\|A^{(2)}(\tau)\|_{L_x^2}\|f(\tau)\|_{L_x^2}\,\mathrm{d}\tau\lesssim\epsilon_1^3\int_1^t\tau^{-\frac53+\delta}(1+\log\tau)\,\mathrm{d}\tau\lesssim\epsilon_1^3.
\end{align*}

\smallskip
\noindent\textit{The frequency derivative falls on the profile.}
Similarly,
\begin{align*}
&\int_1^t\tau^{-1}\int_{\RR^2}\chi_h(\tau^{\frac13+\delta}|\xi-\eta|)|\widehat{A^{(2)}}(\tau,\eta)|\,|\nabla_\eta\widehat f(\tau,\xi-\eta)|\,\mathrm{d}\eta\,\mathrm{d}\tau\\
&\quad\lesssim\int_1^t\tau^{-1}\|A^{(2)}(\tau)\|_{L_x^2}\|xf(\tau)\|_{L_x^2}\,\mathrm{d}\tau\lesssim\epsilon_1^3\int_1^t\tau^{-2+\gamma}(1+\log\tau)\,\mathrm{d}\tau\lesssim\epsilon_1^3.
\end{align*}

The only contribution from \eqref{ehatf4} which is not integrable directly is therefore
\begin{equation}\label{ehatf6}
-\frac{1}{2i}\int_1^t\tau^{-1}\int_{\RR^2}e^{i\tau\Phi_2(\xi,\eta)}\chi_h(\tau^{\frac13+\delta}|\xi-\eta|)\nabla_\eta\cdot\widehat{A^{(2)}}(\tau,\eta)\widehat f(\tau,\xi-\eta)\,\mathrm{d}\eta\,\mathrm{d}\tau.
\end{equation}
Indeed, by \eqref{ehatf2}, estimating this term before the second integration by parts would only give
\begin{equation*}
\int_1^t\tau^{-1}\|xA^{(2)}(\tau)\|_{L_x^2}\|f(\tau)\|_{L_x^2}\,\mathrm{d}\tau\lesssim\epsilon_1^3\int_1^t\tau^{-1}\,\mathrm{d}\tau,
\end{equation*}
which is not uniformly bounded.

\subsection{The second integration by parts for the remaining $xA^{(2)}$-term}

We apply the second frequency integration by parts only to
\eqref{ehatf6}. Using
\[
  e^{i\tau\Phi_2(\xi,\eta)}=\frac{\xi-\eta}{2i\tau|\xi-\eta|^2}\cdot\nabla_\eta e^{i\tau\Phi_2(\xi,\eta)},
\]
we obtain, for each fixed $\tau$,
\begin{align*}
  &-\frac{1}{2i\tau}\int_{\RR^2}e^{i\tau\Phi_2(\xi,\eta)}\chi_h(\tau^{\frac{1}{3}+\delta}|\xi-\eta|)\nabla_\eta\cdot\widehat{A^{(2)}}(\tau,\eta)\widehat f(\tau,\xi-\eta)\,\mathrm{d}\eta\notag\\
  &\quad=\frac{1}{(2i\tau)^2}\int_{\RR^2}e^{i\tau\Phi_2(\xi,\eta)}\nabla_\eta\cdot\Bigg[\chi_h(\tau^{\frac{1}{3}+\delta}|\xi-\eta|)\frac{\xi-\eta}{|\xi-\eta|^2}\nabla_\eta\cdot\widehat{A^{(2)}}(\tau,\eta)\widehat f(\tau,\xi-\eta)\Bigg]\,\mathrm{d}\eta.
\end{align*}

\smallskip
\noindent\textit{The second derivative falls on the symbol or on the cutoff.}
On the support of a derivative of the cutoff, $|\xi-\eta|\simeq\tau^{-\frac{1}{3}-\delta}$.
Hence,
\[
  \left|\nabla_\eta\left[\chi_h(\tau^{\frac{1}{3}+\delta}|\xi-\eta|)\frac{\xi-\eta}{|\xi-\eta|^2}\right]\right|\lesssim\tau^{\frac{2}{3}+2\delta}.
\]
By \eqref{booassum} and \eqref{ehatf2},
\begin{align*}
  &\int_1^t\tau^{-2}\int_{\RR^2}\left|\nabla_\eta\left[\chi_h(\tau^{\frac{1}{3}+\delta}|\xi-\eta|)\frac{\xi-\eta}{|\xi-\eta|^2}\right]\right||\nabla_\eta\widehat{A^{(2)}}(\tau,\eta)||\widehat f(\tau,\xi-\eta)|\,\mathrm{d}\eta\,\mathrm{d}\tau\\
  &\quad\lesssim\int_1^t\tau^{-2+\frac{2}{3}+2\delta}\|xA^{(2)}(\tau)\|_{L_x^2}\|f(\tau)\|_{L_x^2}\,\mathrm{d}\tau\lesssim\epsilon_1^3\int_1^t\tau^{-\frac{4}{3}+2\delta}\,\mathrm{d}\tau\lesssim\epsilon_1^3.
\end{align*}

\smallskip
\noindent\textit{The second derivative falls on the profile.}
By the same bounds,
\begin{align*}
  &\int_1^t\tau^{-2}\int_{\RR^2}\chi_h(\tau^{\frac{1}{3}+\delta}|\xi-\eta|)|\xi-\eta|^{-1}|\nabla_\eta\widehat{A^{(2)}}(\tau,\eta)||\nabla_\eta\widehat f(\tau,\xi-\eta)|\,\mathrm{d}\eta\,\mathrm{d}\tau\\
  &\quad\lesssim\int_1^t\tau^{-2+\frac{1}{3}+\delta}\|xA^{(2)}(\tau)\|_{L_x^2}\|xf(\tau)\|_{L_x^2}\,\mathrm{d}\tau\lesssim\epsilon_1^3\int_1^t\tau^{-\frac{5}{3}+\delta+\gamma}\,\mathrm{d}\tau\lesssim\epsilon_1^3.
\end{align*}

\smallskip
\noindent\textit{The second derivative falls on $A^{(2)}$.}
Up to a fixed numerical factor, the only remaining contribution is
\begin{equation}
  \int_1^t\tau^{-2}\int_{\RR^2}e^{i\tau\Phi_2(\xi,\eta)}\chi_h(\tau^{\frac{1}{3}+\delta}|\xi-\eta|)\frac{\xi-\eta}{|\xi-\eta|^2}\cdot\nabla_\eta\bigl(\nabla_\eta\cdot\widehat{A^{(2)}}(\tau,\eta)\bigr)\widehat f(\tau,\xi-\eta)\,\mathrm{d}\eta\,\mathrm{d}\tau. \label{ehatf7}
\end{equation}
We now insert \eqref{eAcom2} into this term.

\subsection{The two frequency derivatives of the quadratic connection}

Omitting fixed coefficients and component indices, it suffices to estimate the following contribution to $\widehat{A^{(2)}}(\tau,\eta)$:
\begin{equation}
  \int_{\RR^2}e^{i\tau(|\eta-\zeta|^2-|\zeta|^2)}\frac{\zeta}{|\eta-\zeta|^2+|\zeta|^2}\overline{\widehat f(\tau,\zeta-\eta)}\widehat f(\tau,\zeta)\,\mathrm{d}\zeta. \label{ehatf8}
\end{equation}
The conjugate term is estimated in the same way. For the conjugate term, after a change of the inner integration variable, the numerator in \eqref{ehatf8} is replaced by $\eta-\zeta$, up to a fixed factor. Under the inner $\zeta$ integral, each $\eta$ derivative may be replaced by $\nabla_\eta+\nabla_\zeta$, since the additional $\zeta$ derivative integrates to zero. This operator fixes $\eta-\zeta$ and differentiates the inner phase with gradient $-2\zeta$, so the following estimates apply after interchanging the two internal factors.

The two \(\eta\)-derivatives in \eqref{ehatf7} can fall on the inner time phase, the heat symbol, or the internal profile \(\overline{\widehat f(\zeta-\eta)}\).  Since \(\widehat f(\zeta)\) is independent of \(\eta\), neither derivative falls on that factor. Direct expansion gives
\begin{align}
  &\nabla_\eta^2\Bigg[e^{i\tau(|\eta-\zeta|^2-|\zeta|^2)}\frac{\zeta}{|\eta-\zeta|^2+|\zeta|^2}\overline{\widehat f(\zeta-\eta)}\widehat f(\zeta)\Bigg]\notag\\
  &=e^{i\tau(|\eta-\zeta|^2-|\zeta|^2)}\Bigg[-4\tau^2(\eta-\zeta)\otimes(\eta-\zeta)\frac{\zeta}{|\eta-\zeta|^2+|\zeta|^2}\overline{\widehat f(\zeta-\eta)}\widehat f(\zeta)\notag\\
  &\quad+2i\tau\,\mathrm{Id}\frac{\zeta}{|\eta-\zeta|^2+|\zeta|^2}\overline{\widehat f(\zeta-\eta)}\widehat f(\zeta)\notag\\
  &\quad+2i\tau(\eta-\zeta)\otimes\nabla_\eta\frac{\zeta}{|\eta-\zeta|^2+|\zeta|^2}\overline{\widehat f(\zeta-\eta)}\widehat f(\zeta)\notag\\
  &\quad+2i\tau\nabla_\eta\frac{\zeta}{|\eta-\zeta|^2+|\zeta|^2}\otimes(\eta-\zeta)\overline{\widehat f(\zeta-\eta)}\widehat f(\zeta)\notag\\
  &\quad+2i\tau\frac{\zeta}{|\eta-\zeta|^2+|\zeta|^2}(\eta-\zeta)\otimes\nabla_\eta\overline{\widehat f(\zeta-\eta)}\widehat f(\zeta)\notag\\
  &\quad+2i\tau\frac{\zeta}{|\eta-\zeta|^2+|\zeta|^2}\nabla_\eta\overline{\widehat f(\zeta-\eta)}\otimes(\eta-\zeta)\widehat f(\zeta)\notag\\
  &\quad+\nabla_\eta^2\frac{\zeta}{|\eta-\zeta|^2+|\zeta|^2}\overline{\widehat f(\zeta-\eta)}\widehat f(\zeta)\notag\\
  &\quad+\nabla_\eta\frac{\zeta}{|\eta-\zeta|^2+|\zeta|^2}\otimes\nabla_\eta\overline{\widehat f(\zeta-\eta)}\widehat f(\zeta)\notag\\
  &\quad+\nabla_\eta\overline{\widehat f(\zeta-\eta)}\otimes\nabla_\eta\frac{\zeta}{|\eta-\zeta|^2+|\zeta|^2}\widehat f(\zeta)\notag\\
  &\quad+\frac{\zeta}{|\eta-\zeta|^2+|\zeta|^2}\nabla_\eta^2\overline{\widehat f(\zeta-\eta)}\widehat f(\zeta)\Bigg]. \label{twde}
\end{align}

\subsubsection{Prior inner \(L^2\) estimates}

The first two inner operators considered below have symbols of order \(-1\). Although they could be treated by a dyadic multiplier decomposition, such an argument would require a separate summation over the internal frequency scales, especially near the joint frequency origin. Instead, we recover the heat representation of the elliptic denominator. This converts the discrete frequency analysis into a continuous integration over the heat scale and allows us to combine heat smoothing with the physical time decay of \(\psi\), yielding the desired decay directly. When the resulting symbol has order zero, the heat scale integral becomes critical and no longer provides a derivative gain. In that case, we preserve the order zero structure and use a multilinear multiplier estimate after frequency localization, thereby avoiding any loss of derivatives.

\smallskip
\noindent\textit{No \(\eta\) derivative is placed on the heat symbol or on an internal profile.}
We keep the heat representation
\[
  \frac{1}{|\eta-\zeta|^2+|\zeta|^2}=\int_0^\infty e^{-s(|\eta-\zeta|^2+|\zeta|^2)}\,\mathrm{d}s.
\]
Recall \eqref{eheatoperator} and \eqref{eUfree1},
\begin{align}
  \|e^{s\Delta}\psi(\tau)\|_{L_x^\infty}&\lesssim\epsilon_1(\tau+s)^{-1},&\|\nabla e^{s\Delta}\psi(\tau)\|_{L_x^2}&\lesssim\epsilon_1s^{-\frac{1}{2}},\notag\\
  \|\nabla^2e^{s\Delta}\psi(\tau)\|_{L_x^\infty}&\lesssim\epsilon_1s^{-1}(\tau+s)^{-1}. \label{ehatf9}
\end{align}
By \eqref{ehatf9}, and Plancherel inequality,
\begin{align}
  &\left\|\int e^{i\tau(|\eta-\zeta|^2-|\zeta|^2)}\frac{\zeta}{|\eta-\zeta|^2+|\zeta|^2}\overline{\widehat f(\tau,\zeta-\eta)}\widehat f(\tau,\zeta)\,\mathrm{d}\zeta\right\|_{L_\eta^2}\notag\\
  &\quad\lesssim\int_0^\infty\|e^{s\Delta}\psi(\tau)\|_{L_x^\infty}\|\nabla e^{s\Delta}\psi(\tau)\|_{L_x^2}\,\mathrm{d}s\notag\\
  &\quad\lesssim\epsilon_1^2\int_0^\infty(\tau+s)^{-1}s^{-\frac{1}{2}}\,\mathrm{d}s\lesssim\epsilon_1^2\tau^{-\frac{1}{2}}. \label{ehatf10}
\end{align}

\smallskip
\noindent\textit{One derivative is placed on the inner time phase and one derivative is placed on the heat symbol.}
We use the exact heat representation
\[
  \frac{1}{\bigl(|\eta-\zeta|^2+|\zeta|^2\bigr)^2}=\int_0^\infty se^{-s(|\eta-\zeta|^2+|\zeta|^2)}\,\mathrm{d}s.
\]
Then,
\begin{equation*}
  (\eta_i-\zeta_i)\partial_{\eta_j}\left(\frac{\zeta}{|\eta-\zeta|^2+|\zeta|^2}\right)=-2(\eta_i-\zeta_i)(\eta_j-\zeta_j)\zeta\int_0^\infty se^{-s(|\eta-\zeta|^2+|\zeta|^2)}\,\mathrm{d}s.
\end{equation*}
Therefore, by \eqref{ehatf9}, and Plancherel inequality,
\begin{align}
  &\left\|\int_{\RR^2}e^{i\tau(|\eta-\zeta|^2-|\zeta|^2)}(\eta_i-\zeta_i)\partial_{\eta_j}\left(\frac{\zeta}{|\eta-\zeta|^2+|\zeta|^2}\right)\overline{\widehat f(\tau,\zeta-\eta)}\widehat f(\tau,\zeta)\,\mathrm{d}\zeta
  \right\|_{L_\eta^2}\notag\\
  &\quad\lesssim\int_0^\infty s\|\nabla^2e^{s\Delta}\psi(\tau)\|_{L_x^\infty}\|\nabla e^{s\Delta}\psi(\tau)\|_{L_x^2}\,\mathrm{d}s\notag\\
  &\quad\lesssim\epsilon_1^2\int_0^\infty(\tau+s)^{-1}s^{-\frac{1}{2}}\,\mathrm{d}s\lesssim\epsilon_1^2\tau^{-\frac{1}{2}}. \label{ehatf11}
\end{align}

\smallskip

\noindent\textit{One derivative is placed on the inner time phase and one derivative is placed on the internal profile.}

Here it suffices to bound
\[
\left\|\int e^{i\tau(|\eta-\zeta|^2-|\zeta|^2)}\frac{(\eta_i-\zeta_i)\zeta_j}{|\eta-\zeta|^2+|\zeta|^2}\nabla_\eta\overline{\widehat f(\tau,\zeta-\eta)}\widehat f(\tau,\zeta)\,\mathrm{d}\zeta\right\|_{L_\eta^2},
\]
recovering the heat integral leads to homogeneity at the joint frequency origin and does not yield any decay gain. Instead, we localize both internal frequencies and preserve the order zero structure of the multiplier
\[
  m(\eta,\zeta)=\frac{(\eta_i-\zeta_i)\zeta_j}{|\eta-\zeta|^2+|\zeta|^2}.
\]
The localized multipliers satisfy
\[
  \left\|m(\eta,\zeta)\varphi_{l_1}(\eta-\zeta)\varphi_{l_2}(\zeta)\right\|_{\mathcal S^\infty}\lesssim
  \begin{cases}
    2^{l_1-l_2}, & l_1\leq l_2-5,\\
    2^{l_2-l_1}, & l_2\leq l_1-5,\\
    1, & |l_1-l_2|\leq4.
  \end{cases}
\]

Suppose first that \(l_1\leq l_2-5\), since the output frequency is \(2^{l_2+O(1)}\), Plancherel equality, Lemma~\ref{l multilinear estimates}, Littlewood-Paley square theorem  and Young inequality give
\begin{align*}
  &\left\|\sum_{l_2}\sum_{l_1\leq l_2-5}\int m(\eta,\zeta)\overline{\widehat{e^{i\tau\Delta}P_{l_1}(xf)}}(\tau,\zeta-\eta)\widehat{P_{l_2}\psi}(\tau,\zeta)\,\mathrm{d}\zeta\right\|_{L_\eta^2}\\
  &\quad\lesssim\left(\sum_{l_2}\left[\sum_{l_1\leq l_2-5}2^{l_1-l_2}\|P_{l_1}(xf)\|_{L^2}\|P_{l_2}\psi\|_{L^\infty}\right]^2\right)^{\frac12}\\
  &\quad\lesssim\|\psi\|_{L^\infty}\|xf\|_{L^2}.
\end{align*}

If \(l_2\leq l_1-5\), the output frequency is \(2^{l_1+O(1)}\). Therefore,
\begin{align*}
  &\left\|\sum_{l_1}\sum_{l_2\leq l_1-5}\int m(\eta,\zeta)\overline{\widehat{e^{i\tau\Delta}P_{l_1}(xf)}}(\tau,\zeta-\eta)\widehat{P_{l_2}\psi}(\tau,\zeta)\,\mathrm{d}\zeta\right\|_{L_\eta^2}\\
  &\quad\lesssim\left(\sum_{l_1}\|P_{l_1}(xf)\|_{L^2}^2\left[\sum_{l_2\leq l_1-5}2^{l_2-l_1}\|P_{l_2}\psi\|_{L^\infty}\right]^2\right)^{\frac12}\\
  &\quad\lesssim\|\psi\|_{L^\infty}\|xf\|_{L^2}.
\end{align*}

It remains to consider $|l_1-l_2|\leq4$. By \eqref{bs5},
\begin{align*}
  &\left\|\sum_{|l_1-l_2|\leq 4}\int m(\eta,\zeta)\overline{\widehat{e^{i\tau\Delta}P_{l_1}(xf)}}(\tau,\zeta-\eta)\widehat{P_{l_2}\psi}(\tau,\zeta)\,\mathrm{d}\zeta\right\|_{L_\eta^2}\\
  &\quad\lesssim\sum_{|l_1-l_2|\leq 4}\|P_{l_1}(xf)\|_{L^2}\|P_{l_2}\psi\|_{L^\infty}\\
  &\quad\lesssim\epsilon_1^2(1+\log\tau)\tau^{-1+\gamma}.
\end{align*}

Combining the three frequency cases, we obtain
\begin{equation}
  \left\|\int e^{i\tau(|\eta-\zeta|^2-|\zeta|^2)}\frac{(\eta_i-\zeta_i)\zeta_j}{|\eta-\zeta|^2+|\zeta|^2}\nabla_\eta\overline{\widehat f(\tau,\zeta-\eta)}\widehat f(\tau,\zeta)\,\mathrm{d}\zeta\right\|_{L_\eta^2}\lesssim\epsilon_1^2(1+\log\tau)\tau^{-1+\gamma}. \label{ehatf12}
\end{equation}

\subsubsection{At least one derivative does not fall on the inner time phase}

All subsequent bounds for the displayed oscillatory integrals are understood in $L_\xi^\infty$. The total oscillation after substituting \eqref{ehatf8} into \eqref{ehatf7} is always written in full as $e^{i\tau\Phi_3(\xi,\eta,\zeta)}$.

\smallskip
\noindent\textit{The two derivatives form the Hessian of the inner time phase.}
The expression to estimate is
\begin{align}\label{ehatf13}
  \int_1^t\tau^{-1}\iint&e^{i\tau\Phi_3(\xi,\eta,\zeta)}\chi_h(\tau^{\frac13+\delta}|\xi-\eta|)\frac{\xi-\eta}{|\xi-\eta|^2}\frac{\zeta}{|\eta-\zeta|^2+|\zeta|^2}\notag\\
  &\qquad\times\overline{\widehat f(\tau,\zeta-\eta)}\widehat f(\tau,\zeta)\widehat f(\tau,\xi-\eta)\,\mathrm{d}\zeta\,\mathrm{d}\eta\,\mathrm{d}\tau.
\end{align}
Using \eqref{ehatf10} for the inner two factors and then \(L^2*L^2\to L^\infty\) in \(\eta\),
\begin{equation*}
  \eqref{ehatf13}\lesssim\int_1^t\tau^{-1}\tau^{\frac13+\delta}(\epsilon_1^2\tau^{-\frac12})\|f(\tau)\|_{L_x^2}\,\mathrm{d}\tau\lesssim\epsilon_1^3\int_1^t\tau^{-\frac76+\delta}\,\mathrm{d}\tau\lesssim\epsilon_1^3.
\end{equation*}

\smallskip
\noindent\textit{One derivative falls on the inner time phase and one derivative falls on the heat symbol.}
The expression to estimate is
\begin{align}\label{ehatf14}
  \int_1^t\tau^{-1}\iint&e^{i\tau\Phi_3(\xi,\eta,\zeta)}\chi_h(\tau^{\frac13+\delta}|\xi-\eta|)\frac{\xi-\eta}{|\xi-\eta|^2}(\eta_i-\zeta_i)\partial_{\eta_j}\frac{\zeta}{|\eta-\zeta|^2+|\zeta|^2}\notag\\
  &\qquad\times\overline{\widehat f(\tau,\zeta-\eta)}\widehat f(\tau,\zeta)\widehat f(\tau,\xi-\eta)\,\mathrm{d}\zeta\,\mathrm{d}\eta\,\mathrm{d}\tau.
\end{align}
By \eqref{ehatf11} and \(L^2*L^2\to L^\infty\),
\begin{equation*}
  \eqref{ehatf14}\lesssim\int_1^t\tau^{-1}\tau^{\frac13+\delta}(\epsilon_1^2\tau^{-\frac12})\|f(\tau)\|_{L_x^2}\,\mathrm{d}\tau\lesssim\epsilon_1^3\int_1^t\tau^{-\frac76+\delta}\,\mathrm{d}\tau\lesssim\epsilon_1^3.
\end{equation*}

\smallskip
\noindent\textit{One derivative falls on the inner time phase and one derivative falls on the internal profile.}
The expression to estimate is
\begin{align}\label{ehatf15}
  \int_1^t\tau^{-1}\iint&e^{i\tau\Phi_3(\xi,\eta,\zeta)}\chi_h(\tau^{\frac13+\delta}|\xi-\eta|)\frac{\xi-\eta}{|\xi-\eta|^2}\frac{(\eta_i-\zeta_i)\zeta_j}{|\eta-\zeta|^2+|\zeta|^2}\notag\\
  &\qquad\times\nabla_\eta\overline{\widehat f(\tau,\zeta-\eta)}\widehat f(\tau,\zeta)\widehat f(\tau,\xi-\eta)\,\mathrm{d}\zeta\,\mathrm{d}\eta\,\mathrm{d}\tau.
\end{align}
By \eqref{ehatf12} and \(L^2*L^2\to L^\infty\),
\begin{equation*}
  \eqref{ehatf15}\lesssim\int_1^t\tau^{-1}\tau^{\frac13+\delta}\bigl(\epsilon_1^2(1+\log\tau)\tau^{-1+\gamma}\bigr)\|f(\tau)\|_{L^2}\,\mathrm{d}\tau\lesssim\epsilon_1^3\int_1^t(1+\log\tau)\tau^{-\frac53+\delta+\gamma}\,\mathrm{d}\tau\lesssim\epsilon_1^3.
\end{equation*}

The remaining three terms contain no derivative on the inner time phase. Direct differentiation gives
\[
\begin{aligned}
  \left|\nabla_\eta^2\frac{\zeta}{|\eta-\zeta|^2+|\zeta|^2}\right|&\lesssim(|\eta-\zeta|^2+|\zeta|^2)^{-\frac32},\\
  \left|\nabla_\eta\frac{\zeta}{|\eta-\zeta|^2+|\zeta|^2}\right|&\lesssim(|\eta-\zeta|^2+|\zeta|^2)^{-1},\\
  \left|\frac{\zeta}{|\eta-\zeta|^2+|\zeta|^2}\right|&\lesssim(|\eta-\zeta|^2+|\zeta|^2)^{-\frac12}.
\end{aligned}
\]

\smallskip
\noindent\textit{Both derivatives fall on the heat symbol.}
The expression to estimate is
\begin{align}\label{e2derivativeheat}
  \int_1^t\tau^{-2}\iint&e^{i\tau\Phi_3(\xi,\eta,\zeta)}\chi_h(\tau^{\frac13+\delta}|\xi-\eta|)\frac{\xi-\eta}{|\xi-\eta|^2}\nabla_\eta^2\frac{\zeta}{|\eta-\zeta|^2+|\zeta|^2}\notag\\
  &\qquad\times\overline{\widehat f(\tau,\zeta-\eta)}\widehat f(\tau,\zeta)\widehat f(\tau,\xi-\eta)\,\mathrm{d}\zeta\,\mathrm{d}\eta\,\mathrm{d}\tau.
\end{align}
If \(|\eta-\zeta|^2+|\zeta|^2\leq1\), then all profiles are placed in \(L^\infty\), and the four-dimensional singularity is integrable:
\begin{align*}
  &\iint_{|\eta-\zeta|^2+|\zeta|^2\leq1}(|\eta-\zeta|^2+|\zeta|^2)^{-\frac32}|\widehat f(\zeta-\eta)||\widehat f(\zeta)||\widehat f(\xi-\eta)|\,\mathrm{d}\zeta\,\mathrm{d}\eta\notag\\
  &\quad\lesssim\|\widehat f\|_{L^\infty}^3\iint_{|\eta-\zeta|^2+|\zeta|^2\leq1}(|\eta-\zeta|^2+|\zeta|^2)^{-\frac32}\,\mathrm{d}\zeta\,\mathrm{d}\eta\lesssim\epsilon_1^3.
\end{align*}
If \(|\eta-\zeta|^2+|\zeta|^2\geq1\), we put the outer profile in \(L^\infty\) and use Cauchy-Schwarz in \((\eta-\zeta,\zeta)\in\RR^4\):
\begin{align*}
  &\iint_{|\eta-\zeta|^2+|\zeta|^2\geq1}(|\eta-\zeta|^2+|\zeta|^2)^{-\frac32}|\widehat f(\zeta-\eta)||\widehat f(\zeta)||\widehat f(\xi-\eta)|\,\mathrm{d}\zeta\,\mathrm{d}\eta\notag\\
  &\quad\lesssim\|\widehat f\|_{L^\infty}\bigl\|(|\eta-\zeta|^2+|\zeta|^2)^{-\frac32}\bigr\|_{L^2_{\eta,\zeta}(|\eta-\zeta|^2+|\zeta|^2\geq1)}\|f\|_{L^2}^2\lesssim\epsilon_1^3.
\end{align*}
Therefore,
\begin{align*}
  &\eqref{e2derivativeheat}\lesssim\int_1^t\tau^{-2}\tau^{\frac13+\delta}\epsilon_1^3\,\mathrm{d}\tau\lesssim\epsilon_1^3\int_1^t\tau^{-\frac53+\delta}\,\mathrm{d}\tau\lesssim\epsilon_1^3.
\end{align*}

\smallskip
\noindent\textit{One derivative falls on the heat symbol and one derivative falls on the internal profile.}
The expression to estimate is
\begin{align}\label{e1symbol1profile}
  \int_1^t\tau^{-2}\iint&e^{i\tau\Phi_3(\xi,\eta,\zeta)}\chi_h(\tau^{\frac13+\delta}|\xi-\eta|)\frac{\xi-\eta}{|\xi-\eta|^2}\nabla_\eta\left(\frac{\zeta}{|\eta-\zeta|^2+|\zeta|^2}\right)\notag\\
  &\qquad\times\partial_{\eta_j}\overline{\widehat f(\tau,\zeta-\eta)}\widehat f(\tau,\zeta)\widehat f(\tau,\xi-\eta)\,\mathrm{d}\zeta\,\mathrm{d}\eta\,\mathrm{d}\tau.
\end{align}

By the paraproduct decomposition, the symbol bound above, Cauchy-Schwarz, and \eqref{bspco1},
\begin{align*}
&\sum_l2^{-2l}\Bigl(\|\widehat{P_{\leq l-5}(xf)}(\tau)\|_{L_\xi^1}\|\widehat{P_lf}(\tau)\|_{L_\xi^1}+\|\widehat{P_l(xf)}(\tau)\|_{L_\xi^1}\|\widehat{P_{\leq l+4}f}(\tau)\|_{L_\xi^1}\Bigr)\\
&\quad\lesssim\epsilon_1\|xf(\tau)\|_{L_x^2}\left(\sum_{l\leq0}2^l+\sum_{l>0}2^{-\epsilon l}\right)\lesssim\epsilon_1^2\tau^\gamma.
\end{align*}
Placing the outer profile in Fourier $L^\infty$, we obtain
\[
  \eqref{e1symbol1profile}\lesssim\epsilon_1^3\int_1^t\tau^{-\frac{5}{3}+\delta+\gamma}\,\mathrm{d}\tau\lesssim\epsilon_1^3.
\]

\smallskip
\noindent\textit{Both derivatives fall on the same internal profile.}

The expression to estimate is
\begin{align}\label{e2derivativesprofile}
  \int_1^t\tau^{-2}\iint&e^{i\tau\Phi_3(\xi,\eta,\zeta)}\chi_h(\tau^{\frac13+\delta}|\xi-\eta|)\frac{\xi-\eta}{|\xi-\eta|^2}\frac{\zeta}{|\eta-\zeta|^2+|\zeta|^2}\notag\\
  &\qquad\times\nabla_\eta^2\overline{\widehat f(\tau,\zeta-\eta)}\widehat f(\tau,\zeta)\widehat f(\tau,\xi-\eta)\,\mathrm{d}\zeta\,\mathrm{d}\eta\,\mathrm{d}\tau.
\end{align}
We use the paraproduct decomposition
\begin{align*}
  &\nabla_\eta^2\overline{\widehat f(\tau,\zeta-\eta)}\widehat f(\tau,\zeta)\\
  &=-\sum_l\overline{\widehat{P_{\leq l-5}(x^2f)}(\tau,\zeta-\eta)}\widehat{P_lf}(\tau,\zeta)-\sum_l\overline{\widehat{P_l(x^2f)}(\tau,\zeta-\eta)}\widehat{P_{\leq l+4}f}(\tau,\zeta).
\end{align*}

For \(l\leq0\), we place the outer profile in Fourier \(L^\infty\) and estimate the two paraproducts separately. If the undifferentiated internal profile has the larger frequency, by
\eqref{bspco1},
\begin{align*}
  &\sum_{l\leq0}2^{-l}\|\widehat{P_{\leq l-5}(x^2f)}(\tau)\|_{L_\xi^1}\|\widehat{P_lf}(\tau)\|_{L_\xi^1} \\
  &\lesssim\sum_{l\leq0}2^{-l}\bigl(2^l\|\widehat{P_{\leq l-5}(x^2f)}(\tau)\|_{L_\xi^2}\bigr)\bigl(\epsilon_1 2^{2l}\bigr)\lesssim\epsilon_1\|x^2f(\tau)\|_{L_x^2}\sum_{l\leq0}2^{2l}\lesssim\epsilon_1^2\tau^\beta.
\end{align*}
If the twice differentiated internal profile has the larger frequency, by \eqref{bspco1} we also have
\begin{equation*}
  \sum_{l\leq 0}2^{-l}\|\widehat{P_l(x^2f)}(\tau)\|_{L_\xi^1}\|\widehat{P_{\leq l+4}f}(\tau)\|_{L_\xi^1}\lesssim\sum_{l\leq0}2^{-l}\bigl(2^l\|\widehat{P_l(x^2f)}(\tau)\|_{L_\xi^2}\bigr)\bigl(\epsilon_1 2^{2l}\bigr)\lesssim\epsilon_1^2\tau^\beta.
\end{equation*}

For \(l>0\), the localized heat multiplier has order \(-1\), and
\[
  \|\frac{\zeta}{|\eta-\zeta|^2+|\zeta|^2}\varphi_{\leq l-5}(\eta-\zeta)\varphi_l(\zeta)\|_{\mathcal{S}^\infty}+\|\frac{\zeta}{|\eta-\zeta|^2+|\zeta|^2}\varphi_l(\eta-\zeta)\varphi_{\leq l+4}(\zeta)\|_{\mathcal{S}^\infty}\lesssim2^{-l}.
\]
At this point, by Plancherel inequality and Lemma~\ref{l multilinear estimates},
\begin{align*}
  &\sum_{l>0}2^{-l}\Bigl(\bigl(\|P_l\psi(\tau)\|_{L_x^\infty}\bigr)\bigl(\|P_{\leq l-5}(x^2f)(\tau)\|_{L_x^2}\bigr)+\bigl(\|P_{\leq l+4}\psi(\tau)\|_{L_x^\infty}\bigr)\bigl(\|P_l(x^2f)(\tau)\|_{L_x^2}\bigr)\Bigr)\\
  &\lesssim\epsilon_1\tau^{-1}\sum_{l>0}2^{-l}\Bigl(\|P_{\leq l-5}(x^2f)(\tau)\|_{L_x^2}+\|P_l(x^2f)(\tau)\|_{L_x^2}\Bigr)\lesssim\epsilon_1^2\tau^{-1+\beta}.
\end{align*}

Combining the two paraproducts, we obtain
\[
  \sum_{l\leq0}2^{-l}\Bigl(\|\widehat{P_{\leq l-5}(x^2f)}(\tau)\|_{L_\xi^1}\|\widehat{P_lf}(\tau)\|_{L_\xi^1}+\|\widehat{P_l(x^2f)}(\tau)\|_{L_\xi^1}\|\widehat{P_{\leq l+4}f}(\tau)\|_{L_\xi^1}\Bigr)\lesssim\epsilon_1^2\tau^\beta
\]
while the positive frequency contribution is bounded by $\epsilon_1^2\tau^{-1+\beta}$. Hence, the full contribution is bounded by
\begin{equation*}
  \eqref{e2derivativesprofile}\lesssim\epsilon_1^2\int_1^t\tau^{-2+\frac13+\delta+\beta}\|\hat{f}(\tau)\|_{L_\xi^\infty}\,\mathrm{d}\tau+\epsilon_1^2\int_1^t\tau^{-3+\frac13+\delta+\beta}\|\hat{f}(\tau)\|_{L_\xi^2}\,\mathrm{d}\tau\lesssim\epsilon_1^3.
\end{equation*}
Thus all terms in \eqref{twde} are controlled except the first one, in which both derivatives fall on the inner time phase.

\subsection{Both derivatives on the inner time phase}

The two factors of \(\tau\) in the first line of \eqref{twde} cancel the outer factor \(\tau^{-2}\).  The remaining term is
\begin{align}\label{complexterm}
  \int_1^t\iint&e^{i\tau\Phi_3(\xi,\eta,\zeta)}\chi_h(\tau^{\frac13+\delta}|\xi-\eta|)\frac{\xi-\eta}{|\xi-\eta|^2}\frac{(\eta-\zeta)\otimes(\eta-\zeta)\,\zeta}{|\eta-\zeta|^2+|\zeta|^2}\notag\\
  &\qquad\times\overline{\widehat f(\tau,\zeta-\eta)}\widehat f(\tau,\zeta)\widehat f(\tau,\xi-\eta)\,\mathrm{d}\zeta\,\mathrm{d}\eta\,\mathrm{d}\tau.
\end{align}

\subsubsection{The near resonance region}

In view of the formulas for \(\nabla_\eta\Phi_3\) and \(\nabla_\zeta\Phi_3\) above, we insert
\[
  \chi_l\left(\tau^{\frac12-\delta}(|\eta|^2+|\xi-\zeta|^2)^{\frac12}\right).
\]

\smallskip
\noindent\textit{The output frequency satisfies \(|\xi|\gg\tau^{-\frac12+\delta}\).}
On this part of the near resonance region,
\begin{equation*}
  |\eta|+|\xi-\zeta|\lesssim\tau^{-\frac12+\delta}\ll|\xi|, \qquad |\xi-\eta|\simeq |\zeta|\simeq |\eta-\zeta|\simeq |\xi|.
\end{equation*}
On this region,
\[
  \|m(\xi,\eta,\zeta)\|_{L^\infty}=\|\frac{\xi-\eta}{|\xi-\eta|^2}\frac{(\eta-\zeta)\otimes(\eta-\zeta)\,\zeta}{|\eta-\zeta|^2+|\zeta|^2}\|_{L^\infty}\lesssim1.
\]
For fixed \(\xi\), the change of variables \((\eta,\zeta)\mapsto(\eta,\xi-\zeta)\) has Jacobian one, and the four dimensional volume of the support of the near resonance cutoff is bounded by \(\tau^{-2+4\delta}\).  Therefore the contribution at time \(\tau\) is bounded by
\[
  \epsilon_1^3\tau^{-2+4\delta}.
\]

\smallskip
\noindent\textit{The output frequency satisfies \(|\xi|\lesssim\tau^{-\frac12+\delta}\).}
On the support of the near resonance cutoff,
\[
  |\zeta|\lesssim|\xi|+|\xi-\zeta|\lesssim\tau^{-\frac12+\delta}, \qquad |\xi-\eta|\lesssim|\xi|+|\eta|\lesssim\tau^{-\frac12+\delta}.
\]
Hence the pointwise symbol satisfies
\[
  \frac{|\eta-\zeta|^2|\zeta|}{(|\eta-\zeta|^2+|\zeta|^2)|\xi-\eta|}\lesssim\frac{|\zeta|}{|\xi-\eta|}\lesssim\frac{\tau^{-\frac12+\delta}}{|\xi-\eta|}.
\]
Placing all three profiles in Fourier \(L^\infty\), we obtain at fixed \(\tau\)
\begin{align*}
  &\iint\chi_h(\tau^{\frac13+\delta}|\xi-\eta|)\chi_l\left(\tau^{\frac12-\delta}(|\eta|^2+|\xi-\zeta|^2)^{\frac12}\right)\notag\\
  &\qquad\times\frac{|\eta-\zeta|^2|\zeta|}{(|\eta-\zeta|^2+|\zeta|^2)|\xi-\eta|}|\widehat f(\tau,\zeta-\eta)||\widehat f(\tau,\zeta)||\widehat f(\tau,\xi-\eta)|\,\mathrm{d}\zeta\,\mathrm{d}\eta\notag\\
  &\lesssim\epsilon_1^3\tau^{-\frac12+\delta}\int_{|\xi-\eta|\lesssim\tau^{-\frac12+\delta}}\frac{\mathrm{d}\eta}{|\xi-\eta|}\int_{|\xi-\zeta|\lesssim\tau^{-\frac12+\delta}}\mathrm{d}\zeta\notag\\
  &\lesssim\epsilon_1^3\tau^{-2+4\delta}.
\end{align*}
Combining the two output frequency cases, the near resonance contribution is bounded by
\begin{align*}
  &\int_1^t\iint\chi_h(\tau^{\frac13+\delta}|\xi-\eta|)\chi_l\left(\tau^{\frac12-\delta}(|\eta|^2+|\xi-\zeta|^2)^{\frac12}\right)\frac{|\eta-\zeta|^2|\zeta|}{(|\eta-\zeta|^2+|\zeta|^2)|\xi-\eta|}\notag\\
  &\qquad\times|\widehat f(\tau,\zeta-\eta)||\widehat f(\tau,\zeta)||\widehat f(\tau,\xi-\eta)|\,\mathrm{d}\zeta\,\mathrm{d}\eta\,\mathrm{d}\tau\notag\\
  &\quad\lesssim\epsilon_1^3\int_1^t\tau^{-2+4\delta}\,\mathrm{d}\tau\lesssim\epsilon_1^3,
\end{align*}
because \(\delta<\frac14\).

\subsubsection{The far resonance region and the exact integration by parts}

It remains to insert
\[
  \chi_h\left(\tau^{\frac12-\delta}(|\eta|^2+|\xi-\zeta|^2)^{\frac12}\right),
\]
which is supported where
\[
  |\eta|^2+|\xi-\zeta|^2\gtrsim\tau^{-1+2\delta}.
\]
On this region we use the exact identity
\begin{equation}
  e^{i\tau\Phi_3(\xi,\eta,\zeta)}=\frac{(\xi-\zeta)\cdot\nabla_\eta-\eta\cdot\nabla_\zeta}{2i\tau(|\eta|^2+|\xi-\zeta|^2)}e^{i\tau\Phi_3(\xi,\eta,\zeta)}. \label{eIBPfrequency}
\end{equation}
After inserting the outer and far resonance cutoffs, two integrations by parts in $(\eta,\zeta)$ give the following expression, with fixed coefficients and tensor contractions omitted:
\begin{align}
  &\frac{1}{(2i\tau)^2}\iint e^{i\tau\Phi_3(\xi,\eta,\zeta)}\bigl((\xi-\zeta)\cdot\nabla_\eta-\eta\cdot\nabla_\zeta\bigr)\Bigg\{\frac{1}{|\eta|^2+|\xi-\zeta|^2}\notag\\
  &\quad\times\bigl((\xi-\zeta)\cdot\nabla_\eta-\eta\cdot\nabla_\zeta\bigr)\Bigg[\frac{1}{|\eta|^2+|\xi-\zeta|^2}\chi_h(\tau^{\frac13+\delta}|\xi-\eta|)\chi_h\left(\tau^{\frac12-\delta}(|\eta|^2+|\xi-\zeta|^2)^{\frac12}\right)\notag\\
  &\quad\times\frac{\xi-\eta}{|\xi-\eta|^2}\frac{(\eta-\zeta)\otimes(\eta-\zeta)\,\zeta}{|\eta-\zeta|^2+|\zeta|^2}\overline{\widehat f(\tau,\zeta-\eta)}\widehat f(\tau,\zeta)\widehat f(\tau,\xi-\eta)\Bigg]\Bigg\}\,\mathrm{d}\zeta\,\mathrm{d}\eta. \label{eIBPtwice}
\end{align}
Here we applied \eqref{eIBPfrequency} twice. A derivative of the outer cutoff is supported where \(|\xi-\eta|\simeq\tau^{-\frac13-\delta}\) and satisfies
\[
  \left|\nabla_{\eta,\zeta}\chi_h(\tau^{\frac13+\delta}|\xi-\eta|)\right|\lesssim|\xi-\eta|^{-1}, \qquad \left|\nabla_{\eta,\zeta}^2\chi_h(\tau^{\frac13+\delta}|\xi-\eta|)\right|\lesssim|\xi-\eta|^{-2}.
\]
A derivative of the far resonance cutoff is supported where \(|\eta|^2+|\xi-\zeta|^2\simeq\tau^{-1+2\delta}\) and satisfies
\[
  \left|\nabla_{\eta,\zeta}\chi_h\left(\tau^{\frac12-\delta}(|\eta|^2+|\xi-\zeta|^2)^{\frac12}\right)\right|\lesssim(|\eta|^2+|\xi-\zeta|^2)^{-\frac12},
\]
\[
  \left|\nabla_{\eta,\zeta}^2\chi_h\left(\tau^{\frac12-\delta}(|\eta|^2+|\xi-\zeta|^2)^{\frac12}\right)\right|\lesssim(|\eta|^2+|\xi-\zeta|^2)^{-1}.
\]
Thus a derivative falling on either cutoff has the same symbol bound as a derivative falling on the corresponding displayed denominator, and these placements are combined below.

After expanding \eqref{eIBPtwice} and using the bounds for derivatives of the far resonance cutoff, the integrand is bounded by
\begin{align}
  \tau^{-2}\Bigg[&(|\eta|^2+|\xi-\zeta|^2)^{-1}\left|\nabla_{\eta,\zeta}^2\left(\chi_h(\tau^{\frac{1}{3}+\delta}|\xi-\eta|)\frac{\xi-\eta}{|\xi-\eta|^2}\frac{(\eta-\zeta)\otimes(\eta-\zeta)\,\zeta}{|\eta-\zeta|^2+|\zeta|^2}\overline{\widehat f(\zeta-\eta)}\widehat f(\zeta)\widehat f(\xi-\eta)\right)\right|\notag\\
  &+(|\eta|^2+|\xi-\zeta|^2)^{-\frac{3}{2}}\left|\nabla_{\eta,\zeta}\left(\chi_h(\tau^{\frac{1}{3}+\delta}|\xi-\eta|)\frac{\xi-\eta}{|\xi-\eta|^2}\frac{(\eta-\zeta)\otimes(\eta-\zeta)\,\zeta}{|\eta-\zeta|^2+|\zeta|^2}\overline{\widehat f(\zeta-\eta)}\widehat f(\zeta)\widehat f(\xi-\eta)\right)\right|\notag\\
  &+(|\eta|^2+|\xi-\zeta|^2)^{-2}\chi_h(\tau^{\frac{1}{3}+\delta}|\xi-\eta|)\frac{|\eta-\zeta|^2|\zeta|}{|\xi-\eta|(|\eta-\zeta|^2+|\zeta|^2)}|\widehat f(\zeta-\eta)||\widehat f(\zeta)||\widehat f(\xi-\eta)|\Bigg]. \label{ehatf16}
\end{align}
The support restrictions of the outer and far resonance cutoffs are maintained throughout.

\subsubsection{Dyadic geometry and localized symbol bounds}

We localize
\begin{align}
  |\xi-\eta|&\simeq2^k,\notag\\
  (|\eta-\zeta|^2+|\zeta|^2)^{\frac12}&\simeq2^l,\notag\\
  (|\eta|^2+|\xi-\zeta|^2)^{\frac12}&\simeq2^j, \label{ehatf17}
\end{align}
where
\begin{equation}
  2^k\gtrsim\tau^{-\frac{1}{3}-\delta}, \qquad 2^j\gtrsim\tau^{-\frac{1}{2}+\delta}. \label{ehatf18}
\end{equation}
The relations among \(k,l,j\) are as follows.

\smallskip
\noindent\textit{The outer frequency is below the internal frequency: \(k\leq l-10\).}
If \(|\eta|\gtrsim2^l\), then immediately \(2^j\gtrsim2^l\).  If \(|\eta|\ll2^l\), then \(|\eta-\zeta|\simeq2^l\), and hence
\begin{align*}
  |\xi-\zeta|=|(\xi-\eta)+(\eta-\zeta)|\simeq2^l.
\end{align*}
Since the reverse inequality \(2^j\lesssim2^l\) follows from the triangle inequality, we have
\[
  k\leq l-10 \quad\Longrightarrow\quad j=l+O(1).
\]

\smallskip
\noindent\textit{The internal frequency is below the outer frequency: \(l\leq k-10\).}
In this case \(|\eta|\lesssim2^l\) and
\begin{align*}
  |\xi-\zeta|=|(\xi-\eta)+(\eta-\zeta)|\simeq 2^k.
\end{align*}
Therefore
\[
  l\leq k-10 \quad\Longrightarrow\quad j=k+O(1).
\]

\smallskip
\noindent\textit{The outer and internal frequencies are comparable: \(|k-l|\leq10\).}
Cancellation may occur, and the only general relation is
\[
  |k-l|\leq10 \quad\Longrightarrow\quad j\leq k+O(1).
\]

On the localization \eqref{ehatf17}, the absolute value of the multiplier
\[
  \frac{\xi-\eta}{|\xi-\eta|^2}\frac{(\eta-\zeta)\otimes(\eta-\zeta)\,\zeta}{|\eta-\zeta|^2+|\zeta|^2}
\]
is bounded by \(2^{l-k}\). If one or two derivatives fall on this multiplier, the corresponding pointwise bounds are
\begin{align}
  2^{l-k}(2^{-k}+2^{-l}), \qquad 2^{l-k}(2^{-k}+2^{-l})^2, \label{ehatf19}
\end{align}
respectively.  Indeed, a derivative of \(\frac{\xi-\eta}{|\xi-\eta|^2}\), or of its localization, costs \(2^{-k}\), while a derivative of the heat symbol, or of its localization, costs \(2^{-l}\). The cutoff derivative bounds proved above show that differentiating the time dependent outer cutoff is covered by the same \(2^{-k}\) loss, whereas differentiating the far resonance cutoff is covered by the same \(2^{-j}\) losses as differentiating the resonance denominators in \eqref{ehatf16}.

We now estimate separately the three relations among \(k,l,j\).  In each case, the two derivatives created by \eqref{eIBPtwice} are first assigned to profiles, to the outer or heat multiplier, or to the resonance denominators and localizations, and only then are the frequency sums performed.  The lower bounds in all three cases are those in \eqref{ehatf18}. All the estimates in the three cases below are obtained after taking absolute values in frequency space and therefore require only the displayed pointwise \(L^\infty\) multiplier bounds, not a bilinear or trilinear \(\mathcal S^\infty\) estimate.

\subsubsection{The case \(k\leq l-10\)}

Here \(j=l+O(1)\). Combining \eqref{ehatf16} and \eqref{ehatf19}, the possible coefficients, with the common factor \(\tau^{-2}\) omitted, are
\begin{align*}
  &2^{-3l-k}&&\text{if no derivative falls on a profile or on the multiplier},\notag\\
  &2^{-2l-2k}&&\text{if no derivative falls on a profile and one falls on the multiplier},\notag\\
  &2^{-l-3k}&&\text{if no derivative falls on a profile and two fall on the multiplier},\notag\\
  &2^{-2l-k}&&\text{if one derivative falls on a profile and none on the multiplier},\notag\\
  &2^{-l-2k}&&\text{if one derivative falls on a profile and one on the multiplier},\notag\\
  &2^{-l-k}&&\text{if both derivatives fall on profiles}.
\end{align*}
Since \(k\leq l-10\), it remains to consider only four cases with the large multiplier bound: no derivative on a profile, only one derivative on a profile, both derivatives on the same profile, and the two derivatives on two different profiles.

\smallskip
\noindent\textit{No derivative falls on a profile.}
In this case we use
\[
  \|m(\xi,\eta,\zeta)\|_{L^\infty}\lesssim \tau^{-2}2^{-l-3k}.
\]
Suppose first that \(k,l\leq0\). We put the outer profile in Fourier \(L^1\), put the two internal profiles in Fourier \(L^\infty\), by \eqref{bspco1},
\begin{align*}
  &\eqref{eIBPtwice}\lesssim\tau^{-2}\sum_{\substack{k\leq l-10\leq-10\\2^k\gtrsim\tau^{-\frac13-\delta}}}2^{-l-3k}\|\widehat{P_kf}\|_{L^1}\|\widehat f\|_{L^\infty}\|\widehat f\|_{L^\infty}2^{2l}\lesssim\tau^{-2}\sum_{\substack{k\leq l-10\leq-10\\
  2^k\gtrsim\tau^{-\frac13-\delta}}}2^{l-3k}\bigl(\epsilon_1^3 2^{2k}\bigr)\\
  &\quad\lesssim\epsilon_1^3\tau^{-2}\sum_{\substack{k\leq0\\2^k\gtrsim\tau^{-\frac13-\delta}}}\sum_{l=k+10}^{0}2^{l-k}\lesssim\epsilon_1^3\tau^{-2}\sum_{\substack{k\leq0\\
  2^k\gtrsim\tau^{-\frac13-\delta}}}2^{-k}\lesssim\epsilon_1^3\tau^{-\frac53+\delta}.
\end{align*}

We next consider $l>0$. At least one internal profile has frequency comparable to $2^l$. By \eqref{booassum} and \eqref{bspco1}, for $k\leq0$,
\begin{align*}
  \eqref{eIBPtwice}&\lesssim\tau^{-2}\sum_{\substack{l>0,\ k\leq0,\ k\leq l-10\\2^k\gtrsim\tau^{-\frac13-\delta}}}2^{-l-3k}\|\widehat{P_kf}\|_{L^1}\|P_lf\|_{L^2}\|f\|_{L^2}\\
  &\lesssim\epsilon_1^3\tau^{-2}\sum_{l>0}2^{-(1+\epsilon)l}\sum_{\substack{k\leq0\\2^k\gtrsim\tau^{-\frac13-\delta}}}2^{-k}\lesssim\epsilon_1^3\tau^{-\frac53+\delta}.
\end{align*}
If $0<k\leq l-10$, then
\begin{align*}
  \eqref{eIBPtwice}&\lesssim\tau^{-2}\sum_{l>0}\sum_{0<k\leq l-10}2^{-l-3k}\|\widehat{P_kf}\|_{L^1}\|P_lf\|_{L^2}\|f\|_{L^2}\\
  &\lesssim\epsilon_1^3\tau^{-2}\sum_{l>0}2^{-(1+\epsilon)l}\sum_{0<k\leq l-10}2^{-(2+\epsilon)k}\lesssim\epsilon_1^3\tau^{-2}.
\end{align*}

\smallskip
\noindent\textit{Only one derivative falls on a profile.}
The remaining derivative may fall on the multiplier, a resonance denominator, or a localization. All these terms are bounded by
\[
  \|m(\xi,\eta,\zeta)\|_{L^\infty}\lesssim\tau^{-2}2^{-l-2k}.
\]

We first suppose that the derivative falls on the outer profile. For \(k,l\leq0\), we put the differentiated outer profile in Fourier \(L^1\) and the two internal profiles in \(L^2\). Since at least one internal profile has frequency comparable to \(2^l\),
\begin{align*}
  \eqref{eIBPtwice}&\lesssim\tau^{-2}\sum_{\substack{k\leq l-10\leq-10\\2^k\gtrsim\tau^{-\frac13-\delta}}}2^{-l-2k}\bigl(2^k\|\widehat{P_k(xf)}\|_{L^2}\bigr)\|P_lf\|_{L^2}\|f\|_{L^2}\\
  &\lesssim\tau^{-2}\sum_{\substack{k\leq l-10\leq-10\\2^k\gtrsim\tau^{-\frac13-\delta}}}2^{-l-2k}\bigl(2^k\|P_k(xf)\|_{L^2}\bigr)\bigl(\epsilon_1^2 2^l\bigr)\\
  &\lesssim\epsilon_1^2\tau^{-2}\sum_{\substack{k\leq0\\2^k\gtrsim\tau^{-\frac13-\delta}}}(1+|k|)2^{-k}\|P_k(xf)\|_{L^2}\\
  &\lesssim\epsilon_1^3(1+\log\tau)\tau^{-\frac53+\delta+\gamma}.
\end{align*}

We next suppose that the derivative falls on an internal profile of frequency comparable to $2^l$. For $k,l\leq0$, we use Fourier $L^1$ for the outer profile and $L^2$ for the two internal profiles:
\begin{align*}
  \eqref{eIBPtwice}&\lesssim\tau^{-2}\sum_{\substack{k\leq l-10\leq-10\\2^k\gtrsim\tau^{-\frac13-\delta}}}2^{-l-2k}\|\widehat{P_kf}\|_{L^1}\|P_l(xf)\|_{L^2}\|f\|_{L^2}\\
  &\lesssim\tau^{-2}\epsilon_1\sum_{\substack{k\leq l-10\leq-10\\2^k\gtrsim\tau^{-\frac13-\delta}}}2^{-l-2k}\bigl(\epsilon_1 2^{2k}\bigr)\|P_l(xf)\|_{L^2}\\
  &\lesssim\epsilon_1^2\tau^{-2}(1+\log\tau)\sum_{\substack{l\leq0\\2^l\gtrsim\tau^{-\frac13-\delta}}}2^{-l}\|P_l(xf)\|_{L^2}\\
  &\lesssim\epsilon_1^3(1+\log\tau)\tau^{-\frac53+\delta+\gamma}.
\end{align*}

It remains to consider the positive-frequency parts. If the derivative falls on the outer profile, $l>0$, and $k\leq0$, then
\begin{align*}
\eqref{eIBPtwice}&\lesssim\tau^{-2}\sum_{\substack{l>0,\ k\leq0,\ k\leq l-10\\2^k\gtrsim\tau^{-\frac13-\delta}}}2^{-l-2k}\bigl(2^k\|P_k(xf)\|_{L^2}\bigr)\|P_lf\|_{L^2}\|f\|_{L^2}\\
&\lesssim\epsilon_1^2\tau^{-2}\sum_{l>0}2^{-(1+\epsilon)l}\sum_{\substack{k\leq0\\2^k\gtrsim\tau^{-\frac13-\delta}}}2^{-k}\|P_k(xf)\|_{L^2}\lesssim\epsilon_1^3\tau^{-\frac53+\delta+\gamma}.
\end{align*}
For $0<k\leq l-10$,
\begin{align*}
\eqref{eIBPtwice}&\lesssim\tau^{-2}\sum_{l>0}\sum_{0<k\leq l-10}2^{-l-2k}\bigl(2^k\|P_k(xf)\|_{L^2}\bigr)\|P_lf\|_{L^2}\|f\|_{L^2}\\
&\lesssim\epsilon_1^2\tau^{-2}\sum_{l>0}\sum_{0<k\leq l-10}2^{-(1+\epsilon)l-k}\|P_k(xf)\|_{L^2}\lesssim\epsilon_1^3\tau^{-2+\gamma}.
\end{align*}

If the derivative falls on the internal profile of frequency \(2^l\), for \(l>0\) and \(k\leq0\),
\begin{align*}
  &\eqref{eIBPtwice}\lesssim\tau^{-2}\sum_{\substack{l>0,\ k\leq0,\ k\leq l-10\\2^k\gtrsim\tau^{-\frac13-\delta}}}2^{-l-2k}\|\widehat{P_kf}\|_{L^1}\|P_l(xf)\|_{L^2}\|f\|_{L^2}\\
  &\quad\lesssim\epsilon_1^2\tau^{-2}(1+\log\tau)\sum_{l>0}2^{-l}\|P_l(xf)\|_{L^2}\lesssim\epsilon_1^3(1+\log\tau)\tau^{-2+\gamma}.
\end{align*}
For $0<k\leq l-10$,
\begin{align*}
\eqref{eIBPtwice}&\lesssim\tau^{-2}\sum_{l>0}\sum_{0<k\leq l-10}2^{-l-2k}\|\widehat{P_kf}\|_{L^1}\|P_l(xf)\|_{L^2}\|f\|_{L^2}\\
&\lesssim\epsilon_1^2\tau^{-2}\sum_{l>0}2^{-l}\|P_l(xf)\|_{L^2}\sum_{0<k\leq l-10}2^{-(1+\epsilon)k}\lesssim\epsilon_1^3\tau^{-2+\gamma}.
\end{align*}
If the derivative falls on the lower frequency internal profile, we place it in $L^2$ and use the low frequency or $H^\epsilon$ bound for the undifferentiated profile of frequency $2^l$. Splitting at $k=0$ and $l=0$, we obtain
\begin{align*}
\eqref{eIBPtwice}&\lesssim\tau^{-2}\sum_{\substack{k\leq l-10\\2^k\gtrsim\tau^{-\frac13-\delta}}}2^{-l-2k}\|\widehat{P_kf}\|_{L^1}\|P_lf\|_{L^2}\|P_{\leq l+4}(xf)\|_{L^2}\\
&\lesssim\epsilon_1^3\tau^{-2+\gamma}\left[\sum_{\substack{k\leq l-10\leq-10\\2^k\gtrsim\tau^{-\frac13-\delta}}}1+(1+\log\tau)\sum_{l>0}2^{-(1+\epsilon)l}+\sum_{l>0}\sum_{0<k\leq l-10}2^{-(1+\epsilon)(k+l)}\right]\\
&\lesssim\epsilon_1^3(1+\log\tau)^2\tau^{-2+\gamma}.
\end{align*}

\smallskip
\noindent\textit{Both derivatives fall on the same profile.}
The pointwise multiplier bound is
\[
  \|m(\xi,\eta,\zeta)\|_{L^\infty}\lesssim\tau^{-2}2^{-l-k}.
\]

Suppose first that both derivatives fall on the outer profile. For \(k\leq l-10\leq-10\), 
\begin{align*}
  \eqref{eIBPtwice}&\lesssim\tau^{-2}\sum_{\substack{k\leq l-10\leq-10\\2^k\gtrsim\tau^{-\frac13-\delta}}}2^{-l-k}\|\widehat{P_k(x^2f)}\|_{L^1}\|\widehat f\|_{L^\infty}^22^{2l}\\
  &\lesssim\epsilon_1^2\tau^{-2}\sum_{\substack{k\leq0\\2^k\gtrsim\tau^{-\frac13-\delta}}}\|P_k(x^2f)\|_{L^2}\sum_{l=k+10}^{0}2^l\lesssim\epsilon_1^3(1+\log\tau)^{\frac12}\tau^{-2+\beta}.
\end{align*}
For $l>0$,
\begin{align*}
\eqref{eIBPtwice}&\lesssim\tau^{-2}\sum_{\substack{l>0,\ k\leq l-10\\2^k\gtrsim\tau^{-\frac13-\delta}}}2^{-l-k}\|\widehat{P_k(x^2f)}\|_{L^1}\|P_lf\|_{L^2}\|f\|_{L^2}\\
&\lesssim\epsilon_1^2\tau^{-2}\sum_{l>0}2^{-(1+\epsilon)l}\sum_{\substack{k\leq l-10\\2^k\gtrsim\tau^{-\frac13-\delta}}}\|P_k(x^2f)\|_{L^2}\lesssim\epsilon_1^3(1+\log\tau)^{\frac12}\tau^{-2+\beta}.
\end{align*}

We next suppose that both derivatives fall on the same internal profile of frequency comparable to $2^l$. If $k\leq0$, we put the outer profile in Fourier $L^1$ and the two internal profiles in $L^2$:
\begin{align*}
  \eqref{eIBPtwice}&\lesssim\tau^{-2}\sum_{\substack{k\leq l-10,\ k\leq0\\2^k\gtrsim\tau^{-\frac13-\delta}}}2^{-l-k}\|\widehat{P_kf}\|_{L^1}\|P_l(x^2f)\|_{L^2}\|f\|_{L^2}\\
  &\lesssim\epsilon_1^2\tau^{-2}\sum_l\|P_l(x^2f)\|_{L^2}\sum_{\substack{k\leq l-10,\ k\leq0\\2^k\gtrsim\tau^{-\frac13-\delta}}}2^{k-l}\lesssim\epsilon_1^3(1+\log\tau)^{\frac12}\tau^{-2+\beta}.
\end{align*}
If $k>0$, \eqref{bspco1} and \eqref{booassum} give
\begin{align*}
\eqref{eIBPtwice}&\lesssim\tau^{-2}\sum_{l>0}\sum_{0<k\leq l-10}2^{-l-k}\|\widehat{P_kf}\|_{L^1}\|P_l(x^2f)\|_{L^2}\|f\|_{L^2}\\
&\lesssim\epsilon_1^2\tau^{-2}\sum_{l>0}\sum_{0<k\leq l-10}2^{-l-\epsilon k}\|P_l(x^2f)\|_{L^2}\lesssim\epsilon_1^3\tau^{-2+\beta}.
\end{align*}
If both derivatives fall on the same lower frequency internal profile, we place it in $L^2$ and use the low frequency or $H^\epsilon$ bound for the undifferentiated profile of frequency $2^l$. Splitting at $k=0$ and $l=0$, we obtain
\begin{align*}
\eqref{eIBPtwice}&\lesssim\tau^{-2}\sum_{\substack{k\leq l-10\\2^k\gtrsim\tau^{-\frac13-\delta}}}2^{-l-k}\|\widehat{P_kf}\|_{L^1}\|P_lf\|_{L^2}\|P_{\leq l+4}(x^2f)\|_{L^2}\\
&\lesssim\epsilon_1^3\tau^{-2+\beta}\left[\sum_{k\leq l-10\leq-10}2^k+\sum_{l>0}\sum_{k\leq0}2^{k-(1+\epsilon)l}+\sum_{l>0}\sum_{0<k\leq l-10}2^{-\epsilon k-(1+\epsilon)l}\right]\lesssim\epsilon_1^3\tau^{-2+\beta}.
\end{align*}

\smallskip
\noindent\textit{The two derivatives fall on two different profiles.}
The pointwise multiplier bound is again
\[
  \|m(\xi,\eta,\zeta)\|_{L^\infty}\lesssim\tau^{-2}2^{-l-k}.
\]

Suppose first that one derivative falls on the outer profile and the other on an internal profile. For $k\leq l-10\leq-10$, we use $\|\widehat{P_{\leq l+4}(xf)}\|_{L_\xi^1} \lesssim2^l\|P_{\leq l+4}(xf)\|_{L_x^2}$,
\begin{align*}
  \eqref{eIBPtwice}&\lesssim\tau^{-2}\sum_{\substack{k\leq l-10\leq-10\\2^k\gtrsim\tau^{-\frac13-\delta}}}2^{-l-k}2^l\|P_k(xf)\|_{L^2}\|P_{\leq l+4}(xf)\|_{L_x^2}\|f\|_{L^2}\\
  &\lesssim\epsilon_1\tau^{-2+\frac13+\delta}(1+\log\tau)\|xf\|_{L^2}^2\lesssim\epsilon_1^3(1+\log\tau)\tau^{-\frac53+\delta+2\gamma}.
\end{align*}

If $l>0$ and the undifferentiated internal profile has frequency comparable to $2^l$, its $H^\epsilon$ norm gives
\begin{align*}
\eqref{eIBPtwice}&\lesssim\tau^{-2}\sum_{\substack{l>0,\ k\leq l-10\\2^k\gtrsim\tau^{-\frac13-\delta}}}2^{-l-k}2^l\|P_k(xf)\|_{L^2}\|P_{\leq l+4}(xf)\|_{L^2}\|P_lf\|_{L^2}\\
&\lesssim\epsilon_1\tau^{-2}\sum_{\substack{l>0,\ k\leq l-10\\2^k\gtrsim\tau^{-\frac13-\delta}}}2^{-k-\epsilon l}\|P_k(xf)\|_{L^2}\|P_{\leq l+4}(xf)\|_{L^2}\lesssim\epsilon_1^3\tau^{-\frac53+\delta+2\gamma}+\epsilon_1^3\tau^{-2+2\gamma}.
\end{align*}
For the complementary case, the differentiated internal profile has frequency $2^l>1$. By \eqref{bspco1}, \(\|\widehat{P_{\leq l+4}f}\|_{L^1}\lesssim\epsilon_12^{(1-\epsilon)l}\). Keeping the factor $2^{-l}$ in the multiplier, we obtain
\begin{align*}
\eqref{eIBPtwice}&\lesssim\tau^{-2}\sum_{\substack{l>0,\ k\leq l-10\\2^k\gtrsim\tau^{-\frac13-\delta}}}2^{-l-k}\|P_k(xf)\|_{L^2}\|P_l(xf)\|_{L^2}\|\widehat{P_{\leq l+4}f}\|_{L^1}\\
&\lesssim\epsilon_1\tau^{-2}\sum_{\substack{l>0,\ k\leq l-10\\2^k\gtrsim\tau^{-\frac13-\delta}}}2^{-k-\epsilon l}\|P_k(xf)\|_{L^2}\|P_l(xf)\|_{L^2}\lesssim\epsilon_1^3\tau^{-\frac53+\delta+2\gamma}+\epsilon_1^3\tau^{-2+2\gamma}.
\end{align*}

Finally, suppose that the two derivatives fall on the two internal profiles, whose frequencies are \(2^l\) and \(2^{l'}\), with \(l'\leq l\). If \(k\leq0\), then by \(\|P_kf\|_{L^2}\lesssim\epsilon_1 2^k\),
\begin{align*}
  \eqref{eIBPtwice}&\lesssim\tau^{-2}\sum_{\substack{k\leq l-10,\ k\leq0\\2^k\gtrsim\tau^{-\frac13-\delta}}}\sum_{l'\leq l}2^{-l-k}\|P_kf\|_{L^2}\|P_l(xf)\|_{L^2}
  \|\widehat{P_{l'}(xf)}\|_{L^1}\\
  &\lesssim\epsilon_1\tau^{-2}\sum_{\substack{k\leq l-10,\ k\leq0\\2^k\gtrsim\tau^{-\frac13-\delta}}}\sum_{l'\leq l}2^{l'-l}\|P_l(xf)\|_{L^2}\|P_{l'}(xf)\|_{L^2}\\
  &\lesssim\epsilon_1\tau^{-2}(1+\log\tau)\sum_l\sum_{l'\leq l}2^{l'-l}\|P_l(xf)\|_{L^2}\|P_{l'}(xf)\|_{L^2}\lesssim\epsilon_1^3(1+\log\tau)\tau^{-2+2\gamma}.
\end{align*}

If $k>0$, necessarily $l>0$. Using \eqref{booassum}, we obtain
\begin{align*}
\eqref{eIBPtwice}&\lesssim\tau^{-2}\sum_{l>0}\sum_{0<k\leq l-10}\sum_{l'\leq l}2^{-l-k}\|P_kf\|_{L^2}\|P_l(xf)\|_{L^2}\|\widehat{P_{l'}(xf)}\|_{L^1}\\
&\lesssim\epsilon_1\tau^{-2}\sum_{l>0}\sum_{0<k\leq l-10}\sum_{l'\leq l}2^{-(1+\epsilon)k}2^{l'-l}\|P_l(xf)\|_{L^2}\|P_{l'}(xf)\|_{L^2}\lesssim\epsilon_1^3\tau^{-2+2\gamma}.
\end{align*}
Combining the four cases, all the contributions of the region \(k\leq l-10\) are integrable on \([1,\infty)\).

\subsubsection{The case \(l\leq k-10\)}

Here \(j=k+O(1)\). Moreover,
\[
  |\eta|+|\eta-\zeta|+|\zeta|\lesssim2^l\ll 2^k
\]
implies
\[
  |\xi|\simeq|\xi-\eta|\simeq 2^k.
\]
Since the Fourier \(L^\infty_\xi\) estimate is carried out for each fixed output frequency \(\xi\), the frequency \(k\) is determined up to \(O(1)\) dyadic values. Thus, there is no summation over \(k\) in this subsection. Combining \eqref{ehatf16} and \eqref{ehatf19}, the possible coefficients, with the common factor \(\tau^{-2}\) omitted, are
\begin{align*}
  &2^{l-5k}&&\text{if no derivative falls on a profile or on the multiplier},\notag\\
  &2^{-4k}&&\text{if no derivative falls on a profile and one falls on the multiplier},\notag\\
  &2^{-l-3k}&&\text{if no derivative falls on a profile and two fall on the multiplier},\notag\\
  &2^{l-4k}&&\text{if one derivative falls on a profile and none on the multiplier},\notag\\
  &2^{-3k}&&\text{if one derivative falls on a profile and one on the multiplier},\notag\\
  &2^{l-3k}&&\text{if both derivatives fall on profiles}.
\end{align*}
Since \(l\leq k-10\), we still only need to consider four cases, recall that \(2^k\gtrsim \tau^{-\frac{1}{3}-\delta}\).

\smallskip
\noindent\textit{No derivative falls on a profile.}

In this case,
\[
  \|m(\xi,\eta,\zeta)\|_{L^\infty}\lesssim\tau^{-2}2^{-l-3k}.
\]
Arguing as in the previous subsection, and using that \(k\) is fixed by the output frequency, we obtain
\[
  \eqref{eIBPtwice}\lesssim\epsilon_1^3\tau^{-2}.
\]

\smallskip
\noindent\textit{Only one derivative falls on a profile.}

In this case,
\[
  \|m(\xi,\eta,\zeta)\|_{L^\infty}\lesssim \tau^{-2}2^{-3k}.
\]
The same Fourier \(L^1\times L^2\times L^2\) and Fourier \(L^\infty\times L^1\times L^1\) estimates used above give
\[
  \eqref{eIBPtwice}\lesssim\epsilon_1^3\tau^{-\frac53+\delta+\gamma}+\epsilon_1^3\tau^{-2+\gamma}.
\]

\smallskip
\noindent\textit{Both derivatives fall on the same profile.}

In this case,
\[
  \|m(\xi,\eta,\zeta)\|_{L^\infty}\lesssim \tau^{-2}2^{l-3k}.
\]
The estimates for \(\widehat{P_k(x^2f)}\) and \(\widehat{P_l(x^2f)}\) are the same as in the previous subsection, with the factors \(2^{l-k}\) and \(2^{2(l-k)}\), respectively.
Hence,
\[
  \eqref{eIBPtwice}\lesssim \epsilon_1^3\tau^{-2+\beta}.
\]

\smallskip
\noindent\textit{The two derivatives fall on two different profiles.}

The pointwise multiplier bound is again
\[
  \|m(\xi,\eta,\zeta)\|_{L^\infty}\lesssim\tau^{-2}2^{l-3k}.
\]
Suppose first that one derivative falls on the outer profile and the other on an internal profile. If the differentiated internal profile has the larger internal frequency \(2^l\), then by
\(
  \sum_{l'\leq l}\|P_{l'}f\|_{L^2}\lesssim\epsilon_1,
\)
\begin{align*}
  \eqref{eIBPtwice}&\lesssim\tau^{-2}\sum_{l\leq k-10}2^{l-3k}\|\widehat{P_k(xf)}\|_{L^1}\|P_l(xf)\|_{L^2}\sum_{l'\leq l}\|P_{l'}f\|_{L^2}\\
  &\lesssim\epsilon_1\tau^{-2}2^{-k}\|P_k(xf)\|_{L^2}\sum_{l\leq k-10}2^{l-k}\|P_l(xf)\|_{L^2}\lesssim\epsilon_1^3\tau^{-\frac53+\delta+2\gamma}.
\end{align*}

If the differentiated internal profile has the lower frequency \(2^{l'}\), \(l'\leq l\), we put it in Fourier \(L^1\). By \eqref{booassum},
\begin{align*}
\eqref{eIBPtwice}&\lesssim\tau^{-2}\sum_{l\leq k-10}\sum_{l'\leq l}2^{l-3k}\|P_k(xf)\|_{L_x^2}\|P_lf\|_{L_x^2}\|\widehat{P_{l'}(xf)}\|_{L_\xi^1}\\
&\lesssim\epsilon_1\tau^{-2}\|P_k(xf)\|_{L_x^2}\|xf\|_{L_x^2}\left(2^{-3k}\sum_{\substack{l\leq k-10\\l\leq0}}2^{3l}+2^{-3k}\sum_{0<l\leq k-10}2^{(2-\epsilon)l}\right)\lesssim\epsilon_1^3\tau^{-2+2\gamma}.
\end{align*}

Finally, suppose that the derivatives fall on the two internal profiles, whose frequencies are \(2^l\) and \(2^{l'}\), with \(l'\leq l\). We put the outer profile in Fourier $L^\infty$ and both differentiated internal profiles in Fourier \(L^1\). Then
\begin{align*}
  \eqref{eIBPtwice}&\lesssim\tau^{-2}\sum_{l\leq k-10}\sum_{l'\leq l}2^{l-3k}\|\widehat{P_kf}\|_{L^\infty}\|\widehat{P_l(xf)}\|_{L^1}\|\widehat{P_{l'}(xf)}\|_{L^1}\\
  &\lesssim\epsilon_1\tau^{-2}\sum_{l\leq k-10}\sum_{l'\leq l}2^{3(l-k)}2^{l'-l}\|P_l(xf)\|_{L^2}\|P_{l'}(xf)\|_{L^2}\lesssim\epsilon_1^3\tau^{-2+2\gamma}.
\end{align*}

Combining the four cases, the contribution of this subsection is bounded by
\[
  \epsilon_1^3\left(\tau^{-\frac{5}{3}+\delta+2\gamma}+\tau^{-2+\beta}+\tau^{-2+2\gamma}\right),
\]
all these powers are integrable on \([1,\infty)\).

\subsubsection{The case \(|k-l|\leq 10\)}

In this case \(j\leq k+O(1)\). We separate
\[
  j\leq k-20 \qquad\text{and}\qquad j=k+O(1).
\]
In the first region,
\[
  |\eta|+|\xi-\zeta|\ll2^k,
\]
and hence
\[
  |\xi-\eta|\simeq|\eta-\zeta|\simeq|\zeta|\simeq|\xi|.
\]
so, for fixed \(\xi\), only \(O(1)\) values of \(k\) occur in the region \(j\leq k-20\), and only the summation in \(j\) remains. In the region \(j=k+O(1)\), there is only one diagonal summation, which may be indexed by either \(j\) or \(k\). Moreover, if one profile has frequency \(2^{l'}\leq2^{k-10}\), this frequency comparability rules out the region \(j\leq k-20\). Hence \(j=k+O(1)\), and the remaining diagonal summation will be indexed by \(k\).

Since \(\lvert k-l\rvert\leq10\), the base multiplier has order zero. Combining \eqref{ehatf16} and \eqref{ehatf19}, the relevant coefficients, with the common factor \(\tau^{-2}\) omitted, are
\[
  2^{-4j}, \qquad 2^{-3j}, \qquad 2^{-2j},
\]
when zero, one, or two derivatives fall on profiles, respectively. A derivative falling on the outer or heat multiplier costs \(2^{-k}\lesssim2^{-j}\), and is therefore covered by the same bounds.

\smallskip
\noindent\textit{No derivative falls on a profile.}

In this case,
\[
  \|m(\xi,\eta,\zeta)\|_{L^\infty}\lesssim\tau^{-2}2^{-4j}.
\]
For \(k\leq0\), note that
\(
  \sum_{\substack{j\leq k+O(1)\\2^j\gtrsim\tau^{-\frac12+\delta}}}1\lesssim1+\log\tau.
\)
Hence
\begin{align*}
  \eqref{eIBPtwice}&\lesssim\tau^{-2}\sum_{\substack{j\leq k+O(1)\\2^j\gtrsim\tau^{-\frac12+\delta}}}2^{-4j}2^{4j}\|\widehat f\|_{L^\infty}^3\lesssim\epsilon_1^3(1+\log\tau)\tau^{-2}.
\end{align*}
For \(k>0\), we place the two profiles of frequencies \(2^{k+O(1)}\) in Fourier \(L^2\). Then
\begin{align*}
  \eqref{eIBPtwice}&\lesssim\epsilon_1^3\tau^{-2}\sum_{\substack{j\leq k+O(1)\\2^j\gtrsim\tau^{-\frac12+\delta}}}2^{-2j}\lesssim\epsilon_1^3\tau^{-1-2\delta}.
\end{align*}

\smallskip
\noindent\textit{Only one derivative falls on a profile.}

In this case,
\[
  \|m(\xi,\eta,\zeta)\|_{L^\infty}\lesssim\tau^{-2}2^{-3j}.
\]
If the differentiated profile has frequency \(2^{k+O(1)}\), then
\begin{align*}
  \eqref{eIBPtwice}&\lesssim\epsilon_1^2\tau^{-2}\sum_{\substack{j\leq k+O(1)\\2^j\gtrsim\tau^{-\frac12+\delta}}}2^{-j}\|\widehat{P_k(xf)}\|_{L^2}\lesssim\epsilon_1^3\tau^{-\frac32-\delta+\gamma}.
\end{align*}

If the differentiated profile has frequency \(2^{l'}\leq2^{k-10}\), then \(j=k+O(1)\), and hence there is no \(j\)-summation. Hence,
\begin{align*}
  \eqref{eIBPtwice}&\lesssim\sum_{2^k\gtrsim\tau^{-\frac13-\delta}}\tau^{-2}2^{-3k}\|\widehat{P_kf}\|_{L^2}^2\sum_{l'\leq k-10}\|\widehat{P_{l'}(xf)}\|_{L^1}\\
  &\lesssim\sum_{2^k\gtrsim\tau^{-\frac13-\delta}}\tau^{-2}2^{-2k}\|P_kf\|_{L^2}^2\|xf\|_{L^2}\lesssim\epsilon_1^3(1+\log\tau)\tau^{-2+\gamma}.
\end{align*}

\smallskip
\noindent\textit{Both derivatives fall on the same profile.}

In this case,
\[
  \|m(\xi,\eta,\zeta)\|_{L^\infty}\lesssim\tau^{-2}2^{-2j}.
\]
Suppose first that the twice differentiated profile has frequency \(2^{k+O(1)}\). For \(k\leq0\),
\begin{align*}
  \eqref{eIBPtwice}&\lesssim\epsilon_1^2\tau^{-2}\sum_{\substack{j\leq k+O(1)\\2^j\gtrsim\tau^{-\frac12+\delta}}}2^{-2j}2^{2j}\|\widehat{P_k(x^2f)}\|_{L^2}\lesssim\epsilon_1^3(1+\log\tau)\tau^{-2+\beta}.
\end{align*}

For \(k>0\), one undifferentiated profile has frequency \(2^{k+O(1)}\). Using its \(H^\epsilon\) bound and \eqref{booassum}, we obtain
\begin{align*}
\eqref{eIBPtwice}&\lesssim\epsilon_1^2\tau^{-2}\sum_{k>0}\sum_{\substack{j\leq k+O(1)\\2^j\gtrsim\tau^{-\frac12+\delta}}}2^{-\epsilon k}\|\widehat{P_k(x^2f)}\|_{L_\xi^2}\\
&\lesssim\epsilon_1^3\tau^{-2+\beta}\sum_{k>0}(1+\log\tau+k)2^{-\epsilon k}\lesssim\epsilon_1^3(1+\log\tau)\tau^{-2+\beta}.
\end{align*}

If the twice differentiated profile has frequency \(2^{l'}\leq2^{k-10}\), then \(j=k+O(1)\). We obtain
\begin{align*}
  \eqref{eIBPtwice}&\lesssim\tau^{-2}\sum_{2^k\gtrsim\tau^{-\frac13-\delta}}2^{-2k}\|\widehat{P_kf}\|_{L^2}^2\sum_{l'\leq k-10}\|\widehat{P_{l'}(x^2f)}\|_{L^1} \\
  &\lesssim\tau^{-2}\|x^2f\|_{L^2}\sum_{2^k\gtrsim\tau^{-\frac13-\delta}}2^{-k}\|P_kf\|_{L^2}^2\lesssim\epsilon_1^3\tau^{-2+\beta}.
\end{align*}

\smallskip
\noindent\textit{The two derivatives fall on two different profiles.}

The pointwise multiplier bound is again
\[
  \|m(\xi,\eta,\zeta)\|_{L^\infty}\lesssim\tau^{-2}2^{-2j}.
\]
Suppose first that both differentiated profiles have frequency \(2^{k+O(1)}\). For \(k\leq0\), 
\[
  \eqref{eIBPtwice}\lesssim\epsilon_1\tau^{-2}\sum_{\substack{j\leq k+O(1)\\2^j\gtrsim\tau^{-\frac12+\delta}}}2^{-2j}2^{2j}\|\widehat{P_k(xf)}\|_{L^2}\|\widehat{P_k(xf)}\|_{L^2}\lesssim\epsilon_1^3(1+\log\tau)\tau^{-2+2\gamma}.
\]
For \(k>0\), 
\[
  \eqref{eIBPtwice}\lesssim\tau^{-2}\sum_{\substack{j\leq k+O(1),\ k>0\\2^j\gtrsim\tau^{-\frac12+\delta}}}2^{-j}\|\widehat{P_k(xf)}\|_{L^2}\|\widehat{P_k(xf)}\|_{L^2}\|f\|_{L^2} \lesssim\epsilon_1^3\tau^{-\frac32-\delta+2\gamma}.
\]

If one differentiated profile has frequency \(2^{l'}\leq2^{k-10}\), then \(j=k+O(1)\). Placing the lower-frequency differentiated profile in Fourier \(L^1\), we obtain
\begin{align*}
  \eqref{eIBPtwice}&\lesssim\tau^{-2}\sum_{2^k\gtrsim\tau^{-\frac13-\delta}}2^{-2k}\|\widehat{P_kf}\|_{L^2}\|\widehat{P_k(xf)}\|_{L^2}\sum_{l'\leq k-10}\|\widehat{P_{l'}(xf)}\|_{L^1}\\
  &\lesssim\tau^{-2}\|xf\|_{L^2}\sum_{2^k\gtrsim\tau^{-\frac13-\delta}}2^{-k}\|P_kf\|_{L^2}\|P_k(xf)\|_{L^2} \\
  &\lesssim\epsilon_1^3(1+\log\tau)\tau^{-2+2\gamma}.
\end{align*}
This covers both the outer-internal and the two-internal placements.

Combining the four cases, the contribution of this subsection is bounded by
\[
  C\epsilon_1^3(1+\log\tau)\left(\tau^{-1-2\delta}+\tau^{-\frac32-\delta+\gamma}+\tau^{-2+\beta}+\tau^{-2+2\gamma}+\tau^{-\frac32-\delta+2\gamma}\right).
\]
All these powers are integrable on \([1,\infty)\).

\subsection{Conclusion}

Combining the above estimates proves \eqref{ehatfimprove} and completes the proof.
\end{proof}

\section{Decay estimate $\mathrm{II}$: improvement of $\|xf\|_{L_x^2}$}\label{sec:decayestimate2}

We now improve the first weighted estimate for the profile. At $s=0$ we write $J=J_0$. By \eqref{eJs4}, $\|J\psi(t)\|_{L_x^2}=\|xf(t)\|_{L_x^2}$.

\begin{proposition}\label{xfimp}
Under the bootstrap assumption \eqref{booassum}, there holds
\begin{equation}\label{exf1}
\|xf(t)\|_{L_x^2}\lesssim\|xf(1)\|_{L_x^2}+\epsilon_1^3t^\gamma,\qquad t\geq1.
\end{equation}
\end{proposition}

\begin{proof}
Since \([i\partial_t+\Delta,J_i]=0\), we use
\begin{equation}\label{exf2}
  J_i(A_l\partial_l\psi_m)=A_l\partial_lJ_i\psi_m+2it(\partial_iA_l)\partial_l\psi_m-A_i\psi_m.
\end{equation}
Applying \(J_i\) to \eqref{ecaloric7} gives
\begin{equation}
\begin{aligned}
  (i\partial_t+\Delta)J_i\psi_m&=-2iA_l\partial_lJ_i\psi_m+(A_3+A_lA_l-i\partial_lA_l)J_i\psi_m\\
  &\quad+2iA_i\psi_m+4t(\partial_iA_l)\partial_l\psi_m+2it\partial_i(A_3+A_lA_l)\psi_m\\
  &\quad+2t(\partial_i\partial_lA_l)\psi_m-iJ_i\bigl(\psi_l\operatorname{Im}(\overline{\psi_l}\psi_m)\bigr).
\end{aligned}
\label{exf3}
\end{equation}
Since \(A\) is real-valued, the magnetic terms in \eqref{exf3} satisfy, by integration by parts,
\[
  \operatorname{Im}\int_{\RR^2}\overline{J_i\psi_m}\bigl[-2iA_l\partial_lJ_i\psi_m-i(\partial_lA_l)J_i\psi_m\bigr]\,\mathrm{d}x=0.
\]
The real-valued potential \(A_3+A_lA_l\) also contributes zero to the energy identity. Hence H\"older inequality and \eqref{eJs4} yield
\begin{align*}
\frac{\mathrm{d}}{\mathrm{d}t}\|xf(t)\|_{L_x^2}&\lesssim\|A\|_{L_x^2}\|\psi\|_{L_x^\infty}+t\|(\nabla A)\nabla\psi\|_{L_x^2}+t\|A\|_{L_x^\infty}\|\nabla A\|_{L_x^2}\|\psi\|_{L_x^\infty}+t\|(\nabla^2A)\psi\|_{L_x^2}\\
&\quad+\left\|J\bigl(\psi_l\operatorname{Im}(\overline{\psi_l}\psi_m)\bigr)\right\|_{L_x^2}+t\|(\nabla A_3)\psi\|_{L_x^2}.
\end{align*}
Since $2it\nabla\psi=J\psi-x\psi$, we have
\[
t\|(\nabla A)\nabla\psi\|_{L_x^2}\lesssim\|\nabla A\|_{L_x^\infty}\|xf\|_{L_x^2}+\|(x\nabla A)\psi\|_{L_x^2}.
\]
\smallskip
\noindent\textit{The quadratic connection.}
We first establish the bilinear estimates that will also be used in Section~\ref{sec:decayestimate3}. By \eqref{ehatf8}, differentiating $\eta\widehat{A^{(2)}}$ produces
\begin{equation}
\begin{aligned}
&\int e^{it(|\eta-\zeta|^2-|\zeta|^2)}\frac{\zeta}{|\eta-\zeta|^2+|\zeta|^2}\overline{\widehat f(t,\zeta-\eta)}\widehat f(t,\zeta)\,\mathrm{d}\zeta\\
&\quad+2it\int e^{it(|\eta-\zeta|^2-|\zeta|^2)}\eta(\eta-\zeta)\frac{\zeta}{|\eta-\zeta|^2+|\zeta|^2}\overline{\widehat f(t,\zeta-\eta)}\widehat f(t,\zeta)\,\mathrm{d}\zeta\\
&\quad+\int e^{it(|\eta-\zeta|^2-|\zeta|^2)}\eta\nabla_\eta\left(\frac{\zeta}{|\eta-\zeta|^2+|\zeta|^2}\right)\overline{\widehat f(t,\zeta-\eta)}\widehat f(t,\zeta)\,\mathrm{d}\zeta\\
&\quad+\int e^{it(|\eta-\zeta|^2-|\zeta|^2)}\frac{\eta\zeta}{|\eta-\zeta|^2+|\zeta|^2}\nabla_\eta\overline{\widehat f(t,\zeta-\eta)}\widehat f(t,\zeta)\,\mathrm{d}\zeta.
\end{aligned}
\label{exf4}
\end{equation}
Since
\[
2it\eta e^{it(|\eta-\zeta|^2-|\zeta|^2)}=-\nabla_\zeta e^{it(|\eta-\zeta|^2-|\zeta|^2)},
\]
integration by parts in $\zeta$ gives
\begin{align}
&2it\int e^{it(|\eta-\zeta|^2-|\zeta|^2)}\eta\frac{(\eta-\zeta)\zeta}{|\eta-\zeta|^2+|\zeta|^2}\overline{\widehat f(\zeta-\eta)}\widehat f(\zeta)\,\mathrm{d}\zeta\notag\\
&=\int e^{it(|\eta-\zeta|^2-|\zeta|^2)}\nabla_\zeta\left(\frac{(\eta-\zeta)\zeta}{|\eta-\zeta|^2+|\zeta|^2}\right)\overline{\widehat f(\zeta-\eta)}\widehat f(\zeta)\,\mathrm{d}\zeta\notag\\
&\quad-\int e^{it(|\eta-\zeta|^2-|\zeta|^2)}\frac{(\eta-\zeta)\zeta}{|\eta-\zeta|^2+|\zeta|^2}\nabla_\eta\overline{\widehat f(\zeta-\eta)}\widehat f(\zeta)\,\mathrm{d}\zeta\notag\\
&\quad+\int e^{it(|\eta-\zeta|^2-|\zeta|^2)}\frac{(\eta-\zeta)\zeta}{|\eta-\zeta|^2+|\zeta|^2}\overline{\widehat f(\zeta-\eta)}\nabla_\zeta\widehat f(\zeta)\,\mathrm{d}\zeta.
\label{exf5}
\end{align}
To estimate the contribution of $t\nabla^2A^{(2)}\psi$ in \eqref{exf3}, we also need to control $t\|\nabla^2A^{(2)}(0,t)\|_{L_x^2}$. Similarly, integrating by parts in $\zeta$, we obtain
\begin{align*}
&t\int e^{it(|\eta-\zeta|^2-|\zeta|^2)}\eta^2\frac{\zeta}{|\eta-\zeta|^2+|\zeta|^2}\overline{\widehat f(\zeta-\eta)}\widehat f(\zeta)\,\mathrm{d}\zeta\\
&\qquad=\frac{1}{2i}\int e^{it(|\eta-\zeta|^2-|\zeta|^2)}\nabla_\zeta\left[\frac{\eta\zeta}{|\eta-\zeta|^2+|\zeta|^2}\overline{\widehat f(\zeta-\eta)}\widehat f(\zeta)\right]\,\mathrm{d}\zeta.
\end{align*}
The conjugate contribution satisfies the same estimates.

Thus the terms without a weighted input have one of the order $-1$ multipliers
\[
\frac{\zeta}{|\eta-\zeta|^2+|\zeta|^2},\qquad
\eta\nabla_\eta\left(\frac{\zeta}{|\eta-\zeta|^2+|\zeta|^2}\right),\qquad
\nabla_\zeta\left(\frac{(\eta-\zeta)\zeta}{|\eta-\zeta|^2+|\zeta|^2}\right),\qquad
\nabla_\zeta\left(\frac{\eta\zeta}{|\eta-\zeta|^2+|\zeta|^2}\right).
\]
The terms with one weighted input have one of the order zero multipliers
\[
\frac{\eta\zeta}{|\eta-\zeta|^2+|\zeta|^2},\qquad
\frac{(\eta-\zeta)\zeta}{|\eta-\zeta|^2+|\zeta|^2}.
\]
We estimate these two classes separately.

\smallskip
\noindent\textit{The terms with an order $-1$ multiplier and no weighted input.}
For each multiplier above, direct differentiation gives
\[
|\nabla_{\eta-\zeta}^{\mu}\nabla_\zeta^\nu m(\eta,\zeta)|\lesssim(|\eta-\zeta|+|\zeta|)^{-1-|\mu|-|\nu|}.
\]
Hence, Lemma~\ref{lmet} gives
\begin{equation*}
\|m(\eta,\zeta)\varphi_{l_1}(\eta-\zeta)\varphi_{l_2}(\zeta)\|_{\mathcal S^\infty}\lesssim2^{-\max\{l_1,l_2\}}.
\end{equation*}
The same estimate holds with a low frequency cutoff in the lower frequency input.

For $p\in\{2,\infty\}$, Lemma~\ref{l multilinear estimates} bounds the Low-High, High-Low, and High-High contributions, respectively, by
\begin{align*}
&\sum_{l_2}2^{-l_2}\|P_{\leq l_2-5}\psi\|_{L_x^p}\|P_{l_2}\psi\|_{L_x^\infty}\lesssim\|\psi\|_{L_x^p}\sum_l2^{-l}\|P_l\psi\|_{L_x^\infty},\\
&\sum_{l_1}2^{-l_1}\|P_{l_1}\psi\|_{L_x^\infty}\|P_{\leq l_1-5}\psi\|_{L_x^p}\lesssim\|\psi\|_{L_x^p}\sum_l2^{-l}\|P_l\psi\|_{L_x^\infty},\\
&\sum_{|l_1-l_2|\leq4}2^{-\max\{l_1,l_2\}}\|P_{l_1}\psi\|_{L_x^p}\|P_{l_2}\psi\|_{L_x^\infty}\lesssim\|\psi\|_{L_x^p}\sum_l2^{-l}\|P_l\psi\|_{L_x^\infty}.
\end{align*}
Using \eqref{bso7} and \eqref{booassum}, we obtain
\begin{align}
\left\|\mathcal F_\eta^{-1}\left[\int m(\eta,\zeta)\overline{\widehat\psi(t,\zeta-\eta)}\widehat\psi(t,\zeta)\,\mathrm{d}\zeta\right]\right\|_{L_x^\infty}&\lesssim\epsilon_1^2t^{-\frac32},\label{exf6}\\
\left\|\mathcal F_\eta^{-1}\left[\int m(\eta,\zeta)\overline{\widehat\psi(t,\zeta-\eta)}\widehat\psi(t,\zeta)\,\mathrm{d}\zeta\right]\right\|_{L_x^2}&\lesssim\epsilon_1^2t^{-\frac12}.\label{exf7}
\end{align}

\smallskip
\noindent\textit{The terms with an order zero multiplier and one weighted internal profile.}
For the two order zero multipliers above, direct differentiation gives, for all multi-indices $\mu,\nu$,
\begin{equation*}
\left|\nabla_{\eta-\zeta}^{\mu}\nabla_\zeta^\nu\left(\eta\frac{\zeta}{|\eta-\zeta|^2+|\zeta|^2}\right)\right|+\left|\nabla_{\eta-\zeta}^{\mu}\nabla_\zeta^\nu\left(\frac{(\eta-\zeta)\zeta}{|\eta-\zeta|^2+|\zeta|^2}\right)\right|\lesssim(|\eta-\zeta|+|\zeta|)^{-|\mu|-|\nu|}.
\end{equation*}
Thus Lemma~\ref{lcoifman-meyer estimates} applies to both multipliers and either position of the weighted internal profile. By \eqref{booassum}, for either multiplier we obtain
\begin{equation}
\begin{aligned}
&\left\|\mathcal F_\eta^{-1}\left[\int_{\mathbb R^2}m(\eta,\zeta)\overline{\widehat{e^{it\Delta}(xf)}(t,\zeta-\eta)}\,\widehat\psi(t,\zeta)\,\mathrm d\zeta\right]\right\|_{L_x^2}\\
&\quad+\left\|\mathcal F_\eta^{-1}\left[\int_{\mathbb R^2}m(\eta,\zeta)\overline{\widehat\psi(t,\zeta-\eta)}\,\widehat{e^{it\Delta}(xf)}(t,\zeta)\,\mathrm d\zeta\right]\right\|_{L_x^2}\\
&\lesssim\|\psi(t)\|_{L_x^\infty}\|e^{it\Delta}(xf(t))\|_{L_x^2}\lesssim\epsilon_1^2t^{-1+\gamma}.
\end{aligned}
\label{exf8}
\end{equation}
Similarly, using \eqref{bsc10}, we obtain
\begin{equation}
\begin{aligned}
&\left\|\mathcal F_\eta^{-1}\left[\int_{\mathbb R^2}m(\eta,\zeta)\overline{\widehat{e^{it\Delta}(xf)}(t,\zeta-\eta)}\,\widehat\psi(t,\zeta)\,\mathrm d\zeta\right]\right\|_{L_x^4}\\
&\quad+\left\|\mathcal F_\eta^{-1}\left[\int_{\mathbb R^2}m(\eta,\zeta)\overline{\widehat\psi(t,\zeta-\eta)}\,\widehat{e^{it\Delta}(xf)}(t,\zeta)\,\mathrm d\zeta\right]\right\|_{L_x^4}\\
&\lesssim\|\psi(t)\|_{L_x^\infty}\|e^{it\Delta}(xf(t))\|_{L_x^4}\lesssim\epsilon_1^2t^{-\frac32+\frac{\beta}{2}}.
\end{aligned}
\label{exf9}
\end{equation}

Combining \eqref{exf4}, \eqref{exf5}, \eqref{exf7}, and \eqref{exf8} with the integration by parts above for $t\nabla^2A^{(2)}$, we obtain
\begin{equation}
\|x\nabla A^{(2)}(t)\|_{L_x^2}+t\|\nabla^2A^{(2)}(t)\|_{L_x^2}\lesssim\epsilon_1^2\left(t^{-\frac12}+t^{-1+\gamma}\right)\lesssim\epsilon_1^2t^{-\frac12}.
\label{exf10}
\end{equation}
Together with \eqref{eAhigh2}, this gives
\begin{equation}
\|\nabla^2A(t)\|_{L_x^2}\lesssim\epsilon_1^2t^{-\frac32}.
\label{exf11}
\end{equation}

By \eqref{exf10}, \eqref{eAhigh4}, and \eqref{exf11},
\[
\|(x\nabla A)\psi\|_{L_x^2}+t\|(\nabla^2A)\psi\|_{L_x^2}\lesssim\epsilon_1^3t^{-\frac32}+\epsilon_1^5t^{-\frac52+\gamma}(1+\log t)\lesssim\epsilon_1^3t^{-\frac32}.
\]
Therefore, \eqref{eAonederivativelinfty} implies
\[
t\|(\nabla A)\nabla\psi\|_{L_x^2}+t\|(\nabla^2A)\psi\|_{L_x^2}\lesssim\epsilon_1^3\left[t^{-\frac32}+t^{-2+\gamma}(1+\log t)\right].
\]
By \eqref{eAl2decay}, \eqref{eAlinftydecay}, and \eqref{eAonederivativel21} with $\sigma=0$,
\[
\|A\|_{L_x^2}\|\psi\|_{L_x^\infty}+t\|A\|_{L_x^\infty}\|\nabla A\|_{L_x^2}\|\psi\|_{L_x^\infty}\lesssim\epsilon_1^3t^{-2}(1+\log t)+\epsilon_1^5t^{-\frac52}.
\]
For the cubic term, the \(J\) product rule \eqref{eJproduct} and
\eqref{booassum} imply
\[
  \left\|J\bigl(\psi_l\operatorname{Im}(\overline{\psi_l}\psi_m)\bigr)\right\|_{L_x^2}\lesssim\|\psi(t)\|_{L_x^\infty}^2\|J\psi(t)\|_{L_x^2}\lesssim\epsilon_1^2t^{-2}\|xf(t)\|_{L_x^2}.
\]
By \eqref{eA34}, integrating by parts in $\zeta$ as in \eqref{exf5}, and applying \eqref{exf7} and \eqref{exf8}, we obtain
\begin{equation*}
t\|\nabla A_3^{(2)}(t)\|_{L_x^2}\lesssim\|\psi(t)\|_{L_x^\infty}\|xf(t)\|_{L_x^2}+\epsilon_1^2t^{-\frac{1}{2}}.
\end{equation*}
Together with \eqref{eA35} and \eqref{booassum}, this gives
\begin{equation*}
t\|(\nabla A_3(t))\psi(t)\|_{L_x^2}\lesssim\epsilon_1^2t^{-2}\|xf(t)\|_{L_x^2}+\epsilon_1^3t^{-\frac{3}{2}}+\epsilon_1^5t^{-2}\lesssim\epsilon_1^3t^{-\frac{3}{2}}.
\end{equation*}

Combining these estimates, we obtain
\[
\frac{\mathrm{d}}{\mathrm{d}t}\|xf(t)\|_{L_x^2}\lesssim\epsilon_1^2t^{-2}\|xf(t)\|_{L_x^2}+\epsilon_1^3t^{-1}.
\]
Gronwall inequality on $[1,t]$ gives
\[
\|xf(t)\|_{L_x^2}\lesssim\|xf(1)\|_{L_x^2}+\epsilon_1^3(1+\log t).
\]
Hence,
\[
\sup_{1\leq\tau\leq t}\tau^{-\gamma}\|xf(\tau)\|_{L_x^2}\lesssim\|xf(1)\|_{L_x^2}+\epsilon_1^3,
\]
which proves \eqref{exf1}.
\end{proof}

\section{Decay estimate $\mathrm{III}$: improvement of $\|x^2f\|_{L_x^2}$} \label{sec:decayestimate3}

In this section we improve the second weighted profile estimate for $t\geq1$. As in Section~\ref{sec:decayestimate2}, we use the vector field $J$. By \eqref{eJs4},
\[
J_iJ_j\psi(t)=e^{it\Delta}(x_ix_jf(t)),\qquad \|J^2\psi(t)\|_{L_x^2}\simeq\|x^2f(t)\|_{L_x^2}.
\]
\begin{proposition}\label{xx2}
Under the bootstrap assumption \eqref{booassum}, there holds
\begin{equation}\label{ex2f1}
\|x^2f(t)\|_{L_x^2}\lesssim\|x^2f(1)\|_{L_x^2}+\epsilon_1^3t^\beta,\qquad t\geq1.
\end{equation}
\end{proposition}

\begin{proof}
The proof of this proposition is split into the following subsections.

\subsection{The twice-commuted equation and the energy reduction}

Since \([i\partial_t+\Delta,J_i]=0\), applying \eqref{exf2} twice to \eqref{ecaloric7}, we obtain
\begin{equation}
\begin{aligned}
(i\partial_t+\Delta)J_iJ_j\psi_m&=-2iA_l\partial_lJ_iJ_j\psi_m\\&\quad+4t(\partial_iA_l)\partial_lJ_j\psi_m+4t(\partial_jA_l)\partial_lJ_i\psi_m\\
  &\quad+8it^2(\partial_i\partial_jA_l)\partial_l\psi_m+2iA_iJ_j\psi_m+2iA_jJ_i\psi_m\\
  &\quad-4t(\partial_iA_j+\partial_jA_i)\psi_m+J_iJ_j\bigl((A_3+A_lA_l-i\partial_lA_l)\psi_m\bigr)\\
  &\quad-iJ_iJ_j\bigl(\psi_l\operatorname{Im}(\overline{\psi_l}\psi_m)\bigr).
\end{aligned} \label{ex2f2}
\end{equation}

Since \(A_3\) is real-valued, its undifferentiated contribution to the energy identity vanishes. Moreover,
\[
  J_iJ_j(A_3\psi_m)=A_3J_iJ_j\psi_m+2it(\partial_iA_3)J_j\psi_m+2it(\partial_jA_3)J_i\psi_m-4t^2(\partial_i\partial_jA_3)\psi_m.
\]
Integration by parts in the transport term, H\"older inequality, and \eqref{eJs4} therefore give
\begin{equation*}
\begin{aligned}
  \frac{\mathrm{d}}{\mathrm{d}t}\|x^2f(t)\|_{L_x^2}&\lesssim\|\nabla A(t)\|_{L_x^\infty}\|x^2f(t)\|_{L_x^2}\\&\quad+t\|(\nabla A(t))\nabla e^{it\Delta}(xf(t))\|_{L_x^2}
  +t^2\|(\nabla^2A(t))\nabla\psi(t)\|_{L_x^2}\\
  &\quad+\|A(t)e^{it\Delta}(xf(t))\|_{L_x^2}+t\|(\nabla A(t))\psi(t)\|_{L_x^2}\\
  &\quad+\max_F\|J^2(F\psi)(t)\|_{L_x^2}+t\|(\nabla A_3(t))J\psi(t)\|_{L_x^2}\\
  &\quad+t^2\|(\nabla^2A_3(t))\psi(t)\|_{L_x^2},
\end{aligned}
\end{equation*}
where \(F\) ranges over the finite collection in \eqref{eFdefinition}, evaluated at heat time zero. By \eqref{eAonederivativelinfty}, Gronwall inequality gives
\begin{equation}
\begin{aligned}
  \|x^2f(t)\|_{L_x^2}&\lesssim\|x^2f(1)\|_{L_x^2}\\
  &\quad+\int_1^t\left[\tau\|(\nabla A(\tau))\nabla e^{i\tau\Delta}(xf(\tau))\|_{L_x^2}+\tau^2\|(\nabla^2A(\tau))\nabla\psi(\tau)\|_{L_x^2}\right]\,\mathrm{d}\tau\\
  &\quad+\int_1^t\left[\|A(\tau)e^{i\tau\Delta}(xf(\tau))\|_{L_x^2}+\tau\|(\nabla A(\tau))\psi(\tau)\|_{L_x^2}\right]\,\mathrm{d}\tau\\
  &\quad+\int_1^t\max_F\|J^2(F\psi)(\tau)\|_{L_x^2}\,\mathrm{d}\tau\\
  &\quad+\int_1^t\left[\tau\|(\nabla A_3(\tau))J\psi(\tau)\|_{L_x^2}+\tau^2\|(\nabla^2A_3(\tau))\psi(\tau)\|_{L_x^2}\right]\,\mathrm{d}\tau.
\end{aligned} \label{ex2f3}
\end{equation}

By \eqref{eAlinftydecay} and \eqref{booassum},
\begin{equation*}
  \int_1^t\left\|A(\tau)e^{i\tau\Delta}(xf(\tau))\right\|_{L_x^2}\,\mathrm{d}\tau\lesssim\int_1^t\|A(\tau)\|_{L_x^\infty}\|xf(\tau)\|_{L_x^2}\,\mathrm{d}\tau\lesssim\epsilon_1^3\int_1^t\tau^{-\frac32+\gamma}\,\mathrm{d}\tau\lesssim\epsilon_1^3t^\beta.
\end{equation*}

By \eqref{eAonederivativel21} and \eqref{booassum},
\begin{equation*}
  \int_1^t\tau\|(\nabla A(\tau))\psi(\tau)\|_{L_x^2}\,\mathrm{d}\tau\lesssim\int_1^t\tau\|\nabla A(\tau)\|_{L_x^2}\|\psi(\tau)\|_{L_x^\infty}\,\mathrm{d}\tau\lesssim\epsilon_1^3\int_1^t \tau^{-1}\,\mathrm{d}\tau\lesssim\epsilon_1^3t^\beta.
\end{equation*}
Thus, both lower order terms in \eqref{ex2f3} are compatible with the second weighted bootstrap.

We next estimate the two main contributions:
\[
  \int_1^t\tau\left\|(\nabla A(\tau))\nabla e^{i\tau\Delta}(xf(\tau))\right\|_{L_x^2}\,\mathrm{d}\tau
\]
and
\[
  \int_1^t\tau^2\left\|(\nabla^2A(\tau))\nabla\psi(\tau)\right\|_{L_x^2}\,\mathrm{d}\tau.
\]
A direct physical space H\"older estimate is not available for them due to the low regularity. We therefore return to frequency space by Plancherel. Up to harmless constants, for every fixed \(\tau\geq1\),
\begin{align*}
&\tau\left\|(\partial_iA(\tau))\cdot\nabla e^{i\tau\Delta}(x_jf(\tau))\right\|_{L_x^2}\\
&\qquad\simeq\tau\left\|\int_{\mathbb R^2}e^{i\tau(|\xi|^2-|\xi-\eta|^2)}\eta_i(\xi-\eta)\cdot\widehat A(\tau,\eta)\partial_{\xi_j}\widehat f(\tau,\xi-\eta)\,d\eta\right\|_{L_\xi^2},
\end{align*}
while
\begin{align*}
&\tau^2\left\|(\partial_i\partial_jA(\tau))\cdot\nabla\psi(\tau)\right\|_{L_x^2}\\
&\qquad\simeq\tau^2\left\|\int_{\mathbb R^2}e^{i\tau(|\xi|^2-|\xi-\eta|^2)}\eta_i\eta_j(\xi-\eta)\cdot\widehat A(\tau,\eta)\widehat f(\tau,\xi-\eta)\,d\eta\right\|_{L_\xi^2}.
\end{align*}

For every fixed $\tau\geq1$, we decompose the frequency interactions and integrate by parts in the frequency variables. We then integrate the resulting $L_\xi^2$ bounds over $[1,t]$.

\subsection{The term containing \(\tau(\nabla A)\nabla e^{i\tau\Delta}(xf)\)}

We first exploit the factor \(\xi-\eta\) coming from the derivative of the outer weighted profile. Since
\[
  \nabla_\eta\bigl(|\xi|^2-|\xi-\eta|^2\bigr)=2(\xi-\eta),
\]
integration by parts in \(\eta\) gives, up to harmless fixed tensor contractions,
\begin{align}
  &\tau\int_{\mathbb R^2}e^{i\tau(|\xi|^2-|\xi-\eta|^2)}\eta(\xi-\eta)\widehat A(\tau,\eta)\nabla_\xi\widehat f(\tau,\xi-\eta)\,\mathrm{d}\eta\notag\\
  &\quad=-\frac{1}{2i}\int_{\mathbb R^2}e^{i\tau(|\xi|^2-|\xi-\eta|^2)}\nabla_\eta\bigl(\eta\widehat A(\tau,\eta)\bigr)\nabla_\xi\widehat f(\tau,\xi-\eta)\,\mathrm{d}\eta\notag\\
  &\qquad+\frac{1}{2i}\int_{\mathbb R^2}e^{i\tau(|\xi|^2-|\xi-\eta|^2)}\eta\widehat A(\tau,\eta)\nabla_\xi^2\widehat f(\tau,\xi-\eta)\,\mathrm{d}\eta. \label{ex2f4}
\end{align}
Therefore, by Plancherel inequality,
\begin{align*}
  &\tau\left\|(\nabla A(\tau))\nabla e^{i\tau\Delta}(xf(\tau))\right\|_{L_x^2}\\
  &\quad\lesssim\left\|(x\nabla A(\tau))e^{i\tau\Delta}(xf(\tau))\right\|_{L_x^2}+\left\|(\nabla A(\tau))e^{i\tau\Delta}(x^2f(\tau))\right\|_{L_x^2}.
\end{align*}

The second term is controlled directly by \eqref{eAonederivativelinfty} and \eqref{booassum}:
\begin{align*}
  &\int_1^t\left\|(\nabla A(\tau))e^{i\tau\Delta}(x^2f(\tau))\right\|_{L_x^2}\,\mathrm{d}\tau\\
  &\quad\lesssim\int_1^t\|\nabla A(\tau)\|_{L_x^\infty}\|x^2f(\tau)\|_{L_x^2}\,\mathrm{d}\tau\lesssim\epsilon_1^3\int_1^t(1+\log\tau)\tau^{-2+\beta}\,\mathrm{d}\tau\lesssim\epsilon_1^3.
\end{align*}

It remains to estimate
\[
  \int_1^t\left\|(x\nabla A(\tau))e^{i\tau\Delta}(xf(\tau))\right\|_{L_x^2}\,\mathrm{d}\tau.
\]
We use the decomposition \eqref{eAcom1}. For the higher order connection, by \eqref{eAhigh4} and \eqref{booassum},
\begin{align*}
  &\int_1^t\left\|(x\nabla A^{(\geq4)}(\tau))e^{i\tau\Delta}(xf(\tau))\right\|_{L_x^2}\,\mathrm{d}\tau\\
  &\quad\lesssim\int_1^t\|x\nabla A^{(\geq4)}(\tau)\|_{L_x^\infty}\|xf(\tau)\|_{L_x^2}\,\mathrm{d}\tau\lesssim\epsilon_1^5\int_1^t\tau^{-\frac52+2\gamma}(1+\log\tau)\,\mathrm{d}\tau\lesssim\epsilon_1^5.
\end{align*}

For the quadratic contribution, we use the decomposition established in \eqref{exf4} and \eqref{exf5}. By \eqref{exf6}, \eqref{exf9}, and \eqref{bsc10},
\begin{align*}
\|(x\nabla A^{(2)}(\tau))e^{i\tau\Delta}(xf(\tau))\|_{L_x^2}&\lesssim\epsilon_1^2\tau^{-\frac32}\|xf(\tau)\|_{L_x^2}+\epsilon_1^2\tau^{-\frac32+\frac\beta2}(1+\log\tau)\|e^{i\tau\Delta}(xf(\tau))\|_{L_x^4}\\
&\lesssim\epsilon_1^3\left[\tau^{-\frac32+\gamma}+\tau^{-2+\beta}(1+\log\tau)\right].
\end{align*}
Both terms are integrable on $[1,\infty)$. Together with the higher order contribution, this yields
\begin{equation*}
\int_1^t\|(x\nabla A(\tau))e^{i\tau\Delta}(xf(\tau))\|_{L_x^2}\,\mathrm{d}\tau\lesssim\epsilon_1^3.
\end{equation*}
Returning to \eqref{ex2f4}, we conclude that
\begin{equation}
\int_1^t\tau\|(\nabla A(\tau))\nabla e^{i\tau\Delta}(xf(\tau))\|_{L_x^2}\,\mathrm{d}\tau\lesssim\epsilon_1^3. \label{ex2f5}
\end{equation}

\subsection{The term containing \(\tau^2(\nabla^2A)\nabla\psi\)}

It remains to prove
\[
  \int_1^t\tau^2\left\|(\nabla^2A(\tau))\nabla\psi(\tau)\right\|_{L_x^2}\,\mathrm{d}\tau\lesssim\epsilon_1^3t^\beta.
\]
We use the decomposition \eqref{eAcom1} and estimate \(A^{(\geq4)}\) and \(A^{(2)}\) separately.

\smallskip
\noindent\textit{The higher order connection.}

By Plancherel inequality, the contribution of \(A^{(\geq4)}\) is equivalent, up to harmless tensor contractions, to
\[
\begin{aligned}
  \tau^2\left\|\int_{\mathbb R^2}e^{i\tau(|\xi|^2-|\xi-\eta|^2)}\eta^2(\xi-\eta)\widehat{A^{(\geq4)}}(\tau,\eta)\widehat f(\tau,\xi-\eta)\,\mathrm{d}\eta\right\|_{L_\xi^2}.
\end{aligned}
\]
One integration by parts in \(\eta\) gives
\[
\begin{aligned}
  &\tau^2\int_{\mathbb R^2}e^{i\tau(|\xi|^2-|\xi-\eta|^2)}\eta^2(\xi-\eta)\widehat{A^{(\geq4)}}(\tau,\eta)\widehat f(\tau,\xi-\eta)\,\mathrm{d}\eta\\
  &\quad=-\frac{\tau}{2i}\int_{\mathbb R^2}e^{i\tau(|\xi|^2-|\xi-\eta|^2)}\nabla_\eta\left[\eta^2\widehat{A^{(\geq4)}}(\tau,\eta)\widehat f(\tau,\xi-\eta)\right]\,\mathrm{d}\eta .
\end{aligned}
\]
Hence,
\[
\begin{aligned}
  &\tau^2\left\|(\nabla^2A^{(\geq4)}(\tau))\nabla\psi(\tau)\right\|_{L_x^2}\\
  &\quad\lesssim\tau\left\|(x\nabla^2A^{(\geq4)}(\tau))\psi(\tau)\right\|_{L_x^2}+\tau\left\|(\nabla^2A^{(\geq4)}(\tau))e^{i\tau\Delta}(xf(\tau))\right\|_{L_x^2}.
\end{aligned}
\]
By \eqref{eAhigh5} and \eqref{booassum},
\[
\int_1^t\tau\|(x\nabla^2A^{(\geq4)}(\tau))\psi(\tau)\|_{L_x^2}\,\mathrm{d}\tau\lesssim\epsilon_1^5\int_1^t\tau^{-2+\gamma}(1+\log\tau)\,\mathrm{d}\tau\lesssim\epsilon_1^5.
\]
For the second term, Sobolev embedding and \eqref{eAhigh2} give
\[
\|\nabla^2A^{(\geq4)}(\tau)\|_{L_x^4}\lesssim\|\nabla^2A^{(\geq4)}(\tau)\|_{H^1}\lesssim\epsilon_1^4\tau^{-\frac52}.
\]
Together with \eqref{bsc10}, this yields
\[
\int_1^t\tau\|(\nabla^2A^{(\geq4)}(\tau))e^{i\tau\Delta}(xf(\tau))\|_{L_x^2}\,\mathrm{d}\tau\lesssim\epsilon_1^5\int_1^t\tau^{-2+\frac\beta2}\,\mathrm{d}\tau\lesssim\epsilon_1^5.
\]
Hence,
\begin{equation}\label{ex2f6}
  \int_1^t\tau^2\left\|(\nabla^2A^{(\geq 4)}(\tau))\nabla\psi(\tau)\right\|_{L_x^2}\,\mathrm{d}\tau\lesssim\epsilon_1^5.
\end{equation}

\smallskip
\noindent\textit{The quadratic connection.}

We consider the expression of \eqref{ehatf8} whose numerator is \(\zeta\). The conjugate symmetric term has the same localized multiplier bounds. After substituting \eqref{ehatf8}, the quadratic contribution is
\begin{equation*}
  \tau^2\Bigg\|\iint_{\mathbb R^2\times\mathbb R^2}e^{i\tau\Phi_3(\xi,\eta,\zeta)}\eta^2(\xi-\eta)\frac{\zeta}{|\eta-\zeta|^2+|\zeta|^2}\overline{\widehat f(\tau,\zeta-\eta)}\widehat f(\tau,\zeta)\widehat f(\tau,\xi-\eta)\,\mathrm{d}\zeta\,\mathrm{d}\eta\Bigg\|_{L_\xi^2},
\end{equation*}
where
\[
  \Phi_3(\xi,\eta,\zeta)=|\xi|^2-|\xi-\eta|^2+|\eta-\zeta|^2-|\zeta|^2.
\]
Before integrating by parts in $\zeta$, we split the internal frequencies using
\[
\chi_l\left(\tau\sqrt{|\eta-\zeta|^2+|\zeta|^2}\right)+\chi_h\left(\tau\sqrt{|\eta-\zeta|^2+|\zeta|^2}\right)=1.
\]
On the support of the low frequency cutoff, $|\eta-\zeta|+|\zeta|\lesssim\tau^{-1}$. We use
\begin{equation}\label{ex2f7}
(\xi-\eta)e^{i\tau\Phi_3}=\frac{1}{2i\tau}\nabla_\eta e^{i\tau\Phi_3}-(\eta-\zeta)e^{i\tau\Phi_3}
\end{equation}
and integrate the first term by parts in $\eta$. The coefficient $\eta^2\zeta(|\eta-\zeta|^2+|\zeta|^2)^{-1}\chi_l$ is $O(\tau^{-1})$, its $\eta$ derivative is $O(1)$, and the remaining coefficient $(\eta-\zeta)\eta^2\zeta(|\eta-\zeta|^2+|\zeta|^2)^{-1}\chi_l$ is $O(\tau^{-2})$. Minkowski and Cauchy-Schwarz give
\begin{align*}
&\tau^2\left\|\iint e^{i\tau\Phi_3}\eta^2(\xi-\eta)\frac{\zeta}{|\eta-\zeta|^2+|\zeta|^2}\chi_l\left(\tau\sqrt{|\eta-\zeta|^2+|\zeta|^2}\right)\overline{\widehat f(\zeta-\eta)}\widehat f(\zeta)\widehat f(\xi-\eta)\,\mathrm{d}\zeta\,\mathrm{d}\eta\right\|_{L_\xi^2}\\
&\quad\lesssim\tau^{-3}\|f\|_{L_x^2}\|\widehat f\|_{L_\xi^\infty}^2+\tau^{-3}\|f\|_{L_x^2}\|xf\|_{L_x^2}\|\widehat f\|_{L_\xi^\infty}+\tau^{-4}\|xf\|_{L_x^2}\|\widehat f\|_{L_\xi^\infty}^2+\tau^{-4}\|f\|_{L_x^2}\|\widehat f\|_{L_\xi^\infty}^2\lesssim\epsilon_1^3\tau^{-3+\gamma}.
\end{align*}
This contribution is integrable in $\tau$. We apply the following two integrations by parts in $\zeta$ only to the high-frequency contribution.
Since $\nabla_\zeta\Phi_3=-2\eta$, we have
\[
  \tau^2\eta_i\eta_j e^{i\tau\Phi_3}=-\frac{1}{4}\partial_{\zeta_i}\partial_{\zeta_j}e^{i\tau\Phi_3}.
\]
Integrating by parts twice in \(\zeta\), it suffices to estimate
\begin{align}
\Bigg\|\iint e^{i\tau\Phi_3}(\xi-\eta)\nabla_\zeta^2\left[\frac{\zeta}{|\eta-\zeta|^2+|\zeta|^2}\chi_h\left(\tau\sqrt{|\eta-\zeta|^2+|\zeta|^2}\right)\overline{\widehat f(\zeta-\eta)}\widehat f(\zeta)\right]\widehat f(\xi-\eta)\,\mathrm{d}\zeta\,\mathrm{d}\eta\Bigg\|_{L_\xi^2}.
\label{ex2f8}
\end{align}
For fixed $\tau$, we write
\begin{align}\label{99}
&m_{-1}(\eta,\zeta)=\frac{\zeta}{|\eta-\zeta|^2+|\zeta|^2}\chi_h\left(\tau\sqrt{|\eta-\zeta|^2+|\zeta|^2}\right)\nonumber\\
&m_{-2}=\nabla_\zeta m_{-1}, m_{-3}=\nabla_\zeta^2m_{-1}, m_{-4}=\nabla_\eta m_{-3}.
\end{align}
All the following symbol bounds hold uniformly in $\tau$. Derivatives of the cutoff are supported where $|\eta-\zeta|+|\zeta|\simeq\tau^{-1}$ and obey the same bounds as derivatives of the corresponding homogeneous symbols. Direct differentiation gives, for \(|\mu|+|\nu|\leq4\),
\[
  \left|\nabla_{\eta-\zeta}^{\mu}\nabla_\zeta^\nu m_{-j}(\eta,\zeta)\right|\lesssim\bigl(|\eta-\zeta|+|\zeta|\bigr)^{-j-|\mu|-|\nu|}, \qquad 1\leq j\leq4.
\]
The expansion in \eqref{ex2f8} consists of the following four classes:
\[
\begin{aligned}
  &m_{-1}(\eta,\zeta)\overline{\nabla_\zeta^2\widehat f(\zeta-\eta)}\widehat f(\zeta), \qquad m_{-1}(\eta,\zeta)\overline{\widehat f(\zeta-\eta)}\nabla_\zeta^2\widehat f(\zeta),\\
  &m_{-1}(\eta,\zeta)\overline{\nabla_\zeta\widehat f(\zeta-\eta)}\nabla_\zeta\widehat f(\zeta),\\
  &m_{-2}(\eta,\zeta)\overline{\nabla_\zeta\widehat f(\zeta-\eta)}\widehat f(\zeta), \qquad m_{-2}(\eta,\zeta)\overline{\widehat f(\zeta-\eta)}\nabla_\zeta\widehat f(\zeta),\\
  &m_{-3}(\eta,\zeta)\overline{\widehat f(\zeta-\eta)}\widehat f(\zeta).
\end{aligned}
\]
Thus we consider successively the cases in which the two derivatives fall on the same internal profile, on two different internal profiles, on one internal profile and the heat multiplier, or both on the heat multiplier.

We first record two localized multiplier estimates used repeatedly below.

\begin{lemma}\label{luu1}
There holds
\[
\begin{aligned}
  &\left\|\mathcal F_\eta^{-1}\left[\int_{\mathbb R^2}e^{i\tau(|\eta-\zeta|^2-|\zeta|^2)}m_{-1}(\eta,\zeta)\overline{\widehat{x^2f}(\tau,\zeta-\eta)}\widehat{xf}(\tau,\zeta)\,\mathrm{d}\zeta\right]\right\|_{L_x^4}\\
  &\qquad\lesssim\|x^2f(\tau)\|_{L_x^2}\|e^{i\tau\Delta}(xf(\tau))\|_{L_x^4}.
\end{aligned}
\]
The same estimate holds with the two input positions exchanged and with \((x^2f,xf)\) replaced by any of
\[
  (x^2f,f),\qquad (xf,xf),\qquad (xf,f),
\]
with the \(L_x^2\) norm placed on the first input and the \(L_x^4\) norm placed on the corresponding Schr\"odinger evolution of the second input.
\end{lemma}

\begin{proof}
Lemma~\ref{lmet} gives
\[
  \left\|m_{-1}(\eta,\zeta)\varphi_k(\eta)\varphi_{k_1}(\eta-\zeta)\varphi_{k_2}(\zeta)\right\|_{\mathcal S^\infty}\lesssim2^{-\max\{k_1,k_2\}}.
\]

If \(k_1\leq k_2-10\), then \(k=k_2+O(1)\). By Lemma~\ref{l multilinear estimates} and Bernstein inequality,
\[
\begin{aligned}
  &\left\|P_k\mathcal F_\eta^{-1}\left[\int e^{i\tau(|\eta-\zeta|^2-|\zeta|^2)}m_{-1}(\eta,\zeta)\overline{\widehat{P_{k_1}(x^2f)}(\zeta-\eta)}\widehat{P_{k_2}(xf)}(\zeta)\,\mathrm{d}\zeta\right]\right\|_{L_x^4}\\
  &\quad\lesssim2^{-k_2}\|P_{k_1}e^{i\tau\Delta}(x^2f)\|_{L_x^\infty}\|P_{k_2}e^{i\tau\Delta}(xf)\|_{L_x^4}\\
  &\quad\lesssim 2^{k_1-k_2}\|P_{k_1}(x^2f)\|_{L_x^2}\|P_{k_2}e^{i\tau\Delta}(xf)\|_{L_x^4}.
\end{aligned}
\]
Hence, Young inequality gives
\[
\begin{aligned}
  &\left(\sum_k\left[\sum_{\substack{k_1\leq k_2-10\\|k_2-k|\leq5}}2^{k_1-k_2}\|P_{k_1}(x^2f)\|_{L_x^2}\|P_{k_2}e^{i\tau\Delta}(xf)\|_{L_x^4}\right]^2\right)^{\frac12}\\
  &\quad\lesssim\|e^{i\tau\Delta}(xf)\|_{L_x^4}\left(\sum_k\left[\sum_{k_1\leq k-5}2^{k_1-k}\|P_{k_1}(x^2f)\|_{L_x^2}\right]^2\right)^{\frac12}\\
  &\quad\lesssim\|x^2f\|_{L_x^2}\|e^{i\tau\Delta}(xf)\|_{L_x^4}.
\end{aligned}
\]

If \(k_2\leq k_1-10\), then \(k=k_1+O(1)\), and
\[
\begin{aligned}
  &\left\|P_k\mathcal F_\eta^{-1}\left[\int e^{i\tau(|\eta-\zeta|^2-|\zeta|^2)}m_{-1}(\eta,\zeta)\overline{\widehat{P_{k_1}(x^2f)}(\zeta-\eta)}\widehat{P_{k_2}(xf)}(\zeta)\,\mathrm{d}\zeta\right]\right\|_{L_x^4}\\
  &\quad\lesssim 2^{-k_1}\|P_{k_1}e^{i\tau\Delta}(x^2f)\|_{L_x^4}\|P_{k_2}e^{i\tau\Delta}(xf)\|_{L_x^\infty}\\
  &\quad\lesssim 2^{\frac{k_2-k_1}{2}}\|P_{k_1}(x^2f)\|_{L_x^2}\|P_{k_2}e^{i\tau\Delta}(xf)\|_{L_x^4}.
\end{aligned}
\]
It follows that
\[
\begin{aligned}
  &\left(\sum_k\left[\sum_{\substack{|k_1-k|\leq5\\k_2\leq k_1-10}}2^{\frac{k_2-k_1}{2}}\|P_{k_1}(x^2f)\|_{L_x^2}\|P_{k_2}e^{i\tau\Delta}(xf)\|_{L_x^4}\right]^2\right)^{\frac12}\\
  &\quad\lesssim\|e^{i\tau\Delta}(xf)\|_{L_x^4}\left(\sum_{k}\|P_{k}(x^2f)\|_{L_x^2}^2\right)^{\frac12}\\
  &\quad\lesssim\|x^2f\|_{L_x^2}\|e^{i\tau\Delta}(xf)\|_{L_x^4}.
\end{aligned}
\]

Finally, if \(|k_1-k_2|\leq10\), then \(k\leq k_1+15\). Bernstein inequality from \(L_x^{\frac{4}{3}}\) to \(L_x^4\) gives
\[
\begin{aligned}
  &\left\|P_k\mathcal F_\eta^{-1}\left[\int e^{i\tau(|\eta-\zeta|^2-|\zeta|^2)}m_{-1}(\eta,\zeta)\overline{\widehat{P_{k_1}(x^2f)}(\zeta-\eta)}\widehat{P_{k_2}(xf)}(\zeta)\,\mathrm{d}\zeta\right]\right\|_{L_x^4}\\
  &\quad\lesssim 2^{k-k_1}\|P_{k_1}(x^2f)\|_{L_x^2}\|P_{k_2}e^{i\tau\Delta}(xf)\|_{L_x^4}.
\end{aligned}
\]
Therefore,
\[
\begin{aligned}
  &\left(\sum_k\left[\sum_{\substack{k_1\geq k-15\\|k_2-k_1|\leq10}}2^{k-k_1}\|P_{k_1}(x^2f)\|_{L_x^2}\|P_{k_2}e^{i\tau\Delta}(xf)\|_{L_x^4}\right]^2\right)^{\frac12}\\
  &\quad\lesssim\|e^{i\tau\Delta}(xf)\|_{L_x^4}\left(\sum_k\left[\sum_{k_1\geq k-15}2^{k-k_1}\|P_{k_1}(x^2f)\|_{L_x^2}\right]^2\right)^{\frac12}\\
  &\quad\lesssim\|x^2f\|_{L_x^2}\|e^{i\tau\Delta}(xf)\|_{L_x^4}.
\end{aligned}
\]
Combining the three frequency cases, the Littlewood-Paley square function estimate and Minkowski inequality prove the lemma.
\end{proof}

\begin{lemma}\label{luu2}
There holds
\[
\begin{aligned}
  &\left\|\mathcal F_\eta^{-1}\left[\int_{\mathbb R^2}e^{i\tau(|\eta-\zeta|^2-|\zeta|^2)}m_{-2}(\eta,\zeta)\overline{\widehat{x^2f}(\tau,\zeta-\eta)}\widehat f(\tau,\zeta)\,\mathrm{d}\zeta\right]\right\|_{L_x^\infty}\\
  &\qquad\lesssim\|x^2f(\tau)\|_{L_x^2}\|f(\tau)\|_{L_x^2}.
\end{aligned}
\]
The same estimate holds after interchanging the two input positions or replacing \((x^2f,f)\) by
\[
  (x^2f,xf),\qquad (xf,xf),\qquad (xf,f),\qquad (f,f).
\]
\end{lemma}

\begin{proof}
Lemma \ref{lmet} gives
\[
  \left\|m_{-2}(\eta,\zeta)\varphi_{k_1}(\eta-\zeta)\varphi_{k_2}(\zeta)\right\|_{\mathcal S^\infty}\lesssim 2^{-2\max\{k_1,k_2\}}.
\]
Lemma~\ref{l multilinear estimates} and Bernstein inequality give
\[
\begin{aligned}
  &\left\|\mathcal F_\eta^{-1}\left[\int e^{i\tau(|\eta-\zeta|^2-|\zeta|^2)}m_{-2}(\eta,\zeta)\overline{\widehat{P_{k_1}(x^2f)}(\zeta-\eta)}\widehat{P_{k_2}f}(\zeta)\,\mathrm{d}\zeta\right]\right\|_{L_x^\infty}\\
  &\quad\lesssim 2^{k_1+k_2-2\max\{k_1,k_2\}}\|P_{k_1}(x^2f)\|_{L_x^2}\|P_{k_2}f\|_{L_x^2}.
\end{aligned}
\]
We compute
\[
\begin{aligned}
  &\sum_{k_2}\sum_{k_1\leq k_2-10}2^{k_1-k_2}\|P_{k_1}(x^2f)\|_{L_x^2}\|P_{k_2}f\|_{L_x^2}\lesssim \|x^2f\|_{L_x^2}\|f\|_{L_x^2},\\
  &\sum_{k_1}\sum_{k_2\leq k_1-10}2^{k_2-k_1}\|P_{k_1}(x^2f)\|_{L_x^2}\|P_{k_2}f\|_{L_x^2}\lesssim \|x^2f\|_{L_x^2}\|f\|_{L_x^2},\\
  &\sum_{|k_1-k_2|\leq 10}\|P_{k_1}(x^2f)\|_{L_x^2}\|P_{k_2}f\|_{L_x^2}\lesssim \|x^2f\|_{L_x^2}\|f\|_{L_x^2}.
\end{aligned}
\]
This proves the desired estimate.
\end{proof}
The proofs of Lemma~\ref{luu1} and Lemma~\ref{luu2} use only the corresponding order $-1$ and order $-2$ symbol bounds. Their conclusions therefore remain valid for any multipliers satisfying these bounds, with constants depending only on the corresponding symbol seminorms.

\smallskip
\noindent\textit{Both derivatives fall on the same internal profile.}
We separate the outer frequency from the two internal frequencies to avoid placing the outer derivative on $\psi$. We implement the outer high region by the smooth cutoff
\(\chi\left(2^{20}\frac{|\eta-\zeta|^2+|\zeta|^2}{|\xi-\eta|^2}\right)\), and the complementary region by one minus this cutoff. The region restrictions below denote these smooth localizations. We first consider
\[
  |\xi-\eta|\lesssim|\eta-\zeta|+|\zeta|.
\]
After localizing \(|\xi-\eta|\simeq2^l, |\eta-\zeta|\simeq2^{k_1}, |\zeta|\simeq2^{k_2}\), Lemma~\ref{lmet} gives
\[
\begin{aligned}
  &\left\|(\xi-\eta)m_{-1}(\eta,\zeta)\varphi_l(\xi-\eta)\varphi_{\leq k-5}(\eta-\zeta)\varphi_k(\zeta)\right\|_{\mathcal S^\infty}\lesssim2^{l-k},\\
  &\left\|(\xi-\eta)m_{-1}(\eta,\zeta)\varphi_l(\xi-\eta)\varphi_k(\eta-\zeta)\varphi_{\leq k+4}(\zeta)\right\|_{\mathcal S^\infty}\lesssim2^{l-k},\qquad l\leq k+O(1).
\end{aligned}
\]
By Lemma~\ref{l multilinear estimates}, the contribution with the weighted input at the lower internal frequency is bounded by
\[
  \sum_k\sum_{l\leq k+O(1)}2^{l-k}\|P_{\leq k-5}(x^2f)\|_{L_x^2}\|P_k\psi\|_{L_x^\infty}\|P_l\psi\|_{L_x^\infty}\lesssim\|x^2f\|_{L_x^2}\|\psi\|_{L_x^\infty}\sum_k\|P_k\psi\|_{L_x^\infty}.
\]
For the other internal frequency case, Cauchy-Schwarz, Young inequality, and the Littlewood-Paley square theorem give
\[
\begin{aligned}
  &\sum_k\sum_{l\leq k+O(1)}2^{l-k}\|P_k(x^2f)\|_{L_x^2}\|P_{\leq k+4}\psi\|_{L_x^\infty}\|P_l\psi\|_{L_x^\infty}\\
  &\quad\lesssim\|\psi\|_{L_x^\infty}\left(\sum_k\|P_k(x^2f)\|_{L_x^2}^2\right)^{\frac12}\left(\sum_k\left[\sum_{l\leq k+O(1)}2^{l-k}\|P_l\psi\|_{L_x^\infty}\right]^2\right)^{\frac12}\\
  &\quad\lesssim\|x^2f\|_{L_x^2}\|\psi\|_{L_x^\infty}\sum_l\|P_l\psi\|_{L_x^\infty}.
\end{aligned}
\]
By \eqref{epsihighfre2} at heat time zero,
\begin{equation}
\begin{aligned}
  &\left\|\iint_{|\xi-\eta|\lesssim|\eta-\zeta|+|\zeta|}e^{i\tau\Phi_3}(\xi-\eta)m_{-1}(\eta,\zeta)\overline{\widehat{x^2f}(\zeta-\eta)}\widehat f(\zeta)\widehat f(\xi-\eta)\,\mathrm{d}\zeta\,\mathrm{d}\eta\right\|_{L_\xi^2}\\
  &\quad\lesssim\|x^2f\|_{L_x^2}\|\psi\|_{L_x^\infty}\sum_k\|P_k\psi\|_{L_x^\infty}\\
  &\quad\lesssim\epsilon_1^3\tau^{-2+\beta}(1+\log\tau).
\end{aligned} \label{ex2f9}
\end{equation}
The same argument applies when the weighted internal input is in the other position. We shall use this summation below without repeating the localized symbol calculation.

We next consider
\[
  |\xi-\eta|\gg|\eta-\zeta|+|\zeta|.
\]
When both derivatives fall on \(\widehat f(\zeta-\eta)\), we use \(\nabla_\eta+\nabla_\zeta\), whereas when both derivatives fall on \(\widehat f(\zeta)\), we use \(\nabla_\eta\). Indeed,
\begin{align}
  &(\nabla_\eta+\nabla_\zeta)\widehat f(\zeta-\eta)=0, \qquad (\nabla_\eta+\nabla_\zeta)\Phi_3=2\bigl((\xi-\eta)-\zeta\bigr),\notag\\
  &\nabla_\eta\widehat f(\zeta)=0, \qquad\nabla_\eta\Phi_3=2\bigl((\xi-\eta)+(\eta-\zeta)\bigr). \label{ex2f10}
\end{align}
On the present region,
\[
  |(\xi-\eta)-\zeta|\simeq|\xi-\eta|, \qquad |(\xi-\eta)+(\eta-\zeta)|\simeq|\xi-\eta|.
\]
Thus the corresponding phase integration by parts produces \(\tau^{-1}\), and its denominator cancels the original factor \((\xi-\eta)\). We use
\[
\begin{aligned}
  e^{i\tau\Phi_3}&=\frac{(\xi-\eta)-\zeta}{2i\tau|(\xi-\eta)-\zeta|^2}\cdot(\nabla_\eta+\nabla_\zeta)e^{i\tau\Phi_3},\\
  e^{i\tau\Phi_3}&=\frac{(\xi-\eta)+(\eta-\zeta)}{2i\tau|(\xi-\eta)+(\eta-\zeta)|^2}\cdot\nabla_\eta e^{i\tau\Phi_3}.
\end{aligned}
\]
After integrating by parts, the derivative can fall on the other internal profile, on the outer profile, on \(m_{-1}\), on the integration-by-parts coefficient, or on the frequency cutoff defining the present region. A cutoff derivative has the same bound as a derivative of the integration-by-parts coefficient. Derivatives of \(m_{-1}\) or \(m_{-2}\) in either integration-by-parts direction have the same symbol bounds as the corresponding order \(-2\) or \(-3\) symbols above.

\smallskip
\noindent\textit{The derivative falls on the other internal profile.}
The internal multiplier remains \(m_{-1}\). We localize only the outer frequency:
\[
  |\xi-\eta|\simeq2^l.
\]
The outer high restriction can be implemented by the fixed smooth cutoff \(\chi\left(2^{20}\frac{|\eta-\zeta|^2+|\zeta|^2}{|\xi-\eta|^2}\right)\).
On its support,
\[
  |\eta-\zeta|+|\zeta|\ll2^l,\qquad |(\xi-\eta)-\zeta|\simeq2^l.
\]
Lemma~\ref{lmet} thus yields
\[
  \left\|\frac{(\xi-\eta)\otimes((\xi-\eta)-\zeta)}{|(\xi-\eta)-\zeta|^2}\varphi_{[l-2,l+2]}(\xi-\eta)\chi\left(2^{20}\frac{|\eta-\zeta|^2+|\zeta|^2}{|\xi-\eta|^2}\right)\right\|_{\mathcal S^\infty}\lesssim 1.
\]
The same bound holds with \((\xi-\eta)-\zeta\) replaced by \((\xi-\eta)+(\eta-\zeta)\).

The translation argument in Lemma~\ref{l multilinear estimates}, followed by Lemma~\ref{luu1} and H\"older inequality, gives
\[
\begin{aligned}
  &\Bigg\|\iint_{|\xi-\eta|\gg|\eta-\zeta|+|\zeta|}e^{i\tau\Phi_3}\frac{(\xi-\eta)\otimes((\xi-\eta)-\zeta)}{|(\xi-\eta)-\zeta|^2} m_{-1}(\eta,\zeta)\overline{\widehat{x^2f}(\zeta-\eta)}\widehat{xf}(\zeta)\widehat{P_lf}(\xi-\eta)\,\mathrm{d}\zeta\,\mathrm{d}\eta\Bigg\|_{L_\xi^2}\\
  &\qquad\lesssim\|x^2f\|_{L_x^2}\|e^{i\tau\Delta}(xf)\|_{L_x^4}\|P_l\psi\|_{L_x^4}.
\end{aligned}
\]
Summing over \(l\), we obtain
\begin{equation}
\begin{aligned}
  &\Bigg\|\iint_{|\xi-\eta|\gg|\eta-\zeta|+|\zeta|}e^{i\tau\Phi_3}\frac{(\xi-\eta)\otimes((\xi-\eta)-\zeta)}{|(\xi-\eta)-\zeta|^2}m_{-1}(\eta,\zeta)\overline{\widehat{x^2f}(\zeta-\eta)}\widehat{xf}(\zeta)\widehat f(\xi-\eta)\,\mathrm{d}\zeta\,\mathrm{d}\eta\Bigg\|_{L_\xi^2}\\
  &\qquad\lesssim\|x^2f\|_{L_x^2}\|e^{i\tau\Delta}(xf)\|_{L_x^4}\sum_l\|P_l\psi\|_{L_x^4}.
\end{aligned} \label{ex2f11}
\end{equation}
By \eqref{bsc8} with $p=4$,
\begin{equation}
\sum_l\|P_l\psi(\tau)\|_{L_x^4}\lesssim\epsilon_1\tau^{-\frac12}. \label{ex2f12}
\end{equation}
Therefore, \eqref{ex2f11}, \eqref{ex2f12}, \eqref{booassum}, and \eqref{bsc10} give
\[
  \tau^{-1}\|x^2f(\tau)\|_{L_x^2} \|e^{i\tau\Delta}(xf(\tau))\|_{L_x^4}\sum_l\|P_l\psi(\tau)\|_{L_x^4}\lesssim\epsilon_1^3\tau^{-2+\frac{3}{2}\beta}.
\]

\smallskip
\noindent\textit{The derivative falls on the outer profile.}
The internal multiplier again remains \(m_{-1}\). Using \eqref{ex2f11} with the corresponding input substitution, followed by \eqref{bsc14} and \eqref{booassum}, we obtain
\[
\tau^{-1}\|x^2f(\tau)\|_{L_x^2}\|\psi(\tau)\|_{L_x^4}\sum_l\|P_le^{i\tau\Delta}(xf(\tau))\|_{L_x^4}\lesssim\epsilon_1^3\tau^{-2+\frac32\beta}(1+\log\tau).
\]

\smallskip
\noindent\textit{The derivative falls on the heat multiplier.}
A derivative of \(m_{-1}\) gives an order \(-2\) internal multiplier. The full coefficient is bounded by \((|\eta-\zeta|^2+|\zeta|^2)^{-1}\). Since \(k_1+k_2-2\max\{k_1,k_2\}=-|k_1-k_2|\), Minkowski inequality, translation invariance in \(L_\xi^2\), and Cauchy-Schwarz on each internal frequency annulus give
\begin{equation}
\begin{aligned}
  &\frac{1}{\tau}\Bigg\|\iint_{|\xi-\eta|\gg|\eta-\zeta|+|\zeta|}e^{i\tau\Phi_3}\frac{(\xi-\eta)\otimes((\xi-\eta)-\zeta)}{|(\xi-\eta)-\zeta|^2}m_{-2}(\eta,\zeta)\overline{\widehat{x^2f}(\zeta-\eta)}\widehat f(\zeta)\widehat f(\xi-\eta)\,\mathrm{d}\zeta\,\mathrm{d}\eta\Bigg\|_{L_\xi^2}\\
  &\quad\lesssim\frac{1}{\tau}\|f\|_{L_x^2}\iint\frac{|\widehat{x^2f}(\zeta-\eta)|\,|\widehat f(\zeta)|}{|\eta-\zeta|^2+|\zeta|^2}\,\mathrm{d}\zeta\,\mathrm{d}\eta\\
  &\quad\lesssim\frac{1}{\tau}\|f\|_{L_x^2}\sum_{k_1,k_2}2^{-|k_1-k_2|}\|P_{k_1}(x^2f)\|_{L_x^2}\|P_{k_2}f\|_{L_x^2}\\
  &\quad\lesssim\frac{1}{\tau}\|x^2f\|_{L_x^2}\|f\|_{L_x^2}^2\lesssim\epsilon_1^3\tau^{-1+\beta}.
\end{aligned} \label{ex2f13}
\end{equation}

\smallskip
\noindent\textit{The derivative falls on the integration-by-parts coefficient or the frequency cutoff.}
A derivative of the integration-by-parts coefficient or of the relative frequency cutoff contributes an additional \(|\xi-\eta|^{-1}\). On the present region, after localizing the three frequencies,
\[
  2^{-l-\max\{k_1,k_2\}}=2^{-(l-\max\{k_1,k_2\})}2^{-2\max\{k_1,k_2\}},\qquad l\geq\max\{k_1,k_2\}+10.
\]
The first factor is summable, while the second is the localized bound for an order \(-2\) internal multiplier. The same Fourier-Minkowski argument as in \eqref{ex2f13} therefore gives the bound \(\epsilon_1^3\tau^{-1+\beta}\).

Combining the two outer-frequency regions, we obtain
\begin{equation}
  \|\text{both derivatives on one internal profile}\|_{L_\xi^2}\lesssim\epsilon_1^3\left(\tau^{-2+\frac{3}{2}\beta}(1+\log\tau)+\tau^{-1+\beta}\right). \label{ex2f14}
\end{equation}

\smallskip
\noindent\textit{One derivative falls on each internal profile.}
Recall the smooth cutoff defining the outer high region. Suppose first that
\[
  |\xi-\eta|\lesssim|\eta-\zeta|+|\zeta|.
\]
Use the same two internal paraproducts and localized symbol bounds as in \eqref{ex2f9}. This time both internal inputs are placed in \(L_x^4\), \eqref{bsc10} and \eqref{epsihighfre2} give
\[
\begin{aligned}
  &\left\|\iint_{|\xi-\eta|\lesssim|\eta-\zeta|+|\zeta|}e^{i\tau\Phi_3}(\xi-\eta)m_{-1}(\eta,\zeta)\overline{\widehat{xf}(\zeta-\eta)}\widehat{xf}(\zeta)\widehat f(\xi-\eta)\,\mathrm{d}\zeta\,\mathrm{d}\eta\right\|_{L_\xi^2}\\
  &\qquad\lesssim\|e^{i\tau\Delta}(xf)\|_{L_x^4}^2\sum_l\|P_l\psi\|_{L_x^\infty}\sum_{k\geq l-O(1)}2^{l-k}\\
  &\qquad\lesssim\epsilon_1^3\tau^{-2+\beta}(1+\log\tau).
\end{aligned}
\]

On the region
\[
  |\xi-\eta|\geq2^{10}\bigl(|\eta-\zeta|+|\zeta|\bigr),
\]
we again use \eqref{ex2f10}. If the derivative falls on one internal profile, the same use of \eqref{ex2f11}, \eqref{ex2f12}, and \eqref{bsc10} gives
\[
  \tau^{-1}\|x^2f(\tau)\|_{L_x^2}\|e^{i\tau\Delta}(xf(\tau))\|_{L_x^4}\sum_l\|P_l\psi(\tau)\|_{L_x^4}\lesssim\epsilon_1^3\tau^{-2+\frac{3}{2}\beta}.
\]
If it falls on the outer profile, by \eqref{bsc14} and \eqref{bsc10},
\[
  \tau^{-1}\|xf(\tau)\|_{L_x^2}\|e^{i\tau\Delta}(xf(\tau))\|_{L_x^4}\sum_l\|P_le^{i\tau\Delta}(xf(\tau))\|_{L_x^4}\lesssim\epsilon_1^3\tau^{-2+\beta+\gamma}(1+\log\tau).
\]
If it falls on \(m_{-1}\), the same Fourier-Minkowski argument as in \eqref{ex2f13} gives
\[
  \tau^{-1}\|xf(\tau)\|_{L_x^2}^2\|f(\tau)\|_{L_x^2}\lesssim\epsilon_1^3\tau^{-1+2\gamma}.
\]
A derivative of the integration-by-parts coefficient or of the frequency cutoff again converts the remaining order \(-1\) multiplier into an order \(-2\) multiplier, with a summable outer inner frequency factor. Hence it satisfies the same \(\epsilon_1^3\tau^{-1+2\gamma}\) bound.

Combining these bounds we have proved
\begin{equation}
  \|\text{one derivative on each internal profile}\|_{L_\xi^2}\lesssim\epsilon_1^3\left(\tau^{-2+\frac{3}{2}\beta}(1+\log\tau)+\tau^{-1+2\gamma}\right). \label{ex2f15}
\end{equation}

\smallskip
\noindent\textit{One derivative falls on the heat multiplier and one derivative falls on an internal profile.}
We first consider
\[
  |\xi-\eta|\lesssim|\eta-\zeta|+|\zeta|.
\]
The same rescaling argument gives
\[
  \left\|(\xi-\eta)m_{-2}(\eta,\zeta)\varphi_l(\xi-\eta)\varphi_{k_1}(\eta-\zeta)\varphi_{k_2}(\zeta)\right\|_{\mathcal S^\infty}\lesssim2^{l-2\max\{k_1,k_2\}}.
\]
Then,
\[
\begin{aligned}
  &\sum_{k_1,k_2}\sum_{l\leq\max\{k_1,k_2\}+O(1)}2^{l-2\max\{k_1,k_2\}+\min\{k_1,k_2\}}\|P_{k_1}(xf)\|_{L_x^2}\|P_{k_2}f\|_{L_x^2}\|P_l\psi\|_{L_x^\infty}\\
  &\quad\lesssim\|\psi\|_{L_x^\infty}\sum_{k_1,k_2}2^{-|k_1-k_2|}\|P_{k_1}(xf)\|_{L_x^2}\|P_{k_2}f\|_{L_x^2}\\
  &\quad\lesssim\|xf\|_{L_x^2}\|f\|_{L_x^2}\|\psi\|_{L_x^\infty}.
\end{aligned}
\]
Thus
\[
\begin{aligned}
  &\left\|\iint_{|\xi-\eta|\lesssim|\eta-\zeta|+|\zeta|}e^{i\tau\Phi_3}(\xi-\eta)m_{-2}(\eta,\zeta)\overline{\widehat{xf}(\zeta-\eta)}\widehat f(\zeta)\widehat f(\xi-\eta)\,\mathrm{d}\zeta\,\mathrm{d}\eta\right\|_{L_\xi^2}\\
  &\qquad\lesssim\|xf\|_{L_x^2}\|f\|_{L_x^2}\|\psi\|_{L_x^\infty}\lesssim\epsilon_1^3\tau^{-1+\gamma}.
\end{aligned}
\]

On the region
\[
  |\xi-\eta|\geq2^{10}\bigl(|\eta-\zeta|+|\zeta|\bigr),
\]
we again use \eqref{ex2f10}, choosing the derivative which does not fall on the weighted internal profile. If the integration-by-parts derivative falls on the other internal profile, the Fourier-Minkowski argument in \eqref{ex2f13} gives
\[
  \tau^{-1}\|xf(\tau)\|_{L_x^2}^2\|f(\tau)\|_{L_x^2}\lesssim\epsilon_1^3\tau^{-1+2\gamma}.
\]
If it falls on the outer profile, the same argument gives
\[
  \tau^{-1}\|xf(\tau)\|_{L_x^2}\|f(\tau)\|_{L_x^2}\|xf(\tau)\|_{L_x^2}\lesssim\epsilon_1^3\tau^{-1+2\gamma}.
\]
It remains to estimate the terms in which the derivative falls on \(m_{-2}\), the integration-by-parts coefficient, or the relative frequency cutoff. The full coefficients are bounded by \((|\eta-\zeta|^2+|\zeta|^2)^{-\frac{3}{2}}\), for the last two terms, we use \(|\xi-\eta|^{-1}\lesssim(|\eta-\zeta|^2+|\zeta|^2)^{-\frac{1}{2}}\). Lemma~\ref{lmet} gives
\[
  \left\|m_{-3}(\eta,\zeta)\varphi_{k_1}(\eta-\zeta)\varphi_{k_2}(\zeta)\right\|_{\mathcal S^\infty}\lesssim 2^{-3\max\{k_1,k_2\}}.
\]
Hence,
\[
\begin{aligned}
  &\left\|\mathcal F_\eta^{-1}\left[\int e^{i\tau(|\eta-\zeta|^2-|\zeta|^2)}m_{-3}(\eta,\zeta)\overline{\widehat{P_{k_1}(xf)}(\zeta-\eta)}\widehat{P_{k_2}f}(\zeta)\,\mathrm{d}\zeta\right]\right\|_{L_x^\infty}\\
  &\quad\lesssim 2^{k_1+k_2-3\max\{k_1,k_2\}}\|P_{k_1}(xf)\|_{L_x^2}\|P_{k_2}f\|_{L_x^2}.
\end{aligned}
\]
If $k_1\leq k_2-10$, \eqref{bsc11} and \eqref{booassum} give
\[
\sum_{k_2}\sum_{k_1\leq k_2-10}2^{k_1-2k_2}\|P_{k_1}(xf)\|_{L_x^2}\|P_{k_2}f\|_{L_x^2}\lesssim\epsilon_1^2\tau^{\frac\beta2}\sum_{k_2<0}2^{\frac{k_2}{2}}+\epsilon_1^2\tau^\gamma\sum_{k_2\geq0}2^{-k_2}\lesssim\epsilon_1^2\tau^{\frac\beta2}.
\]
If \(k_2\leq k_1-10\), the same bound holds. If $|k_1-k_2|\leq10$, we have
\[
\sum_{|k_1-k_2|\leq10}2^{-\max\{k_1,k_2\}}\|P_{k_1}(xf)\|_{L_x^2}\|P_{k_2}f\|_{L_x^2}\lesssim\epsilon_1^2\tau^{\frac\beta2}\sum_{k_1<0}2^{\frac{k_1}{2}}+\epsilon_1^2\tau^\gamma\sum_{k_1\geq0}2^{-k_1}\lesssim\epsilon_1^2\tau^{\frac\beta2}.
\]
Combining the three frequency cases, we obtain
\begin{equation}
\begin{aligned}
  &\left\|\mathcal F_\eta^{-1}\left[\int e^{i\tau(|\eta-\zeta|^2-|\zeta|^2)}m_{-3}(\eta,\zeta)\overline{\widehat{xf}(\zeta-\eta)}\widehat f(\zeta)\,\mathrm{d}\zeta\right]\right\|_{L_x^\infty}\lesssim\epsilon_1^2\tau^{\frac{\beta}{2}}.
\end{aligned} \label{ex2f16}
\end{equation}
The same estimate holds with the two input positions exchanged. Using the Fourier-Minkowski argument in \eqref{ex2f13} and the frequency sum used to prove \eqref{ex2f16}, we obtain
\[
\begin{aligned}
  &\frac{1}{\tau}\|f\|_{L_x^2}\sum_{k_1,k_2}2^{k_1+k_2-3\max\{k_1,k_2\}}\|P_{k_1}(xf)\|_{L_x^2}\|P_{k_2}f\|_{L_x^2}\lesssim\epsilon_1^3\tau^{-1+\frac{\beta}{2}}.
\end{aligned}
\]
Therefore,
\begin{equation}
\|\text{one derivative on the multiplier and one on a profile}\|_{L_\xi^2}\lesssim\epsilon_1^3\left(\tau^{-1+2\gamma}+\tau^{-1+\frac\beta2}\right).
\label{ex2f17}
\end{equation}

\smallskip
\noindent\textit{Both derivatives fall on the heat multiplier.}
It remains to estimate
\begin{align*}
  \Bigg\|\iint e^{i\tau\Phi_3(\xi,\eta,\zeta)}(\xi-\eta)m_{-3}(\eta,\zeta)\overline{\widehat f(\tau,\zeta-\eta)}\widehat f(\tau,\zeta)\widehat f(\tau,\xi-\eta)\,\mathrm{d}\zeta\,\mathrm{d}\eta\Bigg\|_{L_\xi^2}.
\end{align*}
By the definition of $m_{-3}$, this contribution is supported where $|\eta-\zeta|^2+|\zeta|^2\gtrsim\tau^{-2}$. We now apply \eqref{ex2f7}. The second term converts $m_{-3}$ into an order $-2$ multiplier. In the first term, derivatives of the heat multiplier and its cutoff are included in $m_{-4}$.

On a dyadic annulus
\[
  |\eta-\zeta|^2+|\zeta|^2\simeq2^{2k}, \qquad\tau^{-1}\lesssim2^k\leq1,
\]
for fixed \(\eta\), the \(\zeta\)-integration region has area \(O(2^{2k})\), and \(|\eta|\lesssim2^k\). Hence
\[
\begin{aligned}
  &\sum_{\substack{k\leq0\\2^k\gtrsim\tau^{-1}}}\left\|\int_{|\eta-\zeta|^2+|\zeta|^2\simeq2^{2k}}e^{i\tau(|\eta-\zeta|^2-|\zeta|^2)}(\eta-\zeta)m_{-3}(\eta,\zeta)\overline{\widehat f(\zeta-\eta)}\widehat f(\zeta)\,\mathrm{d}\zeta\right\|_{L_\eta^2}\\
  &\qquad\lesssim\sum_{\substack{k\leq0\\2^k\gtrsim\tau^{-1}}}2^{-2k}2^{2k}2^k\|\widehat f(\tau)\|_{L_\xi^\infty}^2\lesssim\epsilon_1^2.
\end{aligned}
\]

For \(k\geq0\),
\[
  \left\|(\eta-\zeta)m_{-3}(\eta,\zeta)\varphi_k(\eta-\zeta)\varphi_{\leq k+5}(\zeta)\right\|_{\mathcal S^\infty}+\left\|(\eta-\zeta)m_{-3}(\eta,\zeta)\varphi_{\leq k+5}(\eta-\zeta)\varphi_k(\zeta)\right\|_{\mathcal S^\infty}\lesssim2^{-2k}.
\]
Lemma \ref{lmet} and Lemma \ref{l multilinear estimates} give
\[
\begin{aligned}
  &\sum_{k\geq0}\left\|\int_{|\eta-\zeta|^2+|\zeta|^2\simeq2^{2k}}e^{i\tau(|\eta-\zeta|^2-|\zeta|^2)}(\eta-\zeta)m_{-3}(\eta,\zeta)\overline{\widehat f(\zeta-\eta)}\widehat f(\zeta)\,\mathrm{d}\zeta\right\|_{L_\eta^2}\\
  &\qquad\lesssim\sum_{k\geq0}2^{-2k}\sum_{|k_1-k|\leq5}\|P_{k_1}f(\tau)\|_{L_x^2}\|P_{\leq k+5}\psi(\tau)\|_{L_x^\infty}\\
  &\qquad\lesssim\epsilon_1^2\tau^{-1}.
\end{aligned}
\]

Combining the low and high frequency estimates, the second term in \eqref{ex2f7} is bounded by \(\epsilon_1^3\tau^{-1}\).

It remains to integrate by parts in \(\eta\) in the first term on the right hand side of \eqref{ex2f7}.

When the derivative falls on \(m_{-3}\), the resulting multiplier is \(m_{-4}\). On
\[
  \tau^{-2}\lesssim|\eta-\zeta|^2+|\zeta|^2\lesssim 1,
\]
we place the outer profile in \(L_\xi^2\), the two internal profiles in Fourier \(L^\infty\), and \(m_{-4}\) in four-dimensional \(L^1\):
\[
\begin{aligned}
  \frac{1}{\tau}\|f(\tau)\|_{L_x^2}\|\widehat f(\tau)\|_{L_\xi^\infty}^2\iint_{\tau^{-2}\lesssim|\eta-\zeta|^2+|\zeta|^2\lesssim1}|m_{-4}(\eta,\zeta)|\,\mathrm{d}\zeta\,\mathrm{d}\eta&\lesssim\epsilon_1^3\tau^{-1}\sum_{\tau^{-1}\lesssim2^k\lesssim1}2^{4k}2^{-4k}\\
  &\lesssim\epsilon_1^3\tau^{-1}(1+\log\tau).
\end{aligned}
\]
On
\[
  |\eta-\zeta|^2+|\zeta|^2\gtrsim1,
\]
the order \(-4\) multiplier is in four-dimensional \(L^2\), Cauchy-Schwarz gives
\[
  \frac{1}{\tau}\|f(\tau)\|_{L_x^2}(\iint_{|\eta-\zeta|^2+|\zeta|^2\gtrsim1}|m_{-4}(\eta,\zeta)|^2\,\mathrm{d}\zeta\,\mathrm{d}\eta)^{\frac12}(\iint|\widehat f(\zeta-\eta)|^2|\widehat f(\zeta)|^2\,\mathrm{d}\zeta\,\mathrm{d}\eta)^{\frac12}\lesssim\epsilon_1^3\tau^{-1}.
\]

When the derivative falls on \(\widehat f(\zeta-\eta)\), the internal multiplier remains of order \(-3\). Splitting the integral and applying Cauchy-Schwarz gives
\[
\begin{aligned}
  &\frac{1}{\tau}\|f\|_{L_x^2}\iint_{|\eta-\zeta|^2+|\zeta|^2\gtrsim\tau^{-2}}|m_{-3}(\eta,\zeta)||\widehat{xf}(\zeta-\eta)|\,|\widehat f(\zeta)|\,\mathrm{d}\zeta\,\mathrm{d}\eta\\
  &\quad\lesssim\frac{1}{\tau}\|f\|_{L_x^2}\Bigg[\sum_{\tau^{-1}\lesssim2^k\lesssim1}2^{-3k}\|\widehat{P_{\leq k+5}(xf)}\|_{L_\xi^1}\|\widehat{P_{\leq k+5}f}\|_{L_\xi^1}+\left(\sum_{k\geq0}2^{4k}2^{-6k}\right)^{\frac12}\|xf\|_{L_x^2}\|f\|_{L_x^2}\Bigg]\\
  &\quad\lesssim\epsilon_1^2\tau^{-1}\sum_{\tau^{-1}\lesssim2^k\lesssim1}2^{-3k}2^k2^{2k}\|P_{\leq k+5}(xf)\|_{L_x^2}+\epsilon_1^3\tau^{-1+\gamma}\\
  &\quad\lesssim\epsilon_1^3\tau^{-1+\gamma}(1+\log\tau).
\end{aligned}
\]

When the derivative falls on the outer profile, the same argument bounds this contribution by
\[
\begin{aligned}
  &\frac{1}{\tau}\|xf\|_{L_x^2}\iint_{|\eta-\zeta|^2+|\zeta|^2\gtrsim\tau^{-2}}|m_{-3}(\eta,\zeta)||\widehat f(\zeta-\eta)|\,|\widehat f(\zeta)|\,\mathrm{d}\zeta\,\mathrm{d}\eta\\
  &\quad\lesssim\frac{1}{\tau}\|xf\|_{L_x^2}\left[\|\widehat f\|_{L_\xi^\infty}^2\sum_{\tau^{-1}\lesssim2^k\lesssim1}2^{4k}2^{-3k}+\|f\|_{L_x^2}^2\left(\sum_{k\geq0}2^{4k}2^{-6k}\right)^{\frac12}\right]\\
  &\quad\lesssim\epsilon_1^3\tau^{-1+\gamma}.
\end{aligned}
\]

Combining all contributions in this case, we obtain
\begin{equation}
\left\|\text{both derivatives on the heat multiplier}\right\|_{L_\xi^2}\lesssim\epsilon_1^3\tau^{-1+\gamma}(1+\log\tau). \label{ex2f18}
\end{equation}

Combining \eqref{ex2f14}, \eqref{ex2f15}, \eqref{ex2f17}, \eqref{ex2f18}, and the low frequency estimate above, we obtain
\begin{equation*}
\int_1^t\tau^2\|(\nabla^2A^{(2)}(\tau))\nabla\psi(\tau)\|_{L_x^2}\,\mathrm{d}\tau\lesssim\epsilon_1^3t^\beta.
\end{equation*}
Together with \eqref{ex2f6}, this proves
\begin{equation}
\int_1^t\tau^2\|(\nabla^2A(\tau))\nabla\psi(\tau)\|_{L_x^2}\,\mathrm{d}\tau\lesssim\epsilon_1^3t^\beta. \label{ex2f19}
\end{equation}

\subsection{The remaining nonlinearities}

We estimate the remaining contributions for $\tau\geq1$. By \eqref{eJproduct}, \eqref{bsc10}, and \eqref{booassum},
\begin{equation*}
\|J^2(\psi_l\operatorname{Im}(\overline{\psi_l}\psi_m))\|_{L_x^2}\lesssim\|\psi\|_{L_x^\infty}^2\|J^2\psi\|_{L_x^2}+\|\psi\|_{L_x^\infty}\|J\psi\|_{L_x^4}^2\lesssim\epsilon_1^3\tau^{-2+\beta}.
\end{equation*}
By \eqref{eJproduct}, \eqref{eAlinftydecay}, \eqref{eAonederivativelinfty}, \eqref{eAonederivativel21} with $\sigma=0$, and \eqref{exf11}, we obtain
\begin{align*}
\|J^2(A_lA_l\psi)\|_{L_x^2}&\lesssim\|A\|_{L_x^\infty}^2\|x^2f\|_{L_x^2}+\tau\|A\|_{L_x^\infty}\|\nabla A\|_{L_x^\infty}\|xf\|_{L_x^2}\\
&\quad+\tau^2\left(\|A\|_{L_x^\infty}\|\nabla^2A\|_{L_x^2}+\|\nabla A\|_{L_x^\infty}\|\nabla A\|_{L_x^2}\right)\|\psi\|_{L_x^\infty}\\
&\lesssim\epsilon_1^5\left[\tau^{-3+\beta}+\tau^{-\frac{5}{2}+\gamma}(1+\log\tau)+\tau^{-2}+\tau^{-2}(1+\log\tau)\right]\lesssim\epsilon_1^5\tau^{-\frac{3}{2}}.
\end{align*}
We next estimate the higher order spatial and temporal connections. By Gagliardo-Nirenberg inequality, Sobolev embedding, \eqref{eAhigh1}, \eqref{eAhigh3}, and \eqref{eAhigh2},
\begin{align*}
\|\nabla A^{(\geq4)}(\tau)\|_{L_x^\infty}&\lesssim\|A^{(\geq4)}(\tau)\|_{L_x^2}^{\frac{1}{3}}\|\nabla^3A^{(\geq4)}(\tau)\|_{L_x^2}^{\frac{2}{3}}\lesssim\epsilon_1^4\tau^{-\frac{8}{3}+\frac{2\gamma}{3}}(1+\log\tau)^{\frac{2}{3}}\lesssim\epsilon_1^4\tau^{-\frac{5}{2}},\\
\|\nabla^2A^{(\geq4)}(\tau)\|_{L_x^4}&\lesssim\|\nabla^2A^{(\geq4)}(\tau)\|_{H^1}\lesssim\epsilon_1^4\tau^{-\frac{5}{2}}.
\end{align*}
Together with \eqref{eJproduct}, \eqref{bsc10}, \eqref{booassum}, and \eqref{eAhigh3}, this gives
\begin{align*}
\|J^2((\nabla\cdot A^{(\geq4)})\psi)(\tau)\|_{L_x^2}&\lesssim\|\nabla A^{(\geq4)}\|_{L_x^\infty}\|J^2\psi\|_{L_x^2}+\tau\|\nabla^2A^{(\geq4)}\|_{L_x^4}\|J\psi\|_{L_x^4}+\tau^2\|\nabla^3A^{(\geq4)}\|_{L_x^2}\|\psi\|_{L_x^\infty}\\
&\lesssim\epsilon_1^5\left[\tau^{-\frac{5}{2}+\beta}+\tau^{-2+\frac{\beta}{2}}+\tau^{-2+\gamma}(1+\log\tau)\right].
\end{align*}
By \eqref{eA35}, \eqref{bsc10}, \eqref{booassum}, and Sobolev embedding,
\begin{equation*}
\tau\|(\nabla A_3^{(\geq4)})J\psi\|_{L_x^2}+\tau^2\|(\nabla^2A_3^{(\geq4)})\psi\|_{L_x^2}\lesssim\|A_3^{(\geq4)}\|_{H^2}\left(\tau\|J\psi\|_{L_x^4}+\tau^2\|\psi\|_{L_x^\infty}\right)\lesssim\epsilon_1^5\tau^{-1}.
\end{equation*}
For the quadratic contributions, \eqref{ehatf8} and \eqref{eA34} give the same phase $\Phi_3$. For $J^2((\nabla\cdot A^{(2)})\psi)$, we use $\nabla_\xi\Phi_3=-\nabla_\zeta\Phi_3$ and integrate by parts as in \eqref{ex2f8}. For $\tau(\nabla A_3^{(2)})J\psi$ and $\tau^2(\nabla^2A_3^{(2)})\psi$, we integrate by parts in $\zeta$ once and twice, respectively. Localizing afterwards, Lemma~\ref{lmet} gives
\begin{align*}
&\left\|\nabla_\zeta^a\left(\frac{\eta\zeta}{|\eta-\zeta|^2+|\zeta|^2}\right)\varphi_j(\eta-\zeta)\varphi_k(\zeta)\right\|_{\mathcal S^\infty}\\
&\quad+\left\|\nabla_\zeta^a\left(\frac{\zeta(\eta-\zeta)}{|\eta-\zeta|^2+|\zeta|^2}\right)\varphi_j(\eta-\zeta)\varphi_k(\zeta)\right\|_{\mathcal S^\infty}\lesssim2^{-a\max\{j,k\}},\qquad a=0,1,2.
\end{align*}
The conjugate contributions satisfy the same bounds. We estimate all three quadratic contributions together. By Lemma~\ref{lcoifman-meyer estimates}, the paraproduct estimates above, \eqref{epsihighfre2}, and \eqref{bsc10}, the terms with no derivative on the multiplier are bounded by
\begin{equation*}
\|x^2f\|_{L_x^2}\|\psi\|_{L_x^\infty}\sum_k\|P_k\psi\|_{L_x^\infty}+\|\psi\|_{L_x^\infty}\|J\psi\|_{L_x^4}^2\lesssim\epsilon_1^3\tau^{-2+\beta}(1+\log\tau).
\end{equation*}
By Lemma~\ref{luu1}, \eqref{exf6}, and \eqref{bso7}, the terms with one derivative on the multiplier are bounded by
\begin{equation*}
\|xf\|_{L_x^2}\left(\|\psi\|_{L_x^4}^2+\|\psi\|_{L_x^\infty}\sum_k2^{-k}\|P_k\psi\|_{L_x^\infty}\right)\lesssim\epsilon_1^3\tau^{-1+\gamma}.
\end{equation*}
By Lemma~\ref{l multilinear estimates}, Bernstein inequality, \eqref{epsidecay3}, and \eqref{booassum}, the terms with two derivatives on the multiplier are bounded by
\begin{equation*}
\epsilon_1^2\|\psi\|_{L_x^\infty}\left(\sum_{k\leq0}2^k+\sum_{k>0}2^{-(\epsilon+2)k}\right)\lesssim\epsilon_1^3\tau^{-1}.
\end{equation*}
Hence,
\begin{align*}
&\|J^2((\nabla\cdot A^{(2)})\psi)(\tau)\|_{L_x^2}+\tau\|(\nabla A_3^{(2)})J\psi\|_{L_x^2}+\tau^2\|(\nabla^2A_3^{(2)})\psi\|_{L_x^2}\\
&\quad\lesssim\epsilon_1^3\left[\tau^{-2+\beta}(1+\log\tau)+\tau^{-1+\gamma}\right].
\end{align*}
Combining these estimates and using \eqref{bop}, we obtain
\begin{equation}\label{ex2f20}
\int_1^t\left[\max_F\|J^2(F\psi)(\tau)\|_{L_x^2}+\tau\|(\nabla A_3)J\psi\|_{L_x^2}+\tau^2\|(\nabla^2A_3)\psi\|_{L_x^2}\right]\,\mathrm{d}\tau\lesssim\epsilon_1^3(1+t^\gamma)+\epsilon_1^5(1+\log t)\lesssim\epsilon_1^3t^\beta.
\end{equation}

\subsection{Conclusion}
Finally, inserting \eqref{ex2f5}, \eqref{ex2f19}, and \eqref{ex2f20} into \eqref{ex2f3}, together with the two lower order estimates, gives
\[
\|x^2f(t)\|_{L_x^2}\lesssim\|x^2f(1)\|_{L_x^2}+\epsilon_1^3t^\beta,
\]
which proves \eqref{ex2f1}.
\end{proof}

\section{Proofs of the main results}\label{sec:9}

By \cite[Theorem~1.2]{BIKT2011}, the data in \eqref{ini} determine a unique global solution in the small data solution class satisfying \eqref{ethm1.1-glo} and \eqref{ethm1.1-sobolev}. With $\epsilon_1=\epsilon_0^{\frac{1}{2}}$, the quantitative estimates obtained in closing Proposition~\ref{pmain bootstrap}, together with \eqref{eboot-init}, give the corresponding bounds with $C\epsilon_0$ in place of $\epsilon_1$. Throughout this section, we use these strengthened bounds for $t\geq1$. All decay estimates below are stated for $t\geq1$.

\subsection{Optimal decay of the map gradient}

By \eqref{efield-init2} and \eqref{etimebootstrapimprovement},
\[
  |\nabla\phi(t,x)|=|\psi(t,x)|, \qquad\|\nabla\phi(t)\|_{L_x^\infty}=\|\psi(t)\|_{L_x^\infty}\lesssim\epsilon_0t^{-1}, \qquad t\geq 1.
\]
This proves \eqref{ethm1.1-decay1}.

\subsection{Energy bounds and preliminary decay}

Proposition~\ref{pcaloric gauge} gives, for $s\geq0$ and $0\leq\sigma\leq1+\epsilon$,
\begin{equation}\label{ethmpf1}
\|\tilde\phi(s,t)-Q\|_{\dot H^\sigma}+\|v(s,t)-v_\infty\|_{\dot H^\sigma}+\|w(s,t)-w_\infty\|_{\dot H^\sigma}\lesssim\|\phi(t)-Q\|_{\dot H^\sigma}.
\end{equation}
We first compare the $H^\epsilon$ norms of $\nabla\phi$ and $\psi$ at $s=0$. Taking $\sigma=1$ in \eqref{ethmpf1} and using \eqref{enerconser}, we obtain
\begin{equation*}
\|\nabla v(t)\|_{L_x^2}+\|\nabla w(t)\|_{L_x^2}\lesssim\|\nabla\phi(t)\|_{L_x^2}=\|\nabla\phi_0\|_{L_x^2}\lesssim\epsilon_0.
\end{equation*}
By \eqref{efield2}, fractional Leibniz inequality, and Sobolev embedding,
\begin{equation*}
\begin{aligned}
\|\nabla\phi(t)\|_{H^\epsilon}&\lesssim\|\psi(t)\|_{H^\epsilon}+\|\psi(t)\|_{L_x^{\frac{2}{1-\epsilon}}}\left(\bigl\||\nabla|^\epsilon(v(t)-v_\infty)\bigr\|_{L_x^{\frac{2}{\epsilon}}}+\bigl\||\nabla|^\epsilon(w(t)-w_\infty)\bigr\|_{L_x^{\frac{2}{\epsilon}}}\right)\\
&\lesssim\left(1+\|\nabla v(t)\|_{L_x^2}+\|\nabla w(t)\|_{L_x^2}\right)\|\psi(t)\|_{H^\epsilon}\lesssim\|\psi(t)\|_{H^\epsilon}.
\end{aligned}
\end{equation*}
Applying the same argument to \eqref{efield1} gives the reverse estimate. Hence $\|\nabla\phi(t)\|_{H^\epsilon}\simeq\|\psi(t)\|_{H^\epsilon}$.

We next prove the uniform Sobolev bounds for the heat flow and the frame. Applying \eqref{ethmpf1} with $\sigma=0,1+\epsilon$ and \eqref{ethm1.1-sobolev}, we obtain
\begin{equation}\label{ethmpf2}
\begin{aligned}
&\|\tilde\phi(s,t)-Q\|_{H^{1+\epsilon}}+\|v(s,t)-v_\infty\|_{H^{1+\epsilon}}+\|w(s,t)-w_\infty\|_{H^{1+\epsilon}}\\
&\qquad\lesssim\|\phi(t)-Q\|_{H^{1+\epsilon}}\lesssim\|\phi_0-Q\|_{H^{1+\epsilon}}\lesssim\epsilon_0,\qquad s\geq0,\quad t\in\RR.
\end{aligned}
\end{equation}

To obtain preliminary uniform decay, use \eqref{eAinfty1} and \eqref{epsidecay1} to estimate
\[
  \int_s^\infty\|A(r,t)\psi(r,t)\|_{L_x^\infty}\,\mathrm{d}r\lesssim\epsilon_0^3(t+s)^{-\frac{3}{2}},\qquad t\geq 1.
\]
It follows from \eqref{ecaloric3} and \eqref{epsiinte1} that
\begin{equation}\label{ethmpf3}
  \int_s^\infty\|\psi_0(r,t)\|_{L_x^\infty}\,\mathrm{d}r\lesssim\epsilon_0(t+s)^{-\frac{1}{2}},\qquad t\geq 1.
\end{equation}
The heat direction version of \eqref{efield2}, with $A_0=0$, gives
\[
  \partial_s\tilde\phi=v\operatorname{Re}\psi_0+w\operatorname{Im}\psi_0, \qquad\partial_sv=-\tilde\phi\operatorname{Re}\psi_0, \qquad\partial_sw=-\tilde\phi\operatorname{Im}\psi_0.
\]
Integrating from $s$ to infinity, using \eqref{ethmpf3} and the terminal limits in Proposition~\ref{pcaloric gauge}, yields
\begin{equation}\label{ethmpf4}
  \|\tilde\phi(s,t)-Q\|_{L_x^\infty}+\|v(s,t)-v_\infty\|_{L_x^\infty}+\|w(s,t)-w_\infty\|_{L_x^\infty}\lesssim\epsilon_0(t+s)^{-\frac{1}{2}},\qquad t\geq 1.
\end{equation}

\subsection{A complex coordinate and the reduced equation}

Unless the heat time is displayed explicitly, all fields in this subsection are evaluated at $s=0$.

\smallskip
\noindent\textit{The coordinate and its algebraic identities.}
Define
\begin{equation}\label{ecomplex-coordinate}
  z(t,x)=\frac{(v_\infty+iw_\infty)\cdot\phi(t,x)}{1+Q\cdot\phi(t,x)}.
\end{equation}
By \eqref{ethm1.1-sobolev} and $H^{1+\epsilon}\hookrightarrow L_x^\infty$, taking $\epsilon_0$ sufficiently small gives $\|\phi(t)-Q\|_{L_x^\infty}\leq1$ for every $t\in\RR$. Hence
\begin{equation*}
1+Q\cdot\phi=2-\frac12|\phi-Q|^2\geq\frac32.
\end{equation*}
Since $(Q,v_\infty,w_\infty)$ is orthonormal,
\begin{equation}\label{ezphilinfty}
  |z|^2=\frac{1-Q\cdot\phi}{1+Q\cdot\phi}=\frac{|\phi-Q|^2}{4-|\phi-Q|^2}, \qquad\frac{1}{2}|\phi-Q|\leq|z|\leq\frac{1}{\sqrt{3}}|\phi-Q|.
\end{equation}
Thus \eqref{emass-conservation}, \eqref{ini}, and \eqref{ethmpf4} give
\begin{equation}\label{ezprebound}
  \|z(t)\|_{L_x^2}\lesssim\epsilon_0, \qquad\|z(t)\|_{L_x^\infty}\lesssim\epsilon_0t^{-\frac{1}{2}}, \qquad t\geq 1.
\end{equation}
Solving for the three components of $\phi$ gives
\begin{equation}\label{einversecoor}
  \phi=\frac{(1-|z|^2)Q+2v_\infty\operatorname{Re}z+2w_\infty\operatorname{Im}z}{1+|z|^2}, \qquad 1+Q\cdot\phi=\frac{2}{1+|z|^2}.
\end{equation}
In particular,
\begin{equation}\label{ethmpf5}
\begin{aligned}
  (v_\infty+iw_\infty-zQ)\cdot\phi&=z,\\
  (v_\infty+iw_\infty-zQ)\times\phi&=i(v_\infty+iw_\infty-zQ-z\phi).
\end{aligned}
\end{equation}
Taking the dot product of the last identity with $v$ and using $w=\phi\times v$, we obtain
\begin{equation}\label{ethmpf6}
\begin{aligned}
  (v_\infty+iw_\infty-zQ)\cdot w&=i(v_\infty+iw_\infty-zQ)\cdot v,\\
  \bigl|(v_\infty+iw_\infty-zQ)\cdot v\bigr|&=1.
\end{aligned}
\end{equation}
Indeed, the second identity follows from the first and \eqref{ethmpf5}, since orthogonal decomposition in $(\phi,v,w)$ gives
\[
  2+|z|^2=|v_\infty+iw_\infty-zQ|^2=|z|^2+2\bigl|(v_\infty+iw_\infty-zQ)\cdot v\bigr|^2.
\]
Taking the dot product of the same decomposition with $Q$ yields
\begin{equation}\label{ethmpf7}
  z=-\frac{(v_\infty+iw_\infty-zQ)\cdot v}{1+Q\cdot\phi}(Q\cdot v+iQ\cdot w).
\end{equation}

\smallskip
\noindent\textit{Relation to the differential fields.}
Differentiating \eqref{ecomplex-coordinate} and using \eqref{efield2} and \eqref{ethmpf6}, for $m=1,2,3$,
\begin{equation}\label{ethmpf8}
  \partial_mz=\frac{(v_\infty+iw_\infty-zQ)\cdot v}{1+Q\cdot\phi}\psi_m.
\end{equation}
Hence, \eqref{einversecoor} and \eqref{ethmpf6} imply
\begin{equation}\label{ezgradient}
  |\nabla z|=\frac{1+|z|^2}{2}|\psi|=\frac{1+|z|^2}{2}|\nabla\phi|, \qquad\|\nabla z(t)\|_{L_x^\infty}\lesssim\epsilon_0t^{-1}, \qquad t\geq 1.
\end{equation}

\smallskip
\noindent\textit{Differentiation of the complex coefficient.}
By \eqref{efield2} and \eqref{ethmpf5},
\[
\begin{aligned}
  &\partial_m\left[\frac{(v_\infty+iw_\infty-zQ)\cdot v}{1+Q\cdot\phi}\right]\\
  &=iA_m\frac{(v_\infty+iw_\infty-zQ)\cdot v}{1+Q\cdot\phi}-\frac{1}{1+Q\cdot\phi}\Bigg[(Q\cdot v)\partial_mz+z\operatorname{Re}\psi_m\\
  &\hspace{24mm}+\frac{(v_\infty+iw_\infty-zQ)\cdot v}{1+Q\cdot\phi}\bigl((Q\cdot v)\operatorname{Re}\psi_m+(Q\cdot w)\operatorname{Im}\psi_m\bigr)\Bigg]\\
  &=\left(iA_m-\frac{(Q\cdot v-iQ\cdot w)\psi_m}{1+Q\cdot\phi}\right)\frac{(v_\infty+iw_\infty-zQ)\cdot v}{1+Q\cdot\phi}.
\end{aligned}
\]
The last equality follows by substituting \eqref{ethmpf7} and \eqref{ethmpf8}. Taking the conjugate of \eqref{ethmpf7} and using \eqref{ethmpf6}, we also have
\[
  -\frac{(Q\cdot v-iQ\cdot w)\psi_m}{1+Q\cdot\phi}=\frac{\overline z\,\psi_m}{\overline{(v_\infty+iw_\infty-zQ)\cdot v}}=(1+Q\cdot\phi)\overline z\,\partial_mz=\frac{2\overline z\,\partial_mz}{1+|z|^2}.
\]
Therefore we prove
\begin{equation}\label{ethmpf9}
  \partial_m\left[\frac{(v_\infty+iw_\infty-zQ)\cdot v}{1+Q\cdot\phi}\right]=\left(iA_m+\frac{2\overline z\,\partial_mz}{1+|z|^2}\right)\frac{(v_\infty+iw_\infty-zQ)\cdot v}{1+Q\cdot\phi}.
\end{equation}
By \eqref{eAinfty1}, \eqref{ezprebound}, and \eqref{ezgradient},
\begin{equation}\label{ethmpf10}
  \left\|\nabla\left[\frac{(v_\infty+iw_\infty-zQ)\cdot v}{1+Q\cdot\phi}\right]\right\|_{L_x^\infty}\lesssim\|A(t)\|_{L_x^\infty}+\|z(t)\|_{L_x^\infty}\|\nabla z(t)\|_{L_x^\infty}\lesssim\epsilon_0^2t^{-\frac{3}{2}}, \qquad t\geq 1.
\end{equation}

\smallskip
\noindent\textit{The reduced equation.}
Differentiating \eqref{ethmpf8} once more and using \eqref{ethmpf9}, we obtain
\[
  \partial_l^2z=\frac{(v_\infty+iw_\infty-zQ)\cdot v}{1+Q\cdot\phi}D_l\psi_l+\frac{2\overline z}{1+|z|^2}(\partial_lz)^2.
\]
On the other hand, \eqref{efield4} gives
\[
  i\partial_tz=-\frac{(v_\infty+iw_\infty-zQ)\cdot v}{1+Q\cdot\phi}\sum_{l=1}^2D_l\psi_l.
\]
Thus
\begin{equation}\label{ereducedequation}
  (i\partial_t+\Delta)z=\frac{2\overline z}{1+|z|^2}\sum_{l=1}^2(\partial_lz)^2.
\end{equation}

\begin{remark}
Equation \eqref{einversecoor} is the stereographic formula in the fixed basis $(Q,v_\infty,w_\infty)$. We keep $Q$ arbitrary. Equations \eqref{ethmpf8} and \eqref{ethmpf9} give the explicit relation between this coordinate and the caloric differential fields, the connection coefficients.
\end{remark}

\smallskip
\noindent\textit{Sobolev and initial weighted bounds.}
Note
\[
  z=\frac{(v_\infty+iw_\infty)\cdot(\phi-Q)}{2+Q\cdot(\phi-Q)},
\]
whose denominator is uniformly separated from zero. Applying Moser composition estimate with \eqref{ethmpf2} gives
\begin{equation}\label{sob}
\|z(t)\|_{H^{1+\epsilon}}\lesssim\|\phi(t)-Q\|_{H^{1+\epsilon}}\lesssim\epsilon_0,\qquad t\in\RR.
\end{equation}
At $t=0$, \eqref{ini}, \eqref{ezphilinfty}, and \eqref{ezgradient} yield
\begin{equation}\label{ezinitialbounds}
\|z(0)\|_{H^{1+\epsilon}}+\|\langle x\rangle^2\langle\nabla\rangle z(0)\|_{L_x^2}\lesssim\epsilon_0.
\end{equation}

We state several results of Proposition~\ref{locest} for the weighted estimates below. By \eqref{elocalestimate},
\begin{equation}\label{elocalpsibound}
\sup_{0\leq t\leq1}\left(\|\psi(t)\|_{H^\epsilon}+\|J\psi(t)\|_{L_x^2}+\|J^2\psi(t)\|_{L_x^2}\right)+\left(\int_0^1\|\psi(t)\|_{L_x^\infty}^2\,\mathrm dt\right)^{\frac{1}{2}}\lesssim\epsilon_0.
\end{equation}
Moreover, the dispersive estimate and interpolation give, for $0<t\leq1$,
\begin{equation*}
t\|\psi(t)\|_{L_x^\infty}\lesssim\|f(t)\|_{L_x^1}\lesssim\|f(t)\|_{L_x^2}^{\frac{1}{2}}\||x|^2f(t)\|_{L_x^2}^{\frac{1}{2}}\lesssim\epsilon_0.
\end{equation*}
To estimate the connection at these times, fix $0<t\leq1$. By \eqref{ecaloric8} with $\sigma=1$ and energy conservation,
\begin{align*}
\|\psi(s,t)\|_{L_x^\infty}&=\|\nabla\tilde\phi(s,t)\|_{L_x^\infty}\lesssim\sum_k2^k\gamma_k(1)(1+s2^{2k})^{-20}\lesssim\epsilon_0s^{-\frac{1}{2}},\\
\|\nabla^2\tilde\phi(s,t)\|_{L_x^\infty}&\lesssim\sum_k2^{2k}\gamma_k(1)(1+s2^{2k})^{-20}\lesssim\epsilon_0s^{-1},\\
\|\psi_0(s,t)\|_{L_x^\infty}&=\|\partial_s\tilde\phi(s,t)\|_{L_x^\infty}\leq\|\Delta\tilde\phi(s,t)\|_{L_x^\infty}+\|\nabla\tilde\phi(s,t)\|_{L_x^\infty}^2\lesssim\epsilon_0s^{-1}.
\end{align*}
On the other hand, the heat flow Duhamel formula gives
\begin{equation*}
\nabla\tilde\phi(s,t)=e^{s\Delta}\nabla\phi(t)+\int_0^s\nabla e^{(s-r)\Delta}\left[\tilde\phi(r,t)|\nabla\tilde\phi(r,t)|^2\right]\,\mathrm dr.
\end{equation*}
Using $|\tilde\phi|=1$ and \eqref{eheatoperator}, we obtain
\begin{equation*}
\begin{aligned}
\|\nabla\tilde\phi(s,t)\|_{L_x^\infty}&\leq\|\psi(t)\|_{L_x^\infty}+C\int_0^s(s-r)^{-\frac{1}{2}}\|\nabla\tilde\phi(r,t)\|_{L_x^\infty}^2\,\mathrm dr,\\
\sup_{0\leq r\leq s}\|\nabla\tilde\phi(r,t)\|_{L_x^\infty}&\leq\|\psi(t)\|_{L_x^\infty}+Cs^{\frac{1}{2}}\left(\sup_{0\leq r\leq s}\|\nabla\tilde\phi(r,t)\|_{L_x^\infty}\right)^2.
\end{aligned}
\end{equation*}
For heat times satisfying $s^{\frac{1}{2}}\|\psi(t)\|_{L_x^\infty}>\epsilon_0$, we use the heat flow envelope bounds above. We now consider the case $s^{\frac{1}{2}}\|\psi(t)\|_{L_x^\infty}\leq\epsilon_0$. Fix $s$, if for some $0<s'\leq s$, we have $\sup_{0\leq r\leq s'}\|\nabla\tilde\phi(r,t)\|_{L_x^\infty}\leq2\|\psi(t)\|_{L_x^\infty}$, then this estimate gives
\begin{equation*}
\sup_{0\leq r\leq s'}\|\nabla\tilde\phi(r,t)\|_{L_x^\infty}\leq\|\psi(t)\|_{L_x^\infty}+4C{s'}^{\frac{1}{2}}\|\psi(t)\|_{L_x^\infty}^2\leq\frac{3}{2}\|\psi(t)\|_{L_x^\infty},
\end{equation*}
where we used ${s'}^{\frac{1}{2}}\|\psi(t)\|_{L_x^\infty}\leq s^{\frac{1}{2}}\|\psi(t)\|_{L_x^\infty}\leq\epsilon_0$ and the smallness of $\epsilon_0$. Since $\|\nabla\tilde\phi(0,t)\|_{L_x^\infty}=\|\psi(t)\|_{L_x^\infty}$, a continuity argument in $s'$ yields
\begin{equation*}
\sup_{0\leq r\leq s}\|\nabla\tilde\phi(r,t)\|_{L_x^\infty}\leq2\|\psi(t)\|_{L_x^\infty},\qquad s^{\frac{1}{2}}\|\psi(t)\|_{L_x^\infty}\leq\epsilon_0.
\end{equation*} 
Still assuming that $s>0$ and $s^{\frac{1}{2}}\|\psi(t)\|_{L_x^\infty}\leq\epsilon_0$, we differentiate the heat flow Duhamel formula and use the bound above for $\nabla\tilde\phi$. For every $0<s'\leq s$, this gives
\begin{equation*}
\|\nabla^2\tilde\phi(s',t)\|_{L_x^\infty}\lesssim{s'}^{-\frac{1}{2}}\|\psi(t)\|_{L_x^\infty}+\int_0^{s'}(s'-r)^{-\frac{1}{2}}\left(\|\psi(t)\|_{L_x^\infty}^3+\|\psi(t)\|_{L_x^\infty}\|\nabla^2\tilde\phi(r,t)\|_{L_x^\infty}\right)\,\mathrm dr.
\end{equation*}
Multiplying by ${s'}^{\frac{1}{2}}$, taking the supremum over $0<s'\leq s$, we obtain
\begin{equation*}
\sup_{0<r\leq s}r^{\frac{1}{2}}\|\nabla^2\tilde\phi(r,t)\|_{L_x^\infty}\lesssim\left(1+s\|\psi(t)\|_{L_x^\infty}^2\right)\|\psi(t)\|_{L_x^\infty}+s^{\frac{1}{2}}\|\psi(t)\|_{L_x^\infty}\sup_{0<r\leq s}r^{\frac{1}{2}}\|\nabla^2\tilde\phi(r,t)\|_{L_x^\infty}.
\end{equation*}
Since $s^{\frac{1}{2}}\|\psi(t)\|_{L_x^\infty}\leq\epsilon_0$, absorbing the last term gives
\begin{equation*}
\sup_{0<r\leq s}r^{\frac{1}{2}}\|\nabla^2\tilde\phi(r,t)\|_{L_x^\infty}\lesssim(1+\epsilon_0^2)\|\psi(t)\|_{L_x^\infty}\lesssim\|\psi(t)\|_{L_x^\infty}.
\end{equation*}
By the heat flow equation and the same restriction on $s$,
\begin{equation*}
\begin{aligned}
\|\psi_0(s,t)\|_{L_x^\infty}&=\|\partial_s\tilde\phi(s,t)\|_{L_x^\infty}\lesssim\|\nabla^2\tilde\phi(s,t)\|_{L_x^\infty}+\|\nabla\tilde\phi(s,t)\|_{L_x^\infty}^2\\
&\lesssim s^{-\frac{1}{2}}\|\psi(t)\|_{L_x^\infty}+\|\psi(t)\|_{L_x^\infty}^2\lesssim(1+\epsilon_0)s^{-\frac{1}{2}}\|\psi(t)\|_{L_x^\infty}.
\end{aligned}
\end{equation*}
For $s^{\frac{1}{2}}\|\psi(t)\|_{L_x^\infty}>\epsilon_0$, the heat flow envelope bounds give $\|\psi(s,t)\|_{L_x^\infty}\lesssim\epsilon_0s^{-\frac{1}{2}}$ and $\|\psi_0(s,t)\|_{L_x^\infty}\lesssim\epsilon_0s^{-1}$. Combining the two ranges, we conclude that, for $s>0$ and $0<t\leq1$,
\begin{align*}
\|\psi(s,t)\|_{L_x^\infty}&\lesssim\min\left\{\|\psi(t)\|_{L_x^\infty},\epsilon_0s^{-\frac{1}{2}}\right\},\\
\|\psi_0(s,t)\|_{L_x^\infty}&\lesssim s^{-\frac{1}{2}}\min\left\{\|\psi(t)\|_{L_x^\infty},\epsilon_0s^{-\frac{1}{2}}\right\}.
\end{align*}
\eqref{ecaloric6} therefore implies
\begin{equation*}
\begin{aligned}
\|A(0,t)\|_{L_x^\infty}&\lesssim\int_0^\infty\|\psi_0(s,t)\|_{L_x^\infty}\|\psi(s,t)\|_{L_x^\infty}\,\mathrm ds\\
&\lesssim\int_0^\infty s^{-\frac{1}{2}}\min\left\{\|\psi(t)\|_{L_x^\infty}^2,\epsilon_0^2s^{-1}\right\}\,\mathrm ds\lesssim\epsilon_0\|\psi(t)\|_{L_x^\infty}.
\end{aligned}
\end{equation*}
Using \eqref{ethmpf9}, \eqref{ethmpf8}, and \eqref{sob}, we conclude that
\begin{equation}\label{ethmpf11}
\begin{aligned}
\left\|\nabla\left[\frac{(v_\infty+iw_\infty-zQ)\cdot v}{1+Q\cdot\phi}\right]\right\|_{L_x^\infty}&\lesssim\|A(0,t)\|_{L_x^\infty}+\|z(t)\|_{L_x^\infty}\|\nabla z(t)\|_{L_x^\infty}\lesssim\epsilon_0\|\psi(t)\|_{L_x^\infty},\\
t\left\|\nabla\left[\frac{(v_\infty+iw_\infty-zQ)\cdot v}{1+Q\cdot\phi}\right]\right\|_{L_x^\infty}&\lesssim\epsilon_0^2,\qquad 0<t\leq1.
\end{aligned}
\end{equation}

\subsection{The first and second commuted energy estimates}

\begin{lemma}\label{lJz}
Let $z$ be defined by \eqref{ecomplex-coordinate} with the initial data satisfying \eqref{ezinitialbounds}, then
\begin{equation}\label{eJz}
  \sup_{t\geq0}\|xe^{-it\Delta}z(t)\|_{L_x^2}=\sup_{t\geq0}\|Jz(t)\|_{L_x^2}\lesssim\epsilon_0.
\end{equation}
\end{lemma}

\begin{proof}
Recall that $J_j=x_j+2it\partial_j$ and $[i\partial_t+\Delta,J_j]=0$. The identities
\[
\begin{aligned}
  J_j\left(\sum_{l=1}^2(\partial_lz)^2\right)&=2\sum_{l=1}^2\partial_lz\,\partial_lJ_jz-2z\partial_jz-x_j\sum_{l=1}^2(\partial_lz)^2,\\
  2it\partial_j\left(\frac{2\overline z}{1+|z|^2}\right)&=-\frac{2\overline z^{\,2}J_jz+2\overline{J_jz}}{(1+|z|^2)^2}+\frac{2x_j\overline z}{1+|z|^2}
\end{aligned}
\]
give, upon applying $J_j$ to \eqref{ereducedequation},
\begin{equation}\label{eJz1}
\begin{aligned}
  (i\partial_t+\Delta)J_jz=&\frac{4\overline z}{1+|z|^2}\sum_{l=1}^2\partial_lz\,\partial_lJ_jz-\frac{2\overline z^{\,2}}{(1+|z|^2)^2}\left(\sum_{l=1}^2(\partial_lz)^2\right)J_jz\\
  &-\frac{2}{(1+|z|^2)^2}\left(\sum_{l=1}^2(\partial_lz)^2\right)\overline{J_jz}-\frac{4|z|^2}{1+|z|^2}\partial_jz.
\end{aligned}
\end{equation}
By \eqref{ethmpf8},
\begin{equation}\label{eJz2}
  \partial_lJ_jz=\frac{(v_\infty+iw_\infty-zQ)\cdot v}{1+Q\cdot\phi}J_j\psi_l+2it\partial_j\left[\frac{(v_\infty+iw_\infty-zQ)\cdot v}{1+Q\cdot\phi}\right]\psi_l+\delta_{jl}z.
\end{equation}
For $0<t\leq1$, \eqref{eJz2}, \eqref{ethmpf11}, \eqref{elocalpsibound}, and \eqref{sob} imply
\begin{equation*}
\|\nabla Jz(t)\|_{L_x^2}\lesssim\|J\psi(t)\|_{L_x^2}+\epsilon_0^2\|\psi(t)\|_{L_x^2}+\|z(t)\|_{L_x^2}\lesssim\epsilon_0.
\end{equation*}
The $L_x^2$ energy estimate for \eqref{eJz1} gives
\begin{equation*}
\begin{aligned}
\frac{\mathrm d}{\mathrm dt}\|J_jz(t)\|_{L_x^2}&\lesssim\|z(t)\|_{L_x^\infty}\|\nabla z(t)\|_{L_x^\infty}\|\nabla J_jz(t)\|_{L_x^2}+\|\nabla z(t)\|_{L_x^\infty}^2\|J_jz(t)\|_{L_x^2}\\
&\quad+\|z(t)\|_{L_x^\infty}\|\nabla z(t)\|_{L_x^\infty}\|z(t)\|_{L_x^2}\\
&\lesssim\|\psi(t)\|_{L_x^\infty}^2\|J_jz(t)\|_{L_x^2}+\epsilon_0^2\|\psi(t)\|_{L_x^\infty}.
\end{aligned}
\end{equation*}
By Gronwall inequality, \eqref{elocalpsibound}, and \eqref{ezinitialbounds},
\begin{equation*}
\begin{aligned}
\sup_{0\leq t\leq1}\|Jz(t)\|_{L_x^2}&\lesssim\left(\|xz(0)\|_{L_x^2}+\epsilon_0^2\int_0^1\|\psi(t)\|_{L_x^\infty}\,\mathrm dt\right)\exp\left(C\int_0^1\|\psi(t)\|_{L_x^\infty}^2\,\mathrm dt\right)\\
&\lesssim\left(\|xz(0)\|_{L_x^2}+\epsilon_0^3\right)e^{C\epsilon_0^2}\lesssim\epsilon_0.
\end{aligned}
\end{equation*}
Thus
\begin{equation}\label{eJz3}
\sup_{0<t\leq1}\left(\|Jz(t)\|_{L_x^2}+\|\nabla Jz(t)\|_{L_x^2}\right)\lesssim\epsilon_0.
\end{equation}
For $t\geq1$, \eqref{eJz2}, \eqref{ethmpf10}, \eqref{sob}, and \eqref{etimebootstrapimprovement} give
\begin{equation}\label{eJz4}
\|\nabla Jz(t)\|_{L_x^2}\lesssim\|J\psi(t)\|_{L_x^2}+t\epsilon_0^2t^{-\frac{3}{2}}\|\psi(t)\|_{L_x^2}+\|z(t)\|_{L_x^2}\lesssim\epsilon_0t^\gamma.
\end{equation}
Using \eqref{ezprebound}, \eqref{ezgradient}, and \eqref{eJz4} in the same energy estimate yields
\begin{equation*}
\frac{\mathrm d}{\mathrm dt}\|J_jz(t)\|_{L_x^2}\lesssim\epsilon_0^2t^{-2}\|J_jz(t)\|_{L_x^2}+\epsilon_0^3t^{-\frac{3}{2}+\gamma},\qquad t\geq1.
\end{equation*}
Gronwall inequality gives
\begin{equation*}
\sup_{t\geq1}\|Jz(t)\|_{L_x^2}\lesssim\left(\|Jz(1)\|_{L_x^2}+\epsilon_0^3\int_1^\infty t^{-\frac{3}{2}+\gamma}\,\mathrm dt\right)\exp\left(C\epsilon_0^2\int_1^\infty t^{-2}\,\mathrm dt\right)\lesssim\epsilon_0.
\end{equation*}
Combining the two time intervals proves
\begin{equation}\label{eJz5}
\sup_{t\geq0}\|Jz(t)\|_{L_x^2}\lesssim\|xz(0)\|_{L_x^2}+\epsilon_0^3\lesssim\epsilon_0.
\end{equation}
\end{proof}

For $t\geq1$, the identity \(\partial_j\left(e^{-\frac{i|x|^2}{4t}}z\right)=\frac{1}{2it}e^{-\frac{i|x|^2}{4t}}J_jz\) and Gagliardo-Nirenberg inequality give
\[
  \|z(t)\|_{L_x^p}\lesssim t^{-1+\frac{2}{p}}\|z(t)\|_{L_x^2}^{\frac{2}{p}}\|Jz(t)\|_{L_x^2}^{1-\frac{2}{p}}\lesssim\epsilon_0t^{-1+\frac{2}{p}}, \qquad p=4,6.
\]
Together with \eqref{ezgradient},
\[
  \|z(t)\|_{L_x^\infty}\lesssim\|z(t)\|_{L_x^6}^{\frac{3}{4}}\|\nabla z(t)\|_{L_x^\infty}^{\frac{1}{4}}\lesssim\epsilon_0t^{-\frac{3}{4}}.
\]
Hence,
\begin{equation}\label{eJz6}
  \|z(t)\|_{L_x^4}\lesssim\epsilon_0t^{-\frac{1}{2}},\qquad\|z(t)\|_{L_x^\infty}\lesssim\epsilon_0t^{-\frac{3}{4}},\qquad t\geq1.
\end{equation}

\begin{lemma}\label{lemJ2z}
Let $z$ be defined by \eqref{ecomplex-coordinate} with the initial data satisfying \eqref{ezinitialbounds}, then
\begin{equation}\label{eJ2z}
  \sup_{t\geq0}\|x^2e^{-it\Delta}z(t)\|_{L_x^2}=\sup_{t\geq0}\|J^2z(t)\|_{L_x^2}\lesssim\epsilon_0.
\end{equation}
\end{lemma}

\begin{proof}
Apply $J_i$ to \eqref{eJz1}, using
\(
  2it\partial_i(1+|z|^2)=\overline z\,J_i z-z\overline{J_i z},2it\partial_i\partial_lz=\partial_lJ_i z-\delta_{il}z-x_i\partial_lz,J_i\overline{J_jz}=-\overline{J_iJ_jz}+2x_i\overline{J_jz},
\)
the differentiated coefficients are
\[
  2it\partial_i\left(\frac{4\overline z\,\partial_lz}{1+|z|^2}\right)=\frac{4\overline z}{1+|z|^2}\partial_lJ_i z-\frac{4\partial_lz}{(1+|z|^2)^2}\bigl(\overline{J_i z}+\overline z^{\,2}J_i z\bigr)-\frac{4\delta_{il}|z|^2}{1+|z|^2},
\]
\[
\begin{aligned}
  &2it\partial_i\left[-\frac{2\overline z^{\,2}}{(1+|z|^2)^2}\sum_{l=1}^2(\partial_lz)^2\right]\\
  &=-\frac{4\overline z^{\,2}}{(1+|z|^2)^2}\sum_{l=1}^2\partial_lz\,\partial_lJ_i z+\frac{4|z|^2\overline z\,\partial_i z}{(1+|z|^2)^2}+\frac{4\overline z}{(1+|z|^2)^3}\left(\sum_{l=1}^2(\partial_lz)^2\right)\bigl(\overline{J_i z}+\overline z^{\,2}J_i z\bigr),
\end{aligned}
\]
and
\[
\begin{aligned}
  &2it\partial_i\left[-\frac{2}{(1+|z|^2)^2}\sum_{l=1}^2(\partial_lz)^2\right]-\frac{4x_i}{(1+|z|^2)^2}\sum_{l=1}^2(\partial_lz)^2\\
  &=-\frac{4}{(1+|z|^2)^2}\sum_{l=1}^2\partial_lz\,\partial_lJ_i z+\frac{4z\partial_i z}{(1+|z|^2)^2}+\frac{4}{(1+|z|^2)^3}\left(\sum_{l=1}^2(\partial_lz)^2\right)\bigl(\overline z\,J_i z-z\overline{J_i z}\bigr).
\end{aligned}
\]
The last source satisfies
\[
  J_i\left[-\frac{4|z|^2}{1+|z|^2}\partial_jz\right]=-\frac{4|z|^2}{1+|z|^2}\partial_jJ_i z+\frac{4\delta_{ij}|z|^2z}{1+|z|^2}-\frac{4\partial_jz}{(1+|z|^2)^2}\bigl(\overline z\,J_i z-z\overline{J_i z}\bigr).
\]
Combining these identities cancels every explicit $x_i$ term and gives
\begin{equation}\label{eJ2z1}
\begin{aligned}
  (i\partial_t+\Delta)J_iJ_jz=&\frac{4\overline z}{1+|z|^2}\sum_{l=1}^2\partial_lz\,\partial_lJ_iJ_jz-\frac{2\overline z^{\,2}}{(1+|z|^2)^2}\left(\sum_{l=1}^2(\partial_lz)^2\right)J_iJ_jz\\
  &+\frac{2}{(1+|z|^2)^2}\left(\sum_{l=1}^2(\partial_lz)^2\right)\overline{J_iJ_jz}\\
  &+O\Bigl(|z||\nabla Jz|^2+|Jz||\nabla z||\nabla Jz|+|z||Jz|^2|\nabla z|^2\\
  &\hspace{13mm}+|z|^2|\nabla Jz|+|z||Jz||\nabla z|+|z|^3\Bigr).
\end{aligned}
\end{equation}

We use the energy with weight $(1+|z|^2)^{-2}$. First,
\begin{equation*}
\nabla(1+|z|^2)^{-2}+\frac{4\operatorname{Re}(\overline z\nabla z)}{(1+|z|^2)^3}=0.
\end{equation*}
By \eqref{ereducedequation},
\begin{equation*}
\partial_t|z|^2=-2\operatorname{Im}(\overline z\Delta z)+\frac{4}{1+|z|^2}\operatorname{Im}\left[\overline z^{\,2}\sum_{l=1}^2(\partial_lz)^2\right].
\end{equation*}
Also,
\begin{equation*}
\nabla|z|^2\cdot\operatorname{Im}(\overline z\nabla z)=\operatorname{Im}\left[\overline z^{\,2}\sum_{l=1}^2(\partial_lz)^2\right].
\end{equation*}
Therefore,
\begin{align*}
\partial_t(1+|z|^2)^{-2}&=\frac{4\operatorname{Im}(\overline z\Delta z)}{(1+|z|^2)^3}-\frac{8}{(1+|z|^2)^4}\operatorname{Im}\left[\overline z^{\,2}\sum_{l=1}^2(\partial_lz)^2\right],\\
\nabla\cdot\left[\frac{4\operatorname{Im}(\overline z\nabla z)}{(1+|z|^2)^3}\right]&=\frac{4\operatorname{Im}(\overline z\Delta z)}{(1+|z|^2)^3}-\frac{12}{(1+|z|^2)^4}\operatorname{Im}\left[\overline z^{\,2}\sum_{l=1}^2(\partial_lz)^2\right].
\end{align*}
Subtracting gives
\begin{equation*}
\partial_t(1+|z|^2)^{-2}-\nabla\cdot\left[\frac{4\operatorname{Im}(\overline z\nabla z)}{(1+|z|^2)^3}\right]-\frac{4}{(1+|z|^2)^4}\operatorname{Im}\left[\overline z^{\,2}\sum_{l=1}^2(\partial_lz)^2\right]=0.
\end{equation*}
Multiplying \eqref{eJ2z1} by $(1+|z|^2)^{-2}\overline{J_iJ_jz}$, taking imaginary parts, and integrating by parts, these identities yield
\begin{equation*}
\begin{aligned}
\frac{\mathrm d}{\mathrm dt}\int_{\RR^2}\frac{|J_iJ_jz|^2}{(1+|z|^2)^2}\,\mathrm dx&=4\operatorname{Im}\int_{\RR^2}\frac{\left(\sum_{l=1}^2(\partial_lz)^2\right)\overline{J_iJ_jz}^{\,2}}{(1+|z|^2)^4}\,\mathrm dx\\
&\quad+O\left(\int_{\RR^2}\frac{|J_iJ_jz|}{(1+|z|^2)^2}\left(|z||\nabla Jz|^2+|Jz||\nabla z||\nabla Jz|+|z||Jz|^2|\nabla z|^2\right.\right.\\
&\hspace{55mm}\left.\left.+|z|^2|\nabla Jz|+|z||Jz||\nabla z|+|z|^3\right)\,\mathrm dx\right).
\end{aligned}
\end{equation*}
Since $\|z(t)\|_{L_x^\infty}\lesssim\epsilon_0$, the weighted energy is equivalent to $\|J_iJ_jz(t)\|_{L_x^2}^2$. Cauchy-Schwarz gives
\begin{equation}\label{eJ2z2}
\begin{aligned}
\frac{\mathrm d}{\mathrm dt}\left(\int_{\RR^2}\frac{|J_iJ_jz|^2}{(1+|z|^2)^2}\,\mathrm dx\right)^{\frac{1}{2}}&\lesssim\|\nabla z\|_{L_x^\infty}^2\left(\int_{\RR^2}\frac{|J_iJ_jz|^2}{(1+|z|^2)^2}\,\mathrm dx\right)^{\frac{1}{2}}\\
&\quad+\bigl\||z||\nabla Jz|^2\bigr\|_{L_x^2}+\bigl\||Jz||\nabla z||\nabla Jz|\bigr\|_{L_x^2}+\bigl\||z||Jz|^2|\nabla z|^2\bigr\|_{L_x^2}\\
&\quad+\bigl\||z|^2|\nabla Jz|\bigr\|_{L_x^2}+\bigl\||z||Jz||\nabla z|\bigr\|_{L_x^2}+\bigl\||z|^3\bigr\|_{L_x^2}.
\end{aligned}
\end{equation}
We first apply this estimate on $[0,1]$. By \eqref{eJz3},
\begin{equation*}
\|Jz(t)\|_{L_x^4}\lesssim\|Jz(t)\|_{L_x^2}^{\frac{1}{2}}\|\nabla Jz(t)\|_{L_x^2}^{\frac{1}{2}}\lesssim\epsilon_0.
\end{equation*}
By Gagliardo-Nirenberg inequality and \eqref{elocalpsibound}, we obtain
\begin{equation*}
\|J\psi(t)\|_{L_x^4}^2\lesssim\|\psi(t)\|_{L_x^\infty}\|J^2\psi(t)\|_{L_x^2}\lesssim\epsilon_0\|\psi(t)\|_{L_x^\infty}.
\end{equation*}
Also,
\begin{equation*}
\|\psi(t)\|_{L_x^4}^2\leq\|\psi(t)\|_{L_x^2}\|\psi(t)\|_{L_x^\infty}\lesssim\epsilon_0\|\psi(t)\|_{L_x^\infty}.
\end{equation*}
Thus \eqref{eJz2}, \eqref{ethmpf11}, and \eqref{sob} give
\begin{equation*}
\|\nabla Jz(t)\|_{L_x^4}^2\lesssim\|J\psi(t)\|_{L_x^4}^2+\epsilon_0^4\|\psi(t)\|_{L_x^4}^2+\|z(t)\|_{L_x^4}^2\lesssim\epsilon_0\|\psi(t)\|_{L_x^\infty}+\epsilon_0^2,\qquad 0<t\leq1.
\end{equation*}
The six remainder terms are integrable on $[0,1]$. More precisely, \eqref{elocalpsibound} and the estimates above imply
\begin{align*}
\int_0^1\bigl\||z||\nabla Jz|^2\bigr\|_{L_x^2}\,\mathrm dt&\leq\int_0^1\|z\|_{L_x^\infty}\|\nabla Jz\|_{L_x^4}^2\,\mathrm dt\lesssim\epsilon_0^2\int_0^1\|\psi(t)\|_{L_x^\infty}\,\mathrm dt+\epsilon_0^3\lesssim\epsilon_0^3,\\
\int_0^1\bigl\||Jz||\nabla z||\nabla Jz|\bigr\|_{L_x^2}\,\mathrm dt&\leq\int_0^1\|Jz\|_{L_x^4}\|\nabla z\|_{L_x^\infty}\|\nabla Jz\|_{L_x^4}\,\mathrm dt\\
&\lesssim\epsilon_0^{\frac{3}{2}}\int_0^1\|\psi(t)\|_{L_x^\infty}^{\frac{3}{2}}\,\mathrm dt+\epsilon_0^2\int_0^1\|\psi(t)\|_{L_x^\infty}\,\mathrm dt\lesssim\epsilon_0^3,\\
\int_0^1\bigl\||z||Jz|^2|\nabla z|^2\bigr\|_{L_x^2}\,\mathrm dt&\leq\int_0^1\|z\|_{L_x^\infty}\|Jz\|_{L_x^4}^2\|\nabla z\|_{L_x^\infty}^2\,\mathrm dt\lesssim\epsilon_0^3\int_0^1\|\psi(t)\|_{L_x^\infty}^2\,\mathrm dt\lesssim\epsilon_0^5,\\
\int_0^1\bigl\||z|^2|\nabla Jz|\bigr\|_{L_x^2}\,\mathrm dt&\leq\int_0^1\|z\|_{L_x^\infty}^2\|\nabla Jz\|_{L_x^2}\,\mathrm dt\lesssim\epsilon_0^3,\\
\int_0^1\bigl\||z||Jz||\nabla z|\bigr\|_{L_x^2}\,\mathrm dt&\leq\int_0^1\|z\|_{L_x^\infty}\|Jz\|_{L_x^2}\|\nabla z\|_{L_x^\infty}\,\mathrm dt\lesssim\epsilon_0^2\int_0^1\|\psi(t)\|_{L_x^\infty}\,\mathrm dt\lesssim\epsilon_0^3,\\
\int_0^1\bigl\||z|^3\bigr\|_{L_x^2}\,\mathrm dt&\leq\int_0^1\|z\|_{L_x^\infty}^2\|z\|_{L_x^2}\,\mathrm dt\lesssim\epsilon_0^3.
\end{align*}
Since $\int_0^1\|\nabla z(t)\|_{L_x^\infty}^2\,\mathrm dt\lesssim\epsilon_0^2$, \eqref{eJ2z2}, Gronwall inequality, and \eqref{ezinitialbounds} yield
\begin{equation*}
\sup_{0\leq t\leq1}\|J^2z(t)\|_{L_x^2}\lesssim\left(\|x^2z(0)\|_{L_x^2}+\epsilon_0^3\right)\exp\left(C\int_0^1\|\psi(t)\|_{L_x^\infty}^2\,\mathrm dt\right)\lesssim\epsilon_0.
\end{equation*}
We next consider $t\geq1$. By \eqref{bsc10}, \eqref{eJz2}, \eqref{ethmpf10}, and \eqref{eJz6},
\begin{equation*}
\|\nabla Jz(t)\|_{L_x^4}\lesssim\|J\psi(t)\|_{L_x^4}+t\epsilon_0^2t^{-\frac{3}{2}}\|\psi(t)\|_{L_x^4}+\|z(t)\|_{L_x^4}\lesssim\epsilon_0t^{-\frac{1}{2}+\frac{\beta}{2}}.
\end{equation*}
Also, \eqref{eJz4} and \eqref{eJz5} give
\begin{equation*}
\|Jz(t)\|_{L_x^4}\lesssim\|Jz(t)\|_{L_x^2}^{\frac{1}{2}}\|\nabla Jz(t)\|_{L_x^2}^{\frac{1}{2}}\lesssim\epsilon_0t^{\frac{\gamma}{2}}.
\end{equation*}
Hence, the six terms in the remainder satisfy
\begin{equation*}
\begin{aligned}
\bigl\||z||\nabla Jz|^2\bigr\|_{L_x^2}&\leq\|z\|_{L_x^\infty}\|\nabla Jz\|_{L_x^4}^2\lesssim\epsilon_0^3t^{-\frac{7}{4}+\beta},\\
\bigl\||Jz||\nabla z||\nabla Jz|\bigr\|_{L_x^2}&\leq\|Jz\|_{L_x^4}\|\nabla z\|_{L_x^\infty}\|\nabla Jz\|_{L_x^4}\lesssim\epsilon_0^3t^{-\frac{3}{2}+\frac{\beta+\gamma}{2}},\\
\bigl\||z||Jz|^2|\nabla z|^2\bigr\|_{L_x^2}&\leq\|z\|_{L_x^\infty}\|Jz\|_{L_x^4}^2\|\nabla z\|_{L_x^\infty}^2\lesssim\epsilon_0^5t^{-\frac{11}{4}+\gamma},\\
\bigl\||z|^2|\nabla Jz|\bigr\|_{L_x^2}&\leq\|z\|_{L_x^\infty}^2\|\nabla Jz\|_{L_x^2}\lesssim\epsilon_0^3t^{-\frac{3}{2}+\gamma},\\
\bigl\||z||Jz||\nabla z|\bigr\|_{L_x^2}&\leq\|z\|_{L_x^\infty}\|Jz\|_{L_x^2}\|\nabla z\|_{L_x^\infty}\lesssim\epsilon_0^3t^{-\frac{7}{4}},\\
\bigl\||z|^3\bigr\|_{L_x^2}&\leq\|z\|_{L_x^\infty}^2\|z\|_{L_x^2}\lesssim\epsilon_0^3t^{-\frac{3}{2}}.
\end{aligned}
\end{equation*}
Their sum is bounded by $C\epsilon_0^3t^{-\frac{3}{2}+\frac{\beta+\gamma}{2}}$. Inserting these bounds into \eqref{eJ2z2}, we obtain
\begin{equation*}
\frac{\mathrm d}{\mathrm dt}\left(\int_{\RR^2}\frac{|J_iJ_jz|^2}{(1+|z|^2)^2}\,\mathrm dx\right)^{\frac{1}{2}}\lesssim\epsilon_0^2t^{-2}\left(\int_{\RR^2}\frac{|J_iJ_jz|^2}{(1+|z|^2)^2}\,\mathrm dx\right)^{\frac{1}{2}}\quad+\epsilon_0^3t^{-\frac{3}{2}+\frac{\beta+\gamma}{2}},\qquad t\geq1.
\end{equation*}
Gronwall inequality yields
\begin{equation*}
\sup_{t\geq1}\|J^2z(t)\|_{L_x^2}\lesssim\left(\|J^2z(1)\|_{L_x^2}+\epsilon_0^3\int_1^\infty t^{-\frac{3}{2}+\frac{\beta+\gamma}{2}}\,\mathrm dt\right)\exp\left(C\epsilon_0^2\int_1^\infty t^{-2}\,\mathrm dt\right)\lesssim\epsilon_0.
\end{equation*}
Combining the two time intervals proves
\begin{equation*}
\sup_{t\geq0}\|J^2z(t)\|_{L_x^2}\lesssim\|x^2z(0)\|_{L_x^2}+\epsilon_0^3\lesssim\epsilon_0.
\end{equation*}
\end{proof}

\subsection{Optimal decay of the map near $Q$}

By \eqref{eJz}, \eqref{eJ2z}, and \eqref{sob},
\[
  \|e^{-it\Delta}z(t)\|_{L_x^1}\lesssim\|e^{-it\Delta}z(t)\|_{L_x^2}^{\frac{1}{2}}\||x|^2e^{-it\Delta}z(t)\|_{L_x^2}^{\frac{1}{2}}\lesssim\epsilon_0.
\]
The dispersive estimate, \eqref{ezphilinfty}, and \eqref{ezprebound} imply
\begin{equation}\label{ezdecay1}
  \|\phi(t)-Q\|_{L_x^\infty}\lesssim\|z(t)\|_{L_x^\infty}\lesssim\epsilon_0t^{-1}, \qquad t\geq 1.
\end{equation}
This proves \eqref{th1.1decay2}.

\subsection{Quantitative scattering of the complex coordinate}

We first establish a weighted estimate for the coordinate profile.

\begin{lemma}\label{lzweight}
Let $\sigma$ satisfy
\begin{equation*}
0<\sigma-1<\min\left\{\frac{\beta}{2},\frac{\beta-2\gamma}{2(\beta-\gamma)}\right\}.
\end{equation*}
Then, for every $t\geq1$,
\begin{equation}\label{ethmpf12}
\|\langle x\rangle^\sigma e^{-it\Delta}z(t)\|_{L_x^2}+t^{-\frac{\beta}{2}}\|\langle x\rangle^\sigma\nabla e^{-it\Delta}z(t)\|_{L_x^2}\lesssim\epsilon_0.
\end{equation}
\end{lemma}
\begin{proof}
The range of $\sigma$ is nonempty because $2\gamma<\beta$. Since $\sigma<2$, \eqref{eJ2z} and \eqref{sob} give
\begin{equation*}
\|\langle x\rangle^\sigma e^{-it\Delta}z(t)\|_{L_x^2}\lesssim\|z(t)\|_{L_x^2}+\||x|^2e^{-it\Delta}z(t)\|_{L_x^2}\lesssim\epsilon_0.
\end{equation*}
By \eqref{ethmpf2}, \eqref{sob}, and the Moser estimates,
\begin{equation*}
\left\|\frac{(v_\infty+iw_\infty-zQ)\cdot v}{1+Q\cdot\phi}-\frac{1}{2}\right\|_{H^{1+\epsilon}}\lesssim\epsilon_0.
\end{equation*}
The free kernel formula and Plancherel identity imply
\begin{equation}\label{kform}
\||x|^\sigma e^{-it\Delta}h\|_{L_x^2}=(2t)^\sigma\bigl\||\nabla|^\sigma(e^{-\frac{i|x|^2}{4t}}h)\bigr\|_{L_x^2}.
\end{equation}
Applying this identity to \eqref{ethmpf8}, followed by fractional Leibniz inequality, gives
\begin{equation*}
\begin{aligned}
\||x|^\sigma\nabla e^{-it\Delta}z(t)\|_{L_x^2}&\lesssim\||x|^\sigma f(t)\|_{L_x^2}+t^\sigma\left\||\nabla|^\sigma\left[\frac{(v_\infty+iw_\infty-zQ)\cdot v}{1+Q\cdot\phi}-\frac{1}{2}\right]\right\|_{L_x^2}\|\psi(t)\|_{L_x^\infty}\\
&\lesssim\|xf(t)\|_{L_x^2}^{2-\sigma}\||x|^2f(t)\|_{L_x^2}^{\sigma-1}+\epsilon_0^2t^{\sigma-1}\\
&\lesssim\epsilon_0t^{(2-\sigma)\gamma+(\sigma-1)\beta}+\epsilon_0^2t^{\sigma-1}\lesssim\epsilon_0t^{\frac{\beta}{2}}.
\end{aligned}
\end{equation*}
Here $\sigma-1<\frac{\beta}{2}<\epsilon$ allows us to use the $H^{1+\epsilon}$ bound for the coefficient. The two powers of $t$ satisfy
\begin{equation*}
\sigma-1<\frac{\beta}{2},\qquad(2-\sigma)\gamma+(\sigma-1)\beta=\gamma+(\sigma-1)(\beta-\gamma)<\frac{\beta}{2}.
\end{equation*}
Finally,
\begin{equation*}
\|\langle x\rangle^\sigma\nabla e^{-it\Delta}z(t)\|_{L_x^2}\lesssim\|\nabla z(t)\|_{L_x^2}+\||x|^\sigma\nabla e^{-it\Delta}z(t)\|_{L_x^2}\lesssim\epsilon_0t^{\frac{\beta}{2}}.
\end{equation*}
This proves \eqref{ethmpf12}.
\end{proof}

By the fractional Leibniz and Moser estimates, \eqref{ezdecay1}, \eqref{ezgradient}, and \eqref{sob},
\begin{equation}\label{ethmpf13}
\begin{aligned}
\left\|\frac{2\overline z}{1+|z|^2}\sum_{l=1}^2(\partial_lz)^2\right\|_{H^\epsilon}&\lesssim\|z\|_{L_x^\infty}\|\nabla z\|_{L_x^\infty}\|\nabla z\|_{H^\epsilon}+\|z\|_{H^\epsilon}\|\nabla z\|_{L_x^\infty}^2\\
&\lesssim\epsilon_0^3t^{-2},\qquad t\geq1.
\end{aligned}
\end{equation}
On $[0,1]$, the same product estimate, \eqref{ethmpf8}, and \eqref{elocalpsibound} give
\begin{equation*}
\int_0^1\left\|\frac{2\overline z}{1+|z|^2}\sum_{l=1}^2(\partial_lz)^2\right\|_{H^\epsilon}\,\mathrm{d}t\lesssim\epsilon_0^2\int_0^1\|\psi(t)\|_{L_x^\infty}\,\mathrm{d}t+\epsilon_0\int_0^1\|\psi(t)\|_{L_x^\infty}^2\,\mathrm{d}t\lesssim\epsilon_0^3.
\end{equation*}
Define
\begin{equation}\label{deg+}
\frac{1}{2}g_+=z(0)-i\int_0^\infty e^{-i\tau\Delta}\left[\frac{2\overline z}{1+|z|^2}\sum_{l=1}^2(\partial_lz)^2\right](\tau)\,\mathrm d\tau.
\end{equation}
The integral converges in $H^\epsilon$, and the Duhamel formula gives
\begin{equation*}
e^{-it\Delta}z(t)\longrightarrow\frac{1}{2}g_+\quad\text{strongly in }H^\epsilon.
\end{equation*}
On the other hand, \eqref{eJ2z} and \eqref{sob} imply
\begin{equation*}
\sup_{t\geq0}\left(\|e^{-it\Delta}z(t)\|_{H^{1+\epsilon}}+\|\langle x\rangle^2e^{-it\Delta}z(t)\|_{L_x^2}\right)\lesssim\epsilon_0.
\end{equation*}
Weak lower semicontinuity yields
\begin{equation*}
\|g_+\|_{H^{1+\epsilon}}+\|\langle x\rangle^2g_+\|_{L_x^2}\lesssim\epsilon_0,
\end{equation*}
which proves \eqref{ethm1.1-scat}.
\begin{lemma}\label{lznonlinear}
For $\sigma$ as in Lemma~\ref{lzweight} and $t\geq1$,
\begin{equation}\label{ethmpf14}
\left\|\langle x\rangle^\sigma e^{-it\Delta}\left[\frac{2\overline z}{1+|z|^2}\sum_{l=1}^2(\partial_lz)^2\right](t)\right\|_{L_x^2}\lesssim\epsilon_0^3t^{-2+\frac{\beta}{2}}.
\end{equation}
\end{lemma}
\begin{proof}
Since
\begin{equation*}
e^{-\frac{i|x|^2}{4t}}\frac{2\overline z}{1+|z|^2}\sum_{l=1}^2(\partial_lz)^2=\frac{2\overline{e^{-\frac{i|x|^2}{4t}}z}}{1+|e^{-\frac{i|x|^2}{4t}}z|^2}\sum_{l=1}^2\left(e^{-\frac{i|x|^2}{4t}}\partial_lz\right)^2,
\end{equation*}
by \eqref{kform}, \eqref{ethmpf12}, \eqref{ezdecay1}, \eqref{ezgradient}, homogeneous fractional Leibniz and Moser estimates, we have
\begin{equation*}
\begin{aligned}
\left\||x|^\sigma e^{-it\Delta}\left[\frac{2\overline z}{1+|z|^2}\sum_{l=1}^2(\partial_lz)^2\right]\right\|_{L_x^2}&\lesssim\|\nabla z\|_{L_x^\infty}^2\||x|^\sigma e^{-it\Delta}z\|_{L_x^2}+\|z\|_{L_x^\infty}\|\nabla z\|_{L_x^\infty}\||x|^\sigma\nabla e^{-it\Delta}z\|_{L_x^2}\\
&\lesssim\epsilon_0^3t^{-2}+\epsilon_0^3t^{-2+\frac{\beta}{2}}\lesssim\epsilon_0^3t^{-2+\frac{\beta}{2}}.
\end{aligned}
\end{equation*}
Combining this with \eqref{ethmpf13} proves \eqref{ethmpf14}.
\end{proof}
\begin{proposition}\label{pzscattering}
For every $t\geq1$, the state defined in \eqref{deg+} satisfies
\begin{align}
\|2z(t)-e^{it\Delta}g_+\|_{H^\epsilon}&\lesssim\epsilon_0^3t^{-1},\label{ezscattering1}\\
\|2z(t)-e^{it\Delta}g_+\|_{L_x^\infty}&\lesssim\epsilon_0^3t^{-2+\frac{\beta}{2}}.\label{ezscattering2}
\end{align}
\end{proposition}
\begin{proof}
By \eqref{deg+} and the Duhamel formula,
\begin{equation*}
2e^{-it\Delta}z(t)-g_+=2i\int_t^\infty e^{-i\tau\Delta}\left[\frac{2\overline z}{1+|z|^2}\sum_{l=1}^2(\partial_lz)^2\right](\tau)\,\mathrm d\tau.
\end{equation*}
Thus \eqref{ethmpf13} gives
\begin{equation*}
\|2z(t)-e^{it\Delta}g_+\|_{H^\epsilon}\lesssim\epsilon_0^3\int_t^\infty\tau^{-2}\,\mathrm d\tau\lesssim\epsilon_0^3t^{-1}.
\end{equation*}
Fix $\sigma$ as in Lemma~\ref{lzweight}. By \eqref{ethmpf14},
\begin{equation*}
\|\langle x\rangle^\sigma(2e^{-it\Delta}z(t)-g_+)\|_{L_x^2}\lesssim\epsilon_0^3\int_t^\infty\tau^{-2+\frac{\beta}{2}}\,\mathrm d\tau\lesssim\epsilon_0^3t^{-1+\frac{\beta}{2}}.
\end{equation*}
The integral converges in this weighted space, and its value agrees with the $L_x^2$ limit above. Since $\sigma>1$, dispersive estimate and Cauchy-Schwarz yield
\begin{equation*}
\begin{aligned}
\|2z(t)-e^{it\Delta}g_+\|_{L_x^\infty}&\lesssim t^{-1}\|2e^{-it\Delta}z(t)-g_+\|_{L_x^1}\\
&\lesssim t^{-1}\|\langle x\rangle^\sigma(2e^{-it\Delta}z(t)-g_+)\|_{L_x^2}\lesssim\epsilon_0^3t^{-2+\frac{\beta}{2}}.
\end{aligned}
\end{equation*}
This proves the proposition.
\end{proof}

\subsection{Scattering and pointwise asymptotics of the map}

We first prove \eqref{ethm1.1-scattering}.

By \eqref{einversecoor},
\begin{equation}\label{emapdec}
\begin{aligned}
  &\phi-Q-v_\infty\operatorname{Re}(e^{it\Delta}g_+)-w_\infty\operatorname{Im}(e^{it\Delta}g_+)\\
  &\quad=v_\infty\operatorname{Re}(2z-e^{it\Delta}g_+)+w_\infty\operatorname{Im}(2z-e^{it\Delta}g_+)-\frac{2|z|^2}{1+|z|^2}Q-\frac{2|z|^2}{1+|z|^2}\left(v_\infty\operatorname{Re}z+w_\infty\operatorname{Im}z\right).
\end{aligned}
\end{equation}
The fractional Leibniz and Moser estimates yield
\begin{align*}
\left\|\frac{|z|^2}{1+|z|^2}\right\|_{H^\epsilon}&\lesssim\|z\|_{L_x^\infty}\|z\|_{H^\epsilon}\lesssim\epsilon_0^2t^{-1},\\
\left\|\frac{|z|^2z}{1+|z|^2}\right\|_{H^\epsilon}&\lesssim\|z\|_{L_x^\infty}^2\|z\|_{H^\epsilon}\lesssim\epsilon_0^3t^{-2}.
\end{align*}
Together with \eqref{ezscattering1} and \eqref{emapdec}, these imply
\begin{equation*}
\|\phi(t)-Q-v_\infty\operatorname{Re}(e^{it\Delta}g_+)-w_\infty\operatorname{Im}(e^{it\Delta}g_+)\|_{H^\epsilon}\lesssim\epsilon_0^3t^{-1}+\epsilon_0^2t^{-1}+\epsilon_0^3t^{-2}\lesssim\epsilon_0^2t^{-1}.
\end{equation*}
This proves \eqref{ethm1.1-scattering}. It remains to establish \eqref{ethm1.1-uniform1} to complete the proof of Theorem~\ref{thm:main}. By \eqref{ezscattering2}, \eqref{ezdecay1}, and \eqref{emapdec},
\begin{equation*}
\begin{aligned}
&\|\phi-Q-v_\infty\operatorname{Re}(e^{it\Delta}g_+)-w_\infty\operatorname{Im}(e^{it\Delta}g_+)\|_{L_x^\infty}\\
&\qquad\lesssim\|2z-e^{it\Delta}g_+\|_{L_x^\infty}+\|z\|_{L_x^\infty}^2+\|z\|_{L_x^\infty}^3\lesssim\epsilon_0^2t^{-2}+\epsilon_0^3t^{-2+\frac{\beta}{2}}.
\end{aligned}
\end{equation*}
Since $t\geq1$ and $\epsilon_0\leq1$, the last bound is at most $C\epsilon_0^2t^{-2+\frac{\beta}{2}}$. This proves \eqref{ethm1.1-uniform1} and completes the proof of Theorem~\ref{thm:main}. We turn to Corollary~\ref{cmap asymptotics}, beginning with the normal component. \eqref{einversecoor} gives
\begin{equation*}
  Q\cdot\phi-1+\frac{1}{2}|e^{it\Delta}g_+|^2=\frac{1}{2}\bigl(|e^{it\Delta}g_+|^2-4|z|^2\bigr)+\frac{2|z|^4}{1+|z|^2}.
\end{equation*}
By \eqref{ethm1.1-scat} and dispersive estimate,
\begin{equation*}
\|e^{it\Delta}g_+\|_{L_x^\infty}\lesssim t^{-1}\|g_+\|_{L_x^1}\lesssim t^{-1}\|\langle x\rangle^2g_+\|_{L_x^2}\lesssim\epsilon_0t^{-1}.
\end{equation*}
Therefore,
\begin{equation*}
\begin{aligned}
\left\|Q\cdot\phi-1+\frac{1}{2}|e^{it\Delta}g_+|^2\right\|_{L_x^\infty}&\lesssim\|e^{it\Delta}g_+-2z\|_{L_x^\infty}\left(\|e^{it\Delta}g_+\|_{L_x^\infty}+2\|z\|_{L_x^\infty}\right)+\|z\|_{L_x^\infty}^4\\
&\lesssim\epsilon_0^4t^{-3+\frac{\beta}{2}}+\epsilon_0^4t^{-4}\lesssim\epsilon_0^4t^{-3+\frac{\beta}{2}}.
\end{aligned}
\end{equation*}
This proves \eqref{ethm1.1-uniform3}. Using \eqref{emapdec} once more, we obtain
\begin{equation*}
\begin{aligned}
&\left\|\phi-Q-v_\infty\operatorname{Re}(e^{it\Delta}g_+)-w_\infty\operatorname{Im}(e^{it\Delta}g_+)+\frac{1}{2}|e^{it\Delta}g_+|^2Q\right\|_{L_x^\infty}\\
&\qquad\lesssim\|2z-e^{it\Delta}g_+\|_{L_x^\infty}+\|z\|_{L_x^\infty}^3+\left\|Q\cdot\phi-1+\frac{1}{2}|e^{it\Delta}g_+|^2\right\|_{L_x^\infty}\\
&\qquad\lesssim\epsilon_0^3t^{-2+\frac{\beta}{2}}+\epsilon_0^3t^{-3}+\epsilon_0^4t^{-3+\frac{\beta}{2}}\lesssim\epsilon_0^3t^{-2+\frac{\beta}{2}}.
\end{aligned}
\end{equation*}
This proves \eqref{ethm1.1-uniform2}.

\subsection{Explicit stationary phase asymptotics}

It remains to prove the stationary phase formulas \eqref{ethm1.1-uniform4} and \eqref{ethm1.1-uniform5} in Corollary~\ref{cmap asymptotics}.

For $\langle x\rangle^2h\in L_x^2$, the free kernel formula and Cauchy-Schwarz give, for $t\geq1$,
\[
\begin{aligned}
  \left\|e^{it\Delta}h-\frac{e^{\frac{i|x|^2}{4t}}}{2it}\widehat h\left(\frac{x}{2t}\right)\right\|_{L_x^\infty}&\lesssim t^{-1}\left\|\left(e^{\frac{i|x|^2}{4t}}-1\right)h\right\|_{L_x^1}\\
  &\lesssim t^{-1}\||x|^2h\|_{L_x^2}\left(\int_{\mathbb R^2}\min\left\{1,\frac{|x|^2}{t}\right\}^2|x|^{-4}\,\mathrm dx\right)^{\frac{1}{2}}.
\end{aligned}
\]
The integral equals $\frac{2\pi}{t}$, so
\begin{equation}\label{ephaseremainder}
  \left\|e^{it\Delta}h-\frac{e^{\frac{i|x|^2}{4t}}}{2it}\widehat h\left(\frac{x}{2t}\right)\right\|_{L_x^\infty}\lesssim t^{-\frac{3}{2}}\||x|^2h\|_{L_x^2}.
\end{equation}
Applying \eqref{ephaseremainder} to $h=g_+$ and using \eqref{ethm1.1-scat} and \eqref{ethm1.1-uniform1}, we obtain
\begin{equation*}
\left\|\phi(t)-Q-v_\infty\operatorname{Re}\left[\frac{e^{\frac{i|x|^2}{4t}}}{2it}\widehat g_+\left(\frac{x}{2t}\right)\right]-w_\infty\operatorname{Im}\left[\frac{e^{\frac{i|x|^2}{4t}}}{2it}\widehat g_+\left(\frac{x}{2t}\right)\right]\right\|_{L_x^\infty}\lesssim\epsilon_0t^{-\frac{3}{2}},
\end{equation*}
which proves \eqref{ethm1.1-uniform4}. Moreover, $\|\widehat g_+\|_{L_\xi^\infty}
\lesssim\|g_+\|_{L_x^1}\lesssim\epsilon_0$, and
\[
\begin{aligned}
  \left\||e^{it\Delta}g_+|^2-\frac{1}{4t^2}\left|\widehat g_+\left(\frac{x}{2t}\right)\right|^2\right\|_{L_x^\infty}&\leq\left\|e^{it\Delta}g_+-\frac{e^{\frac{i|x|^2}{4t}}}{2it}\widehat g_+\left(\frac{x}{2t}\right)\right\|_{L_x^\infty}\left(\|e^{it\Delta}g_+\|_{L_x^\infty}+\frac{1}{2t}\|\widehat g_+\|_{L_\xi^\infty}\right)\\
  &\lesssim\epsilon_0^2t^{-\frac{5}{2}}.
\end{aligned}
\]
Together with \eqref{ethm1.1-uniform3}, this gives
\begin{equation*}
\left\|Q\cdot\phi(t)-1+\frac{1}{8t^2}\left|\widehat g_+\left(\frac{x}{2t}\right)\right|^2\right\|_{L_x^\infty}\lesssim\epsilon_0^2t^{-\frac{5}{2}}+\epsilon_0^4t^{-3+\frac{\beta}{2}}\lesssim\epsilon_0^2t^{-\frac{5}{2}}.
\end{equation*}
This proves \eqref{ethm1.1-uniform5} and completes the proof of Corollary~\ref{cmap asymptotics}.

Finally, orthonormality of $(v_\infty,w_\infty)$ and \eqref{ethm1.1-uniform4} imply
\[
  \lim_{t\to\infty}t\|\phi(t)-Q\|_{L_x^\infty}=\frac{1}{2}\|\widehat g_+\|_{L_\xi^\infty}.
\]
By \eqref{emass-conservation} and \eqref{ethm1.1-scattering}, $\|g_+\|_{L_x^2}=\|\phi_0-Q\|_{L_x^2}$. Thus the rate in \eqref{th1.1decay2} is sharp for nonconstant initial data.

\subsection{Proof of Corollary~\ref{cpsi scattering}} \label{subsec:psi-scatter}

We first prove \eqref{epsiscat1}.

Extend the coordinate along the heat flow by
\[
  z(s,t)=\frac{(v_\infty+iw_\infty)\cdot\tilde\phi(s,t)}{1+Q\cdot\tilde\phi(s,t)}.
\]
The algebraic identity \eqref{ezphilinfty} applies with $\phi$ replaced by $\tilde\phi(s,t)$. Thus \eqref{ethmpf4} gives
\[
  \|z(s,t)\|_{L_x^\infty}\lesssim\epsilon_0(t+s)^{-\frac{1}{2}}.
\]
The computations in \eqref{ethmpf8} and \eqref{ethmpf9} also apply to the heat direction. Since $A_0=0$,
\[
\begin{aligned}
  \partial_sz&=\frac{(v_\infty+iw_\infty-zQ)\cdot v}{1+Q\cdot\tilde\phi}\,\psi_0,\\
  \partial_s\left[\frac{(v_\infty+iw_\infty-zQ)\cdot v}{1+Q\cdot\tilde\phi}\right]&=\frac{2\overline z\,\partial_sz}{1+|z|^2}\frac{(v_\infty+iw_\infty-zQ)\cdot v}{1+Q\cdot\tilde\phi}.
\end{aligned}
\]
The coefficient has modulus $\frac{1+|z|^2}{2}$ and, by Proposition~\ref{pcaloric gauge}, tends to $\frac{1}{2}$ as $s\to\infty$. Integrating the second identity and using \eqref{ethmpf3}, we obtain
\begin{equation}\label{ethmpf15}
  \left\|\frac{(v_\infty+iw_\infty-z(s,t)Q)\cdot v(s,t)}{1+Q\cdot\tilde\phi(s,t)}-\frac{1}{2}\right\|_{L_x^\infty}\lesssim\sup_{r\geq s}\|z(r,t)\|_{L_x^\infty}\int_s^\infty\|\psi_0(r,t)\|_{L_x^\infty}\,\mathrm{d}r\lesssim\epsilon_0^2(t+s)^{-1}.
\end{equation}

All fields below are evaluated at $s=0$, and $t\geq1$. By fractional Leibniz rule, \eqref{ethmpf8}, \eqref{ethmpf10}, and \eqref{ethmpf15},
\begin{equation*}
\begin{aligned}
\|\psi(t)-2\nabla z(t)\|_{H^\epsilon}&\lesssim\left(\left\|\frac{(v_\infty+iw_\infty-zQ)\cdot v}{1+Q\cdot\phi}-\frac{1}{2}\right\|_{L_x^\infty}+\left\|\nabla\left[\frac{(v_\infty+iw_\infty-zQ)\cdot v}{1+Q\cdot\phi}\right]\right\|_{L_x^\infty}\right)\|\psi(t)\|_{H^\epsilon}\\
&\lesssim\left(\epsilon_0^2t^{-1}+\epsilon_0^2t^{-\frac{3}{2}}\right)\|\psi(t)\|_{H^\epsilon}\lesssim\epsilon_0^3t^{-1}.
\end{aligned}
\end{equation*}
Together with \eqref{ezscattering1}, this yields
\begin{equation*}
\begin{aligned}
\|\psi(t)-e^{it\Delta}\nabla g_+\|_{H^{\epsilon-1}}&\leq\|\psi(t)-2\nabla z(t)\|_{H^{\epsilon-1}}+\|\nabla(2z(t)-e^{it\Delta}g_+)\|_{H^{\epsilon-1}}\\
&\lesssim\|\psi(t)-2\nabla z(t)\|_{H^\epsilon}+\|2z(t)-e^{it\Delta}g_+\|_{H^\epsilon}\lesssim\epsilon_0^3t^{-1}.
\end{aligned}
\end{equation*}
This proves \eqref{epsiscat1}. Moreover, \eqref{ethm1.1-scat} and \eqref{etimebootstrapimprovement} imply
\begin{equation*}
\|f(t)-\nabla g_+\|_{H^\epsilon}\leq\|\psi(t)\|_{H^\epsilon}+\|g_+\|_{H^{1+\epsilon}}\lesssim\epsilon_0.
\end{equation*}
Interpolation therefore gives
\begin{equation}\label{ethmpf16}
\|f(t)-\nabla g_+\|_{L_x^2}\leq\|f(t)-\nabla g_+\|_{H^{\epsilon-1}}^\epsilon\|f(t)-\nabla g_+\|_{H^\epsilon}^{1-\epsilon}\lesssim\epsilon_0^{1+2\epsilon}t^{-\epsilon}.
\end{equation}
We next prove \eqref{epsiscat2}. For $t\geq1$ and $j\geq0$, \eqref{ethmpf16} and the second weighted profile estimate give
\begin{align*}
\|f(2^{j+1}t)-f(2^jt)\|_{L_x^2}&\lesssim\epsilon_0^{1+2\epsilon}(2^jt)^{-\epsilon},\\
\||x|^2(f(2^{j+1}t)-f(2^jt))\|_{L_x^2}&\lesssim\epsilon_0(2^jt)^\beta.
\end{align*}
Using $\|h\|_{L_x^1}\lesssim\|h\|_{L_x^2}^{\frac{1}{2}}\||x|^2h\|_{L_x^2}^{\frac{1}{2}}$, we obtain
\begin{equation*}
\|f(2^{j+1}t)-f(2^jt)\|_{L_x^1}\lesssim\epsilon_0^{1+\epsilon}(2^jt)^{-\frac{\epsilon-\beta}{2}}.
\end{equation*}
Since $\beta<\epsilon$, summation over $j\geq0$ converges. The $L_x^1$ limit agrees with $\nabla g_+$ by \eqref{ethmpf16}. Hence
\begin{equation}\label{ethmpf17}
\|f(t)-\nabla g_+\|_{L_x^1}\lesssim\epsilon_0^{1+\epsilon}\sum_{j\geq0}(2^jt)^{-\frac{\epsilon-\beta}{2}}\lesssim\epsilon_0^{1+\epsilon}t^{-\frac{\epsilon-\beta}{2}},\qquad t\geq1.
\end{equation}
Since $f(t)\in L_x^1$, this also proves $\nabla g_+\in L_x^1$. The dispersive estimate gives
\begin{equation*}
\|\psi(t)-e^{it\Delta}\nabla g_+\|_{L_x^\infty}\lesssim t^{-1}\|f(t)-\nabla g_+\|_{L_x^1}\lesssim\epsilon_0^{1+\epsilon}t^{-1-\frac{\epsilon-\beta}{2}},
\end{equation*}
which proves \eqref{epsiscat2}. Finally, applying \eqref{ephaseremainder} directly to $h=f(t)$ and using \eqref{ethmpf17}, we obtain
\begin{equation*}
\begin{aligned}
\left\|\psi(t,x)-\frac{e^{\frac{i|x|^2}{4t}}}{2it}\widehat{\nabla g_+}\left(\frac{x}{2t}\right)\right\|_{L_x^\infty}&\leq\left\|e^{it\Delta}f(t)-\frac{e^{\frac{i|x|^2}{4t}}}{2it}\widehat f\left(t,\frac{x}{2t}\right)\right\|_{L_x^\infty}+\frac{1}{2t}\|\widehat f(t)-\widehat{\nabla g_+}\|_{L_\xi^\infty}\\
&\lesssim t^{-\frac{3}{2}}\||x|^2f(t)\|_{L_x^2}+t^{-1}\|f(t)-\nabla g_+\|_{L_x^1}\\
&\lesssim\epsilon_0t^{-\frac{3}{2}+\beta}+\epsilon_0^{1+\epsilon}t^{-1-\frac{\epsilon-\beta}{2}}\lesssim\epsilon_0t^{-1-\frac{\epsilon-\beta}{2}},
\end{aligned}
\end{equation*}
which proves \eqref{epsiscat3}.

Finally, since $|\nabla\phi|=|\psi|$, \eqref{epsiscat3} gives
\[
\lim_{t\to\infty}t\|\nabla\phi(t)\|_{L_x^\infty}=\frac{1}{2}\|\widehat{\nabla g_+}\|_{L_\xi^\infty}.
\]
For nonconstant initial data, $g_+$ is a nonzero $L^2$ function, so $\nabla g_+$ cannot vanish identically. The limit is therefore positive, proving that the rate in \eqref{ethm1.1-decay1} is sharp.

\appendix

\section{Higher order connection estimates}\label{app:A}

In this appendix, we prove Proposition~\ref{phigher order connection} and Lemma~\ref{lA3higherorder}. We first consider the spatial connection coefficients, with $m\in\{1,2\}$. Recall from Definition~\ref{def:connectiondecom} that
\begin{align}
A_m^{(\geq4)}(0,t)&=-\int_0^\infty\operatorname{Im}\left[\overline{\psi_m(s)-e^{s\Delta}\psi_m(0)}\,\partial_l\psi_l(s)\right]\,\mathrm{d}s\notag\\
&\quad-\int_0^\infty\operatorname{Im}\left[\overline{e^{s\Delta}\psi_m(0)}\,\partial_l\left(\psi_l(s)-e^{s\Delta}\psi_l(0)\right)\right]\,\mathrm{d}s\notag\\
&\quad-\int_0^\infty\operatorname{Re}(\overline{\psi_m}A_l\psi_l)(s)\,\mathrm{d}s.\label{eappendixA1}
\end{align}
We shall use both \eqref{ecaloric4} and its nondivergence form
\begin{equation}\label{eappendixA2}
(\partial_s-\Delta)\psi_m=2iA_l\partial_l\psi_m+i(\partial_lA_l)\psi_m-A_lA_l\psi_m+i\operatorname{Im}(\psi_m\overline{\psi_l})\psi_l.
\end{equation}

\subsection{Auxiliary heat estimates}
Throughout the proof of Proposition~\ref{phigher order connection}, we fix $t\geq1$ and omit $t$ in the heat flow expressions.
\begin{lemma}\label{lA1}
Under \eqref{booassum}, one has
\begin{equation}\label{eappendixA3}
\sup_{s>0}\left[\|\psi(s)-e^{s\Delta}\psi(0)\|_{L_x^2}+s^{\frac{1}{2}}\|\nabla(\psi(s)-e^{s\Delta}\psi(0))\|_{L_x^2}\right]\lesssim\epsilon_1^3t^{-\frac{3}{2}},
\end{equation}
\begin{equation}\label{eappendixA4}
\|\psi(s)-e^{s\Delta}\psi(0)\|_{L_x^\infty}\lesssim\epsilon_1^3\min\left\{t^{-\frac{5}{2}}s^{\frac{1}{2}},t^{-2}\right\},
\end{equation}
\begin{equation}\label{eappendixA5}
\left(\sum_k\left[2^{2k}\int_0^\infty\|P_k(\psi(s)-e^{s\Delta}\psi(0))\|_{L_x^2}\,\mathrm{d}s\right]^2\right)^{\frac{1}{2}}\lesssim\epsilon_1^3t^{-1}.
\end{equation}
Moreover, for $a=0,1,2,3$,
\begin{equation}\label{eappendixA6}
\int_0^\infty s^{-\frac{3-a}{2}}\|\nabla^a(\psi(s)-e^{s\Delta}\psi(0))\|_{L_x^\infty}\,\mathrm{d}s\lesssim\epsilon_1^3t^{-\frac{5}{2}}(1+\log t).
\end{equation}
\end{lemma}
\begin{proof}
We first establish the $L_x^2$ bounds. For $F$ as in \eqref{eFdefinition}, we claim that
\begin{equation}\label{eappendixA7}
\max_F\|F(s)\|_{L_x^2}\lesssim\epsilon_1^2(t+s)^{-1}\langle s\rangle^{-\frac{1}{2}}.
\end{equation}
For $F=\nabla\cdot A$ and $s\leq1$, this follows from \eqref{eAonederivativel21} with $\sigma=0$. For $s\geq1$, differentiating \eqref{ecaloricconnection}, using \eqref{epsiAdecay1} with $a=1,2$, \eqref{epsidecay1}, \eqref{eAinfty1}, and \eqref{eAonederivativel21} with $\sigma=0$, together with Gagliardo-Nirenberg inequality, we obtain
\begin{align*}
\|\nabla A(s)\|_{L_x^2}&\lesssim\int_s^\infty\Bigl[\|\psi(r)\|_{L_x^\infty}\|\nabla^2\psi(r)\|_{L_x^2}+\|\nabla A(r)\|_{L_x^2}\|\psi(r)\|_{L_x^\infty}^2\\
&\qquad\qquad+\|A(r)\|_{L_x^\infty}\|\psi(r)\|_{L_x^\infty}\|\nabla\psi(r)\|_{L_x^2}\Bigr]\,\mathrm{d}r\\
&\lesssim\epsilon_1^2\int_s^\infty(t+r)^{-1}r^{-\frac{3}{2}}\,\mathrm{d}r\lesssim\epsilon_1^2(t+s)^{-1}s^{-\frac{1}{2}}.
\end{align*}
The other cases in \eqref{eappendixA7} follow from \eqref{epsiAdecay1}, \eqref{epsiAdecay2}, \eqref{epsidecay1}, and \eqref{eAinfty1}. Hence,
\begin{equation}\label{eappendixA8}
\int_0^\infty\max_F\|(F\psi)(r)\|_{L_x^2}\,\mathrm{d}r\lesssim\epsilon_1^3\int_0^\infty(t+r)^{-2}\langle r\rangle^{-\frac{1}{2}}\,\mathrm{d}r\lesssim\epsilon_1^3t^{-\frac{3}{2}}.
\end{equation}
By \eqref{epsidecay1}, \eqref{epsidecay3}, \eqref{eAinfty1}, \eqref{eAinfty3}, and Bernstein inequality, for $a=0,1,2,3$ and $r>0$,
\begin{align}
\|\nabla^a\psi(r)\|_{L_x^\infty}&\lesssim\epsilon_1r^{-\frac{a}{2}}(t+r)^{-1},\notag\\
\|\nabla^aA(r)\|_{L_x^\infty}&\lesssim\epsilon_1^2r^{-\frac{a}{2}}(t+r)^{-\frac{3}{2}}.\label{eappendixA9}
\end{align}
Combining \eqref{eappendixA2}, \eqref{eappendixA7}, \eqref{eappendixA9}, and \eqref{epsihat2}, we obtain
\begin{align}
\|(\partial_r-\Delta)\psi(r)\|_{L_x^2}&\lesssim\epsilon_1^3r^{-1}(t+r)^{-\frac{3}{2}},\notag\\
\|(\partial_r-\Delta)\psi(r)\|_{L_x^\infty}&\lesssim\epsilon_1^3r^{-\frac{1}{2}}(t+r)^{-\frac{5}{2}}.\label{eappendixA10}
\end{align}
Also, \eqref{eAinfty1} and \eqref{epsihat2} give $\|A(r)\psi(r)\|_{L_x^2}\lesssim\epsilon_1^3(t+r)^{-\frac{3}{2}}r^{-\frac{1}{2}}$. Thus \eqref{eheatduhamelpk}, \eqref{eheatoperator}, and \eqref{eappendixA8} imply
\begin{equation*}
\|\psi(s)-e^{s\Delta}\psi(0)\|_{L_x^2}\lesssim\epsilon_1^3t^{-\frac{3}{2}}\int_0^s(s-r)^{-\frac{1}{2}}r^{-\frac{1}{2}}\,\mathrm{d}r+\int_0^s\max_F\|(F\psi)(r)\|_{L_x^2}\,\mathrm{d}r\lesssim\epsilon_1^3t^{-\frac{3}{2}}.
\end{equation*}
For one derivative, use \eqref{eheatduhamelpk} on $0\leq r\leq\frac{s}{2}$ and \eqref{eappendixA2} on $\frac{s}{2}\leq r\leq s$. By \eqref{eheatoperator}, \eqref{eappendixA8}, and \eqref{eappendixA10},
\begin{align*}
\|\nabla(\psi(s)-e^{s\Delta}\psi(0))\|_{L_x^2}&\lesssim s^{-1}\int_0^{\frac{s}{2}}\|A(r)\psi(r)\|_{L_x^2}\,\mathrm{d}r+s^{-\frac{1}{2}}\int_0^{\frac{s}{2}}\max_F\|(F\psi)(r)\|_{L_x^2}\,\mathrm{d}r\\
&\quad+\epsilon_1^3\int_{\frac{s}{2}}^s(s-r)^{-\frac{1}{2}}r^{-1}(t+r)^{-\frac{3}{2}}\,\mathrm{d}r\lesssim\epsilon_1^3t^{-\frac{3}{2}}s^{-\frac{1}{2}}.
\end{align*}
This proves \eqref{eappendixA3}. By \eqref{eheatoperator} and the second estimate in \eqref{eappendixA10},
\begin{equation*}
\|\psi(s)-e^{s\Delta}\psi(0)\|_{L_x^\infty}\lesssim\epsilon_1^3\int_0^s r^{-\frac{1}{2}}(t+r)^{-\frac{5}{2}}\,\mathrm{d}r\lesssim\epsilon_1^3\min\left\{t^{-\frac{5}{2}}s^{\frac{1}{2}},t^{-2}\right\},
\end{equation*}
which proves \eqref{eappendixA4}.
For \eqref{eappendixA5}, use \eqref{eheatenergyimproved2} with $\sigma=\epsilon$ on $0<r\leq1$, \eqref{epsiAdecay1} with $a=1$ on $r\geq1$, and \eqref{eAinfty1}. Together with \eqref{eappendixA8}, these give
\begin{align}
&\int_0^\infty\left[\|A(r)\|_{L_x^\infty}\|\nabla\psi(r)\|_{L_x^2}+\max_F\|(F\psi)(r)\|_{L_x^2}\right]\,\mathrm{d}r\notag\\
&\quad\lesssim\epsilon_1^3\left[t^{-\frac{3}{2}}\int_0^1r^{-\frac{1-\epsilon}{2}}\,\mathrm{d}r+\int_1^\infty(t+r)^{-\frac{3}{2}}r^{-1}\,\mathrm{d}r+t^{-\frac{3}{2}}\right]\lesssim\epsilon_1^3t^{-\frac{3}{2}}(1+\log t).\label{eappendixA11}
\end{align}
Applying the dyadic Duhamel formula associated with \eqref{eappendixA2}, integrating the heat kernel in $s$, and using Minkowski inequality, we obtain
\begin{align*}
&\left(\sum_k\left[2^{2k}\int_0^\infty\|P_k(\psi(s)-e^{s\Delta}\psi(0))\|_{L_x^2}\,\mathrm{d}s\right]^2\right)^{\frac{1}{2}}\\
&\quad\lesssim\int_0^\infty\|(\partial_r-\Delta)\psi(r)\|_{L_x^2}\,\mathrm{d}r\lesssim\epsilon_1^3t^{-\frac{3}{2}}(1+\log t)\lesssim\epsilon_1^3t^{-1}.
\end{align*}
It remains to prove \eqref{eappendixA6}. By \eqref{eAinfty3}, \eqref{eAonederivativel2} with $\sigma=\epsilon$, \eqref{eenergyenvelope3}, and Bernstein inequality,
\begin{equation*}
\|P_kA(s)\|_{L_x^\infty}\lesssim\epsilon_1^2\min\left\{2^{2k},t^{-\frac{3}{2}},t^{-1}2^{-\epsilon k}\right\}(1+s2^{2k})^{-9}.
\end{equation*}
Splitting the sum at $2^k=t^{-\frac{3}{4}}$ and $2^k=t^{\frac{1}{2\epsilon}}$ yields
\begin{equation*}
\sum_k\sup_{s\geq0}(1+s2^{2k})^9\|P_kA(s)\|_{L_x^\infty}\lesssim\epsilon_1^2t^{-\frac{3}{2}}(1+\log t).
\end{equation*}
Together with \eqref{epsihighfre2} and $2^{bk}\int_0^\infty r^{\frac{b-2}{2}}(1+r2^{2k})^{-9}\,\mathrm{d}r\lesssim1$, this proves, for $b=1,2,3$,
\begin{align}
\int_0^\infty r^{\frac{b-2}{2}}\|\nabla^b\psi(r)\|_{L_x^\infty}\,\mathrm{d}r&\lesssim\epsilon_1t^{-1}(1+\log t),\notag\\
\int_0^\infty r^{\frac{b-2}{2}}\|\nabla^bA(r)\|_{L_x^\infty}\,\mathrm{d}r&\lesssim\epsilon_1^2t^{-\frac{3}{2}}(1+\log t).\label{eappendixA12}
\end{align}
By \eqref{eappendixA2}, \eqref{eappendixA9}, and Leibniz rule, for $b=0,1,2$,
\begin{align}
&\int_0^\infty r^{\frac{b-1}{2}}\|\nabla^b((\partial_r-\Delta)\psi(r))\|_{L_x^\infty}\,\mathrm{d}r\notag\\
&\quad\lesssim\epsilon_1^2t^{-\frac{3}{2}}\sum_{c=1}^{b+1}\int_0^\infty r^{\frac{c-2}{2}}\|\nabla^c\psi(r)\|_{L_x^\infty}\,\mathrm{d}r\notag\\
&\qquad+\epsilon_1t^{-1}\int_0^\infty r^{\frac{b-1}{2}}\|\nabla^{b+1}A(r)\|_{L_x^\infty}\,\mathrm{d}r+\epsilon_1^3\int_0^\infty r^{-\frac{1}{2}}(t+r)^{-3}\,\mathrm{d}r\notag\\
&\quad\lesssim\epsilon_1^3t^{-\frac{5}{2}}(1+\log t),\label{eappendixA13}
\end{align}
where the last inequality follows from \eqref{eappendixA12}.
For $a=1,2,3$, split the Duhamel integral at $r=\frac{s}{2}$, placing all $a$ derivatives on the heat kernel on the first interval and one derivative on the kernel on the second. By \eqref{eheatoperator}, Fubini theorem, and \eqref{eappendixA13},
\begin{align*}
&\int_0^\infty s^{-\frac{3-a}{2}}\|\nabla^a(\psi(s)-e^{s\Delta}\psi(0))\|_{L_x^\infty}\,\mathrm{d}s\\
&\quad\lesssim\int_0^\infty r^{-\frac{1}{2}}\|(\partial_r-\Delta)\psi(r)\|_{L_x^\infty}\,\mathrm{d}r+\int_0^\infty r^{\frac{a-2}{2}}\|\nabla^{a-1}((\partial_r-\Delta)\psi(r))\|_{L_x^\infty}\,\mathrm{d}r\lesssim\epsilon_1^3t^{-\frac{5}{2}}(1+\log t).
\end{align*}
For $a=0$, \eqref{eheatoperator}, Fubini theorem, and \eqref{eappendixA13} with $b=0$ give the same estimate using only the first integral.
\end{proof}

\begin{lemma}\label{lA2}
Under \eqref{booassum}, one has
\begin{equation}\label{eappendixA14}
\sup_{s\geq0}\|xA(s)\|_{L_x^\infty}\lesssim\epsilon_1^2t^{-1+\gamma}(1+\log t),
\end{equation}
\begin{equation}\label{eappendixA15}
\begin{aligned}
&\int_0^\infty\Bigl[\|xA(s)\|_{L_x^\infty}\|\nabla^2\psi(s)\|_{L_x^2}+\|\nabla A(s)\|_{L_x^\infty}\|x\nabla\psi(s)\|_{L_x^2}\\
&\qquad\qquad+\|\nabla^2A(s)\|_{L_x^2}\|x\psi(s)\|_{L_x^\infty}\Bigr]\,\mathrm{d}s\lesssim\epsilon_1^3t^{-1+\gamma}(1+\log t),
\end{aligned}
\end{equation}
\begin{equation}\label{eappendixA16}
\|x(\psi(s)-e^{s\Delta}\psi(0))\|_{L_x^\infty}\lesssim\epsilon_1^3t^{-2+\gamma}(1+\log t)s^{\frac{1}{2}},
\end{equation}
and
\begin{equation}\label{eappendixA17}
\left(\sum_k\left[2^{3k}\int_0^\infty\|P_k(x(\psi(s)-e^{s\Delta}\psi(0)))\|_{L_x^2}\,\mathrm{d}s\right]^2\right)^{\frac{1}{2}}\lesssim\epsilon_1^3t^{-1+\gamma}(1+\log t).
\end{equation}
\end{lemma}
\begin{proof}
By \eqref{eJs1}, \eqref{eJs3}, \eqref{eJimprove}, \eqref{eek2}, \eqref{epsidecay3}, and Bernstein inequality, $\|x\nabla^a\psi(s)\|_{L_x^\infty}\lesssim\epsilon_1t^\gamma s^{-\frac{a+1}{2}}$ for $a=0,1,2$. Applying \eqref{eJs2}, \eqref{eek4}, and \eqref{eUfree2} to the free heat evolution gives
\begin{equation}\label{eappendixA18}
\|x\nabla^a\psi(s)\|_{L_x^\infty}+\|x\nabla^ae^{s\Delta}\psi(0)\|_{L_x^\infty}\lesssim\epsilon_1t^\gamma s^{-\frac{a+1}{2}},\qquad a=0,1,2.
\end{equation}
In \eqref{ecaloricconnection}, place the weight on the undifferentiated $\psi$ in the quadratic term and on $A$ in the remainder. By \eqref{eappendixA18}, \eqref{eappendixA12} with $b=1$, and \eqref{epsiinte4},
\begin{align*}
\sup_{s\geq0}\|xA(s)\|_{L_x^\infty}&\lesssim\epsilon_1t^\gamma\int_0^\infty r^{-\frac{1}{2}}\|\nabla\psi(r)\|_{L_x^\infty}\,\mathrm{d}r+\epsilon_1^2t^{-1}\sup_{s\geq0}\|xA(s)\|_{L_x^\infty}\\
&\lesssim\epsilon_1^2t^{-1+\gamma}(1+\log t)+\epsilon_1^2t^{-1}\sup_{s\geq0}\|xA(s)\|_{L_x^\infty}.
\end{align*}
Absorbing the last term proves \eqref{eappendixA14}.
The first integral in \eqref{eappendixA15} follows from \eqref{eappendixA14} and \eqref{epsiinte5}. For the second, \eqref{eJs1}, \eqref{eJs3}, \eqref{eJimprove}, and \eqref{epsihat2} give
\begin{equation*}
\|x\nabla\psi(s)\|_{L_x^2}\lesssim\epsilon_1t^\gamma s^{-\frac{1}{2}}+(t+s)\|\nabla^2\psi(s)\|_{L_x^2}.
\end{equation*}
Moreover, \eqref{eAonederivativelinfty} on $s\leq t$ and \eqref{eappendixA9} on $s\geq t$ imply
\begin{equation*}
\sup_{s>0}(t+s)\|\nabla A(s)\|_{L_x^\infty}\lesssim\epsilon_1^2t^{-1}(1+\log t).
\end{equation*}
Therefore, by \eqref{eappendixA12} and \eqref{epsiinte5},
\begin{align*}
&\int_0^\infty\|\nabla A(s)\|_{L_x^\infty}\|x\nabla\psi(s)\|_{L_x^2}\,\mathrm{d}s\\
&\quad\lesssim\epsilon_1t^\gamma\int_0^\infty s^{-\frac{1}{2}}\|\nabla A(s)\|_{L_x^\infty}\,\mathrm{d}s+\epsilon_1^2t^{-1}(1+\log t)\int_0^\infty\|\nabla^2\psi(s)\|_{L_x^2}\,\mathrm{d}s\\
&\quad\lesssim\epsilon_1^3t^{-1+\gamma}(1+\log t).
\end{align*}
For the third integral, \eqref{eAonederivativel2} with $\sigma=\epsilon,0$ and \eqref{eenergyenvelope3} give
\begin{equation}\label{eappendixA19}
\|\nabla^2A(s)\|_{L_x^2}\lesssim\epsilon_1^2(t+s)^{-1}
\begin{cases}
s^{-\frac{1-\epsilon}{2}},&0<s\leq1,\\
s^{-\frac{1}{2}},&s\geq1.
\end{cases}
\end{equation}
Together with \eqref{eappendixA18}, this yields
\begin{equation*}
\int_0^\infty\|\nabla^2A(s)\|_{L_x^2}\|x\psi(s)\|_{L_x^\infty}\,\mathrm{d}s\lesssim\epsilon_1^3t^\gamma\left[t^{-1}\int_0^1s^{-1+\frac{\epsilon}{2}}\,\mathrm{d}s+\int_1^\infty\frac{1}{s(t+s)}\,\mathrm{d}s\right]\lesssim\epsilon_1^3t^{-1+\gamma}(1+\log t).
\end{equation*}
This proves \eqref{eappendixA15}.
The weighted heat correction satisfies
\begin{equation}\label{eappendixA20}
(\partial_s-\Delta)\left[x_j(\psi(s)-e^{s\Delta}\psi(0))\right]=x_j(\partial_s-\Delta)\psi(s)-2\partial_j(\psi(s)-e^{s\Delta}\psi(0)).
\end{equation}
By \eqref{eappendixA2}, \eqref{eappendixA14}, \eqref{eappendixA18}, \eqref{eappendixA9}, and \eqref{eAonederivativelinfty},
\begin{equation*}
\|x(\partial_s-\Delta)\psi(s)\|_{L_x^\infty}\lesssim\epsilon_1^3t^{-2+\gamma}(1+\log t)s^{-\frac{1}{2}}.
\end{equation*}
For the commutator term, \eqref{eheatoperator}, Fubini theorem, and \eqref{eappendixA10} give
\begin{equation*}
\int_0^s\|\nabla(\psi(r)-e^{r\Delta}\psi(0))\|_{L_x^\infty}\,\mathrm{d}r\lesssim\int_0^s(s-r)^{\frac{1}{2}}\|(\partial_r-\Delta)\psi(r)\|_{L_x^\infty}\,\mathrm{d}r\lesssim\epsilon_1^3t^{-2}s^{\frac{1}{2}}.
\end{equation*}
The Duhamel formula for \eqref{eappendixA20} therefore proves \eqref{eappendixA16}.
To prove \eqref{eappendixA17}, first note that \eqref{eappendixA11}, \eqref{eAonederivativel21} with $\sigma=0$, and \eqref{epsidecay1} imply
\begin{equation*}
\int_0^\infty\left[\|A(s)\|_{L_x^\infty}\|\nabla\psi(s)\|_{L_x^2}+\|\nabla A(s)\|_{L_x^2}\|\psi(s)\|_{L_x^\infty}+\max_F\|(F\psi)(s)\|_{L_x^2}\right]\,\mathrm{d}s\lesssim\epsilon_1^3t^{-1}.
\end{equation*}
By \eqref{eappendixA18}, \eqref{epsidecay1}, \eqref{eheatenergyimproved2} with $\sigma=\epsilon$, and \eqref{epsiAdecay1} with $a=1$,
\begin{align*}
\int_0^\infty\|x\psi(s)\|_{L_x^\infty}\|\psi(s)\|_{L_x^\infty}\|\nabla\psi(s)\|_{L_x^2}\,\mathrm{d}s&\lesssim\epsilon_1^3t^\gamma\left[t^{-1}\int_0^1s^{-1+\frac{\epsilon}{2}}\,\mathrm{d}s+\int_1^\infty(t+s)^{-1}s^{-\frac{3}{2}}\,\mathrm{d}s\right]\\
&\lesssim\epsilon_1^3t^{-1+\gamma}.
\end{align*}
Similarly, \eqref{eappendixA14}, \eqref{eappendixA11}, \eqref{eAinfty1}, and \eqref{eAonederivativel21} with $\sigma=0$ give
\begin{align*}
\int_0^\infty\|xA(s)\|_{L_x^\infty}\|A(s)\|_{L_x^\infty}\|\nabla\psi(s)\|_{L_x^2}\,\mathrm{d}s&\lesssim\epsilon_1^5t^{-2+\gamma}(1+\log t),\\
\int_0^\infty\|x\psi(s)\|_{L_x^\infty}\|A(s)\|_{L_x^\infty}\|\nabla A(s)\|_{L_x^2}\,\mathrm{d}s&\lesssim\epsilon_1^5t^{-2+\gamma}.
\end{align*}
By \eqref{eappendixA2}, Leibniz rule, \eqref{eappendixA15}, and the product estimates above,
\begin{equation}\label{eappendixA21}
\int_0^\infty\|\nabla[x(\partial_s-\Delta)\psi(s)]\|_{L_x^2}\,\mathrm{d}s\lesssim\epsilon_1^3t^{-1+\gamma}(1+\log t).
\end{equation}
Applying dyadic Duhamel to \eqref{eappendixA20}, integrating the heat kernel, and using \eqref{eappendixA5} and \eqref{eappendixA21}, we obtain
\begin{align*}
&\left(\sum_k\left[2^{3k}\int_0^\infty\|P_k(x(\psi(s)-e^{s\Delta}\psi(0)))\|_{L_x^2}\,\mathrm{d}s\right]^2\right)^{\frac{1}{2}}\\
&\quad\lesssim\int_0^\infty\|\nabla[x(\partial_s-\Delta)\psi(s)]\|_{L_x^2}\,\mathrm{d}s+\left(\sum_k\left[2^{2k}\int_0^\infty\|P_k(\psi(s)-e^{s\Delta}\psi(0))\|_{L_x^2}\,\mathrm{d}s\right]^2\right)^{\frac{1}{2}}\\
&\quad\lesssim\epsilon_1^3t^{-1+\gamma}(1+\log t).
\end{align*}
This proves \eqref{eappendixA17}.
\end{proof}

\subsection{The $L_x^2$ and weighted estimates}
\begin{proof}[Proof of \eqref{eAhigh1}]
By \eqref{eappendixA1}, \eqref{eappendixA3}, \eqref{epsiinte2}, \eqref{eUfree1}, \eqref{eAinfty1}, \eqref{epsidecay1}, and \eqref{epsihat2},
\begin{align*}
\|A_m^{(\geq4)}(0,t)\|_{L_x^2}&\lesssim\sup_{s>0}\|\psi(s)-e^{s\Delta}\psi(0)\|_{L_x^2}\int_0^\infty\|\nabla\psi(s)\|_{L_x^\infty}\,\mathrm{d}s\\
&\quad+\int_0^\infty\|e^{s\Delta}\psi(0)\|_{L_x^\infty}\|\nabla(\psi(s)-e^{s\Delta}\psi(0))\|_{L_x^2}\,\mathrm{d}s\\
&\quad+\int_0^\infty\|A(s)\|_{L_x^\infty}\|\psi(s)\|_{L_x^\infty}\|\psi(s)\|_{L_x^2}\,\mathrm{d}s\\
&\lesssim\epsilon_1^4t^{-2}+\epsilon_1^4t^{-\frac{3}{2}}\int_0^\infty(t+s)^{-1}s^{-\frac{1}{2}}\,\mathrm{d}s+\epsilon_1^4\int_0^\infty(t+s)^{-\frac{5}{2}}s^{-\frac{1}{2}}\,\mathrm{d}s\lesssim\epsilon_1^4t^{-2}.
\end{align*}
\end{proof}
\begin{proof}[Proof of \eqref{eAhigh4}]
For the two quadratic corrections in \eqref{eappendixA1}, place the weight on the full and free heat evolutions, respectively. By Leibniz rule, \eqref{eappendixA18}, and \eqref{eappendixA6},
\begin{align*}
&\left\|x\nabla\int_0^\infty\overline{\psi(s)-e^{s\Delta}\psi(0)}\,\nabla\psi(s)\,\mathrm{d}s\right\|_{L_x^\infty}+\left\|x\nabla\int_0^\infty\overline{e^{s\Delta}\psi(0)}\,\nabla(\psi(s)-e^{s\Delta}\psi(0))\,\mathrm{d}s\right\|_{L_x^\infty}\\
&\quad\lesssim\epsilon_1t^\gamma\sum_{a=0}^2\int_0^\infty s^{-\frac{3-a}{2}}\|\nabla^a(\psi(s)-e^{s\Delta}\psi(0))\|_{L_x^\infty}\,\mathrm{d}s\lesssim\epsilon_1^4t^{-\frac{5}{2}+\gamma}(1+\log t).
\end{align*}
For the remainder, place the weight on an undifferentiated $\psi$ factor. By \eqref{eappendixA18}, \eqref{eAinfty1}, \eqref{epsidecay1}, and \eqref{eappendixA12} with $b=1$,
\begin{align*}
\left\|x\nabla\int_0^\infty A(s)\psi(s)\psi(s)\,\mathrm{d}s\right\|_{L_x^\infty}&\lesssim\epsilon_1t^\gamma\int_0^\infty s^{-\frac{1}{2}}\left[\|A(s)\|_{L_x^\infty}\|\nabla\psi(s)\|_{L_x^\infty}+\|\nabla A(s)\|_{L_x^\infty}\|\psi(s)\|_{L_x^\infty}\right]\,\mathrm{d}s\\
&\lesssim\epsilon_1^4t^{-\frac{5}{2}+\gamma}(1+\log t).
\end{align*}
Combining these estimates proves \eqref{eAhigh4}.
\end{proof}
\begin{proof}[Proof of \eqref{eAhigh5}]
Commuting $x$ through the spatial derivatives in the two quadratic corrections in \eqref{eappendixA1}, the principal terms are
\begin{equation*}
\nabla^2\int_0^\infty\overline{x(\psi(s)-e^{s\Delta}\psi(0))}\,\nabla\psi(s)\,\mathrm{d}s,\qquad \nabla^2\int_0^\infty\overline{e^{s\Delta}\psi(0)}\,\nabla[x(\psi(s)-e^{s\Delta}\psi(0))]\,\mathrm{d}s.
\end{equation*}
The commutator terms contain two spatial derivatives in total and no weight. We first estimate the principal terms by paraproduct decomposition. For the Low-High interaction, with the weighted correction at low frequency, \eqref{eheatenergyimproved2} with $\sigma=0$ and \eqref{eappendixA16} give
\begin{align*}
&\left\|P_k\nabla^2\int_0^\infty P_{\leq k-5}[x(\psi(s)-e^{s\Delta}\psi(0))]\,\nabla P_{\sim k}\psi(s)\,\mathrm{d}s\right\|_{L_x^2}\\
&\quad\lesssim\epsilon_1^4t^{-2+\gamma}(1+\log t)b_k\,2^{3k}\int_0^\infty s^{\frac{1}{2}}(1+s2^{2k})^{-20}\,\mathrm{d}s\lesssim\epsilon_1^4t^{-2+\gamma}(1+\log t)b_k.
\end{align*}
For the High-Low interaction, \eqref{epsidecay1} and Bernstein inequality give
\begin{align*}
&\left\|P_k\nabla^2\int_0^\infty P_{\sim k}[x(\psi(s)-e^{s\Delta}\psi(0))]\,\nabla P_{\leq k-5}\psi(s)\,\mathrm{d}s\right\|_{L_x^2}\\
&\quad\lesssim\epsilon_1t^{-1}2^{3k}\int_0^\infty\|P_{\sim k}[x(\psi(s)-e^{s\Delta}\psi(0))]\|_{L_x^2}\,\mathrm{d}s.
\end{align*}
For the High-High interaction, the corresponding bound is
\begin{equation*}
\epsilon_1t^{-1}\sum_{j\geq k-4}2^{2(k-j)}\left[2^{3j}\int_0^\infty\|P_j[x(\psi(s)-e^{s\Delta}\psi(0))]\|_{L_x^2}\,\mathrm{d}s\right].
\end{equation*}
By \eqref{eUfree1}, \eqref{eUfree2}, \eqref{eappendixA16}, and the same paraproduct decomposition,
\begin{align*}
&\left\|\nabla^2\int_0^\infty\overline{e^{s\Delta}\psi(0)}\,\nabla[x(\psi(s)-e^{s\Delta}\psi(0))]\,\mathrm{d}s\right\|_{L_x^2}\\
&\quad\lesssim\epsilon_1^4t^{-2+\gamma}(1+\log t)+\epsilon_1t^{-1}\left(\sum_k\left[2^{3k}\int_0^\infty\|P_k[x(\psi(s)-e^{s\Delta}\psi(0))]\|_{L_x^2}\,\mathrm{d}s\right]^2\right)^{\frac{1}{2}}\\
&\quad\lesssim\epsilon_1^4t^{-2+\gamma}(1+\log t),
\end{align*}
where the last inequality follows from \eqref{eappendixA17}. Together with the estimates above and \eqref{eenergyenvelope3}, this bounds both principal terms. For the commutator terms, \eqref{eappendixA4}, \eqref{eheatenergyimproved2} with $\sigma=0$, \eqref{epsidecay1}, and the corresponding free heat bounds give
\begin{align*}
&\left\|\nabla\int_0^\infty\overline{\psi(s)-e^{s\Delta}\psi(0)}\,\nabla\psi(s)\,\mathrm{d}s\right\|_{L_x^2}\\
&\quad+\left\|\nabla\int_0^\infty\overline{e^{s\Delta}\psi(0)}\,\nabla(\psi(s)-e^{s\Delta}\psi(0))\,\mathrm{d}s\right\|_{L_x^2}\\
&\quad+\left\|\nabla^2\int_0^\infty\overline{e^{s\Delta}\psi(0)}\,(\psi(s)-e^{s\Delta}\psi(0))\,\mathrm{d}s\right\|_{L_x^2}\\
&\quad\lesssim\epsilon_1^4t^{-2}+\epsilon_1t^{-1}\left(\sum_k\left[2^{2k}\int_0^\infty\|P_k(\psi(s)-e^{s\Delta}\psi(0))\|_{L_x^2}\,\mathrm{d}s\right]^2\right)^{\frac{1}{2}}\\
&\quad\lesssim\epsilon_1^4t^{-2},
\end{align*}
where we used \eqref{eappendixA5} in the last inequality.
For the remainder in \eqref{eappendixA1}, Leibniz rule and Gagliardo-Nirenberg inequality give
\begin{align*}
&\left\|x\nabla^2\int_0^\infty A(s)\psi(s)\psi(s)\,\mathrm{d}s\right\|_{L_x^2}\\
&\quad\lesssim\sup_{s\geq0}\|\psi(s)\|_{L_x^\infty}\int_0^\infty\Bigl[\|xA(s)\|_{L_x^\infty}\|\nabla^2\psi(s)\|_{L_x^2}+\|\nabla A(s)\|_{L_x^\infty}\|x\nabla\psi(s)\|_{L_x^2}\\
&\qquad\qquad+\|\nabla^2A(s)\|_{L_x^2}\|x\psi(s)\|_{L_x^\infty}\Bigr]\,\mathrm{d}s\\
&\quad\lesssim\epsilon_1^4t^{-2+\gamma}(1+\log t),
\end{align*}
where the last inequality follows from \eqref{epsidecay1} and \eqref{eappendixA15}. This proves \eqref{eAhigh5}.
\end{proof}

\subsection{The derivative $L_x^2$ estimates}
\begin{proof}[Proof of \eqref{eAhigh3} and \eqref{eAhigh2}]
We first establish
\begin{equation}\label{eappendixA22}
\left(\sum_k\left[2^{4k}\int_0^\infty\|P_k(\psi(s)-e^{s\Delta}\psi(0))\|_{L_x^2}\,\mathrm{d}s\right]^2\right)^{\frac{1}{2}}\lesssim\epsilon_1^3t^{-2+\gamma}(1+\log t).
\end{equation}
By \eqref{eJs1}, equation \eqref{eappendixA2} becomes
\begin{equation*}
(\partial_s-\Delta)\psi_m=\frac{i}{s+it}\left[A_lJ_{s,l}\psi_m-(x_lA_l)\psi_m\right]+i(\partial_lA_l)\psi_m-A_lA_l\psi_m+i\operatorname{Im}(\psi_m\overline{\psi_l})\psi_l.
\end{equation*}
By the dyadic Duhamel formula and Fubini theorem, integration in heat time gives a factor $2^{-2k}$. We estimate the source terms separately below. For $A_lJ_{s,l}\psi_m$, use \eqref{eAinfty1}, \eqref{eAonederivativel2} with $\sigma=0$, \eqref{eJimprove}, and \eqref{eek2}. The Low-High interaction satisfies
\begin{align*}
2^{2k}\int_0^\infty\frac{\|P_k(P_{\leq k-5}A(s)P_{\sim k}J_s\psi(s))\|_{L_x^2}}{t+s}\,\mathrm{d}s&\lesssim\epsilon_1^3t^\gamma E_k\,2^{2k}\int_0^\infty(t+s)^{-\frac{5}{2}}(1+s2^{2k})^{-10}\,\mathrm{d}s\\
&\lesssim\epsilon_1^3t^{-\frac{5}{2}+\gamma}E_k.
\end{align*}
For the High-Low interaction, Bernstein inequality gives $\|P_{\leq k-5}J_s\psi(s)\|_{L_x^\infty}\lesssim\epsilon_1t^\gamma2^k$, and hence
\begin{align*}
2^{2k}\int_0^\infty\frac{\|P_k(P_{\sim k}A(s)P_{\leq k-5}J_s\psi(s))\|_{L_x^2}}{t+s}\,\mathrm{d}s&\lesssim\epsilon_1^3t^\gamma b_k\,2^{2k}\int_0^\infty(t+s)^{-2}(1+s2^{2k})^{-19}\,\mathrm{d}s\\
&\lesssim\epsilon_1^3t^{-2+\gamma}b_k.
\end{align*}
For the High-High interaction, the corresponding bound is
\begin{equation*}
\epsilon_1^3t^\gamma2^{2k}\sum_{j\geq k-4}b_jE_j\int_0^\infty(t+s)^{-2}(1+s2^{2j})^{-10}\,\mathrm{d}s\lesssim\epsilon_1^3t^{-2+\gamma}\sum_{j\geq k-4}2^{2(k-j)}b_jE_j.
\end{equation*}
These estimates are summable in $\ell_k^2$ by \eqref{eenergyenvelope3} and \eqref{eek2}.
For $(x_lA_l)\psi_m$, Leibniz rule, \eqref{eappendixA15}, \eqref{eappendixA11}, \eqref{eAonederivativel21} with $\sigma=0$, and \eqref{epsidecay1} give
\begin{align*}
&\int_0^\infty\frac{\|\nabla^2[(x_lA_l)\psi](s)\|_{L_x^2}}{t+s}\,\mathrm{d}s\\
&\quad\lesssim t^{-1}\int_0^\infty\left[\|xA(s)\|_{L_x^\infty}\|\nabla^2\psi(s)\|_{L_x^2}+\|\nabla A(s)\|_{L_x^\infty}\|x\nabla\psi(s)\|_{L_x^2}+\|\nabla^2A(s)\|_{L_x^2}\|x\psi(s)\|_{L_x^\infty}\right]\,\mathrm{d}s\\
&\qquad+t^{-1}\int_0^\infty\left[\|A(s)\|_{L_x^\infty}\|\nabla\psi(s)\|_{L_x^2}+\|\nabla A(s)\|_{L_x^2}\|\psi(s)\|_{L_x^\infty}\right]\,\mathrm{d}s\lesssim\epsilon_1^3t^{-2+\gamma}(1+\log t).
\end{align*}
For the remaining source terms, \eqref{eAonederivativel2} with $\sigma=\epsilon,0$ gives
\begin{equation*}
\int_0^\infty\|\nabla^3A(s)\|_{L_x^2}\,\mathrm{d}s\lesssim\epsilon_1^2\left[t^{-1}\int_0^1s^{-1+\frac{\epsilon}{2}}\,\mathrm{d}s+\int_1^\infty(t+s)^{-1}s^{-1}\,\mathrm{d}s\right]\lesssim\epsilon_1^2t^{-1}(1+\log t).
\end{equation*}
Thus the second-order product estimate, \eqref{epsidecay1}, \eqref{eAonederivativelinfty}, and \eqref{epsiinte5} imply
\begin{align*}
\int_0^\infty\|\nabla^2[(\nabla\cdot A)\psi](s)\|_{L_x^2}\,\mathrm{d}s&\lesssim\epsilon_1t^{-1}\int_0^\infty\|\nabla^3A(s)\|_{L_x^2}\,\mathrm{d}s+\sup_{s\geq0}\|\nabla A(s)\|_{L_x^\infty}\int_0^\infty\|\nabla^2\psi(s)\|_{L_x^2}\,\mathrm{d}s\\
&\lesssim\epsilon_1^3t^{-2}(1+\log t).
\end{align*}
The same product estimate, \eqref{eAinfty1}, \eqref{epsidecay1}, \eqref{epsiinte5}, and \eqref{eappendixA19} give
\begin{align*}
\int_0^\infty\|\nabla^2(\psi\overline{\psi}\psi)(s)\|_{L_x^2}\,\mathrm{d}s&\lesssim\epsilon_1^3t^{-2},\\
\int_0^\infty\|\nabla^2(A_lA_l\psi)(s)\|_{L_x^2}\,\mathrm{d}s&\lesssim\epsilon_1^5t^{-3}.
\end{align*}
Combining the source estimates with dyadic Duhamel and Minkowski inequality proves \eqref{eappendixA22}.

We now estimate the two quadratic corrections in \eqref{eappendixA1}. If the heat correction has the higher frequency, \eqref{epsidecay1}, \eqref{eUfree1}, and Bernstein inequality bound the High-Low contribution at output frequency $k$ by
\begin{equation*}
\epsilon_1t^{-1}2^{4k}\int_0^\infty\|P_{\sim k}(\psi(s)-e^{s\Delta}\psi(0))\|_{L_x^2}\,\mathrm{d}s.
\end{equation*}
The High-High contribution is bounded by
\begin{equation*}
\epsilon_1t^{-1}\sum_{j\geq k-4}2^{3(k-j)}\left[2^{4j}\int_0^\infty\|P_j(\psi(s)-e^{s\Delta}\psi(0))\|_{L_x^2}\,\mathrm{d}s\right].
\end{equation*}
By \eqref{eappendixA22}, their $\ell_k^2$ norms are bounded by $\epsilon_1^4t^{-3+\gamma}(1+\log t)$.

For the Low-High interaction, collect the low frequency correction into $P_{\leq k-5}(\psi-e^{s\Delta}\psi(0))$. We first consider the case in which the derivative in the heat integral falls on the correction. By \eqref{eJs1},
\begin{equation*}
2i(s+it)(\xi-\eta)_n\widehat{P_{\sim k}\psi(s)}(\xi-\eta)=\widehat{J_{s,n}P_{\sim k}\psi(s)}(\xi-\eta)+i\partial_{\eta_n}\widehat{P_{\sim k}\psi(s)}(\xi-\eta).
\end{equation*}
Using $\sum_n(\xi-\eta)_n^2=|\xi-\eta|^2$ and integrating by parts in $\eta$ with $\xi$ fixed, we obtain
\begin{align*}
&P_k\partial_i\partial_j\partial_l\left[P_{\sim k}\psi(s)\,\partial_mP_{\leq k-5}(\psi(s)-e^{s\Delta}\psi(0))\right]\\
&\quad=\frac{1}{2i(s+it)}\sum_nP_k\mathcal F_\xi^{-1}\Biggl[\int\frac{\xi_i\xi_j\xi_l\eta_m(\xi-\eta)_n}{|\xi-\eta|^2}\widehat{J_{s,n}P_{\sim k}\psi(s)}(\xi-\eta)\\
&\qquad\qquad\times\widehat{P_{\leq k-5}(\psi(s)-e^{s\Delta}\psi(0))}(\eta)\,\mathrm{d}\eta\Biggr]\\
&\qquad-\frac{1}{2i(s+it)}\sum_nP_k\mathcal F_\xi^{-1}\Biggl[\int\frac{\xi_i\xi_j\xi_l\eta_m(\xi-\eta)_n}{|\xi-\eta|^2}\widehat{P_{\sim k}\psi(s)}(\xi-\eta)\\
&\qquad\qquad\times\widehat{x_nP_{\leq k-5}(\psi(s)-e^{s\Delta}\psi(0))}(\eta)\,\mathrm{d}\eta\Biggr]\\
&\qquad-\frac{1}{2(s+it)}\sum_nP_k\mathcal F_\xi^{-1}\Biggl[\int\partial_{\eta_n}\left[\frac{\xi_i\xi_j\xi_l\eta_m(\xi-\eta)_n}{|\xi-\eta|^2}\right]\widehat{P_{\sim k}\psi(s)}(\xi-\eta)\\
&\qquad\qquad\times\widehat{P_{\leq k-5}(\psi(s)-e^{s\Delta}\psi(0))}(\eta)\,\mathrm{d}\eta\Biggr].
\end{align*}
The commutators with the frequency projections are included in the differentiated symbol below. By Lemma~\ref{lmet},
\begin{align*}
\left\|\frac{\xi_i\xi_j\xi_l\eta_m(\xi-\eta)_n}{|\xi-\eta|^2}\varphi_k(\xi)\varphi_{\sim k}(\xi-\eta)\varphi_{\leq k-5}(\eta)\right\|_{\mathcal S^\infty}&\lesssim2^{3k},\\
\left\|\partial_{\eta_n}\left[\frac{\xi_i\xi_j\xi_l\eta_m(\xi-\eta)_n}{|\xi-\eta|^2}\varphi_k(\xi)\varphi_{\sim k}(\xi-\eta)\varphi_{\leq k-5}(\eta)\right]\right\|_{\mathcal S^\infty}&\lesssim2^{2k}.
\end{align*}
We estimate the heat integral of each of the three terms in $L_x^2$. For the contribution containing $P_{\sim k}J_{s,n}\psi$, Lemma~\ref{l multilinear estimates}, \eqref{eJimprove}, and \eqref{eappendixA4} give
\begin{equation*}
\epsilon_1t^\gamma E_k\,2^{3k}\int_0^\infty\frac{(1+s2^{2k})^{-10}}{t+s}\|\psi(s)-e^{s\Delta}\psi(0)\|_{L_x^\infty}\,\mathrm{d}s\lesssim\epsilon_1^4t^{-\frac{7}{2}+\gamma}E_k.
\end{equation*}
For the contribution containing $P_{\leq k-5}[x_n(\psi-e^{s\Delta}\psi(0))]$, \eqref{eheatenergyimproved2} with $\sigma=0$ and \eqref{eappendixA16} give
\begin{equation*}
\epsilon_1b_k\,2^{3k}\int_0^\infty\frac{(1+s2^{2k})^{-20}}{t+s}\|x(\psi(s)-e^{s\Delta}\psi(0))\|_{L_x^\infty}\,\mathrm{d}s\lesssim\epsilon_1^4t^{-3+\gamma}(1+\log t)b_k.
\end{equation*}
For the combined symbol-derivative term, \eqref{eheatenergyimproved2} with $\sigma=0$ and the uniform bound in \eqref{eappendixA4} give
\begin{equation*}
\epsilon_1b_k\,2^{2k}\int_0^\infty\frac{(1+s2^{2k})^{-20}}{t+s}\|\psi(s)-e^{s\Delta}\psi(0)\|_{L_x^\infty}\,\mathrm{d}s\lesssim\epsilon_1^4t^{-3}b_k.
\end{equation*}
When the derivative in the heat integral falls on the high frequency factor, replace $\eta_m$ by $(\xi-\eta)_m$ throughout, both symbol bounds and the resulting estimates remain unchanged. The free heat factor satisfies the same bounds by \eqref{eJs2}, \eqref{eek4}, and \eqref{eheatenergyimproved2} at $s=0$. For a conjugated high frequency factor, the same computation uses $\overline{J_{s,n}\psi(s)}$ and $s-it$. Finally, \eqref{eenergyenvelope3} and \eqref{eek2} bound the Low-High contributions to both quadratic corrections by $\epsilon_1^4t^{-3+\gamma}(1+\log t)$.

It remains to estimate the remainder in \eqref{eappendixA1}. First suppose that $A$ has frequency comparable to the largest input frequency $2^j$. By \eqref{eAonederivativel2} with $\sigma=0$ and \eqref{epsidecay1}, its contribution at output frequency $k$ is bounded by
\begin{equation*}
\epsilon_1^4\sum_{j\geq k-4}2^{3k-j}b_j\int_0^\infty(t+s)^{-3}(1+s2^{2j})^{-19}\,\mathrm{d}s\lesssim\epsilon_1^4t^{-3}\sum_{j\geq k-4}2^{3(k-j)}b_j.
\end{equation*}
Suppose instead that a $\psi$ factor has frequency $2^j$ and $A$ is supported in frequencies at most $2^{j-5}$. By \eqref{eJs1}, extracting one high frequency derivative and integrating by parts in $\eta$ with $\xi,\zeta$ fixed gives
\begin{align*}
&P_k\partial_i\partial_l\partial_m\left[P_{\sim j}\psi(s)\,P_{\leq j-5}A(s)\,P_{\leq j+4}\psi(s)\right]\\
&\quad=-\frac{1}{2(s+it)}\sum_nP_k\mathcal F_\xi^{-1}\Biggl[\iint\frac{\xi_i\xi_l\xi_m(\xi-\eta)_n}{|\xi-\eta|^2}\widehat{J_{s,n}P_{\sim j}\psi(s)}(\xi-\eta)\\
&\qquad\qquad\times\widehat{P_{\leq j-5}A(s)}(\eta-\zeta)\widehat{P_{\leq j+4}\psi(s)}(\zeta)\,\mathrm{d}\zeta\,\mathrm{d}\eta\Biggr]\\
&\qquad+\frac{1}{2(s+it)}\sum_nP_k\mathcal F_\xi^{-1}\Biggl[\iint\frac{\xi_i\xi_l\xi_m(\xi-\eta)_n}{|\xi-\eta|^2}\widehat{P_{\sim j}\psi(s)}(\xi-\eta)\\
&\qquad\qquad\times\widehat{x_nP_{\leq j-5}A(s)}(\eta-\zeta)\widehat{P_{\leq j+4}\psi(s)}(\zeta)\,\mathrm{d}\zeta\,\mathrm{d}\eta\Biggr]\\
&\qquad+\frac{i}{2(s+it)}\sum_nP_k\mathcal F_\xi^{-1}\Biggl[\iint\partial_{\eta_n}\left[\frac{\xi_i\xi_l\xi_m(\xi-\eta)_n}{|\xi-\eta|^2}\right]\widehat{P_{\sim j}\psi(s)}(\xi-\eta)\\
&\qquad\qquad\times\widehat{P_{\leq j-5}A(s)}(\eta-\zeta)\widehat{P_{\leq j+4}\psi(s)}(\zeta)\,\mathrm{d}\zeta\,\mathrm{d}\eta\Biggr].
\end{align*}
The commutators with the frequency projections are included in the differentiated symbol below. By Lemma~\ref{lmet},
\begin{align*}
\left\|\frac{\xi_i\xi_l\xi_m(\xi-\eta)_n}{|\xi-\eta|^2}\varphi_k(\xi)\varphi_{\sim j}(\xi-\eta)\varphi_{\leq j-5}(\eta-\zeta)\varphi_{\leq j+4}(\zeta)\right\|_{\mathcal S^\infty}&\lesssim2^{3k-j},\\
\left\|\partial_{\eta_n}\left[\frac{\xi_i\xi_l\xi_m(\xi-\eta)_n}{|\xi-\eta|^2}\varphi_k(\xi)\varphi_{\sim j}(\xi-\eta)\varphi_{\leq j-5}(\eta-\zeta)\varphi_{\leq j+4}(\zeta)\right]\right\|_{\mathcal S^\infty}&\lesssim2^{3k-2j}.
\end{align*}
By Lemma~\ref{l multilinear estimates}, \eqref{eJimprove}, \eqref{eheatenergyimproved2} with $\sigma=0$, \eqref{eAinfty1}, \eqref{epsidecay1}, and \eqref{eappendixA14}, the heat integrals of the terms containing $P_{\sim j}J_s\psi$, $P_{\leq j-5}(xA)$, and the differentiated symbol are bounded in $L_x^2$, respectively, by
\begin{align*}
&\epsilon_1^4t^\gamma E_j\,2^{3k-j}\int_0^\infty(t+s)^{-\frac{7}{2}}(1+s2^{2j})^{-10}\,\mathrm{d}s\lesssim\epsilon_1^4t^{-\frac{7}{2}+\gamma}2^{3(k-j)}E_j,\\
&\epsilon_1^4t^{-1+\gamma}(1+\log t)b_j\,2^{3k-j}\int_0^\infty(t+s)^{-2}(1+s2^{2j})^{-20}\,\mathrm{d}s\lesssim\epsilon_1^4t^{-3+\gamma}(1+\log t)2^{3(k-j)}b_j,\\
&\epsilon_1^4b_j\,2^{3k-2j}\int_0^\infty(t+s)^{-\frac{7}{2}}(1+s2^{2j})^{-20}\,\mathrm{d}s\lesssim\epsilon_1^4t^{-3}2^{3(k-j)}b_j.
\end{align*}
For the last inequality, we used
\begin{equation*}
2^j\int_0^\infty(t+s)^{-\frac{7}{2}}(1+s2^{2j})^{-20}\,\mathrm{d}s\lesssim\min\left\{2^jt^{-\frac{5}{2}},2^{-j}t^{-\frac{7}{2}}\right\}\lesssim t^{-3}.
\end{equation*}
Conjugated factors are estimated as above. Combining these bounds with those from the other cases, summing over $j\geq k-4$, and taking the $\ell_k^2$ norm, \eqref{eenergyenvelope3} and \eqref{eek2} give
\begin{equation*}
\|\nabla^3A_m^{(\geq4)}(0,t)\|_{L_x^2}\lesssim\epsilon_1^4t^{-3+\gamma}(1+\log t).
\end{equation*}
This proves \eqref{eAhigh3}.

Finally, interpolation with \eqref{eAhigh1} yields
\begin{align*}
\|\nabla^2A_m^{(\geq4)}(0,t)\|_{H^1}&\lesssim\|A_m^{(\geq4)}(0,t)\|_{L_x^2}^{\frac{1}{3}}\|\nabla^3A_m^{(\geq4)}(0,t)\|_{L_x^2}^{\frac{2}{3}}+\|\nabla^3A_m^{(\geq4)}(0,t)\|_{L_x^2}\\
&\lesssim\epsilon_1^4t^{-\frac{8}{3}+\frac{2\gamma}{3}}	(1+\log t)^{\frac{2}{3}}\lesssim\epsilon_1^4t^{-\frac{5}{2}},
\end{align*}
which proves \eqref{eAhigh2}.
\end{proof}

\subsection{Positive homogeneous estimates for the temporal field} \label{app:psi3sobolev}

We prove Lemma~\ref{lpsi3sobolev}.

\begin{proof}
We fix $t$ below. By \eqref{eheatenergyimproved2}, \eqref{econnectionenve1}, \eqref{edkdef2}, and Bernstein inequality,
\[
\begin{gathered}
  \|P_k\psi(s)\|_{L_x^2}\lesssim\epsilon_1b_k(1+s2^{2k})^{-20}, \qquad \|P_kA(s)\|_{L_x^2}\lesssim\epsilon_1^2b_k^2(1+s2^{2k})^{-19},\\
  \max_F\|P_kF(s)\|_{L_x^2}\lesssim\epsilon_1^2 2^k(1+s2^{2k})^{-19},
\end{gathered}
\]
where $F$ is as in \eqref{eFdefinition}. The estimates for $F$ also hold for $\psi_l\psi_m$.

We first prove
\begin{equation}\label{eappendixA23}
  \sup_{k\in\mathbb Z}\sup_{s\geq0}2^{-k}(1+s2^{2k})^8\|P_k\psi_3(s)\|_{L_x^2}\lesssim\epsilon_1.
\end{equation}
Let $C$ be the best constant in this estimate. The initial identity gives $\|P_k\psi_3(0)\|_{L_x^2}\lesssim\epsilon_1 2^k$. By \eqref{epsi3equation}, Bernstein inequality, \eqref{epolynomialconvolution}, and \eqref{eenergyenvelope3},
\begin{align*}
&\sup_{s\geq0}(1+s2^{2k})^8\|P_k(\psi_3(s)-e^{s\Delta}\psi_3(0))\|_{L_x^2}\\
&\quad\lesssim C\epsilon_1^32^k\left[1+\sum_{j\geq k-O(1)}b_j^2\sum_{l\leq j+4}2^{2(l-j)}+\sum_{j\geq k-O(1)}b_j\sum_{l\leq j+4}2^{l-j}b_l\right]\lesssim C\epsilon_1^32^k.
\end{align*}
The terms containing $A_jA_j\psi_3$ are estimated in the same way and contribute at most $C\epsilon_1^52^k$. Thus $C\lesssim1+\epsilon_1^2C$, proving \eqref{eappendixA23}.

By the same product decomposition, \eqref{eheatenergyimproved2}, \eqref{econnectionenve1}, and \eqref{eenergyenvelope3}, for $0<\sigma\leq\epsilon$,
\begin{equation}\label{eappendixA24}
  \left(\sum_k2^{2(\sigma-1)k}\sup_{s\geq0}(1+s2^{2k})^{16}\max_F\|P_kF(s)\|_{L_x^2}^2\right)^{\frac12}\lesssim\epsilon_1\|\psi(0)\|_{\dot H^\sigma},
\end{equation}
\[
  \left(\sum_k2^{2(\sigma-1)k}\sup_{s\geq0}(1+s2^{2k})^{16}\|P_k(A_l\psi_l)(s)\|_{L_x^2}^2\right)^{\frac12}\lesssim\epsilon_1^2\|\psi(0)\|_{\dot H^\sigma}.
\]
In particular, $\|\psi_3(0)\|_{\dot H^{\sigma-1}}\lesssim\|\psi(0)\|_{\dot H^\sigma}$.

Let $C$ now be the best constant in the $\psi_3$ part of \eqref{epsi3sobolev}. By \eqref{eappendixA23} and \eqref{epolynomialconvolution}, the High-High contribution of $F\psi_3$ satisfies
\[
\begin{aligned}
  &2^{(\sigma-1)k}\sup_{s\geq0}(1+s2^{2k})^8\left\|P_k\int_0^s e^{(s-r)\Delta}\sum_{j\geq k-4}P_jF(r)P_{j+O(1)}\psi_3(r)\,\mathrm{d}r\right\|_{L_x^2}\\
  &\qquad\lesssim\epsilon_1\sum_{j\geq k-4}2^{\sigma(k-j)}2^{(\sigma-1)j}\sup_{r\geq0}(1+r2^{2j})^8\|P_jF(r)\|_{L_x^2}.
\end{aligned}
\]
Taking the $\ell_k^2$ norm and using \eqref{eappendixA24} gives $C_\sigma\epsilon_1^2\|\psi(0)\|_{\dot H^\sigma}$, independently of $C$. The remaining nonlinear terms are bounded by $C_0C\epsilon_1^2\|\psi(0)\|_{\dot H^\sigma}$, where $C_0$ is uniform for $0<\sigma\leq\epsilon<1$. Hence
\begin{equation*}
C\leq C_\sigma+C_0\epsilon_1^2C.
\end{equation*}
Absorbing the last term proves the $\psi_3$ estimate. The estimate for $\partial_l\psi_l+iA_l\psi_l$ follows from \eqref{eheatenergyimproved2} and the bound above for $A_l\psi_l$.
\end{proof}

\subsection{The higher order temporal connection}\label{app:A3remainder}
We prove Lemma~\ref{lA3higherorder}.
\begin{proof}
By \eqref{epsi3equation}, \eqref{epsi31}, \eqref{epsi32}, \eqref{eappendixA9}, \eqref{eappendixA19}, and \eqref{eAonederivativelinfty}, restarting the Duhamel formula on $[\frac{s}{2},s]$ and absorbing the lower order terms give
\begin{align}
\left\|\nabla^a\left(\psi_3(s)-i\partial_l e^{s\Delta}\psi_l(0)\right)\right\|_{L_x^2}&\lesssim\epsilon_1^3t^{-1}s^{-\frac{a}{2}},\qquad a=1,2,\quad s>0,\label{eappendixA25}\\
\left\|\nabla\left(\psi_3(s)-i\partial_l e^{s\Delta}\psi_l(0)\right)\right\|_{L_x^\infty}&\lesssim\epsilon_1^3t^{-\frac{5}{2}}(1+\log t)s^{-\frac{1}{2}},\qquad 0<s\leq1.\label{eappendixA26}
\end{align}
By \eqref{eAL2}, \eqref{eAonederivativel21}, \eqref{eAonederivativelinfty}, \eqref{eheatenergyimproved2}, and \eqref{epsidecay1}, for $0<s\leq1$,
\begin{equation}\label{eappendixA27}
\|A(s)\|_{H^\epsilon}\lesssim\epsilon_1^2t^{-1}(1+\log t),\qquad \max_F\|F(s)\|_{H^\epsilon}\lesssim\epsilon_1^2t^{-1},\qquad \max_F\|F(s)\|_{L_x^\infty}\lesssim\epsilon_1^2t^{-2}(1+\log t).
\end{equation}
The bounds for $F$ also hold for $\psi_l\psi_m$. The fractional Leibniz rule and \eqref{epsi34} give
\begin{equation*}
\|(A_l\psi_l)(0)\|_{H^\epsilon}+\sup_{0<s\leq1}s^{\frac{1}{2}}\left[\|A(s)\partial_l e^{s\Delta}\psi_l(0)\|_{H^\epsilon}+\max_F\|F(s)\partial_l e^{s\Delta}\psi_l(0)\|_{H^\epsilon}\right]\lesssim\epsilon_1^3t^{-1}.
\end{equation*}
Using \eqref{epsi3equation}, \eqref{epsi32}, and $\int_0^s(s-r)^{-\frac{1}{2}}r^{-\frac{1}{2}}\,\mathrm{d}r\lesssim1$, the Duhamel formula and absorption prove the first estimate below. Restarting the Duhamel formula on $[\frac{s}{2},s]$ and using the fractional Leibniz rule, \eqref{eheatoperator}, and \eqref{eappendixA26}, we obtain the second:
\begin{align}
\sup_{0\leq s\leq1}\|\psi_3(s)-i\partial_l e^{s\Delta}\psi_l(0)\|_{H^\epsilon}&\lesssim\epsilon_1^3t^{-1},\notag\\
\left\|\nabla^2\left(\psi_3(s)-i\partial_l e^{s\Delta}\psi_l(0)\right)\right\|_{L_x^2}&\lesssim\epsilon_1^3t^{-1}s^{-1+\frac{\epsilon}{2}},\qquad 0<s\leq1.\label{eappendixA28}
\end{align}
For the first weighted estimate, \eqref{eJs2}, \eqref{eJs3}, and heat energy estimate give
\begin{gather*}
J_{s,m}\left(i\partial_l e^{s\Delta}\psi_l(0)\right)=i\partial_l e^{s\Delta}J_{0,m}\psi_l(0)-ie^{s\Delta}\psi_m(0),\\
\int_0^1\|\nabla e^{s\Delta}J_0\psi(0)\|_{L_x^2}^2\,\mathrm{d}s\leq\frac{1}{2}\|J_0\psi(0)\|_{L_x^2}^2.
\end{gather*}
Hence, \eqref{eAinfty1}, \eqref{eAonederivativel21}, \eqref{eappendixA27}, and \eqref{booassum} imply
\begin{align*}
&\|J_0(A_l\psi_l)(0)\|_{L_x^2}+\left(\int_0^1\left\|A(s)J_s\left(i\partial_l e^{s\Delta}\psi_l(0)\right)\right\|_{L_x^2}^2\,\mathrm{d}s\right)^{\frac{1}{2}}\\
&\quad+\int_0^1\max_F\left\|F(s)J_s\left(i\partial_l e^{s\Delta}\psi_l(0)\right)\right\|_{L_x^2}\,\mathrm{d}s\\
&\quad\lesssim\epsilon_1^3\left[t^{-1}+t^{-\frac{3}{2}}(1+t^\gamma)+t^{-2}(1+t^\gamma)(1+\log t)\right]\lesssim\epsilon_1^3t^{-1}.
\end{align*}
For the curvature term, \eqref{eJs1} gives
\begin{equation*}
J_{s,m}\left[i\operatorname{Im}(\psi_3\overline{\psi_j})\psi_j\right]-i\operatorname{Im}\left[(J_{s,m}\psi_3)\overline{\psi_j}\right]\psi_j=(s+it)\partial_m|\psi_j|^2\psi_3-(s+it)\partial_m(\psi_j^2)\overline{\psi_3}-2it\psi_j^2\partial_m\overline{\psi_3}.
\end{equation*}
By \eqref{eAonederivativelinfty}, \eqref{eappendixA19}, \eqref{eheatenergyimproved2}, and \eqref{epsi34},
\begin{align*}
&\int_0^1(t+s)\left[\|\nabla A(s)\|_{L_x^\infty}\|\nabla^2e^{s\Delta}\psi(0)\|_{L_x^2}+\|\nabla^2A(s)\|_{L_x^2}\|\nabla e^{s\Delta}\psi(0)\|_{L_x^\infty}\right.\\
&\left.+\|\psi(s)\|_{L_x^\infty}\|\nabla\psi(s)\|_{L_x^2}\|\nabla e^{s\Delta}\psi(0)\|_{L_x^\infty}+\|\psi(s)\|_{L_x^\infty}^2\|\nabla^2e^{s\Delta}\psi(0)\|_{L_x^2}\right]\,\mathrm{d}s\lesssim\epsilon_1^3t^{-1}(1+\log t).
\end{align*}
Together with \eqref{epsi31}, \eqref{epsi32}, and \eqref{epsi33}, these estimates control the commutator terms. The $L_x^2$ energy estimate and absorption therefore yield
\begin{equation}\label{eappendixA29}
\begin{aligned}
&\sup_{0\leq s\leq1}\left\|J_s\left(\psi_3(s)-i\partial_l e^{s\Delta}\psi_l(0)\right)\right\|_{L_x^2}\\
&\quad+\left(\int_0^1\left\|\nabla J_s\left(\psi_3(s)-i\partial_l e^{s\Delta}\psi_l(0)\right)\right\|_{L_x^2}^2\,\mathrm{d}s\right)^{\frac{1}{2}}\lesssim\epsilon_1^3t^{-1}(1+\log t).
\end{aligned}
\end{equation}
We now apply \eqref{eA3higherrep}. By \eqref{eappendixA6}, \eqref{epsi31}, \eqref{epsiinte2}, \eqref{eAinfty1}, and \eqref{epsi34},
\begin{equation}\label{eappendixA30}
\|A_3^{(\geq4)}(0,t)\|_{L_x^2}\lesssim\epsilon_1^4t^{-\frac{5}{2}}(1+\log t)+\epsilon_1^4\int_0^\infty s^{-\frac{1}{2}}(t+s)^{-\frac{5}{2}}\,\mathrm{d}s+\epsilon_1^4t^{-2}\lesssim\epsilon_1^4t^{-2}.
\end{equation}
For the two products containing the free temporal factor, the paraproduct estimates in the proof of \eqref{eAhigh3} apply with the same symbol bounds $2^{3k},2^{2k}$ and High-High factor $2^{2(k-j)}$. Thus \eqref{eappendixA22}, \eqref{eenergyenvelope3}, and \eqref{eek2} give
\begin{align*}
&\left\|\int_0^\infty\partial_l(\psi_l(s)-e^{s\Delta}\psi_l(0))\overline{\partial_m e^{s\Delta}\psi_m(0)}\,\mathrm{d}s\right\|_{\dot H^2}+\left\|\int_0^\infty A_l(s)\psi_l(s)\overline{\partial_m e^{s\Delta}\psi_m(0)}\,\mathrm{d}s\right\|_{\dot H^2}\\
&\quad\lesssim\epsilon_1t^{-1}\left(\sum_j\left[2^{4j}\int_0^\infty\|P_j(\psi(s)-e^{s\Delta}\psi(0))\|_{L_x^2}\,\mathrm{d}s\right]^2\right)^{\frac{1}{2}}+\epsilon_1^4t^{-3+\gamma}(1+\log t)\\
&\quad\lesssim\epsilon_1^4t^{-3+\gamma}(1+\log t).
\end{align*}
For the product of the free temporal factor and the temporal correction, we use \eqref{bs5}, \eqref{bsc11}, \eqref{bsc12}, and \eqref{eheatoperator}:
\begin{align}
\left(\int_0^1\|\nabla e^{s\Delta}\psi(0)\|_{L_x^\infty}^2\,\mathrm{d}s\right)^{\frac{1}{2}}+\int_0^1\|\nabla^2e^{s\Delta}\psi(0)\|_{L_x^\infty}\,\mathrm{d}s&\lesssim\epsilon_1t^{-1}(1+\log t),\notag\\
\int_0^1\left[\|\nabla^2e^{s\Delta}J_0\psi(0)\|_{L_x^2}+s^{-\frac{1}{2}}\|\nabla e^{s\Delta}J_0\psi(0)\|_{L_x^2}\right]\,\mathrm{d}s&\lesssim\sum_k\|P_k(xf)\|_{L_x^2}\lesssim\epsilon_1t^{\frac{\beta}{2}}.\label{eaappendixA31}
\end{align}
By \eqref{eJs1}, for complex-valued functions $u,v$, $2it\nabla(u\overline v)=(J_su)\overline v-u\overline{J_sv}-2s\left[(\nabla u)\overline v-u\nabla\overline v\right]$. Apply this identity with $u=\psi_3(s)-i\partial_l e^{s\Delta}\psi_l(0)$ and $v=\partial_m e^{s\Delta}\psi_m(0)$, and differentiate once. By \eqref{eappendixA29}, \eqref{eaappendixA31}, \eqref{eappendixA25}, \eqref{eappendixA26}, and \eqref{epsi32},
\begin{align*}
&\int_0^1\left\|\nabla^2\left[\partial_m e^{s\Delta}\psi_m(0)\overline{\psi_3(s)-i\partial_l e^{s\Delta}\psi_l(0)}\right]\right\|_{L_x^2}\,\mathrm{d}s\\
&\quad\lesssim\epsilon_1^4\left[t^{-3}(1+\log t)^2+t^{-\frac{7}{2}+\frac{\beta}{2}}(1+\log t)+t^{-3}\right]\lesssim\epsilon_1^4t^{-2}.
\end{align*}
On $s\geq1$, \eqref{eappendixA25}, \eqref{epsi32}, \eqref{epsi34}, and \eqref{epsihat2} give
\begin{align*}
&\int_1^\infty\left\|\nabla^2\left[\partial_m e^{s\Delta}\psi_m(0)\overline{\psi_3(s)-i\partial_l e^{s\Delta}\psi_l(0)}\right]\right\|_{L_x^2}\,\mathrm{d}s\\
&\quad\lesssim\epsilon_1^4t^{-1}\int_1^\infty(t+s)^{-1}s^{-\frac{3}{2}}\,\mathrm{d}s+\epsilon_1^4t^{-\frac{5}{2}}(1+\log t)\int_1^\infty s^{-2}\,\mathrm{d}s\lesssim\epsilon_1^4t^{-2}.
\end{align*}
For the remaining products, fractional Leibniz rule, \eqref{epsiinte5}, \eqref{eappendixA19}, \eqref{epsiAdecay2}, \eqref{eAinfty1}, and \eqref{epsidecay1} imply
\begin{equation*}
\int_0^\infty\|\nabla^2(A_l\psi_l)(s)\|_{L_x^2}\,\mathrm{d}s\lesssim\epsilon_1^3\left(t^{-\frac{3}{2}}+t^{-2}+t^{-1}\right)\lesssim\epsilon_1^3t^{-1}.
\end{equation*}
Together with \eqref{eappendixA6}, \eqref{eUimpro1}, \eqref{eappendixA25}, \eqref{eappendixA28}, and \eqref{epsi32}, this gives
\begin{align*}
&\int_0^\infty\left[\left\|\nabla^2\left[\partial_m(\psi_m(s)-e^{s\Delta}\psi_m(0))\overline{\psi_3(s)-i\partial_l e^{s\Delta}\psi_l(0)}\right]\right\|_{L_x^2}\right.\\
&\hspace{24mm}\left.+\left\|\nabla^2\left[A_m(s)\psi_m(s)\overline{\psi_3(s)-i\partial_l e^{s\Delta}\psi_l(0)}\right]\right\|_{L_x^2}\right]\,\mathrm{d}s\\
&\quad\lesssim\epsilon_1^3t^{-1}\sum_{a=1}^3\int_0^\infty s^{-\frac{3-a}{2}}\|\nabla^a(\psi(s)-e^{s\Delta}\psi(0))\|_{L_x^\infty}\,\mathrm{d}s+\epsilon_1^6t^{-\frac{7}{2}}(1+\log t)\left(1+\int_0^1s^{-1+\frac{\epsilon}{2}}\,\mathrm{d}s\right)\\
&\quad\lesssim\epsilon_1^6t^{-\frac{7}{2}}(1+\log t).
\end{align*}
Combining these estimates with \eqref{eA3higherrep}, \eqref{eappendixA30}, and \eqref{bop}, we obtain
\begin{equation*}
\|A_3^{(\geq4)}(0,t)\|_{H^2}\lesssim\epsilon_1^4t^{-2}+\epsilon_1^4t^{-3+\gamma}(1+\log t)+\epsilon_1^6t^{-\frac{7}{2}}(1+\log t)\lesssim\epsilon_1^4t^{-2}.
\end{equation*}
This proves \eqref{eA35}.
\end{proof}

\section{Local estimates for the differential fields}\label{app:B}

In this appendix we prove Proposition~\ref{locest}.

\subsection{The unweighted estimates and the local bootstrap}
The small data global estimates of \cite[Propositions~3.8, 4.5, 6.2 and Section~7]{BIKT2011} imply
\begin{equation}\label{eappendixB1}
\sum_k\langle2^k\rangle^{2\epsilon}\left(\|P_k\psi\|_{L_t^\infty([0,T];L_x^2)}^2+\|P_k\psi\|_{L_t^{\frac{2}{1-\frac{\epsilon}{2}}}([0,T];L_x^{\frac{4}{\epsilon}})}^2\right)\lesssim\epsilon_0^2,\qquad 0<T\leq1.
\end{equation}
Here the second estimate uses the nonendpoint Strichartz extension of their linear estimate. Indeed, $\frac{2}{q}+\frac{2}{r}=1$ for $q=\frac{2}{1-\frac{\epsilon}{2}}>2$ and $r=\frac{4}{\epsilon}<\infty$. For the ordinary forcing terms, the proof of \cite[Lemma~7.3]{BIKT2011} applies by homogeneous Strichartz. For the local smoothing terms, after the angular and Galilean reductions in \cite[Lemma~7.4]{BIKT2011}, dual local smoothing and homogeneous Strichartz give the full time estimate at fixed $x\cdot e=y$, since the input is in $L_{t,x_\perp}^2$ and $q>2$, we apply Christ-Kiselev and then integrate in $y$. Thus \eqref{eappendixB1} does not use any weighted bootstrap assumption. By H\"older, Bernstein and Cauchy-Schwarz inequalities,
\begin{align*}
\sum_{k>0}\|P_k\psi\|_{L_t^2([0,T];L_x^\infty)}&\lesssim T^{\frac{\epsilon}{4}}\sum_{k>0}2^{\frac{\epsilon}{2}k}\|P_k\psi\|_{L_t^{\frac{2}{1-\frac{\epsilon}{2}}}L_x^{\frac{4}{\epsilon}}}\\
&\leq T^{\frac{\epsilon}{4}}\left(\sum_{k>0}2^{-\epsilon k}\right)^{\frac{1}{2}}\left(\sum_{k>0}2^{2\epsilon k}\|P_k\psi\|_{L_t^{\frac{2}{1-\frac{\epsilon}{2}}}L_x^{\frac{4}{\epsilon}}}^2\right)^{\frac{1}{2}}\lesssim\epsilon_0T^{\frac{\epsilon}{4}},\\
\sum_{k\leq0}\|P_k\psi\|_{L_t^2([0,T];L_x^\infty)}&\lesssim T^{\frac{1}{2}}\sum_{k\leq0}2^k\|P_k\psi\|_{L_t^\infty L_x^2}\lesssim\epsilon_0T^{\frac{1}{2}}.
\end{align*}
Hence,
\begin{equation}\label{eappendixB2}
\int_0^T\left(\epsilon_1+\sum_k\|P_k\psi(t)\|_{L_x^\infty}\right)^2\,\mathrm{d}t\lesssim\epsilon_1^2T+\epsilon_0^2\left(T+T^{\frac{\epsilon}{2}}\right).
\end{equation}
With $\epsilon_1=\epsilon_0^{\frac{1}{2}}$, the only additional bootstrap in physical time is
\begin{equation}\label{eappendixB3}
\sup_{0\leq t\leq T}\|x^2f(t)\|_{L_x^2}\leq2\epsilon_1,\qquad 0<T\leq1.
\end{equation}
Proposition~\ref{pinitial differentiate fields} starts this bootstrap. Integration by parts in $\xi$ gives
\begin{align*}
\|J\psi(t)\|_{H^{\frac{\epsilon}{2}}}^2&=\int\langle\xi\rangle^\epsilon|\nabla_\xi\widehat f(t,\xi)|^2\,\mathrm{d}\xi\\
&=-\operatorname{Re}\int\langle\xi\rangle^\epsilon\widehat f\,\overline{\Delta_\xi\widehat f}\,\mathrm{d}\xi+\frac{1}{2}\int\Delta_\xi\langle\xi\rangle^\epsilon|\widehat f|^2\,\mathrm{d}\xi\\
&\lesssim \|\psi(t)\|_{H^\epsilon}\|x^2f(t)\|_{L_x^2}+\|\psi(t)\|_{L_x^2}^2\lesssim\epsilon_1^2.
\end{align*}
Also,
\begin{equation*}
\|\widehat\psi(t)\|_{L_\xi^\infty}=\|\widehat f(t)\|_{L_\xi^\infty}\lesssim\|f(t)\|_{L_x^1}\lesssim\|f(t)\|_{L_x^2}^{\frac{1}{2}}\|x^2f(t)\|_{L_x^2}^{\frac{1}{2}}\lesssim\epsilon_1.
\end{equation*}
For almost every $t\in[0,T]$, the sum in \eqref{eappendixB2} is finite. At such times,
\begin{equation}\label{eappendixB4}
\|J\psi(t)\|_{L_x^4}^2\lesssim\|\psi(t)\|_{L_x^\infty}\|J^2\psi(t)\|_{L_x^2}\lesssim\epsilon_1\sum_k\|P_k\psi(t)\|_{L_x^\infty}.
\end{equation}

\subsection{Unweighted estimates along the heat flow}
Fix $t$ for which the sum in \eqref{eappendixB2} is finite. The fixed time estimates \eqref{eheatenergyimproved2}, \eqref{epsihat2} and \eqref{econnectionenve1} apply since $\|\psi(t)\|_{H^\epsilon}+\|\widehat\psi(t)\|_{L_\xi^\infty}\lesssim\epsilon_1$. By Bernstein inequality and \eqref{ecaloricconnection}, for the derivative orders used below,
\begin{equation}\label{eappendixB5}
\begin{aligned}
\|\nabla^b\psi\|_{L_x^2}&\lesssim\epsilon_1s^{-\frac{(b-\epsilon)_+}{2}},&\|\nabla^b\psi\|_{L_x^\infty}&\lesssim\epsilon_1s^{-\frac{b+1-\epsilon}{2}},&&0<s\leq1,\\
\|\nabla^b\psi\|_{L_x^2}&\lesssim\epsilon_1s^{-\frac{b+1}{2}},&\|\nabla^b\psi\|_{L_x^\infty}&\lesssim\epsilon_1s^{-\frac{b+2}{2}},&&s\geq1,
\end{aligned}
\end{equation}
\begin{equation}\label{eappendixB6}
\begin{aligned}
\|\nabla^b A\|_{L_x^\infty}&\lesssim\epsilon_1^2s^{-\frac{b+1}{2}+\epsilon},&&0<s\leq1,\\
\|\nabla^b A\|_{L_x^2}+s^{\frac{1}{2}}\|\nabla^b A\|_{L_x^\infty}&\lesssim\epsilon_1^2s^{-\frac{b+2}{2}},&&s\geq1.
\end{aligned}
\end{equation}
By \eqref{eenergyenvelope} and \eqref{edkdef2}, $\sum_k(b_k+d_k)\lesssim1$, using the Fourier bound for $k\leq0$ and $H^\epsilon$ for $k>0$. By \eqref{econnectionenve1}, \eqref{eFenve1}, the heat convolution estimates \eqref{epolynomialconvolution} and \eqref{esingularheatconvolution}, we obtain
\begin{equation*}
(1+s2^{2j})^{10}\|P_j(\psi(s,t)-e^{s\Delta}\psi(t))\|_{L_x^\infty}\lesssim\epsilon_1^2(d_j+b_j)\sum_k\sup_{0\leq s'\leq s}(1+s'2^{2k})^{10}\|P_k\psi(s',t)\|_{L_x^\infty}.
\end{equation*}
Take the supremum over $0\leq s'\leq s$, sum in $j$, absorb $C\epsilon_1^2$ times the resulting norm and let $s\to\infty$:
\begin{equation}\label{eappendixB7}
\sum_k\sup_{s\geq0}(1+s2^{2k})^{10}\|P_k\psi(s,t)\|_{L_x^\infty}\lesssim\sum_k\|P_k\psi(t)\|_{L_x^\infty}.
\end{equation}
By \eqref{eappendixB5}, \eqref{eappendixB7} and Minkowski inequality,
\begin{equation}\label{eappendixB8}
\begin{aligned}
&\|\psi\|_{L_s^\infty L_x^\infty}+\|s^{\frac{1}{2}}\nabla\psi\|_{L_s^\infty L_x^\infty}+\|\nabla^2\psi\|_{L_s^1L_x^\infty}+\|\nabla\psi\|_{L_s^2L_x^\infty}\lesssim\sum_k\|P_k\psi(t)\|_{L_x^\infty},\\
&\|\psi\|_{L_s^2L_x^\infty}+\|\nabla\psi\|_{L_s^1L_x^\infty}\lesssim\epsilon_1.
\end{aligned}
\end{equation}
By \eqref{eappendixB8}, the connection formula \eqref{ecaloricconnection} gives $\|A\|_{L_s^\infty L_x^\infty}\lesssim\epsilon_1\sum_k\|P_k\psi(t)\|_{L_x^\infty}$, and its spatial derivative gives $\|\nabla A\|_{L_s^\infty L_x^\infty}\lesssim(\sum_k\|P_k\psi(t)\|_{L_x^\infty})^2$. Interpolating the uniform bound for $\nabla A$ with \eqref{eappendixB6} for $b=1$ yields
\begin{equation}\label{eappendixB9}
\|A\|_{L_s^\infty L_x^\infty}+\sup_{0<s\leq1}s^{\frac{1-\epsilon}{2}}\|\nabla A(s,t)\|_{L_x^\infty}\lesssim\epsilon_1\sum_k\|P_k\psi(t)\|_{L_x^\infty}.
\end{equation}

\subsection{Weighted heat estimates, heat source terms and bilinear heat integrals}
For $(a,\sigma)=(1,\frac{\epsilon}{2}),(2,0)$ and integers $0\leq b\leq4$, we have
\begin{equation}\label{eappendixB10}
\sup_{s>0}s^{\frac{b}{2}}\|\nabla^bJ_s^a\psi\|_{H^\sigma}+\|s^{\frac{b}{2}}\nabla^{b+1}J_s^a\psi\|_{L_s^2H^\sigma}\lesssim\epsilon_1
\end{equation}
and
\begin{equation}\label{eappendixB11}
\begin{aligned}
\|[2(s+it)\nabla]^aA\|_{L_x^2}&\lesssim\epsilon_1^2,&&s\geq0,\\
\|\nabla^b[2(s+it)\nabla]^aA\|_{H^\sigma}&\lesssim\epsilon_1^2s^{-\frac{(b-\epsilon+\sigma)_+}{2}},&&0<s\leq1,\\
\|2(s+it)\nabla A\|_{L_x^\infty}&\lesssim\epsilon_1^2s^{-\frac{1}{2}+\frac{3\epsilon}{4}},&&0<s\leq1.
\end{aligned}
\end{equation}
We induct on $b$, first up to $b=7$ for $a=1$ to obtain the $L_x^4$ and $L_x^\infty$ bounds below, then up to $b=4$ for $a=2$. For $0<s\leq1$, write \eqref{ecaloricconnection} as $A(1,t)$ plus the integral over $[s,1]$. The required derivatives of $A(1,t)$ are bounded by \eqref{eappendixB6}. On $[s,1]$, use $\left|\frac{s+it}{s'+it}\right|\leq1$ and
\begin{equation*}
2(s'+it)\nabla(\overline u v)=\overline uJ_{s'}v-\overline{J_{s'}u}v+4s'(\nabla\overline u)v.
\end{equation*}
Interpolating the $L_x^2$ and $H_x^1$ bounds for multiplication by $\psi(s',t)$ from \eqref{eappendixB5}, then applying Cauchy-Schwarz, gives
\begin{align*}
\int_s^1\|\psi\nabla^{b+1}J_{s'}^a\psi\|_{H^\sigma}\,\mathrm{d}s'
&\lesssim\epsilon_1\left(\int_s^1(s')^{-1-b+\epsilon-\sigma}\,\mathrm{d}s'\right)^{\frac{1}{2}}\|s^{\frac{b}{2}}\nabla^{b+1}J_s^a\psi\|_{L_s^2((0,1);H^\sigma)}\\
&\lesssim\epsilon_1s^{-\frac{(b-\epsilon+\sigma)_+}{2}}\|s^{\frac{b}{2}}\nabla^{b+1}J_s^a\psi\|_{L_s^2((0,1);H^\sigma)}.
\end{align*}
The term with all derivatives on $A$ is absorbed since $\epsilon_1^2\int_0^1(s')^{-1+\epsilon-\frac{\sigma}{2}}\,\mathrm{d}s'\lesssim\epsilon_1^2$. In the $s^b$ weighted energy estimate for \eqref{eJs6} and its second commutation, write
\begin{equation*}
(\partial_l[2(s+it)\nabla]^aA)\psi=\partial_l\big(([2(s+it)\nabla]^aA)\psi\big)-([2(s+it)\nabla]^aA)\partial_l\psi.
\end{equation*}
After applying $s^{\frac{b}{2}}\nabla^b$, the product inside $\partial_l$ is estimated in $L_s^2H^\sigma$ using $s^{-\frac{1}{2}+\frac{\epsilon-\sigma}{2}}\in L_s^2(0,1)$, and the second product in $L_s^1H^\sigma$ using $s^{-1+\frac{\epsilon-\sigma}{2}}\in L_s^1(0,1)$. We insert the connection bound above and absorb the term with coefficient $C\epsilon_1^2$. For $b=0$ use $\|J^a\psi(t)\|_{H^\sigma}\lesssim\epsilon_1$, for $b\geq1$ the term from $\partial_s(s^b)$ is controlled by the dissipation at order $b-1$.

For $a=1$, Bernstein inequality gives
\begin{equation*}
\|\nabla^jJ_s\psi\|_{L_x^4}+s^{\frac{1}{4}}\|\nabla^jJ_s\psi\|_{L_x^\infty}\lesssim\epsilon_1s^{-\frac{j}{2}-\frac{1}{4}+\frac{\epsilon}{4}},\qquad 0<s\leq1,\quad 0\leq j\leq6.
\end{equation*}
Applying $2(s+it)\nabla$ to \eqref{ecaloricconnection} and integrating $(s')^{-\frac{3}{2}+\frac{3\epsilon}{4}}$ over $[s,1]$ gives the $L_x^\infty$ bound in \eqref{eappendixB11}. For $a=2$, placing both factors carrying one weight in $L_x^4$ gives
\begin{align*}
\|\nabla^b(J_s\psi\nabla J_s\psi)\|_{L_x^2}&\lesssim\epsilon_1^2s^{-1-\frac{b}{2}+\frac{\epsilon}{2}},\\
s^{\frac{b}{2}}\|\nabla^b(\psi(J_s\psi)^2)\|_{L_x^2}&\lesssim\epsilon_1^3s^{-1+\epsilon}.
\end{align*}
The remaining products follow similarly from \eqref{eappendixB5}.

For $s\geq1$, \eqref{eappendixB6} gives $\|[2(s+it)\nabla]^aA\|_{L_x^2}\lesssim\epsilon_1^2s^{\frac{a-2}{2}}$. In the energy estimate starting at $s=1$, \eqref{eappendixB5} and \eqref{eappendixB6} bound the lower weight terms after $s^{\frac{b}{2}}\nabla^b$ by $C\epsilon_1^3s^{-\frac{3}{2}}$ in $H^\sigma$, using the case $a=1$ when $a=2$. The principal terms are absorbed since $\|A\|_{L_s^2W_x^{1,\infty}}+\|\nabla A\|_{L_s^1W_x^{1,\infty}}+\|\psi\|_{L_s^2W_x^{1,\infty}}^2\lesssim\epsilon_1^2$ on $(1,\infty)$. This proves \eqref{eappendixB10} and \eqref{eappendixB11}.

For $|\alpha|=b\leq3$ and $0<s\leq1$, expanding $J_0=J_s-2s\nabla$ and using $[\partial_i,J_{s,j}]=\delta_{ij}$ in \eqref{eappendixB5} and \eqref{eappendixB10} gives $\|J_0^2\partial^\alpha\psi\|_{L_x^2}\lesssim\epsilon_1s^{-\frac{b}{2}}$. By Gagliardo-Nirenberg and \eqref{eappendixB7},
\begin{equation*}
\|J_0\partial^\alpha\psi\|_{L_x^4}^2\lesssim\|\partial^\alpha\psi\|_{L_x^\infty}\|J_0^2\partial^\alpha\psi\|_{L_x^2}\lesssim\epsilon_1s^{-b}\sum_k\|P_k\psi(t)\|_{L_x^\infty}.
\end{equation*}
Writing $\partial^\alpha J_s\psi=J_0\partial^\alpha\psi+2s\nabla\partial^\alpha\psi+[\partial^\alpha,x]\psi$ and interpolating between $L_x^2$ and $L_x^\infty$ for the last two terms, we obtain
\begin{equation}\label{eappendixB12}
\sup_{0<s\leq1}s^{\frac{b}{2}}\|\nabla^bJ_s\psi\|_{L_x^4}\lesssim\left(\epsilon_1\sum_k\|P_k\psi(t)\|_{L_x^\infty}\right)^{\frac{1}{2}},\qquad 0\leq b\leq3.
\end{equation}

From \eqref{ecaloric4}, \eqref{eJs6} and its second commutation,
\begin{equation}\label{eappendixB13}
\begin{aligned}
\|(\partial_s-\Delta)\psi\|_{L_s^1L_x^\infty}&\lesssim\epsilon_1^2\sum_k\|P_k\psi(t)\|_{L_x^\infty},\\
\|(\partial_s-\Delta)J_s\psi\|_{L_s^1L_x^4}&\lesssim\epsilon_1^{\frac{5}{2}}\left(\epsilon_1+\sum_k\|P_k\psi(t)\|_{L_x^\infty}\right)^{\frac{1}{2}},\\
\|(\partial_s-\Delta)J_s^2\psi\|_{L_s^1L_x^2}&\lesssim\epsilon_1^3.
\end{aligned}
\end{equation}
For $a=0,1$, Gagliardo-Nirenberg gives $\|\nabla^a[2(s+it)\nabla]A\|_{L_x^4}^2\lesssim\|\nabla^aA\|_{L_x^\infty}\|\nabla^a[2(s+it)\nabla]^2A\|_{L_x^2}$. By \eqref{eappendixB8}-\eqref{eappendixB11} and \eqref{eappendixB12}, both weighted sources have at worst the factor $s^{-1+\frac{\epsilon}{2}}$ on $(0,1)$, and their respective decay rates on $(1,\infty)$ are $s^{-\frac{9}{4}}$ and $s^{-\frac{3}{2}}$.

For the quadratic connection term, let $|\alpha|+|\beta|=2$ and let $a,b$ be nonnegative integers with $1\leq a+b\leq2$. Then
\begin{equation}\label{eappendixB14}
\left\|\int_0^\infty\partial^\alpha e^{s\Delta}J^a\psi(t)\,\partial^\beta e^{s\Delta}J^b\psi(t)\,\mathrm{d}s\right\|_{L_x^{\frac{4}{a+b}}}\lesssim\epsilon_1^{\frac{a+b}{2}}\left(\epsilon_1+\sum_k\|P_k\psi(t)\|_{L_x^\infty}\right)^{2-\frac{a+b}{2}},
\end{equation}
also with either factor conjugated. Since $|\alpha|+|\beta|=2$, integration in $s$ gives a bilinear symbol of order zero, and the estimate follows from Lemma~\ref{lcoifman-meyer estimates}, \eqref{eappendixB3} and \eqref{eappendixB4}. For $a=b=0$, the dyadic heat estimate used in \eqref{eappendixB8} gives the $L_x^\infty$ bound $(\sum_k\|P_k\psi(t)\|_{L_x^\infty})^2$.

\subsection{The magnetic terms}
Write $A=A^{(2)}+A^{(\geq4)}$ as in \eqref{eAcom1}, with the connection evaluated at $s=0$. We first estimate the commutators in \eqref{ex2f2} with $A^{(2)}$. We compute
\begin{equation}\label{eappendixB15}
\begin{aligned}
t\eta_i(\xi-\eta)_le^{it\Phi_3}&=\frac{1}{2i}\left[\eta_i\partial_{\eta_l}-(\zeta-\eta)_l\partial_{\zeta_i}\right]e^{it\Phi_3},\\
t^2\eta_i\eta_j(\xi-\eta)_le^{it\Phi_3}&=\frac{1}{4}\left[\eta_j\partial_{\eta_l}\partial_{\zeta_i}-(\zeta-\eta)_l\partial_{\zeta_i}\partial_{\zeta_j}-\delta_{il}\partial_{\zeta_j}\right]e^{it\Phi_3}.
\end{aligned}
\end{equation}
The first vector field is divergence free, and both identities remove $(\xi-\eta)_l$ by integration by parts without introducing $t^{-1}$.

Since $\eta=\zeta-(\zeta-\eta)$, both coefficients contribute a second internal spatial derivative. For zero, one or two internal $J_s$ weights, apply \eqref{eappendixB8} and \eqref{eappendixB14} in $L_x^\infty$, $L_x^4$ and $L_x^2$, respectively. To pass from $J_0$ to $J_s$, use $J_0=J_s-2s\nabla$ and $J_{0,i}J_{0,j}=J_{s,i}J_{s,j}-2s(\partial_iJ_{s,j}+\partial_jJ_{s,i})+2s\delta_{ij}+4s^2\partial_i\partial_j$. The corrections are integrable, for example, by $J_se^{s\Delta}=e^{s\Delta}J$, for $0\leq j\leq3$,
\begin{equation*}
\int_0^\infty s\|\nabla^j e^{s\Delta}\psi(t)\|_{L_x^\infty}\|\nabla^{3-j}e^{s\Delta}J\psi(t)\|_{L_x^2}\,\mathrm{d}s\lesssim\epsilon_1^2\left(\int_0^1s^{-1+\frac{\epsilon}{2}}\,\mathrm{d}s+\int_1^\infty s^{-\frac{3}{2}}\,\mathrm{d}s\right)\lesssim\epsilon_1^2.
\end{equation*}
The $2s\nabla$ correction in $L_x^4$ is integrable by the same argument. If a frequency derivative hits the linear symbol of $\partial_l$, only one internal spatial derivative remains. When the remaining weight is on a heat field, the product is bounded in $L_x^2$ by $\epsilon_1^2\|\psi(t)\|_{L_x^\infty}$, when it is on the outer field, the bound is $\epsilon_1\|J\psi(t)\|_{L_x^2}\sum_k\|P_k\psi(t)\|_{L_x^\infty}$. Thus H\"older inequality, \eqref{eappendixB3} and \eqref{eappendixB4} give
\begin{equation*}
t\|(\nabla A^{(2)})\nabla J\psi(t)\|_{L_x^2}+t^2\|(\nabla^2A^{(2)})\nabla\psi(t)\|_{L_x^2}\lesssim\epsilon_1\left(\epsilon_1+\sum_k\|P_k\psi(t)\|_{L_x^\infty}\right)^2.
\end{equation*}

The same bound holds for $A^{(\geq4)}$. Apply \eqref{eappendixB15} to \eqref{eappendixA1}. In the first two integrals, expand the $J_s$ terms by Duhamel as in \eqref{eheatduhamelpk}. For source times $s'\leq s''$, apply the bilinear estimate used in \eqref{eappendixB14} with heat variable $s-s''$ and bound the sources by \eqref{eappendixB13}, the other ordering is symmetric. For the $L_x^\infty$ terms and the $2s\nabla$ corrections, estimate the full and free heat fields separately as above. For the last integral, apply the product rule to $A_l\psi_l$ and use \eqref{eappendixB5}--\eqref{eappendixB12}. We omit the remaining product estimates. Combining the bounds for $A^{(2)}$ and $A^{(\geq4)}$ gives
\begin{equation}\label{eappendixB16}
t\|(\nabla A(t))\nabla J\psi(t)\|_{L_x^2}+t^2\|(\nabla^2A(t))\nabla\psi(t)\|_{L_x^2}\lesssim\epsilon_1\left(\epsilon_1+\sum_k\|P_k\psi(t)\|_{L_x^\infty}\right)^2.
\end{equation}

\subsection{Closing the physical time bootstrap}
In the energy identity for \eqref{ex2f2}, the undifferentiated real potentials vanish and
\begin{equation*}
\operatorname{Im}\int_{\mathbb R^2}\overline{J_iJ_j\psi_m}\left[-2iA_l\partial_lJ_iJ_j\psi_m-i(\partial_lA_l)J_iJ_j\psi_m\right]\,\mathrm{d}x=-\int_{\mathbb R^2}\partial_l\left(A_l|J_iJ_j\psi_m|^2\right)\,\mathrm{d}x=0.
\end{equation*}
We collect the two bounds needed for $A_3$ at $s=0$:
\begin{align*}
t\|\nabla A_3(0,t)\|_{L_x^4}&\lesssim\epsilon_1^{\frac{1}{2}}\left(\epsilon_1+\sum_k\|P_k\psi(t)\|_{L_x^\infty}\right)^{\frac{3}{2}},\\
t^2\|\nabla^2 A_3(0,t)\|_{L_x^2}&\lesssim\epsilon_1\left(\epsilon_1+\sum_k\|P_k\psi(t)\|_{L_x^\infty}\right).
\end{align*}
The quadratic term in \eqref{eA33} follows from \eqref{eA34} and \eqref{eappendixB14}. For the two corrections in \eqref{eA3higherrep}, we argue as in the proof of Lemma~\ref{lpsi3estimate} and use the local weighted heat estimates above to obtain, for $(a,p)=(0,\infty),(1,4),(2,2)$ and $0<s\leq1$,
\begin{align*}
&\|J_0^a(D_l\psi_l-\partial_l e^{s\Delta}\psi_l(t))\|_{L_x^p}+\|J_0^a(\psi_3-i\partial_l e^{s\Delta}\psi_l(t))\|_{L_x^p}\\
&\quad\lesssim\epsilon_1^{2+\frac{a}{2}}\left(\epsilon_1+\sum_k\|P_k\psi(t)\|_{L_x^\infty}\right)^{1-\frac{a}{2}}s^{-\frac{1}{2}+\frac{\epsilon}{8}}.
\end{align*}
Indeed, commute the heat equations with $J_s$ and use $J_0=J_s-2s\nabla$. At $s=0$, estimate $J_0^a(iA_l\psi_l)$ and $J_0^a(-A_l\psi_l)$, respectively, in $L_x^q$ with $\frac{1}{q}=\frac{1}{p}+\frac{1}{2}-\frac{\epsilon}{8}$. For the source terms, split the heat integrals at $\frac{s}{2}$ and estimate the potential coefficients in $L_x^2$ on the first half and in $L_x^\infty$ on the second. The resulting products in \eqref{eA3higherrep} have the integrable factor $s^{-1+\frac{\epsilon}{8}}$. For $s\geq1$, differentiate the products directly using $t\leq1$. We omit the details.

For $F=A_3,\partial_lA_l,A_lA_l$, H\"older inequality, \eqref{eJproduct} and \eqref{eappendixB4} give
\begin{align*}
&\|AJ\psi\|_{L_x^2}+t\|(\nabla A)\psi\|_{L_x^2}+\|J^2(\psi_l\operatorname{Im}(\overline{\psi_l}\psi_m))\|_{L_x^2}\\
&\quad+t\|(\nabla F)J\psi\|_{L_x^2}+t^2\|(\nabla^2F)\psi\|_{L_x^2}\lesssim\epsilon_1\left(\epsilon_1+\sum_k\|P_k\psi(t)\|_{L_x^\infty}\right)^2.
\end{align*}
For $\partial_lA_l$, differentiate the formulas for $A^{(2)}$ and $A^{(\geq4)}$ separately. The product rule for $J_0$ distributes the remaining derivatives carrying a factor of $t$ between the two inputs, giving the same heat integrals as in the previous subsection. The estimates for the terms involving $A_lA_l$ and the lower order magnetic terms follow from \eqref{eappendixB9}, \eqref{eappendixB11}, and Gagliardo-Nirenberg inequality.

By \eqref{eappendixB16}, \eqref{eappendixB3} and \eqref{eappendixB2}, for $0\leq t\leq T\leq1$,
\begin{align*}
\|x^2f(t)\|_{L_x^2}&\lesssim\|x^2f(0)\|_{L_x^2}+\epsilon_1\int_0^t\left(\epsilon_1+\sum_k\|P_k\psi(\tau)\|_{L_x^\infty}\right)^2\,\mathrm{d}\tau\\
&\lesssim\epsilon_0+\epsilon_1^3T+\epsilon_1\epsilon_0^2\left(T+T^{\frac{\epsilon}{2}}\right)\lesssim\epsilon_0.
\end{align*}
Since $\epsilon_1=\epsilon_0^{\frac{1}{2}}$, continuity closes \eqref{eappendixB3} up to $T=1$ for sufficiently small $\epsilon_0$. Finally,
\begin{equation*}
\|xf(t)\|_{L_x^2}+\|f(t)\|_{L_x^1}\lesssim\|f(t)\|_{L_x^2}+\|x^2f(t)\|_{L_x^2}\lesssim\epsilon_0.
\end{equation*}
Thus $\|\widehat f(t)\|_{L_\xi^\infty}\lesssim\epsilon_0$, and $\|\psi(1)\|_{L_x^\infty}=\|e^{i\Delta}f(1)\|_{L_x^\infty}\lesssim\|f(1)\|_{L_x^1}\lesssim\epsilon_0$. Together with \eqref{eappendixB1} and the two dyadic estimates used to derive \eqref{eappendixB2}, this proves \eqref{elocalestimate} and \eqref{eboot-init}.

\end{document}